\documentclass[12pt]{article} 

\usepackage{latexsym, amsmath, 
amscd, amsthm, graphicx, epsfig, amssymb}
\usepackage{xypic}
\usepackage{psfrag}
\DeclareGraphicsExtensions{.ps}
   
\newtheorem{thm}{Theorem}[section]
\newtheorem{defn}[thm]{Definition}
\newtheorem{prop}[thm]{Proposition}

\newtheorem{rema}[thm]{Remark}
\newtheorem{lemma}[thm]{Lemma}
\newtheorem{ass}[thm]{Assumption}

\newcommand{\halmos}{\rule{1ex}{1.4ex}}

\newcommand{\ds}{\displaystyle}

\newcommand{\beq}{\begin{equation}}
\newcommand{\eeq}{\end{equation}}
\newcommand{\bnu}{\begin{enumerate}}
\newcommand{\enu}{\end{enumerate}}
 \newcommand{\bea}{\begin{eqnarray}}
\newcommand{\eea}{\end{eqnarray}}
 \newcommand{\nn}{\nonumber \\}
	
	\newcommand{\lbar}{\bigg\vert}

\renewcommand{\hom}{\mbox{\rm Hom}}

\newcommand{\edo}{\mbox{\rm End}\;}
 \newcommand{\pf}{{\it Proof.}\hspace{2ex}}
 \newcommand{\epf}{\hspace*{\fill}\mbox{$\halmos$}}

\newcommand{\id}{\mbox{\rm id}}
  
\newcommand{\one}{\mathbf{1}}

\newcommand{\tr}{\mbox{\rm Tr}}

\newcommand{\A}{\mathcal{A}}
\newcommand{\C}{\mathbb{C}}
 
\newcommand{\HH}{\mathbb{H}}
\newcommand{\I}{\mathcal{I}}
\newcommand{\N}{\mathbb{N}}

\newcommand{\R}{\mathbb{R}}
\newcommand{\Z}{\mathbb{Z}}
\newcommand{\Y}{\mathcal{Y}}
\newcommand{\V}{\mathcal{V}}

\title{ {\bf Cardy condition for open-closed field algebras } }
\date{}
\author{Liang Kong}

\begin{document}

\setlength{\unitlength}{1cm}

\bibliographystyle{alpha}
\maketitle

\begin{abstract} 
Let $V$ be a vertex operator
algebra satisfying certain reductivity and finiteness conditions
such that $\mathcal{C}_V$, the category of $V$-modules, is a 
modular tensor category. 
We study open-closed field algebras over $V$ equipped with
nondegenerate invariant bilinear forms for both open and closed 
sectors. We show that they give 
algebras over certain $\C$-extension of the so-called
Swiss-cheese partial dioperad, and we can obtain
Ishibashi states easily in such algebras.
Cardy condition can be formulated as an additional condition on 
such open-closed field algebras 
in terms of the action of the modular transformation
$S: \tau \mapsto -\frac{1}{\tau}$
on the space of intertwining operators of $V$.
We then derive a graphical representation of $S$ in 
the modular tensor category $\mathcal{C}_V$. 
This result enables us to give a categorical formulation of 
the Cardy condition and 
the modular invariance condition for 1-point correlation functions on torus. 
Then we incorporate these two conditions and the axioms of open-closed field
algebra over $V$ equipped with nondegenerate invariant bilinear 
forms into a tensor-categorical notion called
Cardy $\mathcal{C}_V|\mathcal{C}_{V\otimes V}$-algebra.
In the end, we give a categorical construction of 
Cardy $\mathcal{C}_V|\mathcal{C}_{V\otimes V}$-algebra in Cardy case.
\end{abstract}

\tableofcontents

\renewcommand{\theequation}{\thesection.\arabic{equation}}
\renewcommand{\thethm}{\thesection.\arabic{thm}}
\setcounter{equation}{0}
\setcounter{thm}{0}
\setcounter{section}{0}

\section{Introduction}

This work is a continuation of the works \cite{HKo1}-\cite{HKo3}\cite{Ko1}\cite{Ko2} and a part of an open-string extension of a program on the closed conformal field theory via the theory of vertex operator algebra. This program was initiated by I. Frenkel and largely developed by Huang \cite{H2}-\cite{H12}. Zhu's work \cite{Z} is also very influential in this development.   

Segal defined the (closed) conformal field theory \cite{Se1} as a projective monoidal functor from the category of finite ordered sets with morphisms being the conformal equivalent classes of Riemann surfaces with parametrized boundaries to the category of locally convex complete topological vector spaces. This definition is very difficult to work with directly. Taking advantage of the theory of vertex operator algebra, Huang suggested to first construct all necessary structure on a densed subspace of the relevant complete topological vector space \cite{H3}, then complete it properly later \cite{functional-1}\cite{functional-2}. This idea guided us in all our previous works \cite{HKo1}-\cite{HKo3}\cite{Ko1}\cite{Ko2}, in particular in our formulation of the notion of algebra over the Swiss-cheese partial operad which catches only some genus-zero information of the whole structure on the densed subspace. This seemingly temporary structure on the densed subspace does not necessarily follow from that on the complete topological vector space, thus has its own independent values and is worthwhile to be formulated properly and extended to a theory of all genus. We will call such a theory on the densed subspace as {\it (closed) partial conformal field theory}.

More precisely, a partial conformal field theory is a projective monoidal functor $\mathcal{F}: \mathcal{RS} \rightarrow \mathcal{GV}$ between two partial categories $\mathcal{RS}, \mathcal{GV}$ in which the compositions of morphisms are not always well-defined. $\mathcal{RS}$ is the category of finite ordered sets, and $\text{Mor}_{\mathcal{RS}}(S_1,S_2)$ for any pair of such sets $S_1, S_2$ is the set of the conformal equivalent classes of closed Riemann surfaces with $|S_1|$ positively oriented punctures and local coordinates and $|S_2|$ negatively oriented punctures and local coordinates \cite{H3}, and the compositions of morphisms in $\mathcal{RS}$ are given by the sewing operations on oppositely oriented punctures (\cite{H3}). Such sewing operations are only partially defined.  $\mathcal{GV}$ is the category of graded vector spaces with finite dimensional homogeneous spaces and a weak topology induced from the restricted dual spaces. For any pair of $A,B \in \text{Ob}(\mathcal{GV})$, $\text{Mor}_{\mathcal{GV}}(A,B)$ is the set of continuous linear maps from $A$ to $\overline{B} :=\prod_{n\in G_B} B_{(n)}$, where the abelian group $G_B$ gives the grading on $B$. For any pair of morphisms $A \xrightarrow{g} \overline{B}$ and $B \xrightarrow{f} \overline{C}$ in $\mathcal{GV}$, using the projector $P_n: \overline{B} \rightarrow B_{(n)}$, we define $f\circ g(u) := \sum_{n\in G_B} f(P_ng(u)), \forall u\in A$, and $f\circ g$ is well-defined only when the sum is absolutely convergent for all $u\in A$.  We remark that by replacing the surfaces with parametrized boundaries in Segal's definition by the surfaces with oriented punctures and local coordinates in $\mathcal{RS}$ we have enlarged the morphism sets. Thus a Segal's functor may not be extendable to a functor on $\mathcal{RS}$. 

The above definition can be easily extended to include open strings by adding to the morphism sets of $\mathcal{RS}$ the conformal equivalent classes of Riemann surfaces with (unparametrized!) boundaries and both oriented interior punctures and oriented punctures on the boundaries. We will call this open-string extended theory as {\it open-closed partial conformal field theory}. We denote the graded vector spaces associated to interior punctures (closed strings at infinity) and boundary punctures (open strings at infinity) by $V_{cl}$ and $V_{op}$ respectively.

The sets of all genus-zero closed surfaces (spheres) with arbitrary number of positively and negatively oriented punctures form a structure of sphere partial dioperad $\mathbb{K}$ (see the definition in Section 1.1). It includes as a substructure the sphere partial operad \cite{H3}, which only includes spheres with a single negatively oriented puncture.  Sphere partial dioperad allows all sewing operations as long as the surfaces after sewing are still genus-zero. 
Hence to construct genus-zero partial conformal field theory amounts to construct projective $\mathbb{K}$-algebras or algebras over certain extension of $\mathbb{K}$. In \cite{HKo2}, we introduced the notion of conformal full field algebra over $V^L\otimes V^R$  equipped with a nondegenerate invariant bilinear form, where $V^L$ and $V^R$ are two vertex operator algebras of central charge $c^L$ and $c^R$ respectively and satisfy certain finiteness and reductivity conditions. Theorem 2.7 in \cite{Ko1} can be reformulated as follow: a conformal full field algebra $V_{cl}$ over $V^L\otimes V^R$ equipped with a nondegenerate invariant bilinear form canonically gives on $V_{cl}$ an algebra over $\tilde{\mathbb{K}}^{c^L}\otimes \overline{\tilde{\mathbb{K}}^{\overline{c^R}}}$, which is a partial dioperad extension of $\mathbb{K}$.

The sets of all genus-zero surfaces with one (unparametrized!) boundary component (disks) and arbitrary number of oriented boundary punctures form the so-called disk partial dioperad denoted by $\mathbb{D}$ (see the definition in Section 1.1). $\mathcal{F}$ restricted on $\mathbb{D}$ induces a structure of algebra over certain extension of $\mathbb{D}$ on $V_{op}$. In \cite{HKo1}, Huang and I introduced the notion of open-string vertex operator algebra. We will show in Section 1.3 that an open-string vertex operator algebra of central charge $c$ equipped with 
a nondegenerate invariant bilinear form canonically gives 
an algebra over $\tilde{\mathbb{D}}^c$, 
which is a partial dioperad extension of $\mathbb{D}$. 

\begin{figure}
\begin{center}
\includegraphics[width=0.6\textwidth]{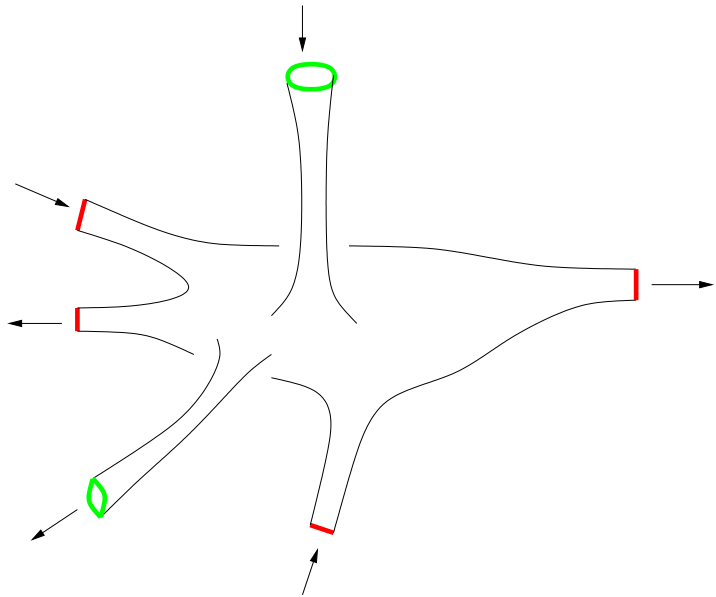}
\end{center}
\caption{\label{sw-diop-fig}
a typical element in Swiss-Cheese partial dioperad $\mathbb{S}$}
\end{figure}

The sets of all genus-zero surfaces with only one (unparametrized!) boundary component and arbitrary number of oriented interior punctures and boundary punctures form the so-called Swiss-Cheese partial dioperad $\mathbb{S}$ (see the definition in Section 2.2). A typical elements in $\mathbb{S}$ is depicted in Figure \ref{sw-diop-fig}, where boundary punctures are drawn as an infinitely long strip (or an open string) and interior punctures are drawn as an infinitely long tube (or a closed string). Let $V$ be a vertex operator algebra of central charge $c$ satisfying the conditions in Theorem \ref{ioa}. We will show in Section 2.3 that an open-closed field algebra over $V$, which contains an open-string vertex operator algebra $V_{op}$ and a conformal full field algebra $V_{cl}$, satisfying a $V$-invariant boundary condition \cite{Ko2} and equipped with nondegenerate invariant bilinear forms on both $V_{op}$ and $V_{cl}$, canonically gives an algebra over $\tilde{\mathbb{S}}^c$, which is an extension of $\mathbb{S}$. 

Note that all surfaces of any genus with arbitrary number of boundary components, interior punctures and boundary punctures can be obtained by applying sewing operations to elements in $\mathbb{S}$. Therefore, except some compatibility conditions coming from different decompositions of the same surfaces with higher genus, an algebra over $\tilde{\mathbb{S}}^c$ if extendable uniquely determine the entire theory of all genus. In particular, the famous Ishibashi states \cite{I} can be obtained in an open-closed field algebra equipped with nondegenerate invariant bilinear forms. An Ishibashi state is a coherent state $\psi$ in $\overline{V_{cl}}$ such that $(L^L(n) - L^R(-n))   \psi =0, \forall n\in \Z$. In physics, Ishibashi states are obtained by solving above equation. In Section 2.4, we show how to obtain Ishibashi states constructively and geometrically from vacuum-like states in $V_{op}$.

\begin{figure}
\begin{center}
\includegraphics[width=0.3\textwidth]{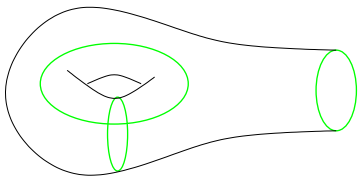}
\end{center}
\caption{\label{mod-inv-top-fig}modular invariance condition for 1-pt correlation functions on torus}
\end{figure}

\begin{figure}
\begin{center}
\includegraphics[width=0.85\textwidth]{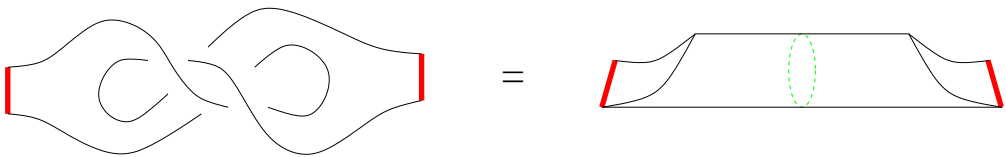}
\end{center}
\caption{\label{cardy-top-fig}Cardy condition}
\end{figure}

Other surfaces which are not included in $\mathbb{S}$ only provide additional compatibility conditions. In 2-d topological field theories, only three additional compatibility conditions are needed to ensure the consistency of a theory of all genus \cite{La}\cite{Mo1}\cite{Se2}\cite{MSeg}\cite{AN}\cite{LP}. The first compatibility condition says that both $V_{op}$ and $V_{cl}$ are finite dimensional. This guarantee the convergence of all higher genus correlation functions. The second condition is the {\it modular invariance condition for 1-point correlation functions on torus}. It is due to two different decompositions of the same torus as depicted in Figure \ref{mod-inv-top-fig}. This condition is automatically satisfied in 2-d topological field theories and but nontrivial in conformal field theories. The third condition is the famous {\it Cardy condition} which is again due to two different decompositions of a single surface as shown in Figure \ref{cardy-top-fig}.

Now we turn to the compatibility conditions in conformal field theory. In this case, both $V_{cl}$ and $V_{op}$ in any nontrivial theory are infinite dimensional. We need require the convergence of all correlation functions of all genus. This is a highly nontrivial condition and not easy to check for examples. So far the only known convergence results are in genus-zero \cite{H8} and genus-one theories \cite{Z}\cite{DLM}\cite{Mi1}\cite{Mi2}\cite{H10}. We recall a theorem by Huang.
\begin{thm}[\cite{H11}\cite{H12}]  \label{ioa}
If $(V, Y, \one, \omega)$ is a simple vertex operator
algebra $V$ satisfying the following conditions: 1. $V$ is $C_2$-cofinite; 2. $V_{n}=0$ for $n<0$, $V_{(0)}=\C \one$, $V'$ is isomorphic to $V$ as $V$-module; 3. all $\N$-gradable weak $V$-modules are completely reducible, then the direct sum of all inequivalent irreducible $V$-modules has a natural structure of intertwining operator algebra \cite{H7}, and the category $\mathcal{C}_V$ of $V$-modules has a structure of vertex tensor category \cite{HL1}-\cite{HL4}\cite{H4} and modular tensor category \cite{RT}\cite{T}. 
\end{thm}


For the intertwining operator algebra given in Theorem \ref{ioa}, 
Huang also proved in \cite{H8}\cite{H10} that the products of 
intertwining operators and $q$-traces of them 
have certain nice convergence and analytic extension properties. 
These properties are sufficient for the construction 
of genus-zero and genus-one correlation
functions in closed partial conformal field theory \cite{HKo2}\cite{HKo3}. 
Since the modular tensor category $\mathcal{C}_V$ supports 
an action of mapping class groups of all genus, it is
reasonable to believe that the above conditions on $V$ 
are also sufficient to guarantee the convergence of correlation functions 
of all genus. 
\begin{ass}  
In this work, we fix a vertex operator algebra $V$ with central charge $c$, which is assumed to satisfy the conditions in
Theorem \ref{ioa} without further announcement. 
\end{ass}

Besides the convergence condition and the axioms of projective $\mathbb{K}$-algebra, Sonoda \cite{So} argued on a physical level of rigor that it is sufficient to check the modular invariance condition in Figure \ref{mod-inv-top-fig} in order to have a consistent partial conformal field theory of all genus. This modular invariance condition was studied in \cite{HKo3}. In the framework of  conformal full field algebra over $V^L\otimes V^R$, it was formulated algebraically as a modular invariance property of an intertwining operator of $V^L\otimes V^R$.

For an open-closed theory, besides the convergence condition, the axioms of algebra over $\tilde{\mathbb{S}}^c$ and the modular invariance condition, Lewellen \cite{Le} argued on a physical level of rigor that the only remaining compatibility condition one needs is the Cardy condition. Cardy condition in open-closed conformal field theory is more complicated than that in topological theory and has never been fully written down by physicists. In Section 3.1, we derive the Cardy condition from the axioms of open-closed partial conformal field theory by writing out two sides of the Cardy condition in Figure \ref{cardy-top-fig} explicitly in terms of the ingredients of open-closed field algebra (see Definition 3.4). Using results in \cite{H10}\cite{H11}, we can show that
the Cardy condition can be reformulated as an invariance condition of 
the modular transformation 
$S: \tau \mapsto -\frac{1}{\tau}$ on intertwining operators.

There are still more compatibility conditions which were not discussed in \cite{So}\cite{Le}. One also need to proved certain algebraic version of uniformization theorems (see \cite{H2} \cite{H3} for the genus-zero case). Such results for genus large than $0$ are still not available. But it seems that no additional assumption on $V$ is needed. They should follow automatically from the properties of Virasoro algebra and intertwining operators. This uniformization problem and convergence problems are not pursued further in this work.

In order to take the advantages of some powerful tools, such as the graphic calculus in tensor category, in the study of the modular invariance condition and the Cardy condition, we would like to obtain categorical formulations of these conditions. It requires us to know the action of the modular transformation $S$ in $\mathcal{C}_V$. Although the action of $SL(2,\Z)$ in a modular tensor category is explicitly known \cite{MSei2}\cite{V}\cite{Ly}\cite{Ki}\cite{BK2}, its relation to the modular transformation of the $q$-trace of the products of intertwining operators of $V$ is not completely clear. This relation was first suggested by I. Frenkel and studied by Moore and Seiberg in \cite{MSei2} but only on a physical level of rigor. Using Huang's results on modular tensor category \cite{H11}\cite{H12}, we derive a graphical representation of $S$ in $\mathcal{C}_V$. This result enables us to give categorical formulations of the modular invariance condition and the Cardy condition. We incorporate them with the categorical formulation of open-closed field algebra over $V$ equipped with nondegenerate invariant bilinear forms into a tensor-categorical notion called {\it Cardy $\mathcal{C}_V|\mathcal{C}_{V\otimes V}$-algebra}. As we discussed in previous paragraphs, it is reasonable to believe that open-closed partial conformal field theories of all genus satisfying the $V$-invariant boundary condition \cite{Ko2} are classified by Cardy $\mathcal{C}_V|\mathcal{C}_{V\otimes V}$-algebras. However, to construct the high-genus theories explicitly is still a hard open problem which is not pursued in this work. In the end, we give an explicit construction of Cardy $\mathcal{C}_V|\mathcal{C}_{V\otimes V}$-algebra in the so-called Cardy case in physics literature (see for example \cite{FFRS}).

Note that this work is somewhat complementary to the works of Fjelstad, Fuchs, Runkel, Schweigert \cite{FS}\cite{FRS1}-\cite{FRS4}\cite{FjFRS1}\cite{FjFRS2}. We will leave a detailed study of the relationship between these two approaches in \cite{cardysew}.

The layout of this work is as follow: in Section 1, 
we introduce the notion of 
sphere partial dioperad and disk partial dioperad
and study algebras over them; 
in Section 2, we introduce the notions of 
Swiss-cheese partial dioperad $\mathbb{S}$ and its $\C$-extension
$\tilde{\mathbb{S}}^c$, and show that
an open-closed conformal field algebra over $V$ equipped with
nondegenerate invariant bilinear forms canonically gives
an algebra over $\tilde{\mathbb{S}}^c$, and we also construct Ishibashi
states in such algebras; in Section 3, we give two
formulations of the Cardy condition; in Section 4, we derive
a graphic representation of the modular transformation $S$; 
in Section 5, we give the categorical formulations of the nondegenerate
invariant bilinear forms, the modular invariance condition and the 
Cardy condition. Then we introduce the 
notion of Cardy $\mathcal{C}_V|\mathcal{C}_{V\otimes V}$-algebra
and give a construction.

Convention of notations: $\N, \Z, \Z_+, \R, \R_+, \C$ 
denote the set of 
natural numbers, integers, positive integers, real numbers, 
positive real numbers, complex numbers, respectively. 
Let $\HH = \{ z\in \C | \text{Im} z >0\}$
and $\overline{\HH} = \{ z\in \C | \text{Im} z < 0\}$.  
Let $\hat{\R}$, $\hat{\C}$ and $\hat{\HH}$ 
be the one point compactification of real line, complex plane and
up-half plane (including the $\R$-boundary) respectively. 
Let $\R_+$ and $\C^{\times}$ be the multiplication groups 
of positive real and nonzero complex numbers respectively. 
The ground field is always chosen to be $\C$. 

Throughout this work, we choose a branch cut 
for logarithm as follow: 
\beq  \label{branch-cut}
\log z = \log |z| + i \text{Arg}\, z, \quad \quad 0\leq 
\text{Arg}\, z  < 2\pi. 
\eeq
We define power functions of two different types of 
complex variables as follow:
\beq  \label{power-f}
z^s := e^{s\log z}, \quad\quad
\bar{z}^s := e^{s\, \overline{\log z}}, \quad\quad \forall s\in \R.
\eeq


\paragraph{Acknowledgment}
The results in Section 2.4 and 3.1 are included in author's thesis. 
I want to thank my advisor Yi-Zhi Huang 
for introducing me to this interesting field and 
for his constant support and many important suggestions 
for improvement. I thank C. Schweigert for 
telling me the meaning of boundary states from
a physical point of view. I also want to 
thank I. Frenkel, J. Fuchs, A. Kirillov, Jr.
and C. Schweigert for some inspiring conversations on
the subject of Section 4.2.

\renewcommand{\theequation}{\thesection.\arabic{equation}}
\renewcommand{\thethm}{\thesection.\arabic{thm}}
\setcounter{equation}{0}
\setcounter{thm}{0}

\section{Partial dioperads}

In Section 1.1, we recall the definition of 
(partial) dioperad and algebra over it, and 
introduce sphere partial dioperad $\mathbb{K}$, 
disk partial dioperad $\mathbb{D}$ and their extensions 
$\tilde{\mathbb{K}}^{c^L}\otimes \overline{\tilde{\mathbb{K}}^{\overline{c^R}}}$,
$\tilde{\mathbb{D}}^c$ as examples. 
In Section 1.2, we discuss an algebra over 
$\tilde{\mathbb{K}}^{c^L}\otimes \overline{\tilde{\mathbb{K}}^{\overline{c^R}}}$
from a conformal full field algebra over $V^L\otimes V^R$. 
In Section 1.3, we discuss an algebra over 
$\tilde{\mathbb{D}}^c$ from an open-string vertex operator algebra.

\subsection{Partial dioperads}

Let us first recall the definition of dioperad given by 
Gan \cite{G}. Let $S_n$ be the automorphism group of the set
$\{1, \dots, n\}$ for $n\in \Z_+$.   
Let $m=m_1+\dots +m_n$ be an ordered partition and $\sigma\in S_n$. 
The block permutation $\sigma_{(m_1,\dots, m_n)} \in S_m$ is
the permutation which permutes $n$ intervals
of lengths $m_1, \dots, m_n$ in the same way as $\sigma$ permutes
$1, \dots, n$. Let $\sigma_i\in S_{m_i}$, $i=1, \dots, k$.
We view the element 
$(\sigma_1, \dots, \sigma_k)\in S_{m_1}\times \dots \times S_{m_k}$
naturally as an element in $S_m$ by the canonical embedding 
$S_{m_1}\times \dots \times S_{m_k}\hookrightarrow S_m$.
For any $\sigma\in S_n$ and $1\leq i\leq n$, we define a map
$\hat{i}: \{ 1, \dots, n-1\}  
\rightarrow \{1, \dots, n\} $ 
by $\hat{i}(j)=j$ if $j<i$ 
and $\hat{i}(j)=j+1$ if $j\geq i$ and an element
$\hat{i}(\sigma)\in S_{n-1}$ by 
$$
\hat{i}(\sigma)(j) := \hat{i}^{-1}
\circ \sigma \circ \widehat{\sigma^{-1}(i)} (j).
$$

\begin{defn}  {\rm
A {\it dioperad} consists of a family of sets 
$\{ \mathcal{P}(m,n) \}_{m,n\in \N}$
with an action of $S_m\times S_n$ on $\mathcal{P}(m,n)$ 
for each pair of $m,n\in \Z_+$, a distinguished 
element $I_{\mathcal{P}}\in \mathcal{P}(1,1)$ and substitution maps 
\bea
 \mathcal{P}(m,n) \times 
\mathcal{P}(k_1,l_1) \times
\dots \times \mathcal{P}(k_n,l_n) 
&\xrightarrow{\gamma_{(i_1,\dots, i_n)}}&
\mathcal{P}(m-n+k_1\dots +k_n, l_1+\dots +l_n)  \nn
(P, P_1, \dots, P_n) &\mapsto & 
\gamma_{(i_1,\dots, i_n)}(P; P_1, \dots, P_n)
\eea
for $m,n,l_1,\dots,l_n\in \N, k_1,\dots,k_n\in \Z_+$
and $1\leq i_j \leq k_j, j=1,\dots, n$,
satisfying the following axioms:
\bnu
\item {\it Unit properties}: for $P\in \mathcal{P}(m,n)$, 
\bnu
\item {\it left unit property}:
$\gamma_{(i)} (I_{\mathcal{P}}; P) = P$ for $1\leq i\leq m$, 
\item {\it right unit property}:
$\gamma_{(1,\dots, 1)}(P; I_{\mathcal{P}}, \dots, I_{\mathcal{P}}) = P$. 
\enu
\item {\it Associativity}: for $P\in \mathcal{P}(m,n)$,
$Q_i\in \mathcal{P}(k_i,l_i), i=1,\dots, n$, 
$R_j\in \mathcal{P}(s_j,t_j), j=1,\dots, l=l_1+\dots +l_n$, we have
\beq
\gamma_{(q_1, \dots, q_l)} 
\big( \gamma_{(p_1,\dots, p_n)}(P; Q_1,\dots, Q_n); 
R_1, \dots, R_l \big)   
= \gamma_{(p_1,\dots, p_n)}(P ;  P_1, \dots, P_n)  \label{diop-asso}
\eeq
where $P_i = \gamma_{(q_{l_1+\dots +l_{i-1}+1}, \dots, q_{l_1+\dots+l_i})}
(Q_i ;  R_{l_1+\dots +l_{i-1}+1}, \dots, R_{l_1+\dots+l_i})$
for $i=1,\dots, n$. 

\item {\it Permutation property}:
For $P\in \mathcal{P}(m,n)$, $Q_i\in P(k_i,l_i), i=1,\dots, n$, 
$(\sigma, \tau) \in S_m\times S_n$, $(\sigma_i, \tau_i)\in 
S_{k_i}\times S_{l_i}, i=1,\dots,n$,
\bea
&&\hspace{-1cm}
\gamma_{(i_1, \dots, i_n)}( (\sigma, \tau) (P) ; Q_1, \dots, Q_n) \nn
&&\hspace{-0.5cm}=  
((\sigma, \tau_{(k_1-1,\dots,k_n-1)}), \tau_{(l_1,\dots, l_n)})  
\gamma_{(i_{\tau(1)}, \dots, i_{\tau(n)})}( P ; Q_{\tau(1)}, \dots, Q_{\tau(n)}), \nn
&&\hspace{-1cm}\gamma_{(i_1, \dots, i_n)}( P;  (\sigma_1, \tau_1)(Q_1), 
\dots, (\sigma_n,\tau_n)(Q_n)) \nn
&&\hspace{-0.5cm}=( (\id, \hat{i}_1(\sigma_1), 
\dots, \hat{i}_n(\sigma_n)), (\tau_1, \dots, \tau_n))
\gamma_{(\sigma_1^{-1}(i_1), \dots, 
\sigma_n^{-1}(i_n))}(P; Q_1, \dots, Q_n). \nonumber
\eea

\enu
}
\end{defn}

We denote such dioperad as 
$(\mathcal{P}, \gamma_{\mathcal{P}}, I_{\mathcal{P}})$ or
simply $\mathcal{P}$.

\begin{rema} {\rm
We define compositions $_{^i}\circ_{^j}$ as follow:
$$
P {}_{^i}\circ_{^j} Q :=
\gamma_{(1,\dots,1,j,1,\dots 1)}(P; I_{\mathcal{P}},\dots, I_{\mathcal{P}},
Q, I_{\mathcal{P}},\dots, I_{\mathcal{P}}). 
$$ 
It is easy to see that $\gamma_{(i_1,\dots, i_n)}$ can be
reobtained from $_{^i}\circ_{^j}$.
In \cite{G}, the definition of dioperad is given in terms of
$_{^i}\circ_{^j}$ instead of $\gamma_{(i_1,\dots, i_n)}$.
}
\end{rema}

\begin{defn}  {\rm 
A {\it partial dioperad} has a similar definition 
as that of dioperad except the map $\gamma_{(i_1,\dots, i_n)}$
or $_{^i}\circ_{^j}$ 
is only partially defined and the same associativity hold
whenever both sides of (\ref{diop-asso}) exist. 
A {\it (partial) nonassociative dioperad} 
consists of the same data as those of (partial) dioperad 
satisfying all the axioms of (partial) dioperad 
except the associativity.
}
\end{defn}

The notion of homomorphism and isomorphism of 
(partial pseudo-) dioperad are naturally defined. 

\begin{rema} {\rm 
In the case of partial dioperad, the definition using 
$\gamma_{(i_1,\dots, i_n)}$
or $_{^i}\circ_{^j}$ may have subtle differences in 
the domains on which $\gamma_{(i_1,\dots, i_n)}$
or $_{^i}\circ_{^j}$ is defined (see appendix $C$ in 
\cite{H4} for more details).
These differences have no effect on those algebras over partial
dioperads considered in this work. So we will simply 
ignore these differences.    
}
\end{rema}

\begin{rema}  {\rm 
Notice that a (partial nonassociative) 
dioperad $\{ \mathcal{P}(m,n) \}_{m,n\in \N}$
naturally contains a (partial nonassociative) operad 
$\{ \mathcal{P}(1, n)\}_{n\in \N}$ as a substructure. 
}
\end{rema}

\begin{defn} {\rm 
A subset $G$ of $\mathcal{P}(1,1)$ is called a 
{\it rescaling group of $\mathcal{P}$} if
the following conditions are satisfied: 
\bnu
\item For any $g, g_1, \dots, g_n\in G, Q\in P(m,n)$, 
$\gamma_{(i)}(g; Q)$ and $\gamma_{(1,\dots,1)}(Q; g_1, \dots, g_n)$
are always well-defined
for $1\leq i\leq m$.
\item $I_{\mathcal{P}}\in G$ 
and $G$ together with the identity element $I_{\mathcal{P}}$
and multiplication map $\gamma_{(1)}: G\times G \rightarrow G$ 
is a group.  
\enu
}
\end{defn}

\begin{defn} {\rm 
A {\it $G$-rescalable partial dioperad} is partial dioperad 
$\mathcal{P}$ such that for any $P\in \mathcal{P}(m,n)$,
$Q_i\in \mathcal{P}(k_i, l_i), i=1,\dots,n$ 
there exist $g_i\in G, i=1,\dots, n$ 
such that $\gamma_{(i_1,\dots, i_n)}\big( P ;  
\gamma_{(i_1)}(g_1; Q_1), \dots, \gamma_{(i_n)}(g_n; Q_n) \big)$
is well-defined. 
}
\end{defn}

The first example of partial dioperad in which 
we are interested in this work comes from 
$\mathbb{K}= \{ \mathbb{K}(n_-,n_+) \}_{n_-,n_+\in \N}$ \cite{Ko2},
a natural extension of sphere partial operad $K$ \cite{H3}. 
More precisely, $\mathbb{K}(n_-,n_+)$ is the set of the conformal
equivalent classes of sphere with $n_-$ ($n_+$) ordered negatively 
(positively) oriented punctures and local coordinate
map around each puncture. In particular, $\mathbb{K}(0,0)$ 
is an one-element set consisting of the conformal equivalent class of 
a sphere with no additional structure. We simply denote this
element as $\hat{\C}$. We use
\bea   \label{ele-Q}
Q &=& (\, (z_{-1}; a_0^{(-1)}, A^{(-1)}), \dots, 
(z_{-n_-}; a_0^{(-n_-)}, A^{(-n_-)})  \, |  \nn
&& \hspace{3cm} (z_{1}; a_0^{(1)}, A^{(1)})\dots 
(z_{n_+}; a_0^{(n_+)}, A^{(n_+)}) \,)_{\mathbb{K}}, 
\eea
where $z_i\in \hat{\C}, a_0^{(i)}\in \C^{\times}, A^{(i)}\in \C^{\infty}$ for
$i=-n_-, \dots, -1, 1, \dots, n_+$, 
to denote a sphere $\hat{\C}$ with positively (negatively) oriented
punctures at $z_i\in \hat{\C}$ for $i=1,\dots, n_+$ 
($i=-1, \dots, -n_-$), and with local coordinate map 
$f_i$ around each punctures $z_i$ given by: 
\bea
f_i(w) &=& e^{\sum_{j=1}^{\infty} A_j^{(i)} x^{j+1}\frac{d}{dx}} 
(a_0^{(i)})^{x\frac{d}{dx}}x 
\lbar_{x=w-z_i} \hspace{1cm} \mbox{if $z_i\in \C$},   
\label{K-case-1}   \\
&=& e^{\sum_{j=1}^{\infty} A_j^{(i)} x^{j+1}\frac{d}{dx} } 
(a_0^{(i)})^{x\frac{d}{dx}}x 
\lbar_{x=\frac{-1}{w}} \hspace{1.3cm}\mbox{if $z_i=\infty$}.
\label{K-case-2}
\eea
We introduce a useful notation $\bar{Q}$ defined as follow:
\bea
\bar{Q} &=& (\, 
(\bar{z}_{-1}; \overline{a_0^{(-1)}}, \overline{A^{(-1)}}), \dots,
(\bar{z}_{-n_-}; \overline{a_0^{(-n_-)}}, \overline{A^{(-n_-)})}|  \nn
&& \hspace{3cm} 
(\bar{z}_{1}; \overline{a_0^{(1)}}, \overline{A^{(1)}})\dots 
(\bar{z}_{n_+}; \overline{a_0^{(n_+)}}, \overline{A^{(n_+)}}) 
\,)_{\mathbb{K}}, 
\eea
where the ``overline'' represents complex conjugations.

We denote the set of all such $Q$ as 
$\mathcal{T}_{\mathbb{K}}(n_-, n_+)$. Let $\mathcal{T}_{\mathbb{K}}: = 
\{ \mathcal{T}_{\mathbb{K}}(n_-, n_+) \}_{n_-, n_+\in \N}$. 
There is an action of $SL(2, \C)$ on 
$\mathcal{T}_{\mathbb{K}}(n_-, n_+)$ as Mobius transformations. 
It is clear that 
\beq
\mathbb{K}(n_-, n_+) = \mathcal{T}_{\mathbb{K}}(n_-, n_+) / SL(2, \C). 
\eeq
We denote the quotient map 
$\mathcal{T}_{\mathbb{K}} \rightarrow \mathbb{K}$ as $\pi_{\mathbb{K}}$. 
The identity $I_{\mathbb{K}}\in \mathbb{K}(1,1)$ is given by 
\beq
I_{\mathbb{K}} = \pi_{\mathbb{K}}(\, (\infty, 1, \mathbf{0}) | (0, 1, \mathbf{0})\, )
\eeq 
where $\mathbf{0}=(0,0,\dots) \in \prod_{n=1}^{\infty} \C$.
The composition $ _{^i}\circ_{^j}$ is provided
by the sewing operation $ _{^i}\infty_{^{-j}}$ \cite{H3}. 
In particular, for
$n_1,m_2\geq 1$, $P\in \mathbb{K}(m_1,n_1)$ and $Q\in \mathbb{K}(m_2,n_2)$,
$P _{^i}\infty_{^{-j}} Q$ is the sphere with punctures obtained by
sewing the $i$-th positively oriented puncture of $P$ 
with the $j$-th negatively oriented puncture of $Q$.  
The $S_{n_-}\times S_{n_+}$-action on $\mathbb{K}(n_-,n_+)$ 
(or $\mathcal{T}(n_-,n_+)$) is the natural one. 
Moreover,  the set 
\beq
\{ (\, (\infty, 1,\mathbf{0}) | (0, a, \mathbf{0})\,) |a\in \C^{\times} 
\}
\eeq  
together with multiplication $_{^1}\infty_{^1}$ is a group which
can be canonically identified with group $\C^{\times}$. 
It is clear that 
$\mathbb{K}$ is a $\C^{\times}$-rescalable partial dioperad.
We call it {\it sphere partial dioperad}.

The $\C$-extensions of $\mathbb{K}$, such as $\mathbb{K}^c$
and $\tilde{\mathbb{K}}^{c^L}\otimes 
\overline{\tilde{\mathbb{K}}^{\overline{c^R}}}$ for $c, c^L, c^R\in \C$,
are trivial line bundles over $\mathbb{K}$ with 
natural $\C^{\times}$-rescalable partial dioperad structures
(see Section 6.8 in \cite{H4}). Moreover, we denote the 
canonical section $\mathbb{K} \rightarrow \tilde{\mathbb{K}}^{c^L}\otimes 
\overline{\tilde{\mathbb{K}}^{\overline{c^R}}}$ as 
$\psi_{\mathbb{K}}$.

The next example of partial dioperad is 
$\mathbb{D} = \{ \mathbb{D}(n_-, n_+) \}_{n_-, n_+\in \N}$, 
which is an extension of the partial operad of disk with strips
$\Upsilon$ introduced \cite{HKo1}. 
More precisely, $\mathbb{D}(n_-, n_+)$ is the set of conformal
equivalent classes of disks with ordered punctures 
on their boundaries and local coordinate map
around each puncture. In particular, $\mathbb{D}(0,0)$ is
an one-element set consisting of the conformal equivalent class of 
a disk with no additional structure. We simply denote this
element as $\hat{\HH}$. 
We use
\bea   \label{ele-Q-D}
Q &=& (\, (r_{-n_-}; b_0^{(-n_-)}, B^{(-n_-)}), \dots, 
(r_{-1}; b_0^{(-1)}, B^{(-1)})|  \nn
&& \hspace{3cm} (r_{1}; b_0^{(1)}, B^{(1)})\dots 
(r_{n_+}; b_0^{(n_+)}, B^{(n_+)}) \,)_{\mathbb{D}}, 
\eea
where $r_i\in \hat{\R}, a_0^{(i)}\in \R_+, B^{(i)}\in \R^{\infty}$ for
$i=-n_-, \dots, -1, 1, \dots, n_+$, 
to denote a disk $\hat{\HH}$ with positively (negatively) oriented
punctures at $r_i\in \hat{\R}$ for $i=1,\dots, n_+$ 
($i=-1, \dots, -n_-$), and with local coordinate map 
$g_i$ around each punctures $r_i$ given by: 
\bea
g_i(w) &=& e^{\sum_{j=1}^{\infty} B_j^{(i)} x^{j+1}\frac{d}{dx}} 
(b_0^{(i)})^{x\frac{d}{dx}} x 
\lbar_{x=w-r_i} \hspace{1cm} \mbox{if $r_i\in \R$},   
\label{D-case-1}   \\
&=& e^{\sum_{j=1}^{\infty} B_j^{(i)} x^{j+1}\frac{d}{dx} } 
(b_0^{(i)})^{x\frac{d}{dx}}x 
\lbar_{x=\frac{-1}{w}} \hspace{1.3cm}\mbox{if $r_i=\infty$}.
\label{D-case-2}
\eea

We denote the set of all such $Q$ as 
$\mathcal{T}_{\mathbb{D}}(n_-, n_+)$. Let 
$\mathcal{T}_{\mathbb{D}}: = 
\{ \mathcal{T}_{\mathbb{D}}(n_-, n_+) \}_{n_-, n_+\in \N}$. 
The automorphism group of $\HH$, $SL(2,\R)$, 
naturally acts on $\mathcal{T}_{\mathbb{D}}(n_-, n_+)$. 
It is clear that 
\beq
\mathbb{D}(n_-, n_+) = \mathcal{T}_{\mathbb{K}}(n_-, n_+) / SL(2, \R). 
\eeq
We denote the quotient map as $\pi_{\mathbb{D}}$. 
The identity $I_{\mathbb{D}} \in \mathbb{D}(1,1)$ is given by 
\beq
I_{\mathbb{D}} = \pi_{\mathbb{D}}
(\, (\infty, 1, \mathbf{0}) | (0, 1, \mathbf{0})\, ).
\eeq 
The composition $ _{^i}\circ_{^j}$ is provided
by the sewing operation $ _{^i}\infty_{^{-j}}$ \cite{HKo1}. 
In particular, for
$n_1,m_2\geq 1$, $P\in \mathbb{D}(m_-, m_+)$ 
and $Q\in \mathbb{D}(n_-,n_+)$,
$P _{^i}\infty_{^{-j}} Q$ is the disk with strips obtained by
sewing the $i$-th positively oriented puncture of $P$ 
with the $j$-th negatively oriented puncture of $Q$.  
The $S_{n_-}\times S_{n_+}$-action on $\mathbb{D}(n_-,n_+)$ 
(or $\mathcal{T}(n_-,n_+)$) is the natural one. 
The set 
\beq
\{ (\, (\infty, 1,\mathbf{0}) | (0, a, \mathbf{0})\,) |a\in \R_+ \}
\eeq
together with multiplication $_{^1}\infty_{^1}$ is a group which
can be canonically identified with group $\R_+$. It is clear that 
$\mathbb{D}$ is a $\R_+$-rescalable partial dioperad. We call it
{\it disk partial dioperad}.

$\mathbb{D}$ can be naturally embedded to $\mathbb{K}$ 
as a sub-dioperad. 
The $\C$-extension $\mathbb{D}^c$ of $\mathbb{D}$ 
for $c\in \C$ is just 
the restriction of the line bundle $\mathbb{K}^c$ on $\mathbb{D}$.
$\mathbb{D}^c$ is also a $\R_+$-rescalable partial dioperad
and a partial sub-dioperad of $\mathbb{K}^c$. 
We denote the 
canonical section on $\mathbb{D} \rightarrow \mathbb{D}^c$ as 
$\psi_{\mathbb{D}}$.

Now we discuss an example of partial nonassociative dioperad 
which is important for us. Let $U = \oplus_{n\in J} U_{(n)}$
be a graded vector space and $J$ an index set.  
We denote the projection $U \rightarrow U_{(n)}$ as $P_n$. 
Now we consider a family of spaces of multi-linear maps 
$\mathbb{E}_{U} = \{ \mathbb{E}_{U}(m,n) \}_{m,n\in \N}$, 
where 
\beq
\mathbb{E}_{U}(m,n) :=
\hom_{\C}(U^{\otimes m}, \overline{U^{\otimes n}}) .
\eeq
For $f\in \mathbb{E}_{U}(m,n)$,  
$g_j \in \mathbb{E}_{U}(k_j,l_j)$
and $u_{p_j}^{(j)} \in U$, $1\leq p_j \leq l_j$, $j=1,\dots, n$, 
we say that 
\bea  \label{f-g-def}
&&\Gamma_{(i_1,\dots, i_n)}(f; g_1, \dots, g_n)
(u_1^{(1)}\otimes \dots \otimes u_{l_n}^{(n)})  \nn
&&\hspace{0.5cm} := \sum_{s_1, \dots, s_n \in J}  f\big( 
P_{s_1}g_1(u_1^{(1)}\otimes \dots \otimes
u_{l_1}^{(1)})\otimes \dots \otimes P_{s_n}
g_n(u_1^{(n)}\otimes \dots \otimes u_{l_n}^{(n)}) \big)\nonumber
\eea
is well-defined if the multiple sum converges absolutely.
This give arise to a partially defined substitution map,
for $1\leq i_j \leq k_j, j=1,\dots, n$,
$$
\Gamma_{i_1, \dots, i_n}: \mathbb{E}_{U}(m,n) \otimes 
\mathbb{E}_{U}(k_1,l_1) \otimes \dots \otimes 
\mathbb{E}_{U}(k_n, l_n) \rightarrow 
\mathbb{E}_{U}(m+k-n, l)
$$
where $k=k_1+\dots+k_n$ and $l=l_1+\dots+l_n$. 
In general, the compositions of three substitution maps
are not associative. 
The permutation groups actions on $\mathbb{E}_{U}$ 
are the usual one. Let $\Gamma=\{ \Gamma_{i_1, \dots, i_n} \}$. 
It is clear that $(\mathbb{E}_{U}, \Gamma, \id_U)$ 
is a partial nonassociative dioperad. We often denote it
simply by $\mathbb{E}_U$.

\begin{defn} {\rm 
Let $(\mathcal{P}, \gamma_{\mathcal{P}}, I_{\mathcal{P}})$ 
be a partial dioperad. 
A {\it $\mathcal{P}$-algebra} $(U, \nu)$ consists of
a graded vector space $U$
and a morphism of partial nonassociative dioperad 
$\nu: \mathcal{P} \rightarrow \mathbb{E}_{U}$.
}
\end{defn}

When $U = \oplus_{n\in J} U_{(n)}$ 
is a completely reducible module for a group $G$, $J$ is the 
set of equivalent classes of irreducible $G$-modules
and $U_{(n)}$ is a direct sum of irreducible $G$-modules of
equivalent class $n\in J$, we denote $\mathbb{E}_{U}$ by
$\mathbb{E}_U^G$.

\begin{defn} {\rm
Let $(\mathcal{P}, \gamma_{\mathcal{P}}, I_{\mathcal{P}})$ 
be a $G$-rescalable partial dioperad. 
A {\it $G$-rescalable $\mathcal{P}$-algebra} 
$(U, \nu)$ is a $\mathcal{P}$-algebra and the morphism
$\nu: \mathcal{P} \rightarrow \mathbb{E}_{U}^G$ is so that
$\nu|_G: G \rightarrow \edo \, U$ coincides with the given
$G$-module structure on $U$. 
}
\end{defn}

\subsection{Conformal full field algebras}

Let $(V^L, Y_{V^L}, \one^L, \omega^L)$ 
and $(V^R, Y_{V^R}, \one^R, \omega^R)$ 
be two vertex operator algebras with 
central charge $c^L$ and $c^R$ respectively, satisfying 
the conditions in Theorem \ref{ioa}. Let 
$(V_{cl}, m_{cl}, \iota_{cl})$ be
a conformal full field algebra over $V^L\otimes V^R$. 
A bilinear form $(\cdot, \cdot)_{cl}$ on 
$V_{cl}$ is invariant \cite{Ko1} if, for any $u, w_1, w_2\in V_{cl}$, 
\bea \label{inv-form-ffa-1}
&&\hspace{-1cm}(w_2, \mathbb{Y}_f(u; x, \bar{x})w_1)_{cl}  \nn
&&\hspace{-0.8cm}=(\mathbb{Y}_f(e^{-xL^L(1)}x^{-2L^L(0)}\otimes 
e^{-\bar{x}L^R(1)}\bar{x}^{-2L^R(0)} \, u; 
e^{\pi i} x^{-1}, e^{-\pi i} \bar{x}^{-1})w_2, w_1)_{cl}.
\eea 
or equivalently, 
\bea \label{inv-form-ffa-2}
&&\hspace{-1cm}(\mathbb{Y}_f(u; e^{\pi i}x, e^{-\pi i}\bar{x})w_2, 
w_1)_{cl}  \nn
&&\hspace{-0.5cm}=(w_2, \mathbb{Y}_f(e^{xL^L(1)}x^{-2L^L(0)}\otimes 
e^{\bar{x}L^R(1)}\bar{x}^{-2L^R(0)} \, u; x^{-1}, \bar{x}^{-1})w_1)_{cl}. 
\eea
We showed in \cite{Ko1} 
that an invariant bilinear form on $V_{cl}$ is 
automatically symmetric. Namely, for $u_1, u_2\in V_{cl}$,
we have 
\beq
(u_1, u_2)_{cl} = (u_2, u_1)_{cl}.
\eeq

$V_{cl}$ has a countable basis. 
We choose it to be $\{ e_i \}_{i\in \N}$. 
Assume that $(\cdot, \cdot)_{cl}$ is also nondegenerate, we also have 
the dual basis $\{ e^i \}_{i\in \N}$. Then we define a 
linear map $\Delta_{cl}: \C \rightarrow \overline{V_{cl}\otimes V_{cl}}$ 
as follow:
\beq
\Delta_{cl} :  1\mapsto \sum_{i\in \N} e_i\otimes e^i. 
\eeq

The correlation functions maps $m_{cl}^{(n)}, n\in \N$ 
of $V_{cl}$ are canonically determined by 
$\mathbb{Y}$ and the identity 
$\one_{cl}:= \iota_{cl}(\one^L \otimes \one^R)$ \cite{Ko2}.

For $Q\in \mathcal{T}_{\mathbb{K}}(n_-,n_+)$ 
given in (\ref{ele-Q}), we define, for $\lambda\in \C$, 
\beq  \label{def-psi-cl-0}
\nu_{cl} \big( \lambda \psi_{\mathbb{K}}(\pi_{\mathbb{K}}(Q)) \big)
( u_1 \otimes \dots \otimes u_{n_+} )
\eeq 
in the following three cases: 
\bnu
\item If $z_k\neq \infty$ for $k=-n_-, \dots, -1, 1, \dots, n_+$,  
(\ref{def-psi-cl-0}) is given by
\bea  \label{def-psi-cl}
&&\hspace{-0.5cm}\lambda \sum_{i_1, \dots, i_{n_-} \in \N} \big( 
\one_{cl}, \, m_{cl}^{(n_-+n_+)}(
e^{-L_+^L(A^{(-1)})-L_+^R(\overline{A^{(-1)}})} 
(a_0^{(-1)})^{-L^L(0)} \overline{a_0^{(-1)}}^{-L^R(0)}e_{i_1}, \nn
&&\hspace{2cm} \dots, 
e^{-L_+^L(A^{(-n_-)})-L_+^R(\overline{A^{(-n_-)}})}
(a_0^{(-n_-)})^{-L^L(0)} \overline{a_0^{(-n_-)}}^{-L^R(0)}e_{i_{n_-}}, \nn  
&&\hspace{3.7cm} 
e^{-L_+^L(A^{(1)})-L_+^R(\overline{A^{(1)}})} 
(a_0^{(1)})^{-L^L(0)} \overline{a_0^{(1)}}^{-L^R(0)}u_1,  \nn
&&\hspace{2cm} \dots,
e^{-L_+^L(A^{(n_+)})-L_+^R(\overline{A^{(n_+)}})}
(a_0^{(n_+)})^{-L^L(0)}\overline{a_0^{(n_+)}}^{-L^R(0)}u_{n_+}; \nn
&&\hspace{1cm}  z_{-1}, \bar{z}_{-1}, \dots, 
z_{-n_-}, \bar{z}_{-n_-},
z_1, \bar{z}_1, \dots, z_{n_+}, \bar{z}_{n_+}) \big)_{cl} \,  
e^{i_{1}}\otimes \dots \otimes e^{i_{n_-}}
\eea
where $L_+^L(A) = \sum_{j=1}^{\infty} A_j L_+^L(j)$
and $L_+^R(A) = \sum_{j=1}^{\infty} A_j L_+^R(j)$ for 
$A= (A_1, A_2, \dots )\in \prod_{j=1}^{\infty} \C$;

\item If $\exists \, k \in \{ -n_-, \dots, -1\}$ such that 
$z_k=\infty$ (recall (\ref{K-case-2})), 
(\ref{def-psi-cl-0}) is given by the formula obtained 
from (\ref{def-psi-cl}) by exchanging $\one_{cl}$ in 
(\ref{def-psi-cl}) with
\beq  \label{def-psi-cl-1}
e^{-L_+^L(A^{(k)})-L_+^R(\overline{A^{(k)}})}(a_0^{(k)})^{-L^L(0)} 
\overline{a_0^{(k)}}^{-L^R(0)} e_{i_{-k}};
\eeq

\item If $\exists \, k \in \{1, \dots, n_+\}$ such that 
$z_k=\infty$ (recall (\ref{K-case-2})), 
(\ref{def-psi-cl-0}) is given by the formula 
obtained from (\ref{def-psi-cl}) by exchanging $\one_{cl}$ 
in (\ref{def-psi-cl}) with
\beq  \label{def-psi-cl-2}
e^{-L_+^L(A^{(k)})-L_+^R(\overline{A^{(k)}})}(a_0^{(k)})^{-L^L(0)} 
\overline{a_0^{(k)}}^{-L^R(0)} u_k.
\eeq

\enu

The following result is proved in \cite{Ko1}.
\begin{prop}  \label{prop-welldef-cl}
The map $\nu_{cl}$ is $SL(2,\C)$-invariant, 
\end{prop}

Hence $\nu_{cl}$ induces a map
$\tilde{\mathbb{K}}^{c^L} \otimes 
\overline{\tilde{\mathbb{K}}^{\overline{c^R}}}
\rightarrow \mathbb{E}_{V_{cl}}^{\C^{\times}}$, which is 
still denoted as $\nu_{cl}$. 
Some interesting special cases are listed below: 
\bea
\nu_{cl}( \psi_{\mathbb{K}} ( \hat{\C} )) &=& 
(\one_{cl}, \one_{cl})_{cl} \id_{\C},  \nn
\nu_{cl}( \psi_{\mathbb{K}}(\pi_{\mathbb{K}}((\, (\infty, 1, \mathbf{0})|\, )_{\mathbb{K}}) ) 
&=& \one_{cl}, \nn
\nu_{cl} (\psi_{\mathbb{K}}( \pi_{\mathbb{K}}( (\,| 
(\infty,1,\mathbf{0}), (0,1,\mathbf{0}) \,)_{\mathbb{K}} )) )
(u\otimes v) 
&=& (u, v)_{cl},  \nn
\nu_{cl} (\psi_{\mathbb{K}}( \pi_{\mathbb{K}}((\, 
(\infty,1,\mathbf{0}), (0,1,\mathbf{0}) | \,)_{\mathbb{K}}) )) 
&=& \Delta_{cl},  \nn
\nu_{cl}(\psi_{\mathbb{K}} ( \pi_{\mathbb{K}}((\, (\infty; 1, \mathbf{0}) | 
(z; 1, \mathbf{0}), 
(0; 1, \mathbf{0})\, )_{\mathbb{K}}) ))(u\otimes v) &=& 
\mathbb{Y}(u; z, \bar{z})v,  \nonumber
\eea
where $\hat{\C}$ is the single element in $\mathbb{K}(0,0)$.

\begin{defn} {\rm
A $\tilde{\mathbb{K}}^{c^L}\otimes 
\overline{\tilde{\mathbb{K}}^{\overline{c^R}}}$-algebra $(U, \nu)$ 
is called {\it smooth} if 
\bnu
\item $U=\oplus_{m,n\in \R}U_{(m,n)}$ is 
a completely reducible $\C^{\times}$-module, where
$z\cdot u = z^m\bar{z}^n u, \forall z\in \C^{\times}, u\in U_{(m,n)}$. 
\item $\dim U_{(m,n)}<\infty, \forall m,n\in \R$
and $\dim U_{(m,n)}=0$ for $m$ or $n$ 
sufficiently small. 
\item $\nu$ is linear on fiber and smooth
on the base space $\mathbb{K}$.  
\enu
}
\end{defn}

The Theorem 2.7 in \cite{Ko1} can be restated as the following theorem.
\begin{thm}  \label{ffa-K-alg}
$(V_{cl},  \nu_{cl})$  is a smooth $\tilde{\mathbb{K}}^{c^L}\otimes 
\overline{\tilde{\mathbb{K}}^{\overline{c^R}}}$-algebra. 
\end{thm}

\subsection{Open-string vertex operator algebras}

Let $(V_{op}, Y_{op}, \one_{op}, \omega_{op})$ be an 
open-string vertex operator algebra. 
For $r>0$ and $v_1,v_2\in V_{op}$, we define $Y_{op}(v_1, -r) v_2$ by
\beq  \label{Y-op-opposite}
Y_{op}(v_1, -r)v_2 := e^{-rL(-1)}Y_{op}(v_2, r)v_1. 
\eeq

\begin{rema} {\rm
Taking the analogy between open-string vertex operator algebra
and associative algebra, $Y_{op}(\cdot, -r)\cdot$ 
corresponds to the opposite product \cite{HKo1}.  }
\end{rema}

An invariant bilinear form on an open-string vertex operator
algebra $V_{op}$ is a bilinear form $(\cdot, \cdot)_{op}$ 
on $V_{op}$ satisfying the following properties: 
\bea 
(v_3, Y_{op}(v_1, r)v_2)_{op} 
&=&(Y_{op}(e^{-rL(1)}r^{-2L(0)}\, v_1, -r^{-1})v_3, v_2)_{op}  
\label{inv-form-osvoa-1}  \\
(Y_{op}(v_1, r) v_3, v_2)_{op}
&=& (v_3, Y_{op}(e^{-rL(1)}r^{-2L(0)}v_1, -r^{-1})v_2)_{op}
\label{inv-form-osvoa-2}
\eea
for $r>0$ and $v_1, v_2, v_3\in V_{op}$. 

\begin{lemma}
\beq
(v_1, v_2)_{op} = (v_2, v_1)_{op}
\eeq
\end{lemma}
\pf
The proof is exactly same as that of Proposition 2.3 in \cite{Ko1}. 
\epf

We further assume that $(\cdot, \cdot)_{op}$ is nondegenerate. 
Let $\{ f_i \}_{i\in \R}$ be a basis of $V_{op}$ 
and $\{ f^i \}_{i\in \R}$ its dual basis. 
We define linear map $\Delta_{op}: \C \rightarrow 
\overline{V_{cl}\otimes V_{cl}}$ as follow: 
\beq
\Delta_{op}: 1 \mapsto \sum_{i\in \R} f_i \otimes f^i. 
\eeq

The open-string vertex operator algebra 
$(V_{op}, Y_{op}, \one_{op}, \omega_{op})$ naturally gives
a boundary field algebra $(V_{op}, m_{op}, \mathbf{d}_{op}, D_{op})$
in which the correlation-function maps
$m_{op}^{(n)}, n\in \N$ are completely determined by $Y_{op}$ and 
$\one_{op}$ \cite{Ko2}. 

For any $Q\in \mathcal{T}_{\mathbb{D}}(n_-, n_+)$ given in (\ref{ele-Q-D}).
Let $\alpha$ be a bijective map 
\beq   \label{alpha-map}
\{ -n_-, \dots, -1, 1, \dots, n_+ \} \xrightarrow{\alpha} 
\{ 1, \dots, n_-+n_+ \}
\eeq
so that $s_1, \dots, s_{n_-+n_+}$, defined by $s_i := r_{\alpha^{-1}(i)}$, 
satisfy $\infty\geq s_1>\dots >s_{n_-+n_+}\geq 0$.  
Then we define, for $\lambda \in \C$,    
\beq  \label{def-psi-op-0}
\nu_{op} ( \lambda \psi_{\mathbb{D}}( \pi_{\mathbb{D}}(Q)))
( v_1\otimes \dots \otimes v_{n_+} )
\eeq
as follow: 
\bnu
\item If $r_k\neq \infty, \forall k=-n_-, \dots, -1, 1, \dots, n_+$, 
(\ref{def-psi-op-0}) is given by 
\beq  \label{def-psi-op}
\lambda \sum_{i_1, \dots, i_{n_-} \in \R} \big( \one_{op}, \, 
m_{op}^{(n_-+n_+)} (w_1, \dots, w_{n_-+n_+}; s_1, \dots, s_{n_-+n_+}) \big)_{op}
f^{i_{1}}\otimes \dots \otimes f^{i_{n_-}}
\eeq
where $w_{\alpha(p)} = e^{-L_+(B^{(p)})} (b_0^{(p)})^{-L(0)} f_{i_{-p}}$
and $w_{\alpha(q)} = e^{-L_+(B^{(q)})} (b_0^{(q)})^{-L(0)} v_q$
for $p=-1, \dots -n_-$ and $q=1, \dots, n_+$;  

\item If $\exists \, k \in \{ -n_-, \dots, -1, 1, \dots, n_+\}$ 
such that $r_k=\infty$, 
(\ref{def-psi-op-0}) is given by the formula obtained 
from (\ref{def-psi-op}) by exchanging the $\one_{op}$ with $w_{\alpha(k)}$.
\enu

\begin{prop}
$\nu_{op}$ is $SL(2, \R)$ invariant. 
\end{prop}
\pf
The proof is same as that of Proposition \ref{prop-welldef-cl}. 
\epf

Hence $\nu_{op}$ induces a map 
$\tilde{\mathbb{D}}^c \rightarrow \mathbb{E}_{V_{op}}^{\R_+}$,
which is still denoted as $\nu_{op}$.  
Some interesting special cases are listed explicitly below: 
\bea
\nu_{op}( \psi_{\mathbb{D}} ( \hat{\HH} )) &=& 
(\one_{op}, \one_{op})_{op} \id_{\C},  \nn
\nu_{op}( \psi_{\mathbb{D}}( \pi_{\mathbb{D}}( 
(\, (\infty, 1, \mathbf{0})|\, )) )) 
&=& \one_{op}, \nn
\nu_{op} (\psi_{\mathbb{D}}( \pi_{\mathbb{D}}( (\,| 
(\infty,1,\mathbf{0}), (0,1,\mathbf{0}) \,) ) ) )(u\otimes v)
&=& (u, v)_{op},   \nn
\nu_{op} (\psi_{\mathbb{D}}( \pi_{\mathbb{D}}( (\, 
(\infty,1,\mathbf{0}), (0,1,\mathbf{0}) | \,)) )) &=& \Delta_{op},  \nn
\nu_{op}(\psi_{\mathbb{D}} ( \pi_{\mathbb{D}}( (\, (\infty; 1, \mathbf{0}) | 
(r; 1, \mathbf{0}), 
(0; 1, \mathbf{0})\, )) ))(u\otimes v) &=& Y_{op}(u, r)v.  \nonumber
\eea
where $\hat{\HH}$ is the single element in $\mathbb{D}(0,0)$
and $r>0$.

\begin{defn} {\rm
A $\tilde{\mathbb{D}}^{c}$-algebra $(U, \nu)$ 
is called {\it smooth} if 
\bnu
\item $U=\oplus_{n\in \R}U_{(n)}$ is 
a completely reducible $\R_+$-module, where
$r\cdot u = r^n u, \forall r\in \R_+, u\in U_{(m,n)}$. 
\item $\dim U_{(n)}<\infty, \forall n\in \R$
and $\dim U_{(n)}=0$ for $n$ sufficiently small. 
\item $\nu$ is linear on fiber and smooth
on the base space $\mathbb{D}$.  
\enu
}
\end{defn}

\begin{thm}
$(V_{op}, \nu_{op})$ is a smooth $\tilde{\mathbb{D}}^c$-algebra.  
\end{thm}
\pf
The proof is same as that of Theorem \ref{ffa-K-alg}.
\epf

\renewcommand{\theequation}{\thesection.\arabic{equation}}
\renewcommand{\thethm}{\thesection.\arabic{thm}}
\setcounter{equation}{0}
\setcounter{thm}{0}

\section{Swiss-cheese partial dioperad}

In Section 3.1, we introduce the notion of 
$2$-colored partial dioperad and algebra over it.
In Section 3.2, we study a special example of 
$2$-colored partial dioperad
called Swiss-cheese partial dioperad $\mathbb{S}$ 
and its $\C$-extension $\tilde{\mathbb{S}}^c$.
In Section 3.3, we show that an open-closed field algebra over $V$
equipped with nondegenerate invariant bilinear forms 
canonically gives an algebra over $\tilde{\mathbb{S}}^c$. 
In Section 3.4, we define boundary states in such algebra
and show that some of the boundary states are Ishibashi states.

\subsection{$2$-colored (partial) dioperads}

\begin{defn} {\rm
A {\it right module over a dioperad} $(\mathcal{Q}, \gamma_{\mathcal{Q}}, 
I_{\mathcal{Q}})$, or a right $\mathcal{Q}$-module, is a 
family of sets $\{ \mathcal{P}(m,n) \}_{m,n\in \N}$ with 
an $S_m\times S_n$-action on each set $\mathcal{P}(m,n)$ and 
substitution maps: 
\bea
 \mathcal{P}(m,n) \times 
\mathcal{Q}(k_1,l_1) \times
\dots \times \mathcal{Q}(k_n,l_n) 
&\xrightarrow{\gamma_{(i_1,\dots, i_n)}}&
\mathcal{P}(m+k_1\dots +k_n-n, l_1+\dots +l_n)  \nn
(P, Q_1, \dots, Q_n) &\mapsto & 
\gamma_{(i_1,\dots, i_n)}(P; Q_1, \dots, Q_n)
\eea
for $m,n,l_1,\dots,l_n\in \N, k_1,\dots,k_n\in \Z_+$
and $1\leq i_j \leq k_j, j=1,\dots, n$,
satisfying the right unit property, 
the associativity and the permutation axioms of 
dioperad but with
the right action of $\mathcal{P}$ on itself in the definition of
dioperad replaced by that of $\mathcal{Q}$. 
}
\end{defn}

Homomorphism and isomorphism between two right 
$\mathcal{Q}$-modules can be naturally defined. 
The right module over a partial dioperad can also be defined
in the usual way.

\begin{defn}  {\rm
Let $\mathcal{Q}$ be a $G$-rescalable partial dioperad. 
A right $\mathcal{Q}$-module is called $G$-rescalable if 
for any $P\in \mathcal{P}(m_-, m_+)$, 
$Q_i \in \mathcal{Q}(n_-^{(i)}, n_+^{(i)})$ and
$1\leq j_i \leq n_-^{(i)}, i=1,\dots,m_+$, 
there exist $g_j\in G, j=1,\dots, m_+$ such that 
\beq
\gamma_{(j_1,\dots, j_{m_+})} \big( P ;  \gamma_{(j_1)}(g_1; Q_1), \dots, 
\gamma_{(j_{m_+})}(g_{m_+}; Q_{m_+}) \big)
\eeq
is well-defined. 
}
\end{defn}

\begin{defn}  {\rm 
A {\it $2$-colored dioperad} consists of a dioperad 
$(\mathcal{Q}, \gamma_{\mathcal{Q}}, I_{\mathcal{Q}})$ and a family of
sets $\mathcal{P}(n_-^B, n_+^B|n_-^I, n_+^I)$ equipped with a 
$S_{n_-^B} \times S_{n_+^B} \times S_{n_-^I} \times S_{n_+^I}$-action
for $n_{\pm}^B, n_{\pm}^I\in \N$, a distinguished element 
$I_{\mathcal{P}}\in \mathcal{P}(1,1|0,0)$, and maps
\bea
&&\hspace{-0.5cm}
\mathcal{P}(m_-, m_+ | n_-, n_+) 
\times \mathcal{P}(k_-^{(1)}, k_+^{(1)} | l_-^{(1)}, l_+^{(1)}) \times
\dots \times 
\mathcal{P}(k_-^{(m_+)}, k_+^{(m_+)} | l_-^{(m_+)}, l_+^{(m_+)})  \nn
&&\hspace{2cm} \xrightarrow{\gamma_{(i_1,\dots, i_{m_+})}^B} 
\mathcal{P}(m_--m_++k_- , k_+ \, | n_-+l_- , n_++l_+), \nn
&&\hspace{-0.5cm} 
\mathcal{P}(m_-, m_+ | n_-, n_+)
\times \mathcal{Q}(p_-^{(1)}, p_+^{(1)})\times \dots \times 
\mathcal{Q}(p_-^{(n_+)}, p_+^{(n)}) \nn
&&\hspace{2cm} \xrightarrow{\gamma_{(j_1,\dots, j_{n_+})}^I}
\mathcal{P}(m_-, m_+ | n_--n_+ + p_-, p_+)
\eea 
where $k_{\pm}=k_{\pm}^{(1)}+\dots + k_{\pm}^{(m_+)}$, 
$l_{\pm}=l_{\pm}^{(1)}+\dots + l_{\pm}^{(m_+)}$ and 
$p_{\pm}=p_{\pm}^{(1)}+\dots + p_{\pm}^{(n_+)}$, 
for $1\leq i_r \leq k_-^{(r)}, r=1,\dots, m_+$, 
$1\leq j_s \leq p_-^{(s)}, s=1,\dots, n_+$, satisfying
the following axioms:
\bnu
\item The family of sets $\mathcal{P}:=\{ 
\mathcal{P}(m_-, m_+) \}_{m_-,m_+\in \N}$, where 
$$\mathcal{P}(m_-,m_+):= \cup_{n_-, n_+\in \N} 
\mathcal{P}(m_-, m_+|n_-, n_+),$$ equipped with the 
natural $S_{m_-}\times S_{m_+}$-action on $\mathcal{P}(m_-, m_+)$, 
 together with 
identity element $I_{\mathcal{P}}$ and the family of maps 
$\gamma^B:=\{ \gamma_{(i_1,\dots, i_{m_+})}^B \}$ gives an dioperad.

\item The family of maps 
$\gamma^I:=\{\gamma_{(j_1, \dots, j_{n_+})}^I \}$ gives each 
$\mathcal{P}(m_-, m_+)$ a right $\mathcal{Q}$-module structure
for $m_-,m_+\in \N$. 

\enu
}
\end{defn}

We denote such $2$-colored dioperad as 
$(\mathcal{P}|\mathcal{Q}, (\gamma^B, \gamma^I))$. 
{\it $2$-colored partial dioperad} can be naturally defined. 
If the associativities of $\gamma_{\mathcal{Q}}, \gamma^B, \gamma^I$
do not hold, then it is called {\it 
$2$-colored (partial) nonassociative dioperad}.

\begin{rema}  {\rm
If we restrict to $\{ \mathcal{P}(1,m|0,n) \}_{m,n\in \N}$ and
$\{ \mathcal{Q}(1,n) \}_{n\in \N}$, they simply gives a
structure of $2$-colored (partial) 
operad \cite{V}\cite{Kont}\cite{Ko2}. }
\end{rema}

Now we discuss an important example of $2$-colored partial 
nonassociative dioperad. Let $J_1, J_2$ be two index sets. 
$U_i=\oplus_{n\in J_i}(U_i)_{(n)}, i=1,2$ be two graded vector spaces. 
Consider two families of sets, 
\bea
\mathbb{E}_{U_2}(m,n) &=& 
\hom_{\C}(U_2^{\otimes m}, \overline{U_2^{\otimes n}}),    \nn
\mathbb{E}_{U_1|U_2}(m_-, m_+|n_-, n_+) &=& 
\hom_{\C}(U_1^{\otimes m_+}\otimes U_2^{\otimes n_+}, 
\overline{U_1^{\otimes m_-}\otimes U_2^{\otimes n_-}}), \nonumber
\eea
for $m,n, m_{\pm}, n_{\pm}\in \N$. 
We denote both of the projection operators 
$U_1 \rightarrow (U_1)_{(n)}$, $U_2\rightarrow (U_2)_{(n)}$ as $P_n$ 
for $n\in J_1$ or $J_2$.  
For $f\in \mathbb{E}_{U_1|U_2}(k_-,k_+|l_-,l_+)$,
$g_i \in \mathbb{E}_{U_1|U_2}(m_-^{(i)}, m_+^{(i)}|n_-^{(i)}, n_+^{(i)}), 
i=1,\dots, k_+$, 
and $h_j \in \mathbb{E}_{U_2}(p_-^{(j)}, p_+^{(j)})$, $j=1,\dots,l_+$, 
we say that 
\bea
&&\hspace{-1cm}\Gamma_{(i_1, \dots, i_{k_+})}^B
(f; g_1, \dots, g_{k_+}) (u_1^{(1)}\otimes \dots \otimes 
v_{n_+^{(k_+)}}^{(k_+)}\otimes v_1\otimes \dots \otimes v_{l_+})  \nn
&&\hspace{0cm}:= \sum_{s_1, \dots, s_k\in J_1} 
f \big( P_{s_1}g_1(u_{1}^{(1)} \otimes \dots \otimes 
u_{m_+^{(1)}}^{(1)}\otimes v_{1}^{(1)} \otimes \dots \otimes 
v_{n_+^{(1)}}^{(1)}) \otimes \nn
&&\hspace{0cm} \dots \otimes P_{s_{k_+}} 
g_{k_+}(u_{1}^{(k_+)}\otimes  \dots \otimes u_{m_+^{(k_+)}}^{(k_+)} \otimes
v_{1}^{(k_+)}\otimes \dots \otimes v_{n_+^{(k_+)}}^{(k_+)})\otimes 
v_1\otimes \dots \otimes v_{l_+} \big) \nn
&&\hspace{-1cm}
\Gamma_{(j_1, \dots, j_{l_+})}^I(f; h_1, \dots, h_{l_+})(u_1\otimes \dots
\otimes u_{k_+}\otimes w_1^{(1)}\otimes \dots \otimes w_{p_+^{(l_+)}}^{(l_+)}) 
\nn
&&\hspace{0cm}:=  \sum_{t_1, \dots, t_l\in J_2} 
f \big( u_1\otimes \dots \otimes u_{k_+}\otimes 
P_{t_1}h_1(w_1^{(1)}\otimes \dots \otimes w_{p_+^{(1)}}^{(1)})   \nn
&&\hspace{4cm} \otimes \dots \otimes 
P_{t_{l_+}}h_{l_+}(w_1^{(l_+)}\otimes \dots \otimes w_{p_+^{(l_+)}}^{(l_+)}) 
\big) \nonumber
\eea
for $u_j^{(i)}\in U_1, v_j^{(i)}, w_j^{(i)}\in U_2$, 
are well-defined if 
each multiple sum is absolutely convergent. These give arise 
to partially defined substitution maps: 
\bea
&&\hspace{-0.5cm} \mathbb{E}_{U_1|U_2}(k_-,k_+ | l_-,l_+) 
\otimes \mathbb{E}_{U_1|U_2}(m_-^{(1)},m_+^{(1)} |n_-^{(1)},n_+^{(1)}) 
\otimes \dots \otimes 
\mathbb{E}_{U_1|U_2}(m_-^{(k_+)},m_+^{(k_+)}|n_-^{(k_+)},n_+^{(k_+)})   \nn
&&\hspace{1.5cm} \xrightarrow{\Gamma_{(i_1,\dots, i_{k_+})}^B}
\mathbb{E}_{U_1|U_2}(k_--k_++m_-, m_+ |l_-+ n_-, n_+).    \nn
&&\mathbb{E}_{U_1|U_2}(k_-,k_+|l_-,l_+) 
\otimes \mathbb{E}_{U_2}(p_-^{(1)}, p_+^{(1)}) 
\otimes \dots \otimes \mathbb{E}_{U_2}(p_-^{(l_+)}, p_+^{(l_-)})   \nn
&&\hspace{1.5cm}\xrightarrow{\Gamma_{(j_1+\dots+j_{l_+})}^I}
\mathbb{E}_{U_1|U_2}(k_-,k_+| l_--l_++p_-^{(1)}+\dots +p_-^{(l_+)}, 
p_+^{(1)} +\dots+p_+^{(l_+)}). 
\nonumber
\eea
where $m_{\pm} = m_{\pm}^{(1)}+\dots+m_{\pm}^{(k_+)}$
and $n_{\pm}=n_{\pm}^{(1)}+\dots+n_{\pm}^{(k_+)}$.  
In general, $\Gamma_{(i_1,\dots, i_{k_+})}^B$ and $\Gamma_{(j_1+\dots+j_{l_+})}^I$ 
do not satisfy the associativity.
Let 
$$
\mathbb{E}_{U_1|U_2}=\{ \mathbb{E}_{U_1|U_2}(m_-,m_+|n_-,n_+) 
\}_{m_{\pm},n_{\pm}\in \N},
$$
$\mathbb{E}_{U_2} = \{ \mathbb{E}_{U_2}(n) \}_{n\in \N}$,
$\Gamma^B := \{ \Gamma_{(i_1,\dots, i_n)}^B \}$ and 
$\Gamma^I := \{ \Gamma_{(j_1, \dots, j_n)}^I \}$. 
It is obvious that 
$(\mathbb{E}_{U_1|U_2}|\mathbb{E}_{U_2}, (\Gamma^B, \Gamma^I))$
is a $2$-colored partial nonassociative operad.

Let $U_1$ be a completely reducible $G_1$-modules and 
$U_2$ a completely reducible $G_2$-modules. Namely, 
$U_1=\oplus_{n_1\in J_1} (U_1)_{(n_1)}, U_2= \oplus_{n_2\in J_2} (U_2)_{(n_2)}$
where $J_i$ is the set of equivalent classes of 
irreducible $G_i$-modules and $(U_i)_{(n_i)}$ 
is a direct sum of irreducible $G_i$-modules of
equivalent class $n_i$ for $i=1,2$. 
In this case, we denote $\mathbb{E}_{U_1|U_2}$ 
by $\mathbb{E}_{U_1|U_2}^{G_1|G_2}$. 

\begin{defn} {\rm
A homomorphism between two $2$-colored (partial) dioperads 
$$(\mathcal{P}_i|\mathcal{Q}_i, (\gamma_i^B, \gamma_i^I)), 
i=1,2$$ consists of two (partial) dioperad homomorphisms: 
\beq
\nu_{\mathcal{P}_1|\mathcal{Q}_1}: \mathcal{P}_1 \rightarrow \mathcal{P}_2,  
\quad\quad \mbox{and} \quad\quad
\nu_{\mathcal{Q}_1}: \mathcal{Q}_1 \rightarrow \mathcal{Q}_2  \nonumber
\eeq 
such that $\nu_{\mathcal{P}_1|\mathcal{Q}_1}: \mathcal{P}_1 \rightarrow 
\mathcal{P}_2 $, where $\mathcal{P}_2$ has a
right $\mathcal{Q}_1$-module structure induced by dioperad
homomorphism $\nu_{\mathcal{Q}_1}$, 
is also a right $\mathcal{Q}_1$-module homomorphism.
}
\end{defn}

\begin{defn} {\rm
An {\it algebra over a $2$-colored partial dioperad 
$(\mathcal{P}| \mathcal{Q}, (\gamma^B, \gamma^I))$}, or a 
{\it $\mathcal{P}|\mathcal{Q}$-algebra} 
consists of two graded vector spaces $U_1, U_2$ and 
a $2$-colored partial dioperad homomorphism 
$(\nu_{\mathcal{P}|\mathcal{Q}}, \nu_{\mathcal{Q}}):
(\mathcal{P}|\mathcal{Q}, (\gamma^B, \gamma^I))\rightarrow 
(\mathbb{E}_{U_1|U_2}| \mathbb{E}_{U_2}, (\Gamma^B, \Gamma^I))$. 
We denote this algebra as 
$(U_1|U_2, \nu_{\mathcal{P}|\mathcal{Q}}, \nu_{\mathcal{Q}})$. 
}
\end{defn}

\begin{defn}  {\rm
If a $2$-colored partial dioperad $(\mathcal{P}| \mathcal{Q}, \gamma)$
is so that $\mathcal{P}$ 
is a $G_1$-rescalable partial operad
and a $G_2$-rescalable right $\mathcal{Q}$-module, then 
it is called {\it $G_1|G_2$-rescalable}. 
}
\end{defn}

\begin{defn} {\rm
A {\it $G_1|G_2$-rescalable $\mathcal{P}|\mathcal{Q}$-algebra}
$(U_1|U_2,\nu_{\mathcal{P}|\mathcal{Q}}, \nu_{\mathcal{Q}})$ 
is a $\mathcal{P}|\mathcal{Q}$-algebra so that 
$\nu_{\mathcal{P}|\mathcal{Q}}: \mathcal{P}\rightarrow 
\mathbb{E}_{U_1|U_2}^{G_1|G_2}$
and $\nu_{\mathcal{Q}}: \mathcal{Q} \rightarrow \mathbb{E}_{U_2}^{G_2}$
are dioperad homomorphisms such that 
$\nu_{\mathcal{P}|\mathcal{Q}}: G_1\rightarrow \edo U_1$
coincides with the $G_1$-module structure on $U_1$ 
and $\nu_{\mathcal{Q}}: G_2\rightarrow \edo U_2$ coincides with
the $G_2$-module structure on $U_2$. 
}
\end{defn}

\subsection{Swiss-cheese partial dioperads}

A disk with strips and tubes of type $(m_-, m_+ ; n_-, n_+)$ 
\cite{Ko2} is a disk consisting of $m_+$ ($m_-$) 
ordered positively (negatively) oriented punctures 
on the boundary of the disk, 
and $n_+$ ($n_-$) ordered positively (negatively) oriented punctures  
in the interior of the disk,
and local coordinate map around each puncture. 
Two disks are conformal equivalent if there exists a 
biholomorphic map between them 
preserving order, orientation and local coordinates. 
We denote the moduli space of disks with strips and tubes of 
type $(m_-, m_+; n_-, n_+)$ as $\mathbb{S}(m_-, m_+ | n_-, n_+)$. 

We use the following notation 
\bea
&&\hspace{-1cm} \big[ \, ( \, (r_{-m_-}, b_0^{-m_-}, B^{(-m_-)}), \dots, 
(r_{-1}, b_0^{-1}, B^{(-1)}) \,\, | \nn
&&\hspace{1cm} 
(r_{1}, b_0^{1}, B^{(1)}), \dots, (r_{m_+}, b_0^{m_+}, B^{(m_+)}) \, ) 
\,\,\,  \big\|  \nn
&&\hspace{2cm}
( \, (z_{-n_-}, a_0^{-n_-}, A^{(-n_-)}), \dots, 
(z_{-1}, a_0^{-1}, A^{(-1)}) \,\, | \nn
&&\hspace{5cm} (z_{1}, a_0^{1}, A^{(1)}), \dots, 
(z_{n_+}, a_0^{n_+}, A^{(n_+)}) \, ) \big]_{\mathbb{S}}   \label{gen-ele-S}
\eea
where $r_i \in \hat{\R}$, $b_0^{(i)}\in \R^{\times}$, 
$B_k^{(i)}\in \R$ 
and $z_j\in \HH$, $a_0^{(j)}\in \C^{\times}$, 
$A_l^{(j)}\in \C$ for all
$i=-m_-, \dots, -1, 1, \dots, m_+$, 
$j=-n_-, \dots, -1,1, \dots, n_+$ and $k,l\in \Z_+$, 
to represent a disk with strips at $r_i$ with local coordinate map
$f_i$ and tubes at $z_j$ with local coordinate map $g_j$ given as follow:
\bea
f_{i}(w) &=& e^{\sum_k B_k^{(i)} x^{k+1}\frac{d}{dx}} (b_0^{(i)})^{x\frac{d}{dx}}x 
|_{x=w-r_i} \quad\, \mbox{if $r_i\in \R$ }, \label{case-1-S} \\
&=& e^{\sum_k B_k^{(i)} x^{k+1}\frac{d}{dx}} (b_0^{(i)})^{x\frac{d}{dx}}x 
|_{x=\frac{-1}{w}}  \quad\quad \mbox{if $r_i=\infty$}, \label{case-2-S}\\
g_j(w) &=& 
e^{\sum_k A_k^{(j)} x^{k+1}\frac{d}{dx}} (a_0^{(j)})^{x\frac{d}{dx}}x |_{x=w-z_j}. 
\label{case-3-S}
\eea

We denote the set of all such disks given in (\ref{gen-ele-S}) 
as $\mathcal{T}_{\mathbb{S}}(m_-, m_+ | n_-, n_+)$. 
The automorphisms  of the upper half plane, which is $SL(2, \R)$, 
change the disk (\ref{gen-ele-S}) to a different but 
conformal equivalent disk. It is clear that we have 
\beq
\mathbb{S}(m_-, m_+ | n_-, n_+) = 
\mathcal{T}_{\mathbb{S}}(m_-, m_+ | n_-, n_+)/SL(2, \R). 
\eeq

Let $\mathbb{S}= 
\{ \mathbb{S}(m_-, m_+ | n_-, n_+)\}_{m_-, m_+, n_-, n_+\in \N }$. 
The permutation groups
$$
(S_{m_-}\times S_{m_+}) \times (S_{n_-}\times S_{n_+})
$$
acts naturally on $\mathbb{S}(m_-, m_+ | n_-, n_+)$.
There are so-called {\em boundary sewing operations}
\cite{HKo1}\cite{Ko2}
on $\mathbb{S}$, denoted as $_{^i}\infty_{^{-j}}^{B}$, which
sews the $i$-th positively oriented boundary puncture of the
first disk with the $j$-th negatively
oriented boundary puncture of the second disk.
Boundary sewing operations naturally induce the following maps:
\bea
&&\hspace{-0.5cm}
\mathbb{S}(m_-, m_+ | n_-, n_+) 
\times \mathbb{S}(k_-^{(1)}, k_+^{(1)} | l_-^{(1)}, l_+^{(1)}) \times
\dots \times 
\mathbb{S}(k_-^{(m_+)}, k_+^{(m_+)} | l_-^{(m_+)}, l_+^{(m_+)})  \nn
&&\hspace{2cm} \xrightarrow{\gamma_{(i_1,\dots, i_{m_+})}^B} 
\mathbb{S}(m_--m_++k_- , k_+ \, | n_-+l_- , n_++l_+), \nonumber
\eea 
where $k_{\pm}=k_{\pm}^{(1)}+\dots + k_{\pm}^{(m_+)}$, 
$l_{\pm}=l_{\pm}^{(1)}+\dots + l_{\pm}^{(m_+)}$ for 
$1\leq i_r \leq k_-^{(r)}, r=1,\dots, m_+$. 
It is easy to see that 
boundary sewing operations or $\gamma_{(i_1,\dots, i_{m_+})}^B$, 
together with permutation group
actions on the order of boundary punctures, provide
$\mathbb{S}$ with a structure of partial dioperad.

There are also so-called {\it interior sewing operations} 
\cite{HKo1}\cite{Ko2} on $\mathbb{S}$,
denoted as $_{^i}\infty_{^{-j}}^{I}$, which 
sew the $i$-th positively oriented interior puncture of a 
disk with the $j$-th negatively oriented puncture of a sphere. 
The interior sewing operations 
define a right action of $\mathbb{K}$ on $\mathbb{S}$: 
\bea
&&\hspace{-0.5cm} 
\mathbb{S}(m_-, m_+ | n_-, n_+)
\times \mathbb{K}(p_-^{(1)}, p_+^{(1)})\times \dots \times 
\mathbb{K}(p_-^{(n_+)}, p_+^{(n_+)}) \nn
&&\hspace{2cm} \xrightarrow{\gamma_{(j_1,\dots, j_{n_+})}^I}
\mathbb{S}(m_-, m_+ | n_--n_+ + p_-, p_+)
\eea
where $p_{\pm}=p_{\pm}^{(1)}+\dots + p_{\pm}^{(n_+)}$, for 
$1\leq j_s \leq p_-^{(s)}, s=1,\dots, n_+$.
Such action gives $\mathbb{S}$ a right $\mathbb{K}$-module structure. 

Let $\gamma^B=\{ \gamma_{(i_1,\dots, i_n)}^B\}$ and 
$\gamma^I=\{ \gamma_{(i_1,\dots, i_n)}^I\}$. 
The following proposition is clear. 
\begin{prop}
$(\mathbb{S}|\mathbb{K}, (\gamma^B, \gamma^I))$ is a 
$\R_+|\C^{\times}$-rescalable $2$-colored partial dioperad. 
\end{prop}

We call $(\mathbb{S}|\mathbb{K}, (\gamma^B, \gamma^I))$
{\it Swiss-cheese partial dioperad}.
When it is restricted on 
$\mathfrak{S}=\{ \mathbb{S}(1,m|0,n) \}_{m,n\in \N}$, it is nothing
but the so-called Swiss-cheese partial operad \cite{HKo1}\cite{Ko2}.

In \cite{HKo1}\cite{Ko2}, we show that the Swiss-cheese partial 
operad $\mathfrak{S}$ can be naturally embedded into the sphere
partial operad $K$ via the so-called doubling map, denoted as 
$\delta: \mathfrak{S} \hookrightarrow K$. Such doubling map 
$\delta$ obviously can be extended to a doubling map 
$\mathbb{S} \hookrightarrow \mathbb{K}$, still
denoted as $\delta$. 
In particular, the general element (\ref{gen-ele-S}) maps under
$\delta$ to
\bea
&&\hspace{-1cm}(\, (z_{-1}, a_0^{(-1)}, A^{(-1)}), \dots, 
(z_{-n_-}, a_0^{(-n_-)}, A^{(-n_-)}), \nn
&&\hspace{1cm}
(\bar{z}_{-1}, \overline{a_0^{(-1)}}, \overline{A^{(-1)}}), \dots 
(\bar{z}_{-n_-}, \overline{a_0^{(-n_-)}}, \overline{A^{(-n_-)}}),\nn
&&\hspace{3cm}(r_{-1}, b_0^{(-1)}, B^{(-1)}), \dots, 
(r_{-m_-}, b_0^{(-m_-)}, B^{(-m_-)}) \, | \nn
&&\hspace{-0.5cm} (z_{1}, a_0^{(1)}, A^{(1)}), 
\dots,  (z_{n_+}, a_0^{(n_+)}, A^{(n_+)}),  \nn
&&\hspace{1cm}  (\bar{z}_{1}, \overline{a_0^{(1)}}, \overline{A^{(1)}}),  
\dots, (\bar{z}_{n_+}, \overline{a_0^{(n_+)}}, \overline{A^{(n_+)}}), \nn
&&\hspace{3cm}
(r_{1}, b_0^{(1)}, B^{(1)}), \dots, (r_{m_+}, b_0^{(m_+)}, B^{(m_+)})\, 
)_{\mathbb{K}}.
\eea

The following proposition is clear. 
\begin{prop}
Let $P_i\in \mathbb{S}(m_-^{(i)}, m_+^{(i)}| n_-^{(i)}, n_+^{(i)}), i=1,2$
and $Q\in \mathbb{K}(m_-,m_+)$. If $P_1 {}_{^i}\infty_{^{-j}}^B P_2$ 
and $P_1 {}_{^k}\infty_{^{-l}}^I Q$ exists
for $1\leq i \leq m_+^{(1)}, 1\leq j\leq m_-^{(2)}$ and
$1\leq k \leq n_+^{(1)}, 1\leq l\leq m_-$,  then we have 
\bea
\delta(P_1 \,\, {}_{^i}\infty_{^{-j}}^B \, P_2) &=& 
\delta(P_1) \,\, {}_{^{2n_+^{(1)} + i}}\infty_{^{-(2n_-^{(2)}+j)}} 
\, \delta(P_2), \nn
\delta( P_1 \,\, {}_{^k}\infty_{^{-l}}^I \, Q )
&=& (\delta(P_1) \,\, {}_{^{k}}\infty_{^{-l}} \, Q \,)
\,\, {}_{^{n_+^{(1)}+m_+-1+k}}\infty_{^{-l}} \, \bar{Q}.   \label{double-sew}
\eea

\end{prop}

By above Proposition, we can identify $\mathbb{S}$ as 
its image under $\delta$ in $\mathbb{K}$ with boundary 
sewing operations replaced by ordinary sewing operations in
$\mathbb{K}$ and interior sewing operations replaced by 
double-sewing operations in $\mathbb{K}$ 
as given in (\ref{double-sew}).

The $\C$-extension $\tilde{\mathbb{S}}^c(m_-, m_+ | n_-, n_+)$ 
of $\mathbb{S}(m_-, m_+ | n_-, n_+)$ is defined to be the pullback 
bundle of $\tilde{\mathbb{K}}^c(2n_-+m_-, 2n_++m_+)$.  
We denote the canonical section on 
$\tilde{\mathbb{S}}^c$, which is induced from that on 
$\tilde{\mathbb{K}}^c$, as $\psi_{\mathbb{S}}$. 
The boundary (interior) 
sewing operations can be 
naturally extended to $\tilde{\mathbb{S}}^c$. 
We denote them as $\widetilde{\infty}^B$ ($\widetilde{\infty}^I$)
and corresponding substitution maps 
as $\tilde{\gamma}^B$ ($\tilde{\gamma}^I$). There is a natural
right action of $\tilde{\mathbb{K}}^c\otimes 
\overline{\tilde{\mathbb{K}}^{\bar{c}}}$ on $\tilde{\mathbb{S}}^c$
defined by $\tilde{\gamma}^I$.  

The following proposition is also clear. 
\begin{prop}
$(\tilde{\mathbb{S}}^c|\tilde{\mathbb{K}}^c\otimes \overline{\tilde{\mathbb{K}}^{\bar{c}}}, (\tilde{\gamma}^B, \tilde{\gamma}^I))$
is a $\R_+|\C^{\times}$-rescalable $2$-colored partial dioperad. 
\end{prop}

We will call the structure 
$(\tilde{\mathbb{S}}^c|\tilde{\mathbb{K}}^c\otimes \overline{\tilde{\mathbb{K}}^{\bar{c}}}, (\tilde{\gamma}^B, \tilde{\gamma}^I))$
as {\it Swiss-cheese
partial dioperad with central charge $c$}.

\begin{defn} {\rm
An algebra over $\tilde{\mathbb{S}}^{c}$ viewed as 
a dioperad, $(U, \nu)$, is called {\it smooth} if 
\bnu
\item $U=\oplus_{n\in \R}U_{(n)}$ is 
a completely reducible $\R_+$-module, where
$r\cdot u = r^n u, \forall r\in \R_+, u\in U_{(m,n)}$. 
\item $\dim U_{(n)}<\infty, \forall n\in \R$
and $\dim U_{(n)}=0$ for $n$ sufficiently small. 
\item $\nu$ is linear on fiber and smooth
on the base space $\mathbb{S}$.  
\enu
}
\end{defn}

\begin{defn} {\rm
An $\tilde{\mathbb{S}}^c|\tilde{\mathbb{K}}^c\otimes \overline{\tilde{\mathbb{K}}^{\bar{c}}}$-algebra 
$(U_1|U_2, \nu_{\tilde{\mathbb{S}}^c|\tilde{\mathbb{K}}^c\otimes \overline{\tilde{\mathbb{K}}^{\bar{c}}}}, 
\nu_{\tilde{\mathbb{K}}^c\otimes \overline{\tilde{\mathbb{K}}^{\bar{c}}}})$ 
is {\it smooth} if both $(U_1, \nu_{\tilde{\mathbb{S}}^c|\tilde{\mathbb{K}}^c\otimes \overline{\tilde{\mathbb{K}}^{\bar{c}}}})$
and $(U_2, \nu_{\tilde{\mathbb{K}}^c\otimes \overline{\tilde{\mathbb{K}}^{\bar{c}}}})$ as algebras over dioperads are smooth. 

}
\end{defn}

\subsection{Open-closed field algebras over $V$}

Let $(V_{op}, Y_{op}, \iota_{op})$ be an open-string
vertex operator algebra over $V$ and 
$(V_{cl}, \mathbb{Y}, \iota_{cl})$ a conformal full field
algebra over $V\otimes V$. Let
\beq
(\, (V_{cl}, \mathbb{Y}, \iota_{cl}), 
(V_{op}, Y_{op}, \iota_{op}), \mathbb{Y}_{cl-op}\, )  \label{opcl-fa}
\eeq
be an open-closed field algebra over $V$ \cite{Ko2}. 
We denote the formal vertex operators
associated with $Y_{op}$ and $\mathbb{Y}$
as $Y_{op}^f$ and $\mathbb{Y}^f$ respectively.
Let $\omega^L = \iota_{cl}(\omega \otimes \one)$, 
$\omega^R = \iota_{cl}(\one \otimes \omega)$ and 
$\omega_{op} = \iota_{op}(\omega)$. We have 
\bea
\mathbb{Y}^f( \omega^L ; x, \bar{x}) &=& \mathbb{Y}^f(\omega^L, x) =
\sum_{n\in \Z} L(n) \otimes 1 x^{-n-2},  \nn
\mathbb{Y}^f( \omega^R ; x, \bar{x}) &=& \mathbb{Y}^f(\omega^R, \bar{x})
= \sum_{n\in \Z} 1 \otimes L(n) \bar{x}^{-n-2},   \nn
Y_{op}^f(\omega_{op}, x) &=& \sum_{n\in \Z} L(n) x^{-n-2}. 
\eea
We also set $L^L(n)=L(n)\otimes 1$ and $L^R(n)=1\otimes L(n)$
for $n\in \Z$. 

For $u\in V$, we showed in \cite{Ko2} that 
$\mathbb{Y}_{cl-op}(\iota_{cl}(u\otimes \one); z, \bar{z})$
and $\mathbb{Y}_{cl-op}(\iota_{cl}(\one\otimes u); z, \bar{z})$
are holomorphic and antiholomorphic respectively. So we also denote
them simply by $\mathbb{Y}_{cl-op}(\iota_{cl}(u\otimes \one), z)$
and $\mathbb{Y}_{cl-op}(\iota_{cl}(\one\otimes u), \bar{z})$ respectively. 

By the $V$-invariant boundary condition \cite{Ko2}, we have 
\beq
\mathbb{Y}_{cl-op}(\omega^L, r) = Y_{op}(\omega_{op}, r) =
\mathbb{Y}_{cl-op}(\omega^R, r). 
\eeq

We assume that both $V_{op}$ and 
$V_{cl}$ are equipped with nondegenerate invariant bilinear forms 
$(\cdot, \cdot)_{op}$ and $(\cdot, \cdot)_{cl}$ respectively.

\begin{lemma}  \label{lemma-1-st}
For any $u\in V_{cl}$ and $v_1, v_2\in V_{op}$ and $z\in \HH$,
we have
\bea  \label{inv-form-bcft-1}
&&\hspace{-1cm}(v_2, \mathbb{Y}_{cl-op}(u; z, \bar{z})v_1)_{op}  \nn
&& \hspace{-0.5cm}=(\mathbb{Y}_{cl-op}(e^{-zL(1)}z^{-2L(0)}\otimes 
e^{-\bar{z}L(1)}\bar{z}^{-2L(0)} u; -z^{-1}, 
-\bar{z}^{-1})v_2,  v_1)_{op}  
\eea
\end{lemma}
\pf
Using (\ref{inv-form-osvoa-1}), for fixed $z\in \HH$, we have, 
\bea
&&(v_2, \mathbb{Y}_{cl-op}(u; z, \bar{z})v_1)_{op}   \nn
&&\hspace{1cm}= (v_2, 
\mathbb{Y}_{cl-op}(u; z, \bar{z})Y_{op}(\one, r)v_1)_{op}  \nn
&&\hspace{1cm}=(v_2, 
Y_{op}(\mathbb{Y}_{cl-op}(u; z-r, \bar{z}-r)\one, r)v_1)_{op}  \nn
&&\hspace{1cm}=(e^{-r^{-1}L(-1)}Y_{op}(v_2, r^{-1})
e^{-rL(1)}r^{-2L(0)}\mathbb{Y}_{cl-op}(u;z-r,\bar{z}-r)\one, v_1)_{op} 
\eea
for $|z|>r>|z-r|>0$.
Notice that $e^{-r^{-1}L(-1)}\in \text{Aut}(\overline{V}_{op})$,
$\forall r\in \C$. By taking $v_2=\one_{op}$, it is easy to see that 
\beq  \label{inv-form-bcft-equ-1}
e^{-rL(1)}r^{-2L(0)}\mathbb{Y}_{cl-op}(u;z-r,\bar{z}-r)\one
\eeq
is a well-defined element in $\overline{V}_{op}$ for $|z|>r>|z-r|>0$.  
Because of the chirality splitting property of $\mathbb{Y}_{cl-op}$
(see (1.72)(1.73) in \cite{Ko2}), it is easy to show that 
(\ref{inv-form-bcft-equ-1}) equals to
\beq
\mathbb{Y}_{cl-op}(e^{-zL(1)}z^{-2L(0)}\otimes e^{-\bar{z}L(1)}\bar{z}^{-2L(0)}u, 
r^{-1}-z^{-1}, r^{-1}-\bar{z}^{-1})\one
\eeq
for $r>|z-r|>0$.
By the commutativity I of analytic open-closed field algebra 
proved in \cite{Ko2}, we know that for fixed $z\in \HH$, 
$$
e^{-r^{-1}L(-1)}Y_{op}(v_2, r^{-1})\mathbb{Y}_{cl-op}(e^{-zL(1)}z^{-2L(0)}
\otimes e^{-\bar{z}L(1)}\bar{z}^{-2L(0)}u, 
r^{-1}-z^{-1}, r^{-1}-\bar{z}^{-1})\one
$$
and 
\bea  \label{inv-form-bcft-equ-2}
&&e^{-r^{-1}L(-1)}\mathbb{Y}_{cl-op}(e^{-zL(1)}z^{-2L(0)}
\otimes e^{-\bar{z}L(1)}\bar{z}^{-2L(0)}u, 
r^{-1}-z^{-1}, r^{-1}-\bar{z}^{-1})Y_{op}(v_2, r^{-1})\one  \nn
&&=e^{-r^{-1}L(-1)}\mathbb{Y}_{cl-op}(e^{-zL(1)}z^{-2L(0)}
\otimes e^{-\bar{z}L(1)}\bar{z}^{-2L(0)}u, 
r^{-1}-z^{-1}, r^{-1}-\bar{z}^{-1})e^{r^{-1}L(-1)}v_2  \nn
\eea
converge in different domains for $r$, but 
are analytic continuation of each other along a path in $r\in \R_+$. 
Moreover, using $L(-1)$ property of intertwining operator
and chirality splitting property of $\mathbb{Y}_{cl-op}$ again, 
the right hand side of (\ref{inv-form-bcft-equ-2}) equals to
$$
\mathbb{Y}_{cl-op}(e^{-zL(1)}z^{-2L(0)}
\otimes e^{-\bar{z}L(1)}\bar{z}^{-2L(0)}u, -z^{-1}, -\bar{z}^{-1})v_2
$$
for $|r^{-1}-z^{-1}|>|r^{-1}|$. Therefore, the both sides of
(\ref{inv-form-bcft-1}) as constant functions of $r$ 
are analytic continuation of each other. 
Hence (\ref{inv-form-bcft-1}) must hold
identically for all $z\in \HH$.  
\epf

For $u_1, \dots, u_l\in V_{cl}, v_1, \dots, v_n\in V_{cl}$, 
$r_1,\dots, r_n\in \R, r_1>\dots >r_n$
and $z_1, \dots, z_l\in \HH$, we define
\bea
&&m_{cl-op}^{(l;n)}(u_1, \dots, u_l; v_1, \dots, v_n; 
z_1, \bar{z}_1, \dots, z_l, \bar{z}_l; r_1, \dots, r_n)  \nn
&&\hspace{1cm}
:= e^{-r_nL(-1)} m_{cl-op}^{(l;n)}(u_1, \dots, u_l; v_1, \dots, v_n; 
z_1-r_n, \bar{z}_1-r_n, \dots, z_l-r_n, \bar{z}_l-r_n; \nn
&&\hspace{5cm} r_1-r_n, 
\dots, r_{n-1}-r_n, 0).
\eea
We simply extend the definition of $m_{cl-op}^{(l;n)}$ to 
a domain where some of $r_i$ can be negative. 
Note that such definition is compatible with $L(-1)$-properties of 
$m_{cl-op}$. 

\begin{lemma}  \label{lemma-2-st}
For $u_1, \dots, u_l\in V_{cl}, v, v_1, \dots, v_n\in V_{cl}$, 
$r_1,\dots, r_n\in \R, r_1>\dots >r_n=0$ and $z_1, \dots, z_l\in \HH$,
\bea  \label{lemma-2}
&&(v, \, m_{cl-op}^{(l;n)}(u_1, \dots, u_l; v_1, \dots, v_n; 
z_1, \bar{z}_1, \dots, z_l, \bar{z}_l; r_1, \dots, r_n))_{op}\nn
&&\hspace{0.2cm}=(v_n,\, m_{cl-op}^{(l;n)}( F_1u_1, \dots, F_l u_l;
v, G_1v_1, \dots, G_{n-1}v_{n-1}; \nn
&&\hspace{2.5cm} 
-z_1^{-1}, -\bar{z}_1^{-1}, \dots, -z_l^{-1}, -\bar{z}_l^{-1}; 0,
-r_1^{-1}, \dots, -r_{n-1}^{-1}))_{op}
\eea
where
\bea
F_i &=& e^{-z_i L(1)}z_i^{-2L(0)}\otimes 
e^{-\bar{z}_i L(1)}\bar{z}_i^{-2L(0)}, \hspace{1cm} i=1,\dots, l,    \nn
G_j &=&  e^{-r_j L(1)}r_j^{-2L(0)}, \hspace{4cm} j=1, \dots, n-1.
\eea
\end{lemma}
\pf
By Lemma \ref{lemma-1-st}, (\ref{lemma-2}) 
is clearly true for $l=0,1;n=0,1$. By (\ref{inv-form-osvoa-1})
and (\ref{inv-form-osvoa-2}), (\ref{lemma-2}) is true for $l=0,n=2$. 
We then prove the Lemma by induction. 
Assume that (\ref{lemma-2}) is true for $l=k\geq 0, n=m\geq 2$
or $l=k\geq 1, n=m\geq 1$.  

Let $l=k$ and $n=m+1$. 
It is harmless to assume that 
$0< r_{n-1}, |z_i| < r_{n-2}, i= l_1+1, \dots, l$ for some
$l_1\leq l$. Using the induction hypothesis, we obtain
\bea  \label{lemma-2-equ-1}
&&(v, \, m_{cl-op}^{(l;n)}(u_1, \dots, u_l; v_1, \dots, v_n; 
z_1, \bar{z}_1, \dots, z_l, \bar{z}_l; r_1, \dots, r_n))_{op}\nn
&&\hspace{0.2cm}=(v, m_{cl-op}^{(l_1; n-1)}
(u_1, \dots, u_{l_1}; v_1, \dots, v_{n-2}, m_{cl-op}^{(l-l_1; 2)}
(u_{l_1+1}, \dots, u_l;  v_{n-1}, v_n;  \nn
&&\hspace{2cm}   z_{l_1+1}, \bar{z}_{l_1+1},
\dots, z_l, \bar{z}_l; r_{n-1},0);
z_1, \bar{z}_1, \dots, z_{l_1}, \bar{z}_{l_1}; 
r_1, \dots, r_{n-2},0))_{op}\nn
&&\hspace{0.2cm}=(m_{cl-op}^{(l_1;n-1)}( F_1u_1, \dots, F_{l_1} u_{l_1};
v, G_1v_1, \dots, G_{n-1}v_{n-1}; \nn
&&\hspace{2.5cm} 
-z_1^{-1}, -\bar{z}_1^{-1}, \dots, -z_{l_1}^{-1}, -\bar{z}_{l_1}^{-1}; 0,
-r_1^{-1}, \dots, -r_{n-1}^{-1}), \nn 
&&\hspace{3cm}m_{cl-op}^{(l-l_1; 2)}
(u_{l_1+1}, \dots, u_l; v_{n-1}, v_n; z_{l_1+1}, \bar{z}_{l_1+1}, 
\dots, z_l, \bar{z}_l; r_{n-1},0) )_{op}  \nn
&&\hspace{0.2cm}=
(m_{cl-op}^{(l-l_1; 2)}(F_{l_1+1}u_{l_1+1}, \dots, F_lu_l; 
m_{cl-op}^{(l_1;n-1)}( F_1u_1, \dots, F_{l_1} u_{l_1}; \nn
&&\hspace{2cm} v, G_1v_1, \dots, G_{n-1}v_{n-1}; 
-z_1^{-1}, -\bar{z}_1^{-1}, \dots, -z_{l_1}^{-1}, -\bar{z}_{l_1}^{-1}; 0,
-r_1^{-1}, \dots, -r_{n-2}^{-1}),   \nn
&&\hspace{3cm}G_{n-1}v_{n-1}; -z_{l_1+1}^{-1}, -\bar{z}_{l_1+1}^{-1}, 
\dots, -z_l^{-1}, -\bar{z}_l^{-1}; 0, -r_{n-1}^{-1}), v_n)_{op}  \nn
&&\hspace{0.2cm}=\sum_{s\in \R} (e^{-r_{n-1}^{-1}L(-1)}
m_{cl-op}^{(l-l_1; 2)}( F_{l_1+1}u_{l_1+1}, \dots, F_lu_l; 
P_s e^{-r_{n-2}^{-1}L(-1)} m_{cl-op}^{(l-l_1; n-1)}( F_1u_1, \nn 
&&\hspace{3cm}\dots, F_{l_1} u_{l_1};
 v, G_1v_1, \dots, G_{n-1}v_{n-1};
-z_1^{-1}+r_{n-2}^{-1}, -\bar{z}_1^{-1}+r_{n-2}^{-1}, \nn 
&&\hspace{2cm}
\dots, -z_{l_1}^{-1}+r_{n-2}^{-1}, -\bar{z}_{l_1}^{-1}+r_{n-2}^{-1}; 
r_{n-2}^{-1}, -r_1^{-1}+r_{n-2}^{-1}, \dots, -r_{n-3}^{-1}+r_{n-2}^{-1},0), \nn
&&\hspace{3cm} G_{n-1}v_{n-1}; -z_{l_1+1}^{-1}+r_{n-1}^{-1}, 
-\bar{z}_{l_1+1}^{-1}+r_{n-1}^{-1}, \nn
&&\hspace{4cm}
\dots, -z_l^{-1}+r_{n-1}^{-1}, -\bar{z}_l^{-1}+r_{n-1}^{-1}; r_{n-1}^{-1}, 0),
v_n)_{op} 
\eea
Note that the position of $P_s$ and $e^{-r_{n-2}^{-1}L(-1)}$
can not be exchanged in general. Because if we exchange their
position, the sum may not converge and then the associativity law
does not hold. 
We want to use analytic continuation to move it to 
a domain such that we can freely apply associativity law.
By our assumption on $V$, both sides of (\ref{lemma-2-equ-1})
are restrictions of analytic function of 
$z_{l_1+1}, \zeta_{l_1+1}, \dots, z_l, \zeta_l, r_{n-1}$ on 
$\zeta_{l_1+1}=\bar{z}_{l_1+1}, \dots, \zeta_l=\bar{z}_l$.  
Let $\tilde{z}_{l_1+1}, \dots, \tilde{z}_l, \tilde{r}_{n-1}$ 
satisfy the following conditions:
\bea  \label{condition-r-n-1}
&& |-\tilde{z}_p^{-1}+r_{n-2}^{-1}|, 
\tilde{r}_{n-1}^{-1}- r_{n-2}^{-1} > 
|-z_i^{-1}+r_{n-2}^{-1}|, r_{n-2}^{-1}, -r_j^{-1}+r_{n-2}^{-1},  
\eea
for all $i=1,\dots, l,  j=1, \dots, n-3$ and $p=l_1+1, \dots, l$. 
Note that such condition define a nonempty open subset on 
$\HH^l\times \R_+$. 
Choose a path $\gamma_1$ in the complement of the diagonal in 
$\HH^l$ from initial point 
$(z_{l_1+1}, \dots, z_l)$ to $(\tilde{z}_{l_1+1}, \dots, \tilde{z}_l)$
and a path $\gamma_2$ in $\R_+$ from $r_{n-1}$ to $\tilde{r}_{n-1}$.
We also denote the complex conjugate of path $\gamma_1$ as
$\bar{\gamma}_1$, which is a path in $\overline{\HH}^l$. 
Combine $\gamma_1, \bar{\gamma}_1$ with $\gamma_2$, we obtain
a path $\gamma$ from  
$(z_{l_1+1}, \bar{z}_{l_1+1}, \dots, z_l, \bar{z}_l, r_{n-1})$
to
$(\tilde{z}_{l_1+1}, \bar{\tilde{z}}_{l_1+1}, \dots, 
\tilde{z}_l, \bar{\tilde{z}}_l, \tilde{r}_{n-1})$. 
Analytically continuating 
the right hand side of (\ref{lemma-2-equ-1})
along the path $\gamma$, we obtain, 
by the properties of $m_{cl-op}^{(l;n)}$ \cite{Ko2}, 
\bea
&&\hspace{0cm}(e^{-\tilde{r}_{n-1}^{-1}L(-1)}
m_{cl-op}^{(l-l_1; 2)}( F_{l_1+1}u_{l_1+1}, \dots, F_lu_l; 
P_s m_{cl-op}^{(l-l_1; n-1)}( F_1u_1, \dots, F_{l_1} u_{l_1};\nn 
&&\hspace{3cm}
 v, G_1v_1, \dots, G_{n-1}v_{n-1};
-z_1^{-1}+r_{n-2}^{-1}, -\bar{z}_1^{-1}+r_{n-2}^{-1}, \nn 
&&\hspace{2cm}
\dots, -z_{l_1}^{-1}+r_{n-2}^{-1}, -\bar{z}_{l_1}^{-1}+r_{n-2}^{-1}; 
r_{n-2}^{-1}, -r_1^{-1}+r_{n-2}^{-1}, \dots, -r_{n-3}^{-1}+r_{n-2}^{-1},0), \nn
&&\hspace{3cm} \tilde{G}_{n-1}v_{n-1}; -\tilde{z}_{l_1+1}^{-1}+r_{n-1}^{-1}, 
-\bar{\tilde{z}}_{l_1+1}^{-1}+r_{n-1}^{-1}, \nn
&&\hspace{4cm}
\dots, -\tilde{z}_l^{-1}+r_{n-1}^{-1}, 
-\bar{\tilde{z}}_l^{-1}+r_{n-1}^{-1}; 
\tilde{r}_{n-1}^{-1}-r_{n-2}^{-1}, 0),
v_n)_{op} 
\label{lemma-2-equ-2}
\eea
where $\tilde{G}_{n-1}=e^{-\tilde{r}_{n-1} L(1)}\tilde{r}_{n-1}^{-2L(0)}$. 
Using the associativity of open-closed field algebra 
and $L(-1)$-properties of $m_{cl-op}$, we see that
(\ref{lemma-2-equ-2}) further equals to
\bea
&&(m_{cl-op}^{(l;n+1)}( F_1u_1, \dots, F_l u_l;
v, G_1v_1, \dots, \tilde{G}_{n-1}v_{n-1};
-z_1^{-1}, -\bar{z}_1^{-1}, \dots, \nn 
&&\hspace{1cm} -z_{l_1}^{-1}, -\bar{z}_{l_1}^{-1}, 
-\tilde{z}_{l_1+1}^{-1}, -\bar{\tilde{z}}_{l_1+1}^{-1}, \dots, 
-z_l^{-1}, -\bar{z}_l^{-1}; 0,
-r_1^{-1}, \dots, -\tilde{r}_{n-1}^{-1}), v_n)_{op}.\nn
\label{lemma-2-equ-3}
\eea
By analytically continuating (\ref{lemma-2-equ-3}) along
the path $-\gamma$, which is $\gamma$ reversed, 
we obtain the right hand side of (\ref{lemma-2}).
Hence (\ref{lemma-2-equ-1}) and the right hand side of 
(\ref{lemma-2}) are analytic continuation of each other 
along path $(-\gamma)\circ (\gamma)$ which is a constant path. 
Hence (\ref{lemma-2}) holds for $l=k, n=m+1$.

Now let $l=k+1, n=m$. The proof is similar to the 
case $l=k, n=m+1$. We only point out the difference. 
Using the smoothness of $m_{cl-op}$,
it is enough to prove the case when $|z_i|\neq |z_j|$ for 
$i,j=1,\dots, l$ and $i\neq j$. Without losing generality, 
we assume that $|z_1|>\dots >|z_l|>0$. 
Let $n_1\leq n$ be the smallest so that $0<r_j<|z_l|$
for $j\geq n_1$. Then we have
\bea
&&(v, \, m_{cl-op}^{(l;n)}(u_1, \dots, u_l; v_1, \dots, v_n; 
z_1, \bar{z}_1, \dots, z_l, \bar{z}_l; r_1, \dots, r_n))_{op}\nn
&&\hspace{0.5cm}=
(v, \, m_{cl-op}^{(l-1; n_1)}(u_1, \dots, u_{l-1}; v_1, \dots, v_{n_1-1}, 
m_{cl-op}^{(1;n-n_1+1)}(u_l; v_{n_1}, \dots, v_n;   \nn
&&\hspace{3cm}  z_l,\bar{z}_l; r_{n_1}, \dots, r_n); 
z_1, \bar{z}_1, \dots, z_{l-1}, \bar{z}_{l-1}; 
r_1, \dots, r_{n_1-1}))_{op}
\eea
We can then apply (\ref{lemma-2}) as in (\ref{lemma-2-equ-1}) for the
case $l\leq k, n\leq m$, which is true by our induction hypothesis. 
The rest of proof is entirely same as that of the case $l=k, n=m+1$. 
\epf

We define a map, for $z, \zeta \in \C$ and $z\neq \zeta$,   
$\iota_{cl-op}(z, \zeta): V_{cl} \rightarrow \overline{V}_{op}$ as
$$
\iota_{cl-op}(z, \zeta)(u) = \mathbb{Y}_{cl-op}(u; z, \zeta)\one_{op}.
$$
We denote its adjoint as $\iota_{cl-op}^*(z,\zeta)$. 
Namely, $\iota_{cl-op}^*(z,\zeta): V_{op} \rightarrow \overline{V_{cl}}$ 
is given by 
\begin{equation}  \label{iota-up-*}
(\iota_{cl-op}^*(z,\zeta)(w), u)_{cl} = 
(w, \iota_{cl-op}(z,\zeta)(u))_{op}
\end{equation}
for any $u\in V_{cl}$ and $w\in V_{op}$.

Let $Q$ be an element in $\mathcal{T}_{\mathbb{S}}(n_-, n_+| m_-, m_+)$ 
of form (\ref{gen-ele-S}). Let $\alpha$ be the map 
(\ref{alpha-map}) so that $s_1, \dots, s_{n_-+n_+}$, 
defined as $s_i:=r_{\alpha^{-1}(i)}$, satisfy 
$\infty \geq s_1>\dots >s_{n_-+n_+}\geq 0$. 
Then we define
\beq   \label{def-psi-S-0}
\nu_{cl-op} ( \lambda \psi_{\mathbb{S}} (Q) )
( u_1\otimes \dots \otimes u_{m_+} \otimes v_1\otimes 
\dots \otimes v_{n_+} )
\eeq
as follow: 
\bnu
\item If $s_1 \neq \infty$, (\ref{def-psi-S-0}) is given by 
\bea \label{def-psi-S}
&&\lambda \,\,\, \sum_{i_1, \dots, i_{m_-}; j_1, \dots, j_{n_-}} 
\big( \one_{op}, \,\, m_{cl-op}^{(m_-+m_+;n_-+n_+)}( u_{-1}, \dots, u_{-m_-}, 
u_1, \dots, u_{m_+}; \nn
&&\hspace{2.5cm} w_1, \dots, w_{n_-+n_+};  
z_{-1}, \bar{z}_{-1}, \dots, z_{-m_-}, \bar{z}_{-m_-}, 
z_1, \bar{z}_1, \dots, z_{m_+}, \bar{z}_{m_+};  \nn
&&\hspace{2.5cm}
s_1, \dots, s_{n_-+n_+})\big)_{op} \cdot \, 
e^{i_{1}}\otimes \dots \otimes e^{i_{n_-}}\otimes
f^{j_1} \otimes \dots \otimes f^{j_{m_-}} . 
\eea
where 
\bea
u_{p} &=& e^{-L_+^L(A^{(p)})-L_+^R(\overline{A^{(p)}})} 
(a_0^{(p)})^{-L^L(0)} \overline{a_0^{(p)}}^{-L^R(0)}e_{i_{-p}},  \nn
u_q &=& e^{-L_+^L(A^{(q)})-L_+^R(\overline{A^{(q)}})} 
(a_0^{(q)})^{-L^L(0)} \overline{a_0^{(q)}}^{-L^R(0)} u_q,   \nn
w_{\alpha(k)} &=& e^{-L_+(B^{(k)})}(b_0^{(k)})^{-L(0)}f_{j_{-k}},  \nn
w_{\alpha(l)} &=& e^{-L_+(B^{(l)})}(b_0^{(l)})^{-L(0)} v_l,
\eea
for $p=-1, \dots, -n_-$, $q=1, \dots, n_+$, 
$k=-1, \dots, -m_-$ and $l=1, \dots, m_+$.

\item When $r_k=\infty$ for some $k=-m_-, \dots, -1, 1, \dots, m_+$. 
(\ref{def-psi-S-0}) is given by the formula obtain from 
(\ref{def-psi-S}) by exchanging $\one_{op}$ with $w_{\alpha(k)}$. 

\enu

\begin{lemma} {\rm  \label{well-def-psi-S}
$\nu_{cl-op}$ is $SL(2, \R)$-invariant. 
} \end{lemma}
\pf
The $SL(2,\R)$ is generated by the following three transformations
1. $w\mapsto aw, \forall a\in \R^+$; 2. $w\mapsto w-b, \forall b\in \R$;
3. $w\mapsto \frac{-1}{w}$. That $\nu_{cl-op}$ is invariant
under the first two transformations simply follows from the 
$L(0)$- and $L(-1)$-properties of $m_{cl-op}$. That $\nu_{cl-op}$ is
invariant under the third transformation is proved in 
Lemma \ref{lemma-2-st}. 
\epf

Hence $\nu_{cl-op}$ induces a map 
$\tilde{\mathbb{S}}^c \rightarrow 
\mathbb{E}_{V_{op}|V_{cl}}^{\R_+|\C^{\times}}$, which is still denoted as
$\nu_{cl-op}$. We list a few interesting cases:
\bea
\nu_{cl-op} \big( \psi_{\mathbb{S}} \big( 
\big[ (\,(\infty, 1, \mathbf{0}) |\,) \big\|  (\, |
(z, 1,\mathbf{0}) ) \big]_{\mathbb{S}} \big) \big) 
&=& \iota_{cl-op}(z,\bar{z}) ,  \nn
\nu_{cl-op} \big( \psi_{\mathbb{S}} \big( 
\big[ (\, |(\infty, 1, \mathbf{0})) \big\| ((z, 1,\mathbf{0})|\,) 
\big]_{\mathbb{S}}
\big) \big)
&=& \iota^{*}(z, \bar{z}) 
\eea
and for $b\in \R_+$, $B\in \R^{\infty}$, $a\in \C^{\times}$,
$A\in \C^{\infty}$ and $v\in V_{op}$, we have  
\bea  \label{iota-*}
&&\hspace{-0.5cm}
\nu_{cl-op}\big( \psi_{\mathbb{S}} \big( \big[ (\, |(\infty, b, B)) \big\| 
((z, a, A)|\,) \big]_{\mathbb{S}} \big) \big)   \nn
&&\hspace{0cm} = a^{-L^L(0)}\bar{a}^{-L^R(0)} 
e^{-\sum_{j=1}^{\infty} (-1)^j [A_j L^L(-j)+\overline{A_j}L^R(-j)]} 
\iota^{*}(z, \bar{z})(e^{-L_+(B)} b^{-L(0)} v).
\eea

\begin{thm} \label{thm-op-cl}
$(V_{op}|V_{cl}, \nu_{cl-op}, \nu_{cl})$ is an $\R_+|\C^{\times}$-rescalable
smooth $\tilde{\mathbb{S}}^c|\tilde{\mathbb{K}}^c \otimes 
\overline{\mathbb{K}^{\bar{c}}}$-algebra.
\end{thm}
\pf
The smoothness is automatic. We showed in \cite{Ko2}
that $(V_{op}|V_{cl}, \nu_{cl-op}, \nu_{cl})$ 
is an $\R_+|\C^{\times}$-rescalable
smooth $\tilde{\mathfrak{S}}^c|\tilde{K}^c \otimes 
\overline{K^{\bar{c}}}$-algebra.
The rest of the proof is similar to 
that of Theorem \ref{ffa-K-alg} in \cite{Ko1}.
We omit the detail here.  
\epf

\subsection{Ishibashi states}

As we mentioned in the introduction, 
an open-closed field algebra over $V$ 
equipped with nondegenerate invariant bilinear forms for
both open theory and closed theory
contains all the data needed to grow to an 
open-closed partial field field theory of all genus. Without
adding more compatibility conditions, itself is already
an interesting object to study. 
We show in this subsection that the famous 
``Ishibashi states'' \cite{I} can be studied 
in the framework of such algebra. 
Throughout this subsection, we fix an open-closed algebra over $V$
given in (\ref{opcl-fa}) and
equipped with nondegenerate invariant bilinear forms
$(\cdot, \cdot)_{cl}$ and $(\cdot, \cdot)_{op}$.

For $u\in V_{op}$ and $z_0 \in \HH$, we 
define the boundary state $B_{z_0}(u) \in \overline{V_{cl}}$
associated with $u$ and $z_0$ by 
\beq  \label{b-state-def}
B_{z_0}(u) = e^{L(-1)} (\bar{z}_0-z_0)^{L(0)} \otimes
e^{L(-1)} \overline{\bar{z}_0-z_0}^{\, L(0)}
\iota_{cl-op}^*(z_0, \bar{z}_0) (u). 
\eeq

\begin{prop}   \label{prop-ishi}
If $u\in V_{op}$ is a {\it vacuum-like vector} \cite{LL}, 
i.e.  $L(-1)u=0$, then, for $z_0\in \HH$,  $B_{z_0}(u)$ is an Ishibashi state, i.e.
\beq  \label{ishi-cond}
(L^L(n) - L^R(-n)) B_{z_0}(u) = 0, \quad\quad \forall n\in \Z. 
\eeq
\end{prop}
\pf
For $v\in V_{cl}$, the following two functions of $z$
\bea
& (u, \,\,\, \mathbb{Y}_{cl-op}(\omega^L, z+z_0)
\mathbb{Y}(v; z_0, \bar{z}_0)\one_{op})_{op},  & \nn
& (u, \,\,\,  \mathbb{Y}_{cl-op}(\omega^R;  \overline{z+z_0})
\mathbb{Y}(v; z_0, \bar{z}_0)\one_{op})_{op}  &
\eea
can be extended to a holomorphic function and 
an antiholomorphic function in 
$\{ z | z+z_0\in \HH, z\neq 0 \}$ respectively 
by our assumption on $V$. 
We denote the extended functions
by $g_1(\omega^L, z)$ and $g_2(\omega^R, \bar{z})$ respectively.

The following two limits
\bea  
&& \lim_{z+z_0\rightarrow r} \left( 1- \frac{z}{\bar{z}_0-z_0}
\right)^4 g_1(\omega^L, z)   \nn
&& \lim_{z+z_0\rightarrow r} \left( 1- \frac{z}{\bar{z}_0-z_0}
\right)^4 g_2(\omega^R, \bar{z})  \label{two-lim-on-b}
\eea
exist for all $r\in \R$.  
Using (\ref{inv-form-bcft-1}), it is easy to see that 
above two limits also exist for $r=\infty\in \hat{\R}$
if and only if $u$ is vacuum-like. Hence, 
by $V$-invariant (or conformal invariant) boundary condition, we have
\beq  \label{ishi-proof-equ--1}
\lim_{z+z_0\rightarrow r} \left( 1- \frac{z}{\bar{z}_0-z_0}
\right)^4  g_1(\omega^L, z)
= \lim_{z+z_0\rightarrow r} \left( 1- \frac{z}{\bar{z}_0-z_0}
\right)^4 g_2(\omega^R, z)
\eeq
for all $r\in \hat{\R}$ when $L(-1)u=0$.

On the other hand, for $|z+z_0|>|z_0|$, we have
\bea  \label{ishi-proof-equ-0}
&&(u, \,\, \mathbb{Y}_{cl-op}(\omega^L, z+z_0)
\mathbb{Y}_{cl-op}(v; z_0, \bar{z}_0)\one_{op})_{op}   \nn
&&\hspace{1cm}= (u, \,\, 
\mathbb{Y}_{cl-op}(\mathbb{Y}(\omega^L, z)v; 
z_0, \bar{z}_0)_{op}  \one_{op} )_{op} \nn
&&\hspace{1cm}= (u,\,\, \iota_{cl-op}(z_0, \bar{z}_0)(
\mathbb{Y}(\omega^L, z) v) )_{op}  \nn
&&\hspace{1cm}= (\iota_{cl-op}^*(z_0, \bar{z}_0)(u), \,\, 
\mathbb{Y}(\omega^L, z) v)_{cl}   \nn
&&\hspace{1cm}= (B_{z_0}(u), \,\,  
e^{L(1)} (\bar{z}_0-z_0)^{-L(0)} \otimes e^{L(1)} 
\overline{\bar{z}_0-z_0}^{\, -L(0)}
\mathbb{Y}(\omega^L, z) v)_{cl}.
\eea
Note that one should check the convergence property of
each step in (\ref{ishi-proof-equ-0}). In particular, 
in the last step, the convergence and equality follow
from the convergence of early steps and the fact that 
$(\bar{z}_0-z_0)^{-L(0)}e^{-L(-1)}\otimes \overline{\bar{z}_0-z_0}^{\, -L(0)}
e^{-L(-1)} \in \text{Aut}\, \overline{V_{cl}}$. 
For $0<|z|<|\text{Im}z_0|$, it is easy to show that 
\bea  \label{ishi-proof-equ-1}
&&e^{L(1)} (\bar{z}_0-z_0)^{-L(0)} \otimes e^{L(1)} 
\overline{\bar{z}_0-z_0}^{\, -L(0)}
\mathbb{Y}(\omega^L, z) v  \nn
&&\hspace{0.5cm}= 
\mathbb{Y}(e^{(1-\frac{z}{\bar{z}_0-z_0})L^L(1)}
\left( 1 - \frac{z}{\bar{z}_0-z_0}\right)^{-2L^L(0)} \cdot  \nn
&&\hspace{3cm} \cdot (\bar{z}_0-z_0)^{-L^L(0)} \omega^L, 
\frac{1}{\bar{z}_0 -z_0} \frac{z}{1- \frac{z}{\bar{z}_0-z_0}} )
v_1 \nn
&&\hspace{0.5cm}= 
\left( 1 - \frac{z}{\bar{z}_0-z_0}\right)^{-4} (\bar{z}_0-z_0)^{-2} 
\mathbb{Y}(\omega^L, \frac{1}{\bar{z}_0 -z_0} 
\frac{z}{1- \frac{z}{\bar{z}_0-z_0}} ) v_1
\eea
where $v_1= e^{L(1)} (\bar{z}_0-z_0)^{-L(0)} \otimes e^{L(1)} 
\overline{\bar{z}_0-z_0}^{\, -L(0)}v$. 
Hence, for all $0<|z|<|\text{Im}z_0|$ and 
$|z+z_0|>|z_0|$, we obtain
\bea
&&(u, \mathbb{Y}_{cl-op}(\omega^L, z+z_0)
\mathbb{Y}_{cl-op}(v; z_0, \bar{z}_0)\one_{op})_{op}   \nn
&&\hspace{0.2cm}=
(B_{z_0}(u), \,\, 
\left( 1 - \frac{z}{\bar{z}_0-z_0}\right)^{-4} (\bar{z}_0-z_0)^{-2} 
\mathbb{Y}(\omega^L, f(z)) v_1 )_{cl}  \label{ishi-proof-equ-2}
\eea
where $f$ is the composition of the following maps:
$$
w \mapsto w+z_0 \mapsto -\frac{(w+z_0)-z_0}{(w+z_0)-\bar{z}_0} 
=\frac{1}{\bar{z}_0 -z_0} 
\frac{w}{1- \frac{w}{\bar{z}_0-z_0}} 
$$
which maps the domain $\HH -z_0$ to the unit disk.
Since $g_1(\omega^L, z)$ is analytic and 
free of singularities for $z+z_0\in \HH 
\backslash \{z_0\}$, the right hand side of 
(\ref{ishi-proof-equ-2}) can also be 
extended to an analytic function in $z\in \HH-z_0 
\backslash \{ 0 \}$.  If we view $f(z)$ as a new variable
$\xi$, then the right hand side of (\ref{ishi-proof-equ-2})
can be extended to an analytic function on $\{ \xi | 1>|\xi|>0 \}$,
which has a Laurent series expansion. By the uniqueness of
Laurent expansion, the right hand side of (\ref{ishi-proof-equ-2})
gives exactly such Laurent expansion and 
thus is absolutely convergent in $\{ \xi | 1>|\xi|>0 \}$. 
Moreover, $\lim_{z+z_0\rightarrow r} g_1(\omega^L, z)$ 
exists for all $r\in \hat{\R}$. 
By the properties of Laurent series,  
the right hand side of (\ref{ishi-proof-equ-2}) 
must converge absolutely for all $f(z) \in 
\{ \xi | |\xi|=1 \}$ to the function given by
$\lim_{z+z_0\rightarrow r} g_1(\omega^L, z), r\in \hat{\R}$.

Similarly, for all $0<|\bar{z}|<|\text{Im} z_0|$
and $|z+z_0|>|z_0|$, we have
\bea
&&(u, \mathbb{Y}_{cl-op}(\omega^R, \overline{z+z_0})
\mathbb{Y}_{cl-op}(v; z_0, \bar{z}_0)\one_{op})_{op}   \nn
&&\hspace{0.2cm}= (B_{z_0}(u),  \,\, 
e^{L(1)} (\bar{z}_0-z_0)^{-L(0)} \otimes e^{L(1)} 
\overline{\bar{z}_0-z_0}^{\, -L(0)}
\mathbb{Y}(\omega^R, \bar{z}) v)_{cl} \nn
&&\hspace{0.2cm} =(B_{z_0}(u), \,\, 
\mathbb{Y}(
e^{(1-\frac{\bar{z}}{z_0-\bar{z}_0})L^R(1)}
\left( 1 - \frac{z}{z_0-\bar{z}_0}\right)^{-2L^R(0)}
(z_0-\bar{z}_0)^{-L^R(0)} \omega^R, g(\bar{z}) )v_1)_{cl}   \nn 
&&\hspace{0.2cm}=
\left( 1 - \frac{\bar{z}}{z_0- \bar{z}_0}\right)^{-4} 
(z_0-\bar{z}_0)^{-2}(B_{z_0}(u),
\mathbb{Y}(\omega^R, g(\bar{z}) )v_1)_{cl}    \label{ishi-proof-equ-3}
\eea
where $g$ is the 
composition of the following maps:
$$
w \mapsto w+\bar{z}_0 \mapsto  
- \frac{(w+\bar{z}_0)-z_0}{(w+\bar{z}_0)-\bar{z}_0} \mapsto 
 - \frac{(w+\bar{z}_0)-\bar{z}_0}{(w+\bar{z}_0)-z_0} = 
\frac{1}{z_0-\bar{z}_0} 
\frac{w}{1- \frac{w}{z_0-\bar{z}_0}}
$$
which maps the domain $-\HH-\bar{z}_0$ to the unit disk.  
Moreover, the right hand side of (\ref{ishi-proof-equ-3}), as 
a Laurent series of $g(\bar{z})$,  is
absolutely convergent for all 
$g(\bar{z}) \in \{ \xi | |\xi| = 1 \}$ to 
$\lim_{z+z_0\rightarrow r} g_2(\omega^R, z), r\in \hat{\R}$.

Also notice that 
\beq
g(r-\bar{z}_0) = \frac{1}{f(r-z_0)} \in \{ e^{i\theta} | 0\leq \theta <2\pi \}
\eeq
for all $r\in \hat{\R}$. 
Using (\ref{ishi-proof-equ--1}) and 
by replacing $z$ in (\ref{ishi-proof-equ-2})
by $r-z_0$ and $\bar{z}$ in (\ref{ishi-proof-equ-3}) by
$r-\bar{z}_0$, we obtain the following identity:
\beq  \label{ishi-proof-equ-4}
(B_{z_0}(u), \,\, \mathbb{Y}(\omega^L, e^{i\theta})v_1)_{cl}
= (B_{z_0}(u), \,\, \mathbb{Y}(\omega^R, e^{-i\theta})v_1)_{cl}\,\, 
e^{-4i\theta}
\eeq
where $e^{i\theta}= f(r-z_0)$, for all $0\leq \theta <2\pi$.
Notice that the existence of both sides of
(\ref{ishi-proof-equ-4}) follows directly from 
(\ref{ishi-proof-equ--1}), which further follows 
from the condition of $u$ being vacuum-like. 
Then we obtain
$$
\sum_{n\in \Z} (B_{z_0}(u), \,\, L^L(n) v_1)_{cl}\,\, e^{i\theta (-n-2)} 
=\sum_{n\in \Z} (B_{z_0}(u), L^R(-n) v_1)_{cl} \,\,  e^{i\theta (-n-2)}
$$ 
for all $0\leq \theta <2\pi$. Notice that $v_1$ can be arbitrary. 
Therefore, we must have (\ref{ishi-cond}) when $L(-1)u=0$. 
\epf

In physics, boundary states are usually obtained by
solving the equation (\ref{ishi-cond}). The solutions of 
such equation was first obtained by Ishibashi \cite{I}. 
They are called Ishibashi states.  The definition of 
boundary states we gives in (\ref{b-state-def}) is more general. 
Boundary conditions are also called ``D-branes'' in string theory. 
If $u$ is not a vacuum-like vector, 
the boundary state (\ref{b-state-def}) associated with $u$
is also very interesting in physics (see for example \cite{FFFS1}). 
Such boundary states are associated to the geometry on D-branes. 
In Section 5.2, we will give a more natural 
definition of D-brane (Definition \ref{def:dbrane}).

\renewcommand{\theequation}{\thesection.\arabic{equation}}
\renewcommand{\thethm}{\thesection.\arabic{thm}}
\setcounter{equation}{0}
\setcounter{thm}{0}

\section{Cardy condition}

In this section, we derive the Cardy condition from the axioms of open-closed partial conformal field theory by writing out the algebraic realizations of the both sides of Figure \ref{cardy-top-fig} explicitly. Then we reformulate the Cardy condition in the framework of intertwining operator algebra.  Throughout this section, we fix an open-closed field algebra over $V$ given in (\ref{opcl-fa}) equipped with nondegenerate invariant bilinear forms $(\cdot, \cdot)_{op}$ and $(\cdot, \cdot)_{cl}$. 

\subsection{The first version}

In the Swiss-cheese dioperad, we exclude an interior sewing operation
between two disks with strips and tubes and a self-sewing operation between two oppositely oriented boundary punctures on a single disk. The surface obtained after these two types of sewing operations can be the same cylinder or annulus. The axioms of open-closed partial conformal field theory require that the algebraic realization of these two sewing operations must coincide. This gives a nontrivial condition called Cardy condition (recall Figure \ref{cardy-top-fig}).

Although the Cardy condition only involves
genus-zero surfaces, its algebraic realization is 
genus-one in nature. This fact is manifest if we 
consider the doubling map $\delta$. 
A double of a cylinder is actually a torus. Hence 
the Cardy condition is a condition on the equivalence of
two algebraic realizations of 
two different decompositions of a torus. 
This is nothing but a condition associated to modularity.

\begin{figure}
\begin{center}
\includegraphics[width=1\textwidth]{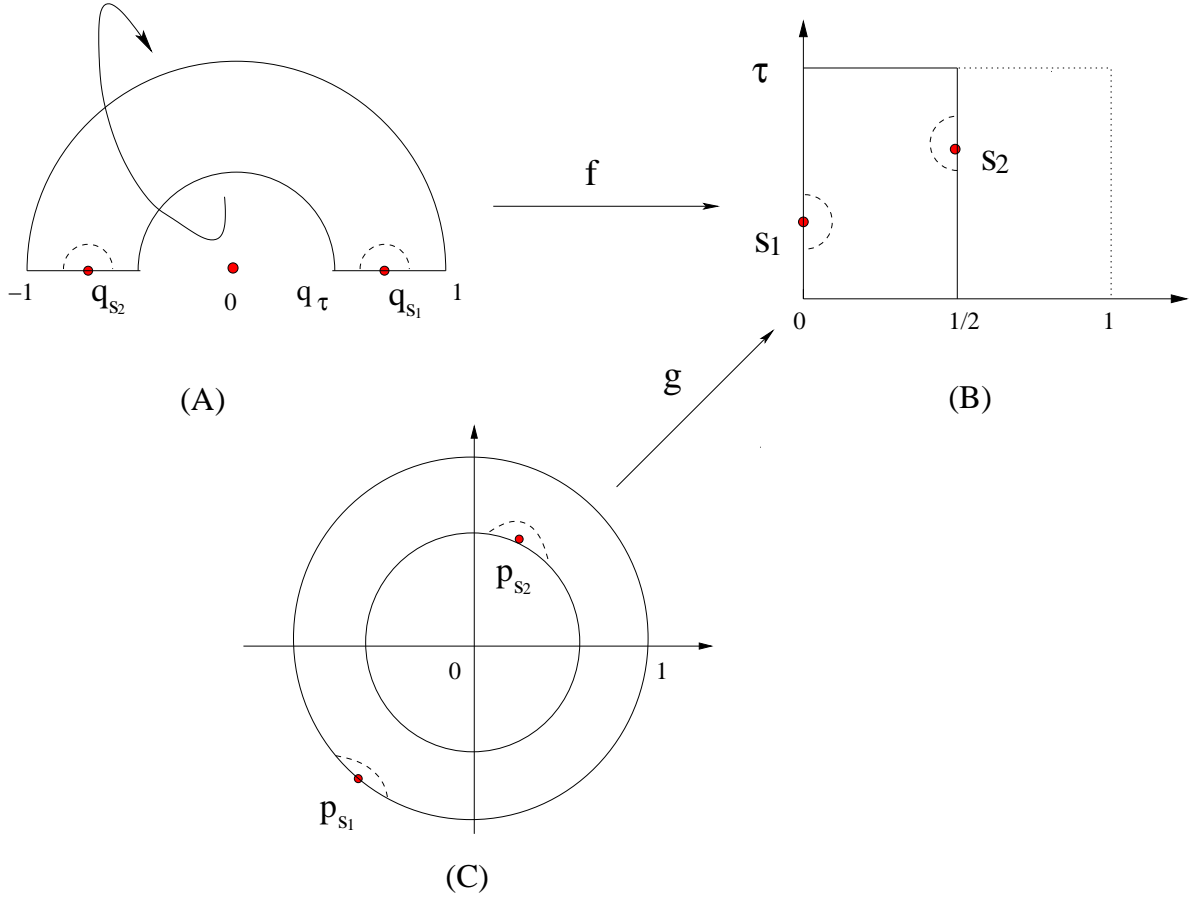}
\end{center}
\caption{Cardy condition: 
two different sewings give same annulus.}
\label{cardy-fig}
\end{figure}

That an annulus can be obtained by two different sewing operations
is also shown in Figure \ref{cardy-fig}. In particular, 
the surface (A) in Figure \ref{cardy-fig} shows how 
an annulus is obtained by
sewing two oppositely oriented boundary punctures 
on the same disk with strips and tubes in $\mathbb{S}(1,3|0,0)$, 
and surface (C) in Figure \ref{cardy-fig}, viewed 
as a propagator of close string, can be obtained by 
sewing an element in $\mathbb{S}(0,1|1,0)$ with an 
element in $\mathbb{S}(0,1|0,1)$ along the interior punctures. 
We only show in Figure \ref{cardy-fig} 
a simple case in which there are only two boundary punctures 
and no interior puncture. 
In general, the number of boundary punctures and interior 
punctures can be arbitrary. However, all general cases can 
be reduced to this simple case by applying associativities. 
Notice that the two boundary punctures in this simple case
can not be reduced further by the associativities. 
We only focus on this case in this work. 

The conformal map $f$ between the surface (A) and (B) and 
$g$ between the surface (C) and (B) in Figure \ref{cardy-fig}
are given by 
\bea
f(w) &=& \frac{1}{2\pi i} \log w,  \nn
g(w) &=& \frac{-\tau}{2\pi i}\log w. 
\eea
It is also useful to know their inverses
$f^{-1}(w)= e^{2\pi i w}$, $g^{-1}(w)=e^{2\pi i (-\frac{1}{\tau}) w}$.
For any $z \in \C$, we set $q_{z}:=e^{2\pi iz}$
and $p_{z}:=e^{2\pi i (-\frac{1}{\tau})z}$. 
The radius of the outer circle of the surface $(C)$
in Figure \ref{cardy-fig} is $|p_{s_1}|=1$ and that of 
inner circle is 
\beq  \label{ab-val-p-s2}
|p_{s_2}|=e^{\pi i (-\frac{1}{\tau})} = q_{-\frac{1}{\tau}}^{1/2}.
\eeq

As we have argued and will show more explicitly later, 
Cardy condition deeply related to modularity. 
In \cite{H10}, Huang introduced the 
so-called geometrically modified intertwining operators, 
which is very convenient for the study of modularity.
He was motivated by the fact that it is much easier to 
study modularity in the global coordinates. Namely, 
one should choose the local coordinates at $s_1, s_2$  in surface (B)
in Figure \ref{cardy-fig} as simple as possible. 
More precisely, we
choose the local coordinates at $s_1,s_2$ as
\bea \label{local-coor-B}
w &\mapsto& e^{\frac{\pi i}{2}}(w-s_1) \nn
w &\mapsto& e^{-\frac{\pi i}{2}}(w-s_2) 
\eea
respectively. 
Correspondingly, the local coordinates at punctures 
$q_{s_1},q_{s_2}$ are
\bea
f_{q_{s_1}}(w) &=& e^{\frac{\pi i}{2}} \frac{1}{2\pi i}\log 
\left( 1 + \frac{1}{q_{s_1}}\, x
\right) \lbar_{x=w-q_{s_1}} \nn
f_{q_{s_2}}(w) &=& e^{-\frac{\pi i}{2}} \frac{1}{2\pi i}\log 
\left( 1 + \frac{1}{q_{s_2}} \, x
\right) \lbar_{x=w-q_{s_2}} 
\eea
respectively.

Notice that both local coordinates
$f_{q_{s_1}}(w)$ and $f_{q_{s_2}}(w)$ are real analytic. 
Hence $\exists B^{(i)}_j\in \R, b_0^{(i)}\in \R_+, i=1,2$ such that 
\bea
f_{q_{s_1}}(w) &=& e^{\sum_{j=1}^{\infty} B_j^{(1)} x^{j+1}\frac{d}{dx}} 
(b_0^{(1)})^{x\frac{d}{dx}} x \lbar_{x=w-q_{s_1}}  \,\, ,\nn
f_{q_{s_2}}(w) &=& e^{\sum_{j=1}^{\infty} B_j^{(2)} x^{j+1}\frac{d}{dx}} 
(b_0^{(2)})^{x\frac{d}{dx}} x \lbar_{x=w-q_{s_2}}  \,\, .
\eea
Then the algebraic realization of the surface (A) gives a map 
$V_{op}\otimes V_{op}\rightarrow \C$. 
We assume $1>|q_{s_1}|>|q_{s_2}|>|q_{\tau}|>0$. By the axiom of open-closed
conformal field theory, this map must be given by 
(recall (\ref{Y-op-opposite}))
\begin{equation} \label{cardy-L-0}
v_1\otimes v_2 \mapsto \text{{\large Tr}}_{V_{op}} 
\left( Y_{op}( T_{q_{s_1}} v_1, q_{s_1})
Y_{op}( T_{q_{s_2}} v_2, q_{s_2}) q_{\tau}^{L(0)} \right),
\end{equation}
where 
\bea
T_{q_{s_1}} &=& e^{-\sum_{j=1}^{\infty} B_j^{(1)}L(j)} (b_0^{(1)})^{-L(0)}, \nn  
T_{q_{s_2}} &=& e^{-\sum_{j=1}^{\infty} B_j^{(2)}L(j)} (b_0^{(2)})^{-L(0)}.
\eea

We need rewrite $T_{q_{s_1}}$ and $T_{q_{s_2}}$. 
Let $A_j, j\in \Z_+$, be the complex numbers defined by
$$
\log (1 + y) = e^{\sum_{j=1}^{\infty} A_j y^{j+1}\frac{d}{dy}}y.
$$
It is clear that $A_j\in \R$. 
Hence we also obtain:
\begin{equation}
\frac{1}{2\pi i} 
\log \left( 1 + \frac{1}{z} x \right)  
= z^{-x\frac{d}{dx}} e^{\sum_{j\in \Z_+} A_j x^{j+1}\frac{d}{dx} } 
(2\pi i)^{-x\frac{d}{dx}} x.
\end{equation}

\begin{lemma} \label{lemma-T-s-1-2}
\bea
T_{q_{s_1}} &=& (q_{s_1})^{L(0)} 
e^{-\sum_{j=1}^{\infty} A_j L(j)} (2\pi i)^{L(0)} e^{-\frac{\pi i}{2} L(0)}, 
\label{T-s-1}  \nn
T_{q_{s_2}} &=& (q_{s_2})^{L(0)} 
e^{-\sum_{j=1}^{\infty} A_j L(j)} (2\pi i)^{L(0)} e^{-2\pi iL(0)} e^{\frac{\pi i}{2} L(0)}. 
\label{T-s-2}
\eea
\end{lemma}
\pf
By results in \cite{BHL}, if
$e^{\sum_{j=1}^{\infty} C_j L(j)} c_0^{L(0)} = e^{\sum_{j=1}^{\infty} D_j L(j)} d_0^{L(0)}$
for any $C_j, D_j \in \C, c_0, d_0\in \C$, we must have
$C_j=D_j$ and $c_0=d_0$. Therefore, 
by moving the factor $(q_{s_1})^{L(0)}$ to the right side of 
$e^{-\sum_{j=1}^{\infty} A_j L(j)}$ in (\ref{T-s-1}) 
and similarly moving the factor $(q_{s_2})^{L(0)}$ in (\ref{T-s-1}), 
we see that it is enough to show
\bea
(b_0^{(1)})^{L(0)} &=& (q_{s_1})^{L(0)}(2\pi i)^{L(0)} e^{-\frac{\pi i}{2} L(0)}, \nn
(b_0^{(1)})^{L(0)} &=& (q_{s_2})^{L(0)} 
(2\pi i)^{L(0)} e^{-2\pi iL(0)} e^{\frac{\pi i}{2} L(0)}.
\eea
Using our conventions (\ref{branch-cut})(\ref{power-f}),
it is easy to check that above identities hold. 
\epf

Let $W$ be a $V$-module. Huang introduced the following operator
in \cite{H10}.
$$
\mathcal{U}(x) := x^{L(0)} e^{-\sum_{j=1}^{\infty} A_jL(j)} (2\pi i)^{L(0)}  
\in \, (\text{End}W)\{ x\}.
$$
Thus (\ref{cardy-L-0}) can be rewritten as follow:
\begin{equation} \label{cardy-L}
v_1\otimes v_2 \mapsto \text{{\large Tr}}_{V_{op}} 
\left( 
Y_{op}(\mathcal{U}(q_{s_1}) v'_1, q_{s_1})
Y_{op}(\mathcal{U}(q_{s_2}) e^{-2\pi iL(0)} v'_2, q_{s_2}) q_{\tau}^{L(0)} \right)
\end{equation}
where $v'_1= e^{-\frac{\pi i}{2}L(0)}v_1$ and 
$v'_2= e^{\frac{\pi i}{2}L(0)}v_2$.

Now we consider the algebraic realization of the 
surface (C) in Figure \ref{cardy-fig} obtained from
an interior sewing operation between an 
element in $\mathbb{S}(0,1|1,0)$ and an element in 
$\mathbb{S}(0,1|0,1)$.

\begin{lemma}
$\forall r\in (0,1)\subset \R_+$, 
the surface $C$ in Figure \ref{cardy-fig} is
conformally equivalent to 
a surface $Q_1\,\, {}_{^1}\infty_{^1}^I \, Q_2$, where 
\bnu
\item $Q_1\in \mathbb{S}(0,1|0,1)$ with 
punctures at $z_1\in \HH, \infty \in \hat{\R}$ 
and local coordinates:
\bea  
f_{z_1}: w &\mapsto& 
\left( \frac{p_{s_1}}{r}\right) \frac{w-z_1}{w-\bar{z}_1}\, , 
\label{lco-z-1}  \\
f_{\infty}^{(1)}: w &\mapsto& 
e^{\frac{\pi i}{2}} \, \left( \frac{-\tau}{2\pi i} \right) \,
\log \frac{w-z_1}{w-\bar{z}_1}\,\, ; \label{lco-infty-1}
\eea
\item $Q_2 \in \mathbb{S}(0,1|1,0)$ with 
punctures at $z_2\in \HH, \infty\in \hat{\R}$ with local coordinates
\bea  
f_{z_2}: w &\mapsto& \left( \frac{-r}{p_{s_2}} \right)
\frac{w-z_2}{w-\bar{z}_2} \, ,
\label{lco-z-2}  \\
f_{\infty}^{(2)}: w &\mapsto& 
e^{-\frac{\pi i}{2}}  \, \left( \frac{-\tau}{2\pi i} \right)\,
\log \frac{w-\bar{z}_2}{w-z_2} \,\, . \label{lco-infty-2}
\eea
\enu
\end{lemma}
\pf
Let us first define two disks $D_1$ and $D_2$. 
$D_1$ is the unit disk, i.e. $D_1:=\{ z\in \C | |z|\leq 1 \}$, which has 
punctures at $p_{s_1}, 0$ and local coordinates:
\bea
g_{p_{s_1}}: w &\mapsto &  
e^{\frac{\pi i}{2}}\left( \frac{-\tau}{2\pi i} \, \log w -s_1 \right),  
\label{g-s-1}  \\
g_0:    w &\mapsto & r^{-1}w.    \label{g-0}
\eea
$D_2$ is the disk $\{ z\in \hat{\C} | |z|\geq |p_{s_2}| \}$ 
which has punctures at $p_{s_2}, \infty$ and local coordinates:
\bea
g_{p_{s_2}}: w &\mapsto & 
e^{-\frac{\pi i}{2}}\left( \frac{-\tau}{2\pi i} \, \log w - s_2 \right),
\label{g-s-2}  \\
g_{\infty}: w &\mapsto & \frac{-r}{w}. \label{g-infty}
\eea
It is not hard to see that 
the surface $C$ in Figure 1 can be obtained by sewing the puncture 
$0\in D_1$ with the puncture $\infty\in D_2$ 
according to the usual definition of interior sewing operation. 
Then it is enough to show that $D_1$ and $D_2$ are conformally equivalent to
$P$ and $Q$ respectively. We define two maps
$h_1: Q_1\rightarrow D_1$ and $h_2: Q_2\rightarrow D_2$ as follow: 
\bea
h_1: w &\mapsto& p_{s_1} \frac{w-z_1}{w-\bar{z}_1},  \nn
h_2: w &\mapsto& p_{s_2} \frac{w-\bar{z}_2}{w-z_2}.  \nonumber
\eea
It is clear that $h_1$ and $h_2$ are both biholomorphic. 
We can check directly 
that $h_1$ and $h_2$ map punctures to punctures and 
preserve local coordinates as well. 
\epf

Using (\ref{iota-*}),  we obtain the 
the algebraic realization the annulus $C$ in Figure 1 as follow:
\beq \label{cardy-R}
v_1\otimes v_2 \mapsto 
( (T_1^L \otimes T_1^R)^{*} \, \iota_{cl-op}^*(z_1, \bar{z}_1) ( T_2 v'_1), 
\,\, (T_3^L \otimes T_3^R)^* \, \iota_{cl-op}^*(z_2, \bar{z}_2) (T_4 v'_2))_{cl}
\eeq
where $T_1^{L,R}, T_2, T_3^{L,R}, T_4$ are conformal transformations
determined by local coordinates $f_{z_1}, f_{\infty}^{(2)},
f_{z_2}^{(1)}, f_{\infty}^{(2)}$, and $(T_i^L \otimes T_i^R)^{*}$
is the adjoint of $T_i^L \otimes T_i^R$ with respect to 
the bilinear form $(\cdot, \cdot)_{cl}$ for $i=1,3$, and  
$v'_1= e^{-\frac{\pi i}{2}L(0)}v_1$ and 
$v'_2= e^{\frac{\pi i}{2}L(0)}v_2$.

\begin{lemma}
Recall the convention (\ref{power-f}), we have
\bea
T_1^L &=& (z_1-\bar{z}_1)^{L(0)} e^{L(1)} 
\left( \frac{p_{s_1}}{r} \right)^{-L(0)}, 
T_1^R =\overline{(z_1-\bar{z}_1)}^{\, L(0)} e^{L(1)} 
\overline{\left( \frac{p_{s_1}}{r} \right)}^{\, -L(0)},  
\label{T-1-L-R}  \\
T_2 &=& e^{\bar{z}_1L(1)}(z_1-\bar{z}_1)^{-L(0)}\mathcal{U}(1) 
\left( -\frac{1}{\tau} \right)^{L(0)}, \label{T-2}  \\
T_3^L &=&  (\bar{z}_2-z_2)^{L(0)} e^{-L(1)} 
\left( \frac{p_{s_2}}{r}\right)^{L(0)}, 
T_3^R=\overline{(\bar{z}_2-z_2)}^{\, L(0)} e^{-L(1)} 
\overline{\left( \frac{p_{s_2}}{r}\right)}^{\, L(0)},
\label{T-3-L-R}  \\
T_4 &=& e^{z_2L(1)} (\bar{z}_2 - z_2)^{-L(0)} \mathcal{U}(1)
\left( -\frac{1}{\tau} \right)^{L(0)}.  \label{T-4}
\eea
\end{lemma}
\pf
From (\ref{lco-z-1}) and (\ref{lco-z-2}), we obtain 
\bea
f_{z_1}: w  &\mapsto&  
(z_1-\bar{z}_1)^{-x\frac{d}{dx}} e^{-x^2\frac{d}{dx}} 
\left( \frac{p_{s_1}}{r} \right)^{x\frac{d}{dx}} x \lbar_{x=w-z_1} \nn
f_{z_2}: w &\mapsto& 
(\bar{z}_2-z_2)^{-x\frac{d}{dx}} e^{x^2\frac{d}{dx}} 
\left( \frac{p_{s_2}}{r} \right)^{-x\frac{d}{dx}} x \lbar_{x=w-z_2}
\eea
Then (\ref{T-1-L-R}) and (\ref{T-3-L-R}) is obvious. 
Notice that the expression 
(\ref{T-1-L-R}) is independent of our choice of branch cut
as long as we keep the convention (\ref{power-f}).

From (\ref{lco-infty-1}) and (\ref{lco-infty-2}), we obtain
\bea
&&\hspace{-0.8cm}f_{\infty}^{(1)}(w) =
e^{-\bar{z}_1x^2\frac{d}{dx}} (z_1-\bar{z}_1)^{x\frac{d}{dx}} 
e^{\sum_{j=1}^{\infty} A_j x^{j+1}\frac{d}{dx}} (2\pi i)^{-x\frac{d}{dx}}
\left( -\frac{1}{\tau} \right)^{-x\frac{d}{dx}} 
e^{\frac{\pi i}{2}x\frac{d}{dx}}  x \lbar_{x=\frac{-1}{w}}  \label{f-infty-1}
\\
&&\hspace{-0.8cm}f_{\infty}^{(2)}(w) =
e^{-z_2x^2\frac{d}{dx}} (\bar{z}_2- z_2)^{x\frac{d}{dx}} 
e^{\sum_{j=1}^{\infty} A_j x^{j+1}\frac{d}{dx}} (2\pi i)^{-x\frac{d}{dx}}
\left( -\frac{1}{\tau} \right)^{-x\frac{d}{dx}} 
e^{-\frac{\pi i}{2}x\frac{d}{dx}}  x \lbar_{x=\frac{-1}{w}}  \label{f-infty-2}
\eea
Recall that $f_{\infty}^{(1)}, f_{\infty}^{(2)}$ are
both real analytic. 
Similar to the proof of Lemma \ref{lemma-T-s-1-2}, 
to show (\ref{T-2}) and (\ref{T-4}) 
is enough to show that
\bea
(b_{\infty}^{(1)})^{L(0)} &=& (z_1-\bar{z}_1)^{-L(0)} 
(2\pi i)^{L(0)} \left( -\frac{1}{\tau} \right)^{L(0)} 
e^{-\frac{\pi i}{2} L(0)} \label{b-finty-1}  \\
(b_{\infty}^{(2)})^{L(0)} &=& (\bar{z}_2- z_2)^{-L(0)} 
(2\pi i)^{L(0)} \left( -\frac{1}{\tau} \right)^{L(0)} 
e^{\frac{\pi i}{2} L(0)} \label{b-finty-2}
\eea
for some $b_{\infty}^{(1)}, b_{\infty}^{(2)}\in \R_+$. 
Using our convention (\ref{branch-cut}) and (\ref{power-f}), 
it is a direct check that (\ref{b-finty-1}) 
and (\ref{b-finty-2}) holds. 
\epf

Combining (\ref{cardy-L}), (\ref{cardy-R}) and additional natural
factors $q_{\tau}^{-\frac{c}{24}}, q_{-\frac{1}{\tau}}^{-\frac{c}{24}}$ 
(see \cite{Z}\cite{H10}), which is due to determinant line bundle on
torus \cite{Se1}\cite{Kr},
we obtain the following formulation of the Cardy condition: 
\begin{defn}  {\rm 
The open-closed field algebra over $V$ given in (\ref{opcl-fa}) 
and equipped with nondegenerate bilinear forms 
$(\cdot, \cdot)_{op}$ and $(\cdot, \cdot)_{cl}$ is said to 
satisfy {\it Cardy condition}
if the left hand sides of the following formula,  
$\forall z_1, z_2\in \HH, v_1, v_2\in V_{op}$,  
\bea  \label{cardy}
&&\hspace{-0.5cm}\text{{\large Tr}}_{V_{op}} 
\left( Y_{op}(\mathcal{U}(q_{s_1}) v_1, q_{s_1}) 
Y_{op}(\mathcal{U}(q_{s_2}) e^{-2\pi iL(0)} v_2, q_{s_2}) 
q_{\tau}^{L(0)-c/24} \right) \nn
&&\hspace{-0.3cm}=
\left( (T_1^L\otimes T_1^R)^* \iota_{cl-op}^*(z_1, \bar{z}_1) ( T_2 v_1), \,
q_{-\frac{1}{\tau}}^{-c/24} (T_3^L\otimes T_3^R)^*   
\iota_{cl-op}^*(z_2, \bar{z}_2) (T_4 v_2) \, \right)_{cl}    
\eea
converge absolutely when $1>|q_{s_1}|>|q_{s_2}|> |q_{\tau}|>0$, 
and the right hand side of (\ref{cardy}) converge absolutely 
for all $s_1, s_2\in \HH$ satisfying 
$\text{Re}\, s_1=0, \text{Re}\, s_2=\frac{1}{2}$. 
Moreover,  the equality (\ref{cardy}) holds
when $1>|q_{s_1}|>|q_{s_2}|> |q_{\tau}|>0$.
}
\end{defn}

\begin{rema}
{\rm The dependence of $z_1,z_2,r$ of the right hand side of 
(\ref{cardy}) is superficial as required by the 
independence of $z_1,z_2,r$ of the left hand side of 
(\ref{cardy}). We will see it more explicitly later. 
}
\end{rema}

Using the definition of boundary states (\ref{b-state-def}), 
(\ref{cardy}) can also be written as follow:
\bea  \label{cardy-b-st}
&&\hspace{-0.5cm}\text{{\large Tr}}_{V_{op}} 
\left( Y_{op}(\mathcal{U}(q_{s_1}) v_1, q_{s_1}) 
Y_{op}(\mathcal{U}(q_{s_2}) e^{-2\pi iL(0)} v_2, q_{s_2}) 
q_{\tau}^{L(0)-c/24} \right) \nn
&&\hspace{-0.3cm}=
\left(  B_{z_1}( T_2 v_1), \,\, 
\left( \frac{p_{s_2}}{-p_{s_1}}\right)^{L(0)} \otimes 
\left( \overline{\frac{p_{s_2}}{-p_{s_1}}} \right)^{L(0)} 
q_{-\frac{1}{\tau}}^{-c/24} \, B_{z_2}(T_4 v_2) \, \right)_{cl} .
\eea

\subsection{The second version}

In this subsection, we rewrite the Cardy condition (\ref{cardy}) 
in the framework of intertwining operator algebra.

Since $V$ satisfies the conditions in Theorem \ref{ioa},
it has only finite number of inequivalent irreducible modules. 
Let $\I$ be the set of equivalence classes of irreducible 
$V$-modules.
We denote the equivalence class of the adjoint module $V$ as $e$,
i.e. $e\in \I$. 
Let $W_a$ be a chosen representative of $a\in \I$.

For any $V$-module $(W, Y_W)$, we denote the graded dual space of 
$W$ as $W'$, i.e. $W'=\oplus_{n\in \C} (W_{(n)})^*$. There is a
contragredient module structure on $W'$ \cite{FHL} 
given by a vertex operator $Y'_W$, which is defined as follow
\beq  \label{contra-module}
\langle Y'_W(u, x) w', w\rangle :=
\langle w', Y_W(e^{-xL(1)}x^{-2L(0)}u, -x^{-1}) w\rangle
\eeq
for $u\in V, w\in W, w'\in W'$. $(W', Y'_W)$ (or simply $W$) 
is the only module structure on $W'$ we use in this work. 
So we can set $Y_{W'}:= Y'_{W}$. We denote the equivalent class
of $W'_{a}$ as $a'$. It is harmless to set $W'_{a}=W_{a'}$. 
Moreover, $W''$ is canonically identified with $W$. 
Hence $a''=a$ for $a\in \I$.

By assumption on $V$, $V'\cong V$, i.e. $e'=e$.   
From \cite{FHL}, there is a nondegenerate invariant bilinear form 
$(\cdot, \cdot)$ on $V$ such that $(\one, \one)=1$. 
This bilinear form specifies a unique isomorphism 
from $V$ to $V'$. In the rest of 
this work, we identify $V'$ with $V$ using this isomorphism without
mentioning it explicitly.

For any triple $V$-modules $W_1, W_2, W_3$, we have isomorphisms 
$$
\Omega_r: \V_{W_1W_2}^{W_3} \rightarrow
\V_{W_2W_1}^{W_3}, \quad \quad \forall r\in \Z
$$
given as follow:
\beq   \label{Omega-r}
\Omega_r(\Y)(w_2, z)w_1 = e^{zL(-1)} \Y(w_1, e^{(2r+1)\pi i}z)w_2,
\eeq
for $\Y \in \V_{W_1W_2}^{W_3}$ and $w_i\in W_i, i=1,2$. 
The following identity
\beq
\Omega_r \circ \Omega_{-r-1} = \Omega_{-r-1} \circ \Omega_r = \id
\eeq
is proved in \cite{HL3}. 

For $\Y \in \V_{W_1W_2}^{W_3}$ and $r\in \Z$,  
a so-called {\em $r$-contragredient operator} $A_r(\Y)$ 
was introduced in \cite{HL3}. 
Here, we use two slightly different operators $\tilde{A}_r(\Y)$
and $\hat{A}_r(\Y)$ introduced in \cite{Ko1} and defined as follow:
\bea  \label{two-new-A-2}
&&\langle \tilde{A}_{r}(\Y)(w_{1}, e^{(2r+1)\pi i}x)w'_{3},w_{2}\rangle 
=\langle w'_{3}, \Y(e^{xL(1)}x^{-2L(0)}w_{1},  x^{-1})w_{2}\rangle,  
\nn
&&\langle \hat{A}_{r}(\Y)(w_{1}, x)w'_{3},w_{2}\rangle 
=\langle w'_{3}, \Y(e^{-xL(1)}x^{-2L(0)}w_{1},  
e^{(2r+1)\pi i}x^{-1})w_{2}\rangle,
\eea
for $\Y \in \V_{W_1W_2}^{W_3}$ and $w_1\in W_1, w_2\in W_2$, 
$w'_3\in W'_3$. 
In particular, when $W_{1}=V$ and $W_{2}=W_{3}=W$, we have
$\tilde{A}_{r}(Y_W)=\hat{A}_r(Y_W)=Y'_W=Y_{W'}, \forall r\in {\Z}$. 
If $\Y \in \V_{W_1W_2}^{W_3}$, 
then $\tilde{A}_r(\Y)$, $\hat{A}_r(\Y) \in \V_{W_1W'_3}^{W'_2}$ 
for $r\in \Z$ and
\begin{equation} \label{inverse-two-A}
\tilde{A}_r \circ \hat{A}_r (\Y) = 
\hat{A}_r \circ \tilde{A}_r (\Y) = \Y.
\end{equation}

Let $\Y \in \V_{a_1a_2}^{a_3}$. We define 
$\sigma_{123} := \Omega_r \circ \tilde{A}_r$.
It is easy to see that 
\bea   \label{omega-r-A-r}
&&\langle w'_{a_3}, \Y(w_{a_1}, x)w_{a_2}\rangle   \nn
&&\hspace{1cm}=\langle e^{-x^{-1}L(-1)} \sigma_{123}(\Y)
(w'_{a_3}, x^{-1}) e^{-xL(1)}x^{-2L(0)} w_{a_1}, w_{a_2}\rangle
\eea
for $w_{a_1}\in W_{a_1}, w_{a_2}\in W_{a_2}, w'_{a_3}\in W'_{a_3}$. 
It is also clear that $\sigma_{123}$ is independent
of $r\in \Z$. It is proved in \cite{Ko1} that 
$\sigma_{123}^3 = \id_{\V_{a_1a_2}^{a_3}}$.
We also denote $\sigma_{123}^{-1}$ as $\sigma_{132}$. 
Clearly, we have $\sigma_{132}= \hat{A}_r \circ \Omega_{-r-1}$ and 
\bea  \label{A-r-omega-r}
&&\langle \sigma_{132}(\Y)(w_1, x)w'_3,  w_2  \rangle  \nn
&&\hspace{1cm}= \langle w'_3, e^{-x^{-1}L(-1)} \Y(w_2, x^{-1}) 
e^{-xL(1)}x^{-2L(0)} w_1 \rangle
\eea
for $w_{a_1}\in W_{a_1}, w_{a_2}\in W_{a_2}, w'_{a_3}\in W'_{a_3}$.

For any $V$-module $W$, we define 
a $V$-module map $\theta_W: W\rightarrow W$ by 
\beq  \label{theta-def}
\theta_W: w \mapsto e^{-2\pi iL(0)}w.
\eeq 
For $W_a$, 
we have $\theta_{W_a} = e^{-2\pi i h_a}\id_{W_a}$
where $h_a\in \C$ is the lowest conformal weight of $W_a$.

We denote the graded dual space of $V_{cl}$ and $V_{op}$
by $V_{cl}'$ and $V_{op}'$ respectively.  
Let $\varphi_{cl}: V_{cl} \rightarrow V_{cl}'$ and 
$\varphi_{op}: V_{op} \rightarrow V_{op}'$ be the isomorphisms 
induced from $(\cdot, \cdot)_{cl}$ and $(\cdot, \cdot)_{op}$
respectively. Namely, we have
\bea
(u_1, u_2)_{cl}  &=& \langle \varphi_{cl}(u_1), u_2\rangle   \nn
(v_1, v_2)_{op}  &=& \langle \varphi_{op}(v_1), v_2\rangle
\label{varphi-form-op}
\eea
for $u_1,u_2\in V_{cl}$ and $v_1,v_2\in V_{op}$. 

$V_{cl}$ as a conformal full field algebra over $V\otimes V$ can be 
expanded as follow:
\beq
V_{cl} = \oplus_{i=1}^{N_{cl}} W_{r_L(i)} \otimes W_{r_R(i)} 
\eeq
where $r_L, r_R: \{ 1, \dots, N_{cl} \} \rightarrow \I$. 
For $a\in \I$, 
we choose a basis $\{ e_{a;\alpha} \}_{\alpha \in \N}$ of $W_a$
and a dual basis $\{ e'_{a;\alpha} \}_{\alpha \in \N}$ of $W'_a$. 
Then 
\beq
\{ e_{r_L(i), \alpha} \otimes  e_{r_R(i), \beta} 
\}_{i=1,\dots, N_{cl}, \alpha, \beta\in \N}
\eeq
is a basis of $V_{cl}$  and 
\beq
\{ \varphi_{cl}^{-1}( e'_{r_L(i);\alpha} \otimes e'_{r_R(i), \beta}) 
\}_{i=1,\dots, N_{cl}, \alpha, \beta\in \N} 
\eeq
is its dual basis with respect to the nondegenerate bilinear form
$(\cdot, \cdot)_{cl}$.

Let $T: \mathcal{C}_{V\otimes V} \rightarrow \mathcal{C}_V$ be the
tensor bifunctor. 
We showed in \cite{Ko2} that there is a morphism 
$\iota_{cl-op}: T(V_{cl}) \rightarrow V_{op}$ in $\mathcal{C}_V$
(see (3.81),(3.82) in \cite{Ko2} for definition). 
We define a morphism $\iota'_{cl-op}: T(V'_{cl}) \rightarrow V_{op}$ 
as a composition of maps as follow:
\beq
\iota'_{cl-op}:  T(V'_{cl}) \xrightarrow{T(\varphi_{cl}^{-1})} 
T(V_{cl}) \xrightarrow{\iota_{cl-op}} V_{op}\, .
\eeq
By the universal property of tensor product $\boxtimes$ 
(\cite{HL1}-\cite{HL4}), 
$\V_{W_1W_2}^{W_3}$ and $\hom_V( W_1\boxtimes W_2, W_3)$ for
any three $V$-modules $W_1,W_2,W_3$ are canonically isomorphic. 
Given a morphism 
$m\in \hom_V( W_1\boxtimes W_2, W_3)$, we denote the 
corresponding intertwining operator as $\Y_m$. 
Conversely, given an intertwining operator $\Y$, we 
denote its corresponding morphism as $m_{\Y}$. 
Therefore, we have two intertwining operators $\Y_{\iota_{cl-op}}$ and 
$\Y_{\iota'_{cl-op}}$ corresponding to morphisms $\iota_{cl-op}$ and 
$\iota'_{cl-op}$ respectively. 

\begin{lemma} 
For $z\in \HH$, we have 
\beq  \label{iota-z-bar-z-1}
\iota_{cl-op}(z,\bar{z}) (e_{r_L(i), \alpha} \otimes  e_{r_R(i), \beta})
= e^{\bar{z} L(-1)} 
\Y_{\iota_{cl-op}}(e_{r_L(i), \alpha}, z-\bar{z} ) e_{r_R(i),\alpha}. 
\eeq
\end{lemma}
\pf
It is proved in \cite{Ko2} that 
\beq  \label{nu-cl-op}
m_{\mathbb{Y}_{cl-op}} = m_{Y_{op}^f} \circ (\iota_{cl-op} \boxtimes \id_{V_{op}}).
\eeq
Using (\ref{nu-cl-op}), 
when $z\in \HH, \zeta\in \overline{\HH}$ 
and $|\zeta|>|z-\zeta|>0$, we have 
\bea  \label{lemma-equ-3}
\iota_{cl-op}(z,\zeta) (e_{r_L(i), \alpha} \otimes  e_{r_R(i), \beta}) &=&
\mathbb{Y}_{cl-op}( e_{r_L(i), \alpha} \otimes  e_{r_R(i), \beta}; 
z,\zeta)\one_{op}   \nn
&=&Y_{op}^f( \Y_{\iota_{cl-op}}( e_{r_L(i), \alpha}, z-\zeta)e_{r_R(i), \beta},
\zeta)\one   \nn
&=&e^{\zeta L(-1)} 
\Y_{\iota_{cl-op}}( e_{r_L(i), \alpha}, z-\zeta) e_{r_R(i), \beta}.
\eea
By the convergence property of 
the iterate of two intertwining operators, 
the right hand side of (\ref{lemma-equ-3})
is a power series of $\zeta$ absolutely convergent 
for $|\zeta|>|z-\zeta|>0$ . 
By the property of power series, 
the right hand side of (\ref{lemma-equ-3})
must converge absolutely for all 
$z\in \HH, \zeta\in \overline{\HH}$. 
Because analytic extension in a simply
connected domain is unique, we obtain that the equality 
(\ref{lemma-equ-3}) holds for all $z\in \HH, \zeta\in \overline{\HH}$.
When $\zeta=\bar{z}$, we obtain (\ref{iota-z-bar-z-1}). 
\epf

Now we consider the both sides of Cardy condition (\ref{cardy})
for an open-closed field algebra over $V$. 
On the left hand side of (\ref{cardy}), we have 
$q_{s_2}<0$. Using (\ref{Y-op-opposite}) 
and (\ref{Omega-r}), we obtain, $\forall v_3\in V_{op}$,
\bea
Y_{op}(\mathcal{U}(q_{s_2})v_2, q_{s_2})v_3 &=&
e^{-|q_{s_2}| L(-1)} Y_{op}(v_3, |q_{s_2}|) \,\, \mathcal{U}(q_{s_2})v_2  \nn
&=& 
\Omega_{-1}(Y_{op}^f)(\mathcal{U}(q_{s_2})v_2 , e^{\pi i} |q_{s_2}|) v_3.
\eea
Hence we can rewrite the left hand side of (\ref{cardy})
as follow:
\beq   \label{cardy-L-2}
\tr_{V_{op}} 
\left( Y_{op}^f(\mathcal{U}(q_{s_1}) v_1, q_{s_1}) 
\Omega_{-1}(Y_{op}^f)(\mathcal{U}(q_{s_2}) e^{-2\pi iL(0)} v_2, 
e^{\pi i} |q_{s_2}|) q_{\tau}^{L(0)-c/24} \right)
\eeq
for $q_{s_1}>|q_{s_2}|>|q_{\tau}|>0$.

We have the following result for the right hand side of 
(\ref{cardy}). 
\begin{prop}
For $s_1, s_2\in \HH, \text{\rm Re}\, s_1=0, \text{\rm Re}\, s_2=0$, 
\bea
&&\hspace{-1cm}( (T_1^L \otimes T_1^R)^{*} \, \iota_{cl-op}^*(z_1, \bar{z}_1) ( T_2 v_1), 
\,\, q_{-\frac{1}{\tau}}^{-\frac{c}{24}}
(T_3^L \otimes T_3^R)^* \, \iota_{cl-op}^*(z_2, \bar{z}_2) 
(T_4 v_2))_{cl}   \nn
&&\hspace{-0.5cm} = \sum_{i=1}^{N_{cl}}  \tr_{W_{r_R(i)}} \,\,
e^{-2\pi iL(0)} \,\, \Y_1 \left( \mathcal{U}(q_{-\frac{1}{\tau}s_1}) 
\left( -\frac{1}{\tau} \right)^{L(0)} v_1, 
q_{-\frac{1}{\tau}s_1} \right)  \nn
&&\hspace{4cm}\Y_2 \left( \mathcal{U}(q_{-\frac{1}{\tau}s_2}) 
\left( -\frac{1}{\tau} \right)^{L(0)} v_2, 
q_{-\frac{1}{\tau}s_2} \right) q_{-\frac{1}{\tau}}^{L(0)-\frac{c}{24}} 
\label{cardy-R-2}
\eea
where $\Y_1$ and $\Y_2$ are intertwining operators of types
$\binom{W_{r_R(i)}}{V_{op}W'_{r_L(i)}}$ and 
$\binom{W'_{r_L(i)}}{V_{op}W_{r_R(i)}}$ respectively and are given by 
\bea
\Y_1 &=& \sigma_{123}(\Y_{\iota'_{cl-op}}) \circ (\varphi_{op} 
\otimes \id_{W'_{r_L(i)}} ),  \label{Y-1} \\
\Y_2 &=& \Omega_0 (\sigma_{132}(\Y_{\iota_{cl-op}})) 
\circ (\varphi_{op} \otimes \id_{W_{r_R(i)}} ). \label{Y-2}
\eea
\end{prop}
\pf
Let $z_3:=z_1-\bar{z}_1$ and $z_4:=z_2-\bar{z}_2$. 
By (\ref{iota-up-*}) and (\ref{iota-z-bar-z-1}), 
the left hand side of (\ref{cardy-R-2}) equals to
\bea
&&\hspace{-1cm}
( (T_1^L \otimes T_1^R)^{*} \, \iota_{cl-op}^*(z_1, \bar{z}_1) ( T_2 v_1), 
\,\, 
(T_3^L \otimes T_3^R)^* \, \iota_{cl-op}^*(z_2, \bar{z}_2)(T_4 v_2))_{cl} 
\nn
&&\hspace{-0.5cm}=
\sum_{i=1}^{N_{cl}} \sum_{\alpha, \beta}
\left( (T_3^L \otimes T_3^R)^* \, 
\iota_{cl-op}^*(z_2, \bar{z}_2) (T_4 v_2),  
e_{r_L(i), \alpha} \otimes  e_{r_R(i), \beta}\right)_{cl}  \nn
&&\hspace{2cm} \left(
\varphi_{cl}^{-1}(e'_{r_L(i), \alpha} \otimes  e'_{r_R(i), \beta}), 
(T_1^L \otimes T_1^R)^{*} \, \iota_{cl-op}^*(z_1, \bar{z}_1) ( T_2 v_1)
\right)_{cl}  \nn
&&\hspace{-0.5cm} = \sum_{i=1}^{N_{cl}}\sum_{\alpha, \beta}
\left( T_4v_2, \iota_{cl-op}(z_2,\bar{z}_2)
(T_3^L e_{r_L(i), \alpha} \otimes T_3^R e_{r_R(i), \beta})  \right)_{op}
\nn
&&\hspace{2cm} \left( \iota_{cl-op}(z_1, \bar{z}_1)(\varphi_{cl}^{-1} 
(T_1^L e'_{r_L(i), \alpha} \otimes  T_1^R e'_{r_R(i), \beta})), 
T_2v_1 \right)_{op}  \nn
&&\hspace{-0.5cm} = \sum_{i=1}^{N_{cl}} \sum_{\alpha, \beta}
\left( T_4 v_2\, ,  \,\, 
e^{\bar{z}_2L(-1)} \Y_{\iota_{cl-op}}( T_3^L e_{r_L(i),\alpha}, 
z_2-\bar{z}_2) T_3^R e_{r_R(i), \beta} \right)_{op}  \nn
&&\hspace{2cm} \left( e^{\bar{z}_1L(-1)} 
\Y_{\iota'_{cl-op}}(T_1^L e'_{r_L(i), \alpha}, z_1-\bar{z}_1) 
T_1^R e'_{r_R(i), \beta}, \,\, \, T_2 v_1 \right)_{op}  \nn
&&\hspace{-0.5cm} = \sum_{i=1}^{N_{cl}} \sum_{\alpha, \beta}
\,\, \langle \, \varphi_{op}(T_4 v_2), \,\,
e^{\bar{z}_2L(-1)} \Y_{\iota_{cl-op}}( T_3^L e_{r_L(i),\alpha}, 
z_2-\bar{z}_2) T_3^R e_{r_R(i), \beta}  \,\rangle \nn
&&\hspace{2cm} \langle \, e^{\bar{z}_1L(-1)} 
\Y_{\iota'_{cl-op}}(T_1^L e'_{r_L(i), \alpha}, \,\, z_1-\bar{z}_1) 
T_1^R e'_{r_R(i), \beta}, \,\, \, \varphi_{op}(T_2 v_1) \,\rangle \nn
&&\hspace{-0.5cm} = \sum_{i=1}^{N_{cl}} \sum_{\alpha, \beta}
\,\, \langle \, 
e^{-z_4^{-1}L(-1)} \sigma_{123}(\Y_{\iota_{cl-op}})
(\varphi_{op}(e^{-\bar{z}_2L(1)} T_4 v_2),  z_4^{-1}) \nn
&&\hspace{4cm} z_4^{-2L(0)} e^{-z_4^{-1}L(1)} T_3^L e_{r_L(i),\alpha}, \,\,
T_3^R e_{r_R(i), \beta}  \,\rangle  \nn
&&\hspace{2cm}
\langle \, T_1^R e'_{r_R(i), \beta}, \,\, e^{-z_3^{-1}L(-1)}
\sigma_{123}(\Y_{\iota'_{cl-op}})(\varphi_{op}(e^{-\bar{z}_1L(1)} T_2 v_1), 
z_3^{-1}) \nn
&&\hspace{4cm} z_3^{-2L(0)} e^{-z_3^{-1} L(1)} T_1^L e'_{r_L(i), \alpha}
\rangle \nn
&&\hspace{-0.5cm} = \sum_{i=1}^{N_{cl}} \sum_{\alpha, \beta}
\langle \, \sigma_{123}(\Y_{\iota_{cl-op}})
(\varphi_{op}(e^{-\bar{z}_2L(1)} T_4 v_2),  z_4^{-1}) 
z_4^{-2L(0)} e^{-z_4^{-1}L(1)} T_3^L e_{r_L(i),\alpha}, \nn
&&\hspace{3cm} 
e^{z_4^{-1}L(1)} T_3^R e_{r_R(i),\beta}\rangle \, 
\langle e'_{r_R(i), \beta},  \,\, (T_1^R)^* e^{-z_3^{-1}L(-1)} \cdot \nn
&&\hspace{1cm}
\sigma_{123}(\Y_{\iota'_{cl-op}})(\varphi_{op}(e^{-\bar{z}_1L(1)} T_2 v_1),
z_3^{-1}) z_3^{-2L(0)} e^{-z_3^{-1} L(1)} T_1^L e'_{r_L(i), \alpha}
\rangle \nn
&&\hspace{-0.5cm} = \sum_{i=1}^{N_{cl}} \sum_{\alpha,\beta}
\langle \, e_{r_L(i),\alpha},  (T_3^L)^* e^{z_4^{-1}L(-1)}z_4^{-2L(0)}
\tilde{A}_0 \circ \sigma_{123} (\Y_{\iota_{cl-op}}) 
(e^{-z_4^{-1}L(1)} z_4^{2L(0)} \cdot \nn 
&&\hspace{1cm} \varphi_{op}(e^{-\bar{z}_2L(1)} T_4 v_2), e^{\pi i}z_4)
e^{z_4^{-1}L(1)} T_3^R e_{r_R(i),\beta}\rangle \, 
\langle e'_{r_R(i), \beta},  \,\,
(T_1^R)^* e^{-z_3^{-1}L(-1)} \cdot \nn
&&\hspace{1cm}
\sigma_{123}(\Y_{\iota'_{cl-op}})(\varphi_{op}(e^{-\bar{z}_1L(1)} T_2 v_1),
z_3^{-1}) z_3^{-2L(0)} e^{-z_3^{-1} L(1)} T_1^L e'_{r_L(i), \alpha}\rangle
\label{cardy-R-pf-equ-1}
\eea
We define two intertwining operators as follow:
\bea
\Y_1^{(0)} &=& \sigma_{123}(\Y_{\iota'_{cl-op}}) 
\circ (\varphi_{op} \otimes \id_{W'_{r_L(i)}}),  \nn
\Y_2^{(0)} &=&
\tilde{A}_0(\sigma_{123}(\Y_{\iota_{cl-op}})) \circ 
(\varphi_{op}\otimes \id_{W'_{r_R(i)}}) .
\eea
Using (\ref{T-1-L-R}),(\ref{T-2}),(\ref{T-3-L-R}) and
(\ref{T-4}), we further obtain that 
the left hand side of (\ref{cardy-R-pf-equ-1}) equals to
\bea
&&\hspace{-1cm} \sum_{i=1}^{N_{cl}} \sum_{\alpha,\beta}
\langle \, e_{r_L(i),\alpha}, \left( \frac{p_{s_2}}{r}\right)^{L(0)}
(\bar{z}_2-z_2)^{L(0)} z_4^{-2L(0)} \cdot \nn
&&\hspace{1.5cm} 
\Y_2^{(0)}(z_4^{2L(0)}(\bar{z}_2-z_2)^{-L(0)}\mathcal{U}(1)
\left( -\frac{1}{\tau} \right)^{L(0)} v_2, e^{\pi i}z_4)  \nn
&&\hspace{1.5cm} 
\overline{ \bar{z}_2 - z_2 }^{L(0)} \overline{\left( 
\frac{p_{s_2}}{r}\right)}^{\, L(0)}
e_{r_R(i),\beta} \rangle  \langle e'_{r_R(i), \beta}, \,\,
\overline{ \left( \frac{p_{s_1}}{r} \right)}^{\, -L(0)}
\overline{ (z_1-\bar{z}_1) }^{L(0)} \nn
&&\hspace{1.5cm} \Y_2^{(0)}(z_3^{-L(0)} \mathcal{U}(1)
\left( -\frac{1}{\tau} \right)^{L(0)}v_1, z_3^{-1})
z_3^{-L(0)}\left( \frac{p_{s_1}}{r} \right)^{-L(0)} e'_{r_L(i),\alpha}\rangle
\nn
&&\hspace{-1cm}= \sum_{i=1}^{N_{cl}} \tr_{W'_{r_L(i)}}
\Y_2^{(0)}(\mathcal{U}(p_{s_2})(-1/\tau)^{L(0)} v_2, E_1)
E_2 \, \Y_1^{(0)}(\mathcal{U}(p_{s_1})(-1/\tau)^{L(0)}v_1, p_{s_1})
\eea
where $E_1=p_{s_2}(\bar{z}_2-z_2) z_4^{-2} e^{\pi i} z_4$ and 
$$
E_2 =p_{s_2}^{L(0)} (\bar{z}_2 -z_2)^{L(0)} z_4^{-2L(0)} 
\overline{ (\bar{z}_2-z_2)}^{L(0)} \overline{p_{s_2}}^{L(0)} 
\overline{p_{s_1}}^{-L(0)} 
\overline{ z_1-\bar{z}_1 }^{L(0)} z_3^{-L(0)} p_{s_1}^{-L(0)} .
$$
For $E_1$, since $z_4^{-1}$ is 
obtained by operations on intertwining operator where
$z_4$ is treated formally, $z_4^{-1}$
really means $|z_4|^{-1} e^{-i\frac{\pi i}{2}}$.  Therefore, we have
\beq
E_1 = p_{s_2} e^{\frac{3\pi i}{2}} e^{-\pi i} e^{\pi i} 
e^{\frac{\pi i}{2}}=  e^{2\pi i} p_{s_2}. 
\eeq
For $E_2$, keep in mind (\ref{branch-cut}) and (\ref{power-f}), 
we have 
\bea
E_2&=&  \lbar \frac{p_{s_2}}{p_{s_1}} \lbar^{2L(0)} 
(\bar{z}_2 - z_2)^{L(0)} (z_2-\bar{z}_2)^{-2L(0)} 
\overline{(\bar{z}_2-z_2)}^{\, L(0)} 
\overline{(z_1-\bar{z}_1)}^{\, L(0)}(z_1-\bar{z}_1)^{-L(0)}\nn
&=& q_{-\frac{1}{\tau}}^{L(0)}\,\, e^{\frac{3\pi i}{2}L(0)} e^{-2\frac{\pi i}{2}L(0)}
e^{-\frac{3\pi i}{2}L(0)} e^{-\frac{\pi i}{2}L(0)} e^{-\frac{\pi i}{2}L(0)}\nn
&=& q_{-\frac{1}{\tau}}^{L(0)} \,\, e^{-2\pi iL(0)}.
\eea
Therefore, we obtain that the left hand 
side of (\ref{cardy-R-pf-equ-1}) further equals to 
\bea   \label{prop-equ-last}
&&\sum_{i=1}^{N_{cl}} \tr_{W'_{r_L(i)}} 
\Omega_0^2(\Y_2^{(0)})\left( \mathcal{U}(p_{s_2}) 
\left( -\frac{1}{\tau}\right)^{L(0)} v_2 , p_{s_2} \right) \nn
&&\hspace{3cm}
q_{-\frac{1}{\tau}}^{L(0)} \, e^{-2\pi iL(0)} \,\, 
\Y_1^{(0)} \left( \mathcal{U}(p_{s_1})
\left( -\frac{1}{\tau}\right)^{L(0)} v_1, p_{s_1} \right),
\eea
where we have used the fact that $\Y(\cdot, e^{2\pi i} x) \cdot =
\Omega_0^2(\Y) (\cdot, x)\cdot$ for any intertwining operator 
$\Y$. By using the property of trace, 
it is easy to see that (\ref{prop-equ-last}) multiplying 
$q_{-\frac{1}{\tau}}^{-\frac{c}{24}}$
is nothing but the right hand side of (\ref{cardy-R-2}).
\epf

\begin{rema}  \label{cardy-R-conv-rema}
{\rm 
It is easy to check that the absolute convergence 
of the left hand side of (\ref{cardy-R-pf-equ-1}) 
by our assumption easily implies the absolute convergence of 
each step in (\ref{cardy-R-pf-equ-1}). 
Notice that the absolute convergence of the right side of 
(\ref{cardy-R-2}) is automatic because $V$ is
assumed to satisfy the conditions in Theorem \ref{ioa}. 
Hence, by tracing back the steps in above proof, we see
that the absolute convergence of the left hand side of 
(\ref{cardy-R-2}) is also automatic. 
}
\end{rema}

Now we recall some results in \cite{H10}\cite{H11}. 
We follow the notations in \cite{HKo3}. 
We denote the unique analytic extension of 
$$
\tr_{W_{a_1}} \Y_{aa_1;i}^{a_1;(1)}(\mathcal{U}(e^{2\pi iz}) w_a , 
e^{2\pi i z}) q_{\tau}^{L(0)-\frac{c}{24}}
$$
in the universal covering space of $1>|q_{\tau}|>0$ as
$$
E(\tr_{W_{a_1}} \Y_{aa_1;i}^{a_1;(1)}(\mathcal{U}(e^{2\pi iz}) w_a , 
e^{2\pi i z}) q_{\tau}^{L(0)-\frac{c}{24}}).
$$
By \cite{Mi2}\cite{H10}, above formula is independent of $z$.
Consider the map: for $w_a\in W_a$,
\beq   \label{Y-z-tau-map}
\Psi_1(\Y_{aa_1;i}^{a_1;(1)}):
w_a \rightarrow E(\tr_{W_{a_1}} 
\Y_{aa_1;i}^{a_1;(1)}(\mathcal{U}(e^{2\pi iz}) w_a , 
e^{2\pi i z}) q_{\tau}^{L(0)-\frac{c}{24}}).
\eeq
We denote the right hand side of (\ref{Y-z-tau-map}) as
$\Psi_1(\Y_{aa_1;i}^{a_1;(1)})(w_a; z, \tau)$. Notice that we
choose to add $z$ in the notation 
even though it is independent of $z$.  
We define an action of $SL(2,\Z)$ on the map 
(\ref{Y-z-tau-map}) as follow: 
\bea
&&\left( \left( \begin{array}{cc} a & b \\ c & d \end{array}\right)
(\Psi_1(\Y_{aa_{1}; i}^{a_{1}; (1)})) \right) (w_a; z, \tau)   \\
&&\hspace{1cm}= 
E\left( \tr_{W_{a_1}} \Y_{aa_1;i}^{a_1;(1)}\left(  
\mathcal{U}(e^{2\pi i z' }) 
\left( \frac{1}{c\tau +d} \right)^{L(0)} w_a, 
e^{2\pi i z'} \right) q_{\tau'}^{L(0)-\frac{c}{24}} \right)
\nonumber
\eea
where $\tau'=\frac{a\tau + b}{c\tau +d}$ and $z'= \frac{z}{c\tau+d}$,
for $\left( \begin{array}{cc} a & b \\ c & d \end{array}\right)
\in SL(2,\Z)$ and $w_a\in W_a$. 
The following Theorem is proved in \cite{Mi2}\cite{H8}. 
\begin{thm}  \label{modular-ioa}
There exists a
unique $A_{a_2a_3}^{ij}\in \C$ for $a_2, a_3\in \I$ such that,
for $w_a\in W_a$, 
\bea
&&\hspace{-1cm}E\left( \tr_{W_{a_1}} \Y_{aa_1;i}^{a_1;(1)}\left(  
\mathcal{U}(e^{2\pi i z' }) 
\left( \frac{1}{c\tau +d} \right)^{L(0)} w_a, 
e^{2\pi i z'} \right) q_{\tau'}^{L(0)-\frac{c}{24}} \right)\nn
&&\hspace{2cm} = \sum_{a_3\in \I} A_{a_1a_2}^{ij} 
E(\tr_{W_{a_1}} \Y_{aa_2;j}^{a_2;(2)}(\mathcal{U}(e^{2\pi iz})w_a, e^{2\pi iz}) 
q_{\tau}^{L(0)-\frac{c}{24}})
\eea
where $\tau'=\frac{a\tau + b}{c\tau +d}$ and $z'= \frac{z}{c\tau+d}$. 
\end{thm}

In particular, the action of 
$S=\left( \begin{array}{cc} 0 & 1 \\ -1 & 0 \end{array}\right)$ 
on (\ref{Y-z-tau-map}) induces, for each $a\in \I$, 
an automorphism on 
$\oplus_{a_1\in \I} \V_{aa_1}^{a_1}$, denoted as $S(a)$. Namely, we have 
\beq  \label{S-act-Y-thm}
S(\Psi_1(\Y_{aa_1; i}^{a_1;(1)}))
= \Psi_1 (S(a)(\Y_{aa_1; i}^{a_1;(1)}))
\eeq
Combining all such $S(a)$, we obtain an automorphism on 
$\oplus_{a,a_1\in \I} \V_{aa_1}^{a_1}$. We
still denote it as $S$, i.e. $S=\oplus_{a\in \I} S(a)$. 
Then $S$ can be further extended to a map on 
$\oplus_{a, a_3\in \I} \V_{aa_3}^{a_3} \otimes \V_{a_1a_2}^{a}$ given
as follow: 
\beq
S(\Y_{aa_3;i}^{a_3;(1)}\otimes \Y_{a_1a_2;j}^{a;(2)}) :=
S(\Y_{aa_3;i}^{a_3;(1)})\otimes \Y_{a_1a_2;j}^{a;(2)}. 
\eeq

There is a fusing isomorphism map \cite{H7}: 
\beq   \label{F-iso}
\mathcal{F}:  
\oplus_{a\in \I} \mathcal{V}_{a_1a}^{a_4}\otimes \mathcal{V}_{a_2a_3}^{a}
\xrightarrow{\cong} 
\oplus_{b\in \I} \mathcal{V}_{ba_3}^{a_4}\otimes \mathcal{V}_{a_1a_2}^{b}
\eeq
for $a_1,a_2,a_3,a_4\in \I$. Using the isomorphism $\mathcal{F}$, 
we obtain a natural action of $S$ 
on $\oplus_{b,a_3} \V_{a_1b}^{a_3}\otimes \V_{a_2a_3}^{b}$.

It is shown in \cite{H8}
that the following 2-points genus-one correlation function,
for $a_1, a_2,a_3,a_4\in \I$, $i=1, \dots, N_{aa_1}^{a_1}$,
$j=1, \dots, N_{a_2a_3}^{a}$ and $w_{a_k}\in W_{a_k}, k=2,3$,
\beq  \label{q-trace}
\tr_{W_{a_{1}}}
\Y_{aa_{1}; i}^{a_{1}; (1)}(\mathcal{U}(q_{z_2})
\Y_{a_{2}a_{3}; j}^{a; (2)}
(w_{a_{2}}, z_{1}-z_{2})w_{a_{3}}, q_{z_2} )
q_{\tau}^{L(0)-\frac{c}{24}}
\eeq
is absolutely convergent when $1>|e^{2\pi iz_2}|>|q_{\tau}|>0$
and $1>|e^{2\pi i(z_1-z_2)}|>0$ and single-valued in the chosen 
branch. It can be extended uniquely to 
a single-valued analytic function on the universal covering space
of 
$$
M_1^2 = \{ (z_1,z_2,\tau)\in \C^{3} |
z_1\neq z_2+p\tau + q, \forall p,q\in Z, \tau\in \HH \}.
$$
This universal covering space is denoted by $\tilde{M}_1^2$. 
We denote this single-valued analytic function on $\tilde{M}_1^2$
as 
$$
E(\tr_{W_{a_1}}
\Y_{aa_{1}; i}^{a_{1}; (1)}(\mathcal{U}(q_{z_2})\Y_{a_{2}a_{3}; j}^{a; (2)}
(w_{a_{2}}, z_{1}-z_{2})w_{a_{3}}, q_{z_2})
q_{\tau}^{L(0)-\frac{c}{24}}).
$$
We denote the space spanned by such functions on $\tilde{M}_1^2$
by $\mathbb{G}_{1; 2}$.

For $\Y_{aa_1;i}^{a_1;(1)} \in \V_{aa_1}^{a_1}$ and 
$\Y_{a_2a_3;j}^{a;(2)} \in \V_{a_2a_3}^{a}$, 
we now define the following linear map:
\begin{eqnarray*}
\Psi_2(\Y_{aa_1; i}^{a_1;(1)}\otimes \Y_{a_2a_3;j}^{a;(2)}): \quad
\oplus_{b_2,b_3\in \I}W_{b_2}\otimes W_{b_3} &\to &  \mathbb{G}_{1; 2}
\end{eqnarray*}
as follow: the map restricted on $W_{b_2}\otimes W_{b_3}$ 
is define by $0$ for $b_2\neq a_2$ or $b_3\neq a_3$, and by 
\beq  \label{ffa-cor-1-1}
E(\tr_{W_{a_{1}}}
\Y_{aa_{1}; i}^{a_{1}; (1)}(\mathcal{U}(q_{z_2})
\Y_{a_{2}a_{3}; j}^{a; (2)}
(w_{a_{2}}, z_{1}-z_{2})w_{a_{3}}, q_{z_2})
q_{\tau}^{L(0)-\frac{c}{24}}),
\eeq
for all $w_{a_k}\in W_{a_k}, k=2,3$. 
The following identity was proved in \cite{H10}.
\bea  \label{S-act-right-0}
&& \left( \Psi_2(\Y_{aa_1; i}^{a_1;(1)}\otimes \Y_{a_2a_3;j}^{a;(2)})
\left( \left( -\frac{1}{\tau}\right)^{L(0)}w_{a_2}\otimes 
\left( -\frac{1}{\tau}\right)^{L(0)}w_{a_3} \right) \right) 
\left( -\frac{1}{\tau}z_{1}, -\frac{1}{\tau}z_{2};
-\frac{1}{\tau} \right)  \nn
&&\hspace{2cm}= 
\big( \Psi_2( S(\Y_{aa_1; i}^{a_1;(1)}\otimes \Y_{a_2a_3;j}^{a;(2)}))
(w_{a_2}\otimes w_{a_3})\big) (z_1, z_2, \tau). 
\eea

One can also produce 2-point genus-one correlation functions from
a product of two intertwining operators. It is proved in \cite{H8} that
$\forall w_{a_k}\in W_{a_k}, k=1,2$, 
\beq  \label{q-trace-2}
\tr_{W_{a_4}} \Y_{a_1a_3;i}^{a_4;(1)}(\mathcal{U}(q_{z_1}) w_{a_1}, q_{z_1})
\Y_{a_2a_4;j}^{a_3;(2)}(\mathcal{U}(q_{z_2}) w_{a_2}, q_{z_2}) 
q_{\tau}^{L(0)- \frac{c}{24}},
\eeq
is absolutely convergent when $1>|q_{z_1}|>|q_{z_2}|>|q_{\tau}|>0$. 
(\ref{q-trace-2}) has a unique extension to the universal 
covering space $\tilde{M}_1^2$, denoted as 
\beq  \label{E-q-trace-2}
E(\tr_{W_{a_4}} \Y_{a_1a_3;i}^{a_4;(1)}(\mathcal{U}(q_{z_1}) w_{a_1}, q_{z_1})
\Y_{a_2a_4;j}^{a_3;(2)}(\mathcal{U}(q_{z_2}) w_{a_2}, q_{z_2}) 
q_{\tau}^{L(0)- \frac{c}{24}}).
\eeq
Such functions on $\tilde{M}_1^2$ also span $\mathbb{G}_{1;2}$. 
We define a map
\beq
\Psi_2(\Y_{a_1a_3;i}^{a_4;(1)}\otimes \Y_{a_2a_4;j}^{a_3;(2)}):\,\,
\oplus_{b_1,b_2\in \I}W_{b_1}\otimes W_{b_2} \to   \mathbb{G}_{1;2}
\eeq
as follow: the map restrict on $W_{b_1}\otimes W_{b_2}$
is defined by $0$ for $b_1\neq a_1, b_2\neq a_2$, and by 
(\ref{E-q-trace-2}) for $w_{a_1}\in W_{a_1}, w_{a_2}\in W_{a_2}$.

It was proved by Huang in \cite{H8} that the fusing isomorphism
(\ref{F-iso}) gives the following associativity: 
\bea  \label{asso-ioa-g-1}
&&E(\tr_{W_{a_4}} \Y_{a_1a_3;i}^{a_4;(1)}(\mathcal{U}(q_{z_1}) w_{a_1}, q_{z_1})
\Y_{a_2a_4;j}^{a_3;(2)}(\mathcal{U}(q_{z_2}) w_{a_2}, q_{z_2}) 
q_{\tau}^{L(0)- \frac{c}{24}} ) \nn
&&\hspace{1cm}=
\sum_{a_5\in \I} \sum_{k,l} 
F(\Y_{a_1a_3;i}^{a_4;(1)}\otimes \Y_{a_2a_4;j}^{a_3;(2)}, 
\Y_{a_5a_4;k}^{a_4;(3)}\otimes \Y_{a_1a_2;l}^{a_5;(4)})  \nn
&&\hspace{2cm} 
E(\tr_{W_{a_4}}
\Y_{a_5a_4;k}^{a_4;(3)}(\mathcal{U}(q_{z_2})\Y_{a_1a_2;l}^{a_5;(4)}
(w_{a_1}, z_{1}-z_{2})w_{a_2}, q_{z_2})
q_{\tau}^{L(0)-\frac{c}{24}}),
\eea
where $F(\Y_{a_1a_3;i}^{a_4;(1)}\otimes \Y_{a_2a_4;j}^{a_3;(2)}, 
\Y_{a_5a_4;k}^{a_4;(3)}\otimes \Y_{a_1a_2;l}^{a_5;(4)})$ is the matrix
representation of $\mathcal{F}$ in the basis 
$\{\Y_{a_1a_3;i}^{a_4;(1)}\otimes \Y_{a_2a_4;j}^{a_3;(2)}\}_{i,j}, 
\{ \Y_{a_5a_4;k}^{a_4;(3)}\otimes \Y_{a_1a_2;l}^{a_5;(4)}\}_{k,l} $.
Therefore, $\forall w_{a_k}\in W_{a_k}, k=1,2$, we obtain
\bea  \label{S-act-right}
&& \left( \Psi_2(\Y_{a_1a; i}^{a_3;(1)}\otimes \Y_{a_2a_3;j}^{a;(2)})
\left( \left( -\frac{1}{\tau}\right)^{L(0)}w_{a_1}\otimes 
\left( -\frac{1}{\tau}\right)^{L(0)}w_{a_2} \right) \right) 
\left( -\frac{1}{\tau}z_{1}, -\frac{1}{\tau}z_{2};
-\frac{1}{\tau} \right)  \nn
&&\hspace{2cm}= 
\big( \Psi_2( S(\Y_{a_1a; i}^{a_3;(1)}\otimes \Y_{a_2a_3;j}^{a;(2)}))
(w_{a_1}\otimes w_{a_2})\big) (z_1, z_2, \tau). 
\eea

Combining (\ref{cardy-L-2}), (\ref{cardy-R-2}), 
(\ref{Y-1}), (\ref{Y-2}) and (\ref{S-act-right}) and 
Remark \ref{cardy-R-conv-rema}, 
we obtain a simpler version of Cardy condition.
\begin{thm}
The Cardy condition can be rewritten as follow: 
\bea  \label{cardy-ioa}
&&  \big( \theta_{W_{r_R(i)}} \circ \sigma_{123}(\Y_{\iota'_{cl-op}}) 
\circ (\varphi_{op}\otimes \id_{W'_{r_L(i)}}) \big) \otimes 
\big( \Omega_0(\sigma_{132}(\Y_{\iota_{cl-op}}))  
\circ (\varphi_{op} \otimes \id_{W_{r_R(i)}}) \big) \nn
&&\hspace{2cm} = S^{-1} \left( Y_{op}^f \otimes 
\big( \Omega_{-1}(Y_{op}^f) \circ 
(\theta_{V_{op}} \otimes \id_{V_{op}})\big) \right).
\eea
\end{thm}

\renewcommand{\theequation}{\thesection.\arabic{equation}}
\renewcommand{\thethm}{\thesection.\arabic{thm}}
\setcounter{equation}{0}
\setcounter{thm}{0}

\section{Modular tensor categories}

This section is independent of the rest of this work. 
The tensor product theory of modules over a vertex operator
algebra is developed by Huang and Lepowsky 
\cite{HL1}-\cite{HL4}\cite{H4}. 
In particular, the notion of vertex tensor category 
is introduced in \cite{HL1}. Huang later proved in \cite{H12} 
that $\mathcal{C}_V$ is a modular tensor category for $V$
satisfying conditions in Theorem \ref{ioa}.  
In Section 4.1, we review some basic ingredients of 
modular tensor category $\mathcal{C}_V$. In Section 4.2, 
we show how to find in $\mathcal{C}_V$ a graphical representation  
of the modular transformation $S: \tau \mapsto -\frac{1}{\tau}$
discussed in Section 3.2. 

\subsection{Preliminaries}
We recall some ingredients of 
vertex tensor category $\mathcal{C}_V$ and those
of modular tensor category structure on $\mathcal{C}_V$
constructed in \cite{H12}.

There is an associativity isomorphism $\mathcal{A}$
$$
\mathcal{A}: W_1 \boxtimes (W_2 \boxtimes W_3) 
\rightarrow (W_1\boxtimes W_2) \boxtimes W_3
$$
for each triple of $V$-modules $W_1, W_2, W_3$.
The relation between the fusing isomorphism $\mathcal{F}$ 
(recall (\ref{F-iso}))
and the associativity isomorphism $\mathcal{A}$ 
in $\mathcal{C}_V$ is described
by the following commutative diagram: 
\beq  \label{nat-F-A}
\xymatrix{
\oplus_{a_5\in \I} \V_{a_1a_5}^{a_4}\otimes \V_{a_2a_3}^{a_5}
\ar[r]^{\hspace{-1.5cm}\cong} \ar[d]^{\mathcal{F}}  &  
\hom_{V}(W_{a_1}\boxtimes (W_{a_2}\boxtimes W_{a_3}), W_{a_4})  
\ar[d]^{(\A^{-1})^*} \\
\oplus_{a_6\in \I} \V_{a_6a_3}^{a_4}\otimes \V_{a_1a_2}^{a_6} 
\ar[r]^{\hspace{-1.5cm} \cong}    &  
\hom_{V}((W_{a_1}\boxtimes W_{a_2})\boxtimes W_{a_3}, W_{a_4}) 
} 
\eeq
where the two horizontal maps are canonical isomorphisms induced
from the universal property of $\boxtimes$.

We recall the braiding structure on $\mathcal{C}_V$.
For each pair of $V$-modules $W_1,W_2$, 
there is also a natural isomorphism, for $z>0$,  
$\mathcal{R}_+^{P(z)}: W_1\boxtimes_{P(z)} W_2 
\rightarrow W_2\boxtimes_{P(z)} W_1$, defined by 
\begin{equation}  \label{R-+-cat}
\overline{\mathcal{R}_+^{P(z)}}(w_1\boxtimes_{P(z)} w_2) = e^{L(-1)} 
\overline{\mathcal{T}}_{\gamma_+} (w_2\boxtimes_{P(-z)} w_1),
\end{equation}
where $\gamma_+$ is a path from $-z$ to $z$ inside the lower half
plane as shown in the following graph
\beq    \label{gamma_--fig}
\begin{picture}(14,2.4)
\put(4, 0){\resizebox{5cm}{2cm}{\includegraphics{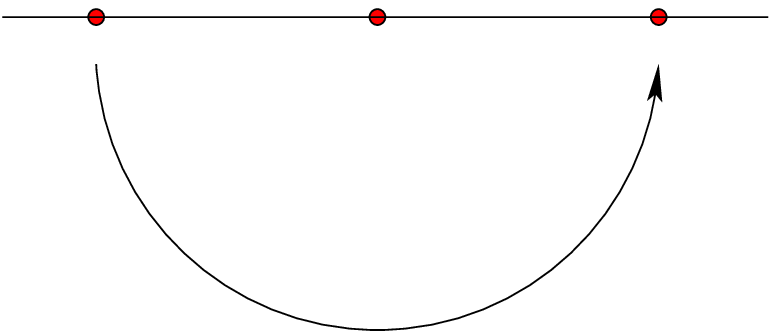}}}
\put(4.3,2.1){$-z$}\put(6.4,2.1){$0$}\put(8.2,2.1){$z$}
\put(10,0){.}
\end{picture}
\eeq
The inverse of $\mathcal{R}_+^{P(z)}$ is denoted by 
$\mathcal{R}_-^{P(z)}$, which is characterized by 
\begin{equation}   \label{R---cat}
\overline{\mathcal{R}_-^{P(z)}}(w_2\boxtimes_{P(z)} w_1) = e^{L(-1)} 
\overline{\mathcal{T}}_{\gamma_-} (w_1\boxtimes_{P(-z)} w_2),
\end{equation}
where $\gamma_-$ is a path 
in the upper half plane as shown in the following graph
\beq        \label{gamma_-+fig}
\begin{picture}(14,2.4)
\put(4, 0.4){\resizebox{5cm}{2cm}{\includegraphics{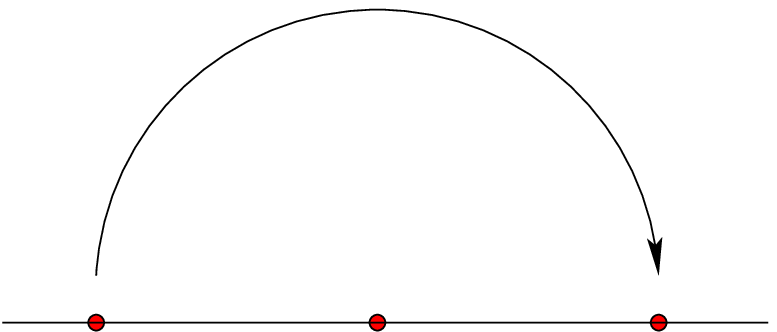}}}
\put(4.3,0){$-z$}\put(6.4,0){$0$}\put(8.2,0){$z$}
\put(10,0){.}
\end{picture}
\eeq 

We denote $\mathcal{R}_{\pm}^{P(1)}$ simply as $\mathcal{R}_{\pm}$. 
The natural isomorphisms $\mathcal{R}_{\pm}$ 
give $\mathcal{C}_V$ two different braiding structures.
We choose $\mathcal{R}_+$ as the default braiding structure on 
$\mathcal{C}_V$. Sometimes, we will denote it by
$(\mathcal{C}_V, \mathcal{R}_+)$ to emphasis our choice of braiding
isomorphisms.

Notice that our choice of $\mathcal{R}_{\pm}$ follows that in 
\cite{Ko2}, which is different 
from that in \cite{H9}\cite{H12}\cite{Ko1}. 
For each $V$-module $W$, (\ref{theta-def}) defines an automorphism 
$\theta_W: W\rightarrow W$ called a {\it twist}. 
A $V$-module $W$ is said to have a trivial twist if $\theta_W=\id_W$. 
The twist $\theta$ and braiding $\mathcal{R}_+$
satisfy the following three balancing axioms 
\bea  
\theta_{W_1\boxtimes W_2} &=& \mathcal{R}_+ \circ \mathcal{R}_+ \circ
(\theta_{W_1} \boxtimes \theta_{W_2}) \label{balanc-axiom-1} \\
\theta_V &=& \id_V  \label{balanc-axiom-2} \\
\theta_{W'} &=& (\theta_W)^*.  \label{balanc-axiom-3}
\eea
for any pair of $V$-modules $W_1, W_2$.

Let $\{ \Y_{ea}^a \}$ be a basis of $\V_{ea}^a$ 
for all $a\in \I$ such
that it coincides with the vertex operator $Y_{W_a}$, which defines
the $V$-module structure on $W_a$, i.e. $\Y_{ea}^a=Y_{W_a}$.
We choose a basis $\{ \Y_{ae}^a \}$ of $\V_{ae}^a$ as follow:
\beq
\Y_{ae}^a = \Omega_{-1} (\Y_{ea}^a).
\eeq
We also choose a basis $\{ \Y_{aa'}^e \}$ of $\V_{aa'}^e$ as 
\beq
\Y_{aa'}^{e}= \Y_{aa'}^{e'} = \hat{A}_{0} (\Y_{ae}^a) = 
\sigma_{132}(\Y_{ea}^a). 
\eeq
Notice that these choices are made 
for all $a\in \I$. In particular, we have
$$
\Y_{a'e}^{a'} = \Omega_{-1} (\Y_{ea'}^{a'}), \quad \quad 
\Y_{a'a}^e = \Y_{a'a}^{e'} = \hat{A}_{0} (\Y_{a'e}^{a'}). 
$$
The following relation was proved in \cite{Ko1}. 
\beq 
\Y_{a'a}^e = e^{2\pi i h_a} \Omega_{0} (\Y_{aa'}^e) =
 e^{-2\pi i h_a} \Omega_{-1} (\Y_{aa'}^e). \label{3-equ-1} 
\eeq

For any $V$-modules $W_1,W_2,W_3$ and $\Y\in \V_{W_1W_2}^{W_3}$, 
we denote by $m_{\Y}$ the morphism 
in $\hom_V(W_1\boxtimes W_2, W_3)$ associated to $\Y$ under
the identification of two spaces induced by 
the universal property of $\boxtimes$.

Now we recall the construction of duality maps \cite{H12}. 
We will follow the convention in \cite{Ko1}. 
Since $\mathcal{C}_V$ is semisimple, we only need to
discuss irreducible modules. For $a\in \I$, the right duality maps 
$e_{a}: W'_a \boxtimes W_a \rightarrow V$
and $i_{a}: V \rightarrow W_a \boxtimes W'_a$
for $a\in \I$ are given by
\beq
e_{a} = m_{\Y_{a'a}^e}, \quad\quad\quad   
m_{\Y_{aa'}^e} \circ i_{a} = \dim a \,\, \id_{V}  \nonumber
\eeq 
where $\dim a \neq 0$ for $a\in \I$ (proved by Huang in \cite{H11}).
The left duality maps 
$e'_{a}: W^a \boxtimes W'_a \rightarrow V$
and $i'_{a}: V \rightarrow W'_a \boxtimes W^a$
are given by 
\beq
e'_{a} = m_{\Y_{aa'}^e}, \quad\quad\quad 
m_{\Y_{a'a}^e} \circ i'_{a} = \dim a \,\, \id_{V}. \nonumber
\eeq

In a ribbon category, 
there is a powerful tool called graphic calculus. One can 
express various morphisms in terms of graphs. In particular, 
the right duality maps $i_a$ and $e_a$ are
denoted by the following graphs:
$$
\begin{picture}(14,1.5)
\put(2,0.7){$i_{a} =$} \put(3.6, 1.3){$a$}\put(5.6,1.3){$a'$}
\put(4,0){\resizebox{1.5cm}{1.5cm}
{\includegraphics{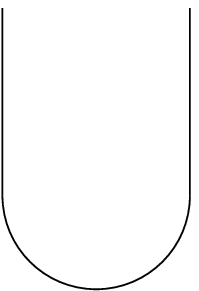}}}

\put(6.5,0){,}

\put(8,0.7){$e_{a} =$}\put(10,0)
{\resizebox{1.5cm}{1.5cm}{\includegraphics{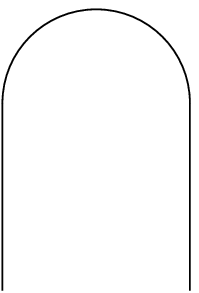}}}
\put(9.6, 0){$a$}\put(11.6, 0){$a'$} \put(12,0){,}
\end{picture}
$$
the left duality maps are denoted by
$$
\begin{picture}(14,1.5)
\put(2,0.7){$i'_{a} =$} \put(3.6, 1.3){$a'$}\put(5.6,1.3){$a$}
\put(4,0){\resizebox{1.5cm}{1.5cm}
{\includegraphics{left-dual-1.eps}}}

\put(6.5,0){,}

\put(8,0.7){$e'_{a} =$}\put(10,0)
{\resizebox{1.5cm}{1.5cm}{\includegraphics{left-dual-3.eps}}}
\put(9.6, 0){$a'$}\put(11.6, 0){$a$} \put(12,0){,}
\end{picture}
$$
and the twist and its inverse, for any object $W$, 
are denoted by
$$
\begin{picture}(14,2)
\put(2,1){$\theta_{W} =$} \put(3.4, 1.7){$W$}\put(3.4,0){$W$}
\put(4,0){\resizebox{0.3cm}{2cm}
{\includegraphics{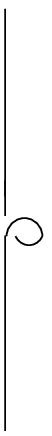}}}

\put(6.5,0){,}

\put(8,1){$\theta_{W}^{-1} =$}\put(9.4,1.7){$W$}\put(9.4,0){$W$}
\put(10,0){\resizebox{0.3cm}{2cm}
{\includegraphics{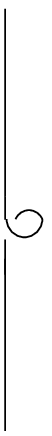}}}\put(11,0){.}

\end{picture}
$$
The identity (\ref{3-equ-1}) proved
in \cite{Ko2} is nothing but the following identity: 
\beq  \label{dual-twist}
\begin{picture}(14,2)
\put(2.6,0){$a'$}\put(4.6,0){$a$}
\put(3,0){\resizebox{1.5cm}{2cm}{\includegraphics{left-dual-3.eps}}}
\put(5,1){$=$}
\put(5.7,0){$a$}\put(7.6,0){$a'$}
\put(6,0){\resizebox{1.5cm}{2cm}{\includegraphics{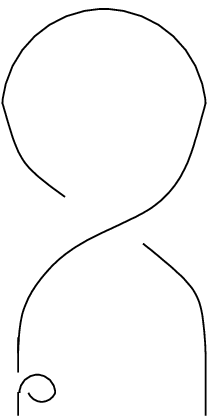}}}
\put(8,1){$=$}
\put(8.4,0){$a$}\put(10.4,0){$a'$}
\put(8.8,0){\resizebox{1.5cm}{2cm}
{\includegraphics{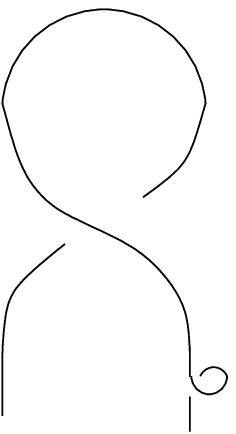}}}  \put(12,0){.}
\end{picture}
\eeq
This formula (\ref{dual-twist}) is implicitly 
used in many graphic calculations in this work.

A basis $\{ \Y_{a_1a_2; i}^{a_3;(1)} \}_{i=1}^{N_{a_1a_2}^{a_3}}$ 
of $\V_{a_1a_2}^{a_3}$ for $a_1,a_2,a_3\in \I$
induces a basis $\{ e_{a_1a_2;i}^{a_3}\}$ 
of $\hom (W_{a_1}\boxtimes W_{a_2}, W_{a_3})$. One can also 
denote $e_{a_1a_2;i}^{a_3}$ as the following graph:
\beq  \label{basis-fig}
\begin{picture}(14, 1.8)
\put(6,0.2){\resizebox{1.5cm}{1.5cm}{\includegraphics{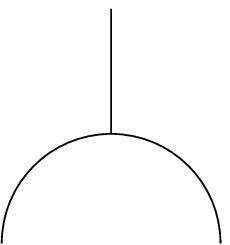}}}
\put(5.5, 0.2){$a_1$}\put(7.6, 0.2){$a_2$}\put(6.3, 1.5){$a_3$}
\put(6.6, 0.5){$i$}\put(8.5,0.2){.}
\end{picture}
\eeq
Note that we will always use $a$ to represent $W_a$ and $a'$ 
to represent $W'_a$ in graph for simplicity.
By the universal property of $\boxtimes_{P(z)}$, 
the map $\Omega_0: \V_{a_1a_2}^{a_3} \rightarrow \V_{a_2a_1}^{a_3}$ 
induces a linear map 
$\Omega_0: \hom_V(W_{a_1}\boxtimes W_{a_2}, W_{a_3}) \rightarrow
\hom_V(W_{a_2}\boxtimes W_{a_1}, W_{a_3})$ given as follow: 
\beq   \label{Omega-0-fig}
\begin{picture}(14, 1.8)
\put(2, 0.8){$ \Omega_0: $}
\put(4,0.2){\resizebox{1.5cm}{1.5cm}{\includegraphics{Y.eps}}}
\put(3.5, 0.2){$a_1$}\put(5.6, 0.2){$a_2$}\put(4.3, 1.5){$a_3$}
\put(4.6, 0.5){$i$}

\put(6.5, 0.8){$\mapsto$}

\put(9,0){\resizebox{1.5cm}{1.8cm}{\includegraphics{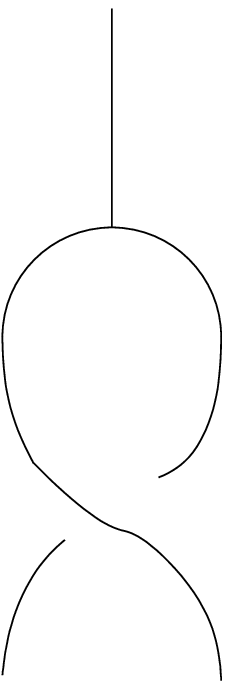}}}
\put(8.5, 0){$a_2$}\put(10.6, 0){$a_1$}\put(9.3, 1.7){$a_3$}
\put(9.7, 0.8){$i$}\put(11.5, 0){.}
\end{picture}
\eeq

Let us choose a basis $\{ f_{a_3;j}^{a_1a_2} \}_{j=1}^{N_{a_1a_2}^{a_3}}$ of 
$\hom_V(W_{a_3}, W_{a_1}\boxtimes W_{a_2})$, denoted by 
\beq  \label{dual-basis-pic}
\begin{picture}(14, 1.5)
\put(4, 0.7){$f_{a_3;j}^{a_1a_2} \, = $}
\put(6.8, 0){\resizebox{1.5cm}{1.5cm}{\includegraphics{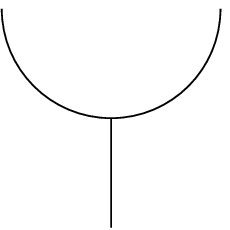}}}
\put(6.4, 1.2){$a_1$}\put(8.3, 1.2){$a_2$}\put(7, 0){$a_3$}
\put(7.4, 1){$j$}\put(8.5, 0){,}
\end{picture} 
\eeq
such that 
\beq \label{Y-dual-Y}
\begin{picture}(14,2)
\put(4, 0){\resizebox{1.5cm}{2cm}{\includegraphics{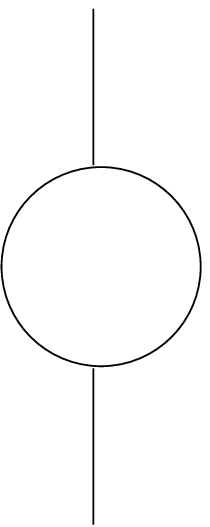}}}
\put(3.5, 0.9){$a_1$}\put(5.6, 0.9){$a_2$}\put(4.8, 1.8){$a_3$}
\put(4.8,0.1){$a_3$} \put(4.6, 1.1){$i$}\put(4.7, 0.8){$j$}
\put(6.9, 0.9){$=$} \put(8, 0.9){$\delta_{ij}$} 
\put(9,0){\resizebox{0.1cm}{2cm}{\includegraphics{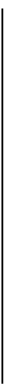}}}
\put(9.2,0){$a_3$}\put(10,0){.}
\end{picture}
\eeq
The following identity is proved in \cite{Ko1}
\beq \label{exp-id}
\begin{picture}(14,2)
\put(2.5, 0){$a_1$}\put(4.6,0){$a_2$}
\put(3, 0){\resizebox{1.5cm}{2cm}{\includegraphics{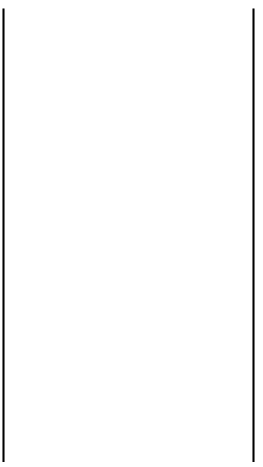}}}

\put(5.5, 1){$\ds =$}
\put(6.5, 1){$\ds \sum_{a_3\in \I} \sum_{i=1}^{N_{a_1a_2}^{a_3}}$}

\put(9, 0){\resizebox{1.5cm}{2cm}{\includegraphics{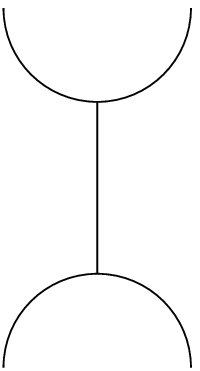}}}
\put(8.5,0){$a_1$}\put(10.6,0){$a_2$}
\put(8.5,1.8){$a_1$}\put(10.6,1.8){$a_2$}
\put(9.9,1){$a_3$}\put(9.7,0.2){$i$}\put(9.7,1.6){$i$}
\put(11.2,0){.}
\end{picture}
\eeq
We prove a similar identity below. 
\begin{lemma}
\beq \label{lemma-dual-equ}
\begin{picture}(14,2)
\put(2, 1){$\ds \sum_{a_4\in \I}\sum_{l} \, \frac{\dim a_4}{\dim b}$}
\put(5.5, 0){\resizebox{2cm}{2cm}
{\includegraphics{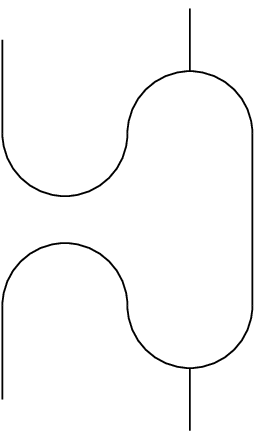}}}
\put(5,0){$a'_3$}\put(5,1.8){$a'_3$}
\put(7.1,0){$b$}\put(7.1,1.8){$b$}\put(7.6,1){$a_4$}
\put(6.9, 0.4){$k$}\put(6.9, 1.3){$k$}
\put(8.3,1){$=$} 
\put(9.5,0){\resizebox{0.1cm}{2cm}{\includegraphics{id.eps}}}
\put(9.7,0){$a'_3$}
\put(10.5,0){\resizebox{0.1cm}{2cm}{\includegraphics{id.eps}}}
\put(10.7,0){$b$}
\end{picture}
\eeq
\end{lemma}
\pf
Using the first balancing axiom (\ref{balanc-axiom-1}), 
we have 
\beq   \label{lemma-dual-Y-equ-1}
\begin{picture}(15,3)
\put(0.5, 0.5){\resizebox{2.5cm}{2cm}
{\includegraphics{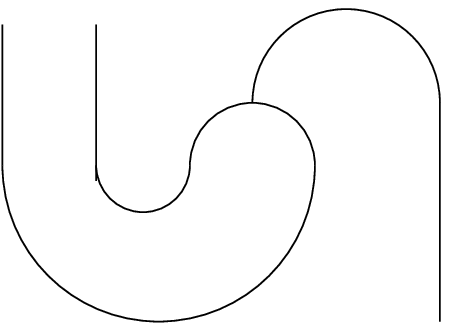}}}
\put(1.1, 1.4){$a_3$}\put(2.4, 1.4){$a_4$}\put(1.7, 2){$b$}
\put(1.8, 1.5){$l$}

\put(3.3,1.4){$=$}

\put(4, 0){\resizebox{3cm}{3cm}
{\includegraphics{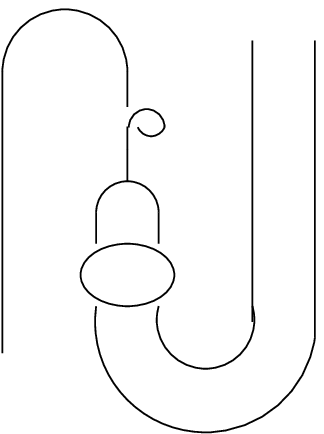}}}
\put(4.4, 1.4){$a_3$}\put(5.6, 1.4){$a_4$}\put(4.9, 2){$b$}
\put(5.1, 1.4){$l$}
\put(4.9, 1){$\theta^{-1}$}

\put(7.3, 1.4){$=$}

\put(8, 0){\resizebox{3cm}{3cm}
{\includegraphics{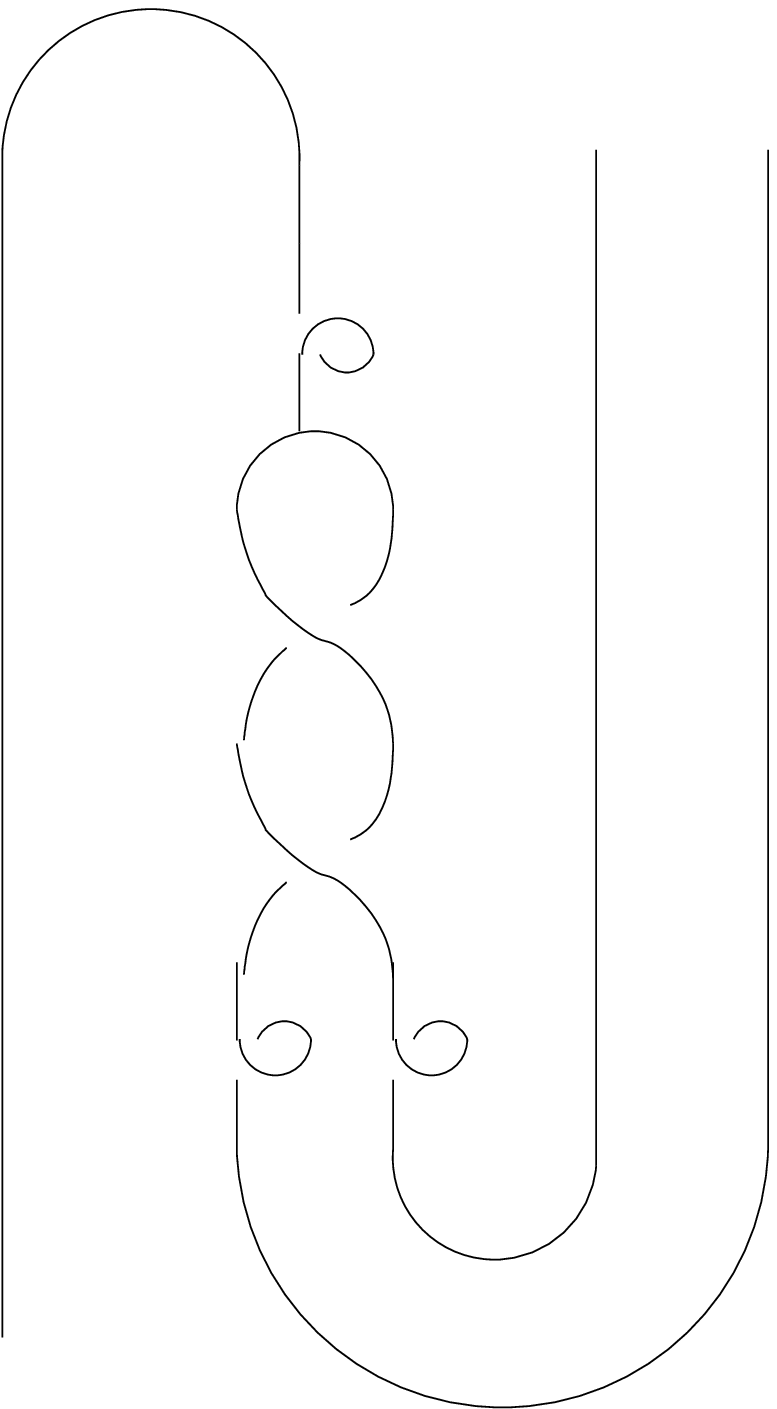}}}
\put(9.7, 1.3){$a_3$}\put(8.5, 1.3){$a_4$}\put(8.8, 2.2){$b$}
\put(9.1, 1.8){$l$}

\put(11.3, 1.4){$=$}

\put(12, 0.5){\resizebox{2.5cm}{2cm}
{\includegraphics{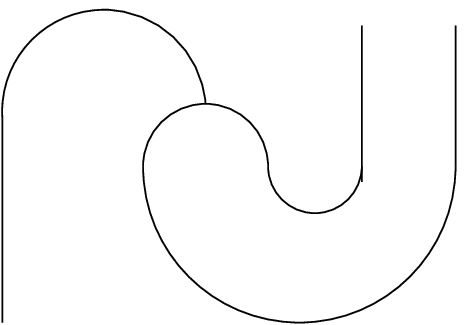}}}    
\put(12.4, 1.4){$a_3$}\put(13.5, 1.4){$a_4$}\put(13.2, 2){$b$}
\put(13.1, 1.5){$l$}

\put(14.8,0){.}

\end{picture}
\eeq
Then the Lemma follows from the following relations: 
\beq
\begin{picture}(14,2.5)
\put(1.7,1.2){$\frac{\dim a_4}{\dim b}$}
\put(3.2,0.2){\resizebox{2cm}{2cm}
{\includegraphics{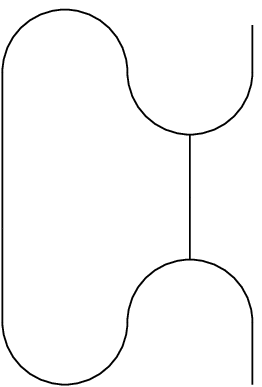}}}
\put(3.4,1.2){$a'_3$}\put(4.8,1.2){$b$}\put(5.3,0.2){$a_4$}
\put(5.3,2){$a_4$}
\put(4.6,1.7){$k$}\put(4.6,0.5){$l$}
\put(6,1.2){$=$}
\put(6.8,1.2){$\frac{1}{\dim b}$}
\put(8.2,0){\resizebox{2.5cm}{2.5cm}
{\includegraphics{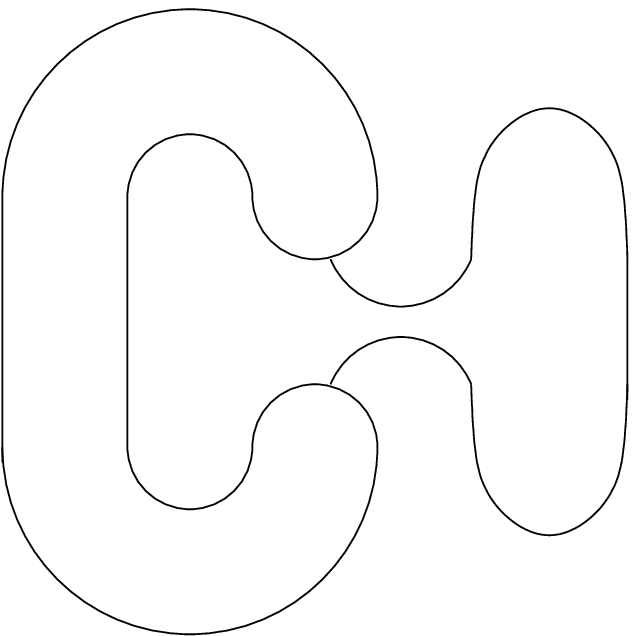}}}
\put(8.8,1.2){$a'_3$}\put(7.7,1.2){$a'_4$}\put(10.8,0.9){$b$}
\put(9.4,1.6){$k$}\put(9.4,0.6){$l$}
\put(11.4,0){\resizebox{0.1cm}{2.5cm}{\includegraphics{id.eps}}}
\put(11.6,0){$a_4$}

\end{picture} 
\eeq
\beq  \label{lemma-dual-Y-equ-2}
\begin{picture}(14,2)
\put(6,1){$=$}
\put(6.8,1){$\frac{1}{\dim b}$}
\put(8.2, 0){\resizebox{2cm}{2cm}
{\includegraphics{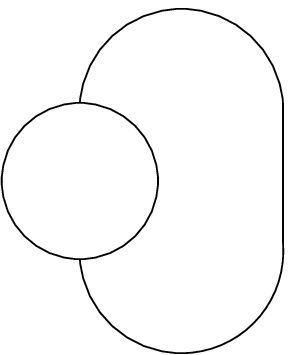}}}
\put(7.7,1){$a'_3$}\put(9.4,1){$a_4$}
\put(8.4, 0.2){$b$}\put(8.6,1.1){$l$}
\put(8.6,0.7){$k$}

\put(10.7,0){\resizebox{0.1cm}{2cm}{\includegraphics{id.eps}}}
\put(10.8,0){$a_4$}
\put(11.2,1){$= \, \delta_{kl}$}
\put(12.5,0){\resizebox{0.1cm}{2cm}{\includegraphics{id.eps}}}
\put(12.7,0){$a_4$} \put(13.2,0){.}
\end{picture}
\eeq
\epf

Similar to $\Omega_0$, 
$\tilde{A}_0$, $\sigma_{123}$ and $\sigma_{132}$ can also be
described graphically as proved in \cite{Ko1}. We recall these 
results below.  
\begin{prop}
\beq \label{tilde-A-0-graph}
\begin{picture}(14, 2)
\put(2, 0.8){$ \tilde{A}_0: $}
\put(4,0.2){\resizebox{1.5cm}{1.5cm}{\includegraphics{Y.eps}}}
\put(3.5, 0.2){$a_1$}\put(5.6, 0.2){$a_2$}\put(4.3, 1.5){$a_3$}
\put(4.6, 0.5){$i$}

\put(7, 0.8){$\mapsto$}
\put(9.5,0){\resizebox{2cm}{2.2cm}{\includegraphics{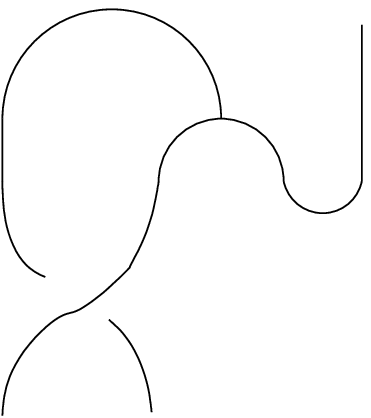}}}
\put(9, 0){$a_1$}\put(10.5, 0){$a'_3$}\put(11.6, 1.9){$a'_2$}
\put(10.6, 1.2){$i$}\put(12, 0){.}
\end{picture}
\eeq
\beq \label{sigma-123-graph}
\begin{picture}(14, 1.5)
\put(2, 0.8){$ \sigma_{123}: $}
\put(4,0.2){\resizebox{1.5cm}{1.5cm}{\includegraphics{Y.eps}}}
\put(3.5, 0.2){$a_1$}\put(5.6, 0.2){$a_2$}\put(4.3, 1.5){$a_3$}
\put(4.6, 0.5){$i$}

\put(7, 0.8){$\mapsto$}
\put(9.5,0){\resizebox{2cm}{1.5cm}{\includegraphics{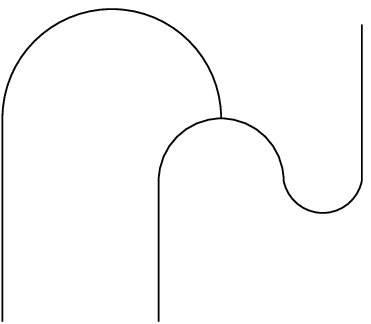}}}
\put(9, 0){$a'_3$}\put(10.5, 0){$a_1$}\put(11.6, 1.3){$a'_2$}
\put(10.6, 0.6){$i$}\put(12, 0){,}
\end{picture}
\eeq
\beq \label{sigma-132-graph}
\begin{picture}(14, 2)
\put(2, 0.8){$ \sigma_{132}: $}
\put(4,0.2){\resizebox{1.5cm}{1.5cm}{\includegraphics{Y.eps}}}
\put(3.5, 0.2){$a_1$}\put(5.6, 0.2){$a_2$}\put(4.3, 1.5){$a_3$}
\put(4.6, 0.5){$i$}

\put(7, 0.8){$\mapsto$}
\put(9.5,0){\resizebox{2cm}{1.5cm}{\includegraphics{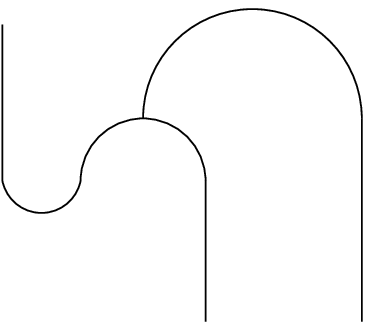}}}
\put(9.1, 1.3){$a'_1$}\put(10.8, 0){$a_2$}\put(11.6, 0){$a'_3$}
\put(10.2, 0.6){$i$}\put(12, 0){.}
\end{picture}
\eeq
\end{prop}

\subsection{Graphical representation of 
$S: \tau \mapsto -\frac{1}{\tau}$}

In \cite{HKo3}, we defined an action of $\alpha, \beta$ on 
$\Psi_2(\Y_{aa_{1}; p}^{a_{1};(1)}\otimes 
\Y_{a_{2}a_{3}; q}^{a; (2)})$. More precisely,
\begin{eqnarray}
\alpha(\Psi_2(\Y_{aa_{1}; p}^{a_{1};(1)}\otimes 
\Y_{a_{2}a_{3}; q}^{a; (2)})):\oplus_{a_2, a_3\in \I}W_{a_{2}}
\otimes W_{a_{3}}&\to &\mathbb{G}_{1;2}\label{alpha-def}\\
\beta(\Psi_2(\Y_{aa_{1}; p}^{a_{1};(1)}\otimes 
\Y_{a_{2}a_{3}; q}^{a; (2)})):\oplus_{a_2, a_3\in \I}W_{a_{2}}
\otimes W_{a_{3}}&\to &\mathbb{G}_{1;2}\label{beta-def}
\end{eqnarray}
are defined by
\bea
&&(\alpha(\Psi_2(\Y_{aa_{1}; p}^{a_{1};(1)}\otimes 
\Y_{a_{2}a_{3}; q}^{a; (2)})))(w_2\otimes w_3)  \nn
&&\hspace{1.5cm} 
=(\Psi_2(\Y_{aa_{1}; p}^{a_{1};(1)}\otimes 
\Y_{a_{2}a'_{2}; q}^{a; (2)}))(w_2 \otimes w_3; z_{1}, z_{2}-1; \tau)  \nn
&&(\beta(\Psi(\Y_{aa_{1}; p}^{a_{1};(1)}\otimes 
\Y_{a_{2}a_{3}; q}^{a; (2)})))
(w_2 \otimes w_3) \nn
&&\hspace{1.5cm} =(\Psi_2(\Y_{aa_{1}; p}^{a_{1};(1)}\otimes 
\Y_{a_{2}a_{3}; q}^{a; (2)}))
(w_2 \otimes w_3; z_{1}, z_{2}+\tau; \tau) 
\eea
if $w_2\otimes w_3 \in W_{a_2}\otimes W_{a_3}$,  
and by $0$ if otherwise.

We also showed in \cite{HKo3} that $\alpha$ induces an automorphism
on $\oplus_{a\in \I} \V_{aa_1}^{a_1}\otimes \V_{a_2a_3}^a$ given
as follow: 
\bea  \label{a-F-F}
&&\hspace{-0.5cm}\Y_{aa_1; i}^{a_1;(1)}\otimes \Y_{a_2a_3;j}^{a;(2)} \mapsto 
\sum_{b,c\in \mathcal{A}}\sum_{k,l,p,q}e^{-2\pi i h_{a_3}}
F^{-1}(\Y_{aa_{1}; i}^{a_{1};(1)}\otimes \Y_{a_{2}a_{3}; j}^{a;(2)};
\Y_{a_{2}b; k}^{a_{1};(3)}\otimes \Y_{a_{3}a_{1}; l}^{b;(4)}) \nn
&&\hspace{5cm}  
F(\Y_{a_{2}b; k}^{a_{1};(3)}\otimes \Omega_{-1}^2(\Y_{a_{3}a_{1}; l}^{b;(4)});
\Y_{ca_{1}; p}^{a_{1};(5)}\otimes \Y_{a_{2}a_{3}; q}^{c;(6)})  \nn
&&\hspace{5cm}\Y_{ca_{1}; p}^{a_{1};(5)}\otimes \Y_{a_{2}a_{3}; q}^{c;(6)}.
\eea
We still denoted this automorphism and its natural extension on 
$\oplus_{a, a_1\in \I} \V_{aa_1}^{a_1}\otimes \V_{a_2a_3}^a$ by $\alpha$. 
The following Lemma follows immediately from (\ref{a-F-F}). 
\begin{lemma}
$\alpha$ can also be expressed graphically as follow:
\beq \label{a-pic}
\begin{picture}(14, 2.5)
\put(1.5,1){$\ds \alpha: $}

\put(3,0){\resizebox{2cm}{2cm}{\includegraphics{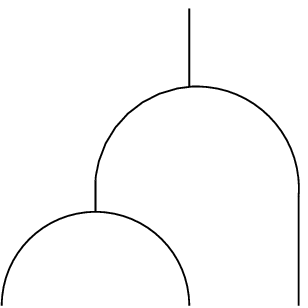}}}
\put(2.5,0){$a_2$}\put(4.4,0){$a_3$}\put(5.1,0){$a_1$}
\put(4.5,1.8){$a_1$}\put(3.5,0.2){$i$}\put(4.2,1.1){$j$}
\put(3.5, 1.2){$a$}

\put(6,1){$\mapsto$}
\put(8,0){\resizebox{2cm}{2.5cm}{\includegraphics{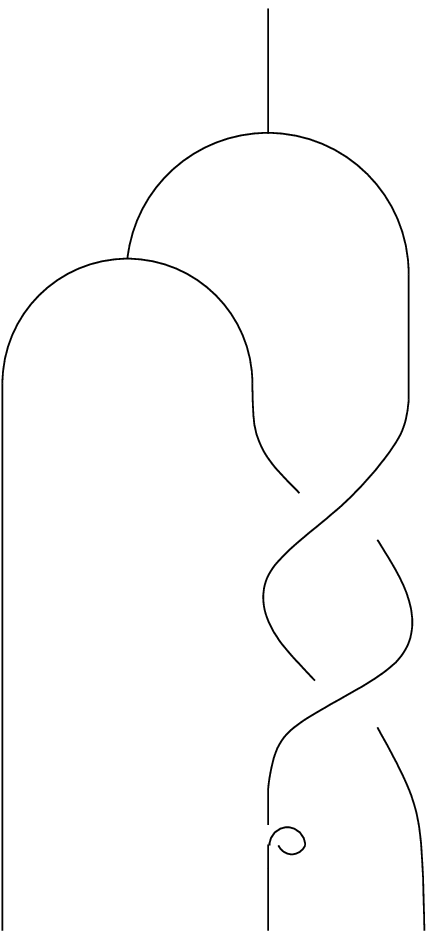}}}
\put(7.5,0){$a_2$}\put(8.6,0){$a_3$}\put(10.1,0){$a_1$}
\put(9.4,2.4){$a_1$}\put(8.5,1.4){$i$}\put(9.2,1.8){$j$}
\put(11,0){.}
\end{picture}
\eeq
\end{lemma}

For $\beta$, we prefer to use maps $\tilde{A}_0$ and $\hat{A}_0$ 
defined in (\ref{two-new-A-2}) instead of the map $A_0$ used in 
\cite{HKo3}. We obtain the following Lemma, which is 
proved in appendix. 
\begin{lemma}  \label{b-F-lemma}
$\beta$ also induces an automorphism
on $\oplus_{a,a_1\in \I} \V_{aa_1}^{a_1}\otimes \V_{a_2a_3}^{a}$ given by 
\bea  \label{b-F-map}
\Y_{aa_1; i}^{a_1;(1)}\otimes \Y_{a_2a_3;j}^{a;(2)} &\mapsto& 
\sum_{b\in \I}\sum_{k,l}\sum_{c\in \I} \sum_{p,q} 
F^{-1}(\Y_{aa_1; i}^{a_1;(1)}\otimes \Y_{a_2a_3;j}^{a;(2)};
\Y_{a_{2}b; k}^{a_{1};(3)}\otimes \Y_{a_{3}a_{1}; l}^{b;(4)} ) \nn
&&
F(\tilde{A}_{0}(\Y_{a_{2}b; k}^{a_{1};(3)})\otimes
  \tilde{A}_{0}(\Y_{a_{3}a_{1}; l}^{b;(4)}), 
  \Y_{cb';p}^{b';(5)} \otimes \Y_{a_2a_3;q}^{c;(6)} )  \nn
&&
\hat{A}_0(\Y_{cb';p}^{b';(5)})  \otimes \Omega_{0}^2(\Y_{a_2a_3;q}^{c;(6)}). 
\eea
We still denote this automorphism by $\beta$. 
\end{lemma}

\begin{lemma}
$\beta$ can be expressed graphically as follow: 
\beq  \label{b-pic}
\begin{picture}(14, 2.5)
\put(1.5,1){$\ds \beta: $}

\put(3,0){\resizebox{2cm}{2cm}{\includegraphics{a-LHS.eps}}}
\put(2.5,0){$a_2$}\put(4.4,0){$a_3$}\put(5.1,0){$a_1$}
\put(3.5, 1.2){$a$}
\put(4.5,1.8){$a_1$}\put(3.5,0.2){$i$}\put(4.2,1.1){$j$}

\put(6,1){$\mapsto$} 
\put(7,1){$\sum_{b\in \I} \sum_{l}$} 
\put(9.5,0){\resizebox{2cm}{2.5cm}{\includegraphics{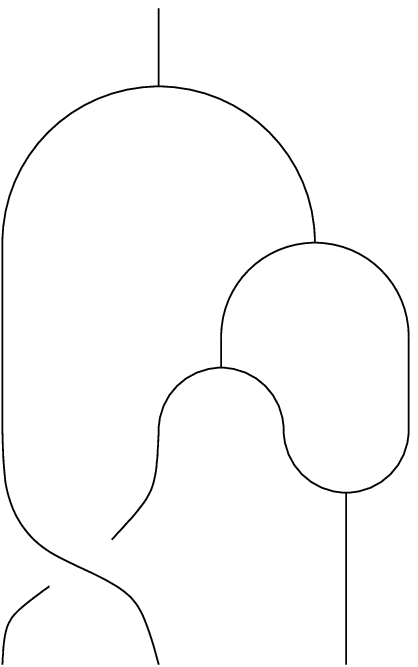}}}
\put(9,0){$a_2$}\put(10.4,0){$a_3$}\put(11.3,0){$b$}
\put(10.4,2.4){$b$}\put(10.5,0.8){$i$}\put(10.9,1.2){$j$}
\put(10.2,1.8){$l$}\put(11.1, 0.8){$l$}\put(10.3, 1.3){$a$}
\put(11.6,1){$a_1$}\put(10.7,0.4){$a_3$}\put(10.9, 2){$a_1$}
\put(12,0){.}
\end{picture}
\eeq
\end{lemma}
\pf
Using (\ref{b-F-map}), we can see that $\beta$ is the 
composition of following maps
\beq  \label{b-comp-maps}
\xymatrix{
\oplus_{a,a_1\in \I} \V_{aa_1}^{a_1}\otimes \V_{a_2a_3}^{a} 
\ar[r]^{\mathcal{F}^{-1}} &   
\oplus_{b,a_1\in \I} \V_{a_2b}^{a_1}\otimes \V_{a_3a_1}^{b}
\ar[r]^{\tilde{A}_0\otimes \tilde{A}_0} &
\oplus_{b,a_1\in \I}\V_{a_2a'_1}^{b'} \otimes \V_{a_3b'}^{a'_1}
\ar[d]^{\mathcal{F}}   \\
&  \oplus_{b, c\in \I} \V_{cb}^{b}\otimes \V_{a_2a_3}^{c}
&  \oplus_{b, c\in \I} \V_{cb'}^{b'}\otimes \V_{a_2a_3}^{c}.
\ar[l]_{\hat{A}_0\otimes \Omega_{0}^2}
} 
\eeq
By the commutative diagram (\ref{nat-F-A}), 
(\ref{b-comp-maps}) can be rewritten graphically as follow: 
\beq  \label{asso-cat-graph}
\begin{picture}(14,2)
\put(1.5,0.7){$\beta:$}
\put(3,0){\resizebox{2cm}{2cm}{\includegraphics{a-LHS.eps}}}
\put(2.5,0){$a_2$}\put(4.4,0){$a_3$}\put(5.1,0){$a_1$}
\put(3.5, 1.2){$a$}
\put(4.5,1.8){$a_1$}\put(3.5,0.2){$i$}\put(4.2,1.1){$j$}

\put(6.4,0.7){$=$} 

\put(7.4,0.7){$\ds \sum_{b\in \I} \sum_{k}$}  
\put(10,0){\resizebox{2cm}{2cm}{\includegraphics{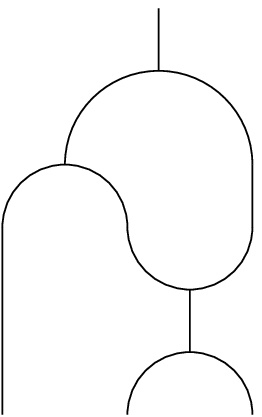}}}
\put(9.5,0){$a_2$}\put(10.5,0){$a_3$}\put(10.6,0.7){$a_3$}
\put(10.4,0.9){$i$}\put(11.1, 1.3){$j$}\put(10.5, 1.6){$a$}
\put(12.1,0){$a_1$}\put(12.1,1){$a_1$}
\put(11.4,1.8){$a_1$}\put(11.4,0.8){$k$}
\put(11.4,0){$k$}\put(11.7,0.35){$b$}
\end{picture}
\eeq
$$
\begin{picture}(14,3.3)
\put(0.5,1.3){$\mapsto$}
\put(1,1.3){$\ds \sum_{b\in \I} \sum_{k}$}  
\put(3,0){\resizebox{3cm}{3cm}{\includegraphics{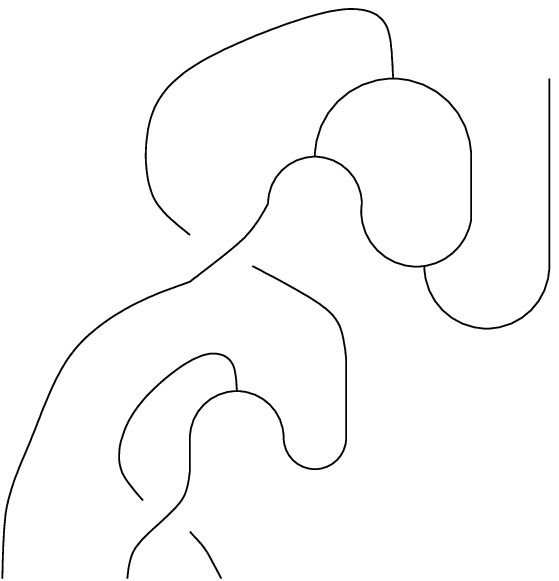}}}
\put(2.5,0){$a_2$}\put(3.2,0){$a_3$}\put(5,1){$a'_1$}
\put(6.1,2.4){$b'$}\put(5.2,1.8){$k$}
\put(4.2,0.6){$k$}\put(4.4,0){$b'$}
\put(4.6,1.8){$i$}\put(5, 2.2){$j$}\put(4.5, 2.4){$a$}

\put(6.4,1.3){$=$}
\put(6.9,1.3){$\ds \sum_{b\in \I} \sum_{k}$}  
\put(9.1,0){\resizebox{3cm}{3cm}{\includegraphics{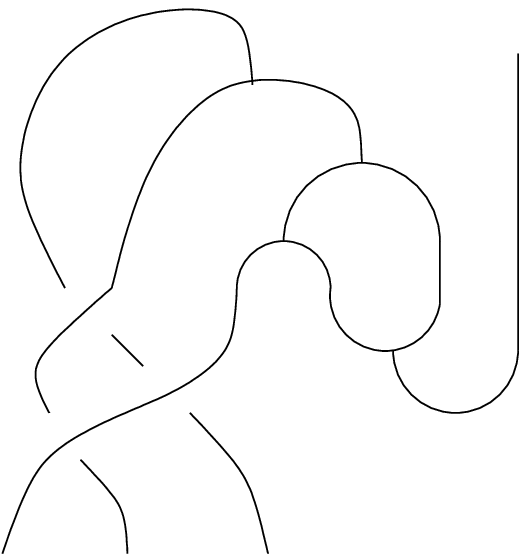}}}
\put(8.7,0){$a_2$}\put(9.9,0){$a_3$}\put(11.2,2.4){$a_1$}
\put(10.8,0){$b'$}\put(11.2, 1.3){$k$}
\put(12.2,2.6){$b'$}\put(10.6,1.4){$i$}\put(11.1,1.8){$j$}
\put(10.4,2.2){$k$}\put(10.5, 1.9){$a$}

\end{picture}
$$


$$
\begin{picture}(14,2.7)
\put(0.5,1.3){$\mapsto$}
\put(1,1.3){$\ds \sum_{b\in \I} \sum_{k}$}
\put(3.5,0){\resizebox{2cm}{2.5cm}{\includegraphics{car-b-RHS.eps}}}
\put(3,0){$a_2$}\put(4.4,0){$a_3$}\put(5.3,0){$b$}
\put(4.4,2.4){$b$}\put(4.5,0.8){$i$}\put(4.9,1.2){$j$}
\put(4.2,1.8){$k$}\put(5.1, 0.8){$k$}\put(4.3, 1.3){$a$}
\put(5.6,1){$a_1$}\put(4.7,0.4){$a_3$}\put(4.9, 2){$a_1$}
\put(6,0){.}

\end{picture}
$$
\epf

We have introduced $S, \alpha, \beta$ all as isomorphisms on 
$\oplus_{a,a_1\in \I} \V_{aa_1}^{a_1}\otimes \V_{a_2a_3}^{a}$.
They satisfy the following well-known equation 
\cite{MSei1}\cite{MSei2}\cite{H11}\cite{HKo3}:
\beq  \label{S-a-b}
S \alpha = \beta S.
\eeq
We proved in \cite{HKo3} that
$S$ is determined by the identity (\ref{S-a-b}) 
up to a constant $S_e^e$. We will solve the equation (\ref{S-a-b})
for $S$ graphically below. 

\begin{prop}
\beq  \label{S-a-map}
\begin{picture}(14,2)
\put(2, 1){$\ds S(a):$}
\put(4,0.2){\resizebox{1.5cm}{1.5cm}{\includegraphics{Y.eps}}}
\put(3.7, 0.2){$a$}\put(5.6, 0.2){$a_1$}\put(4.3, 1.5){$a_1$}
\put(4.7, 0.5){$i$}
\put(6,1){$\ds \mapsto$}\put(7,1){$\ds \sum_{a_2\in \I}$}
\put(7.9,1){$\ds S_e^e$}
\put(8.5, 1){$\ds \dim a_2$}
\put(10.5,0){\resizebox{2.2cm}{2.2cm}{\includegraphics{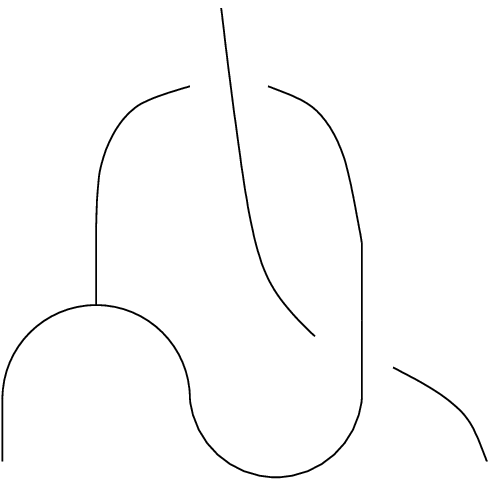}}}
\put(10.2, 0){$\ds a$}\put(10.8,0.5){$i$}
\put(10.4, 1){$a_1$}\put(12.9, 0){$a_2$}
\end{picture}
\eeq
\end{prop}
\pf
Since we know that the equation (\ref{S-a-b}) determine $S$
up to an overall constant $S_e^e$. Hence we only need to
check that (\ref{S-a-map}) gives a solution to (\ref{S-a-b}).

Combining  (\ref{b-pic}) with (\ref{S-a-map}), we obtain that
\beq  \label{b-S-map}
\begin{picture}(14, 3)
\put(1.5,1){$\ds \beta \circ S: $}

\put(3,0){\resizebox{2cm}{2cm}{\includegraphics{a-LHS.eps}}}
\put(2.5,0){$a_2$}\put(4.4,0){$a_3$}\put(5.1,0){$a_1$}
\put(4.5,1.8){$a_1$}\put(3.5,0.2){$i$}\put(4.2,1.1){$j$}
\put(3.5, 1.2){$a$}

\put(6,1){$\mapsto$}
\put(7,1){$\ds \sum_{b, a_4\in \I}\sum_k S_e^e \dim a_4$}
\put(10.8,0){\resizebox{3cm}{3cm}{\includegraphics{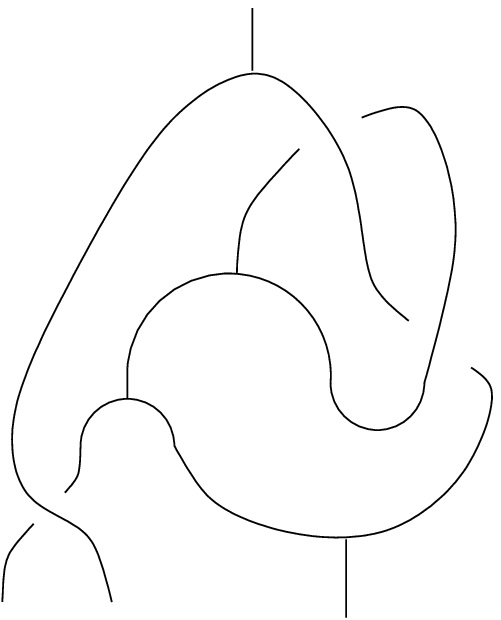}}}
\put(10.3,0){$a_2$}\put(11.6,0){$a_3$}\put(13,0){$b$}
\put(12.5,2.8){$b$}\put(12.2,0.2){$a_3$}\put(13.7,0.5){$a_4$}
\put(12.9,0.5){$k$}\put(11.5,0.7){$i$}\put(12.1,1.3){$j$}
\put(12.3,2.3){$k$}
\put(14,0){.}
\end{picture}
\eeq
The diagram in the right hand side of (\ref{b-S-map}) can be
deformed as follow: 
\beq 
\begin{picture}(14,4)
\put(2,0){\resizebox{4cm}{4cm}{\includegraphics{car-b-comp-S.eps}}}
\put(1.5,0){$a_2$}\put(3,0){$a_3$}\put(4.9,0){$b$}
\put(4.3,3.8){$b$}\put(4,0.3){$a_3$}\put(6,1){$a_4$}
\put(4.7, 0.7){$k$}\put(2.8,1){$i$}\put(3.8,1.8){$j$}
\put(4,3.2){$k$}

\put(6.5,2){$=$}
\put(8,0){\resizebox{4cm}{4cm}{\includegraphics{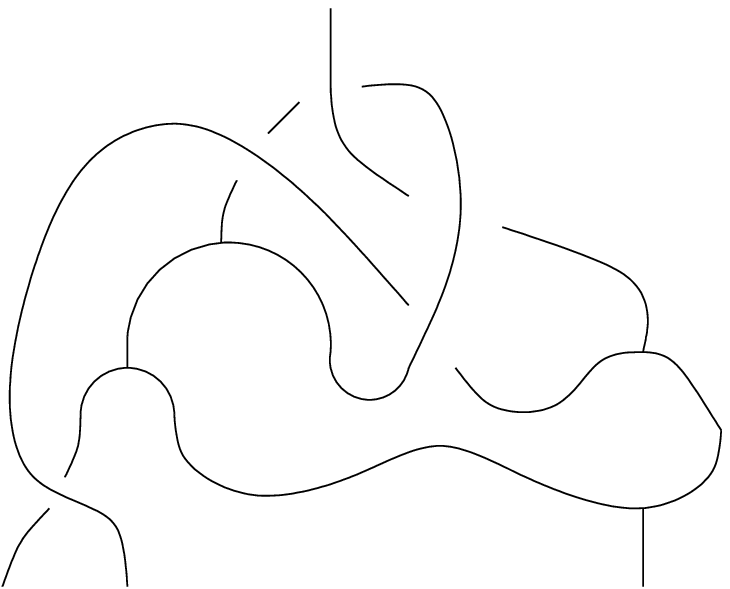}}}
\put(7.5,0){$a_2$}\put(8.8,0){$a_3$}\put(11.6,0){$b$}
\put(9.9,3.8){$b$}\put(10.6,0.5){$a_3$}\put(12.1,1){$a_4$}
\put(11.4, 0.7){$k$}\put(8.7,1){$i$}\put(9.2,1.9){$j$}
\put(11.4, 1.2){$k$}

\end{picture}
\eeq
By (\ref{lemma-dual-equ}), we have 
\beq  \label{b-comp-S-3}
\begin{picture}(14,4)
\put(0.2,2){$\ds \sum_{a_4\in \I}\sum_k S_e^e\dim a_4$} 
\put(3.5,0){\resizebox{4cm}{4cm}{\includegraphics{car-b-comp-S-1.eps}}}
\put(3,0){$a_2$}\put(4.3,0){$a_3$}\put(7.1,0){$b$}
\put(5.4,3.8){$b$}\put(6.1,0.5){$a_3$}\put(7.6,1){$a_4$}
\put(6.9, 0.7){$k$}\put(4.2,1){$i$}\put(4.7,1.9){$j$}
\put(6.9, 1.2){$k$}

\put(7.7,2){$=\, S_e^e\dim b$}
\put(9.8,0){\resizebox{4cm}{4cm}{\includegraphics{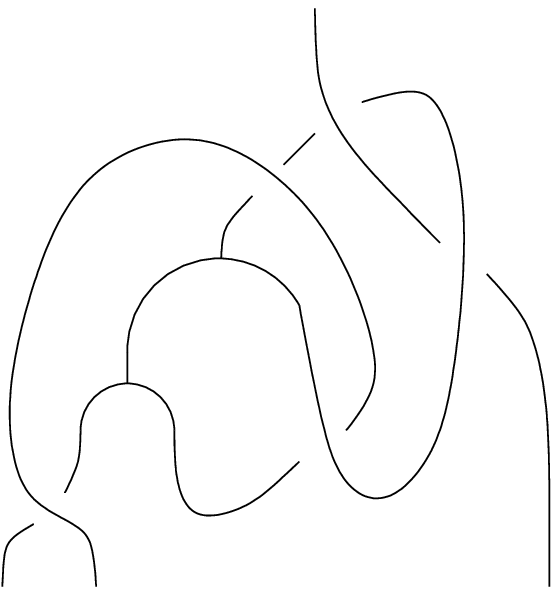}}}
\put(9.3,0){$a_2$}\put(10.7,0){$a_3$}\put(13.9,0){$b$}
\put(10.7,1){$i$}\put(11.4,1.9){$j$}

\end{picture}
\eeq

On the other hand, combining (\ref{a-pic}) with (\ref{S-a-map}), 
we obtain 
\beq \label{S-comp-a-1}
\begin{picture}(14, 3)
\put(1.5,1){$\ds S\circ \alpha: $}

\put(3,0){\resizebox{2cm}{2cm}{\includegraphics{a-LHS.eps}}}
\put(2.5,0){$a_2$}\put(4.4,0){$a_3$}\put(5.1,0){$a_1$}
\put(4.5,1.8){$a_1$}\put(3.5,0.2){$i$}\put(4.2,1.1){$j$}
\put(3.5, 1.2){$a$}

\put(6,1){$\mapsto$}
\put(7,1){$\ds \sum_{b\in \I} \, \, S_e^e \dim b$}
\put(10,0){\resizebox{3cm}{3cm}{\includegraphics{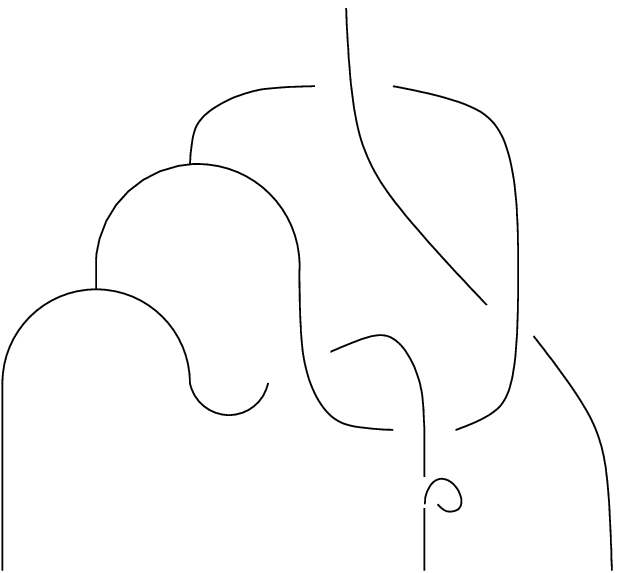}}}
\put(9.5,0){$a_2$}\put(11.5,0){$a_3$}\put(13.1,0){$b$}
\put(10.4,1.1){$i$}\put(10.9,1.8){$j$}\put(10.2, 1.7){$a$}
\put(12.6, 1.9){$a_1$}

\put(13.5,0){.}
\end{picture}
\eeq

Notice the diagram on the right hand side of 
(\ref{S-comp-a-1}) is exactly a deformation 
of the diagram on the right hand side of (\ref{b-comp-S-3}).
Hence we obtain that the map defined as (\ref{S-a-map})
give a solution to the equation $S\alpha=\beta S$.
\epf

To determine $S$ completely, we need determine $S_e^e$. 
This can be done by using other identities satisfied by $S$.

\begin{prop}
\beq  \label{SS-A-0}
S^2(a)(\Y_{aa_1;i}^{a_1;(1)}) = \tilde{A}_0(\Y_{aa_1;i}^{a_1;(1)})
\eeq
\end{prop}
\pf
By the definition of $S$-cation on 
$\Y_{aa_1;i}^{a_1;(1)} \in \V_{aa_1}^{a_1}$, we have
\beq
\Psi_1(S^2(a)(\Y_{aa_1;i}^{a_1;(1)}))(w_a; z, \tau) = 
\Psi_1(\Y_{aa_1;i}^{a_1;(1)})\left(
\tau^{L(0)}\left( \frac{-1}{\tau} \right)^{L(0)} w_a; -z, \tau \right) 
\eeq
for $w_a\in W_a$. Keep in mind of our convention on branch cut for
logarithm. We have 
$$
\tau^{L(0)} \left( \frac{-1}{\tau} \right)^{L(0)} w_a = e^{\pi iL(0)}w_a.
$$
Hence we obtain
\beq \label{SS-equ}
\Psi_1(S^2(a)(\Y_{aa_1;i}^{a_1;(1)}))(w_a; z, \tau) = 
\Psi_1(\Y_{aa_1;i}^{a_1;(1)})( e^{\pi i L(0)}w_a; -z, \tau)
\eeq
By (\ref{A-r-U}), we also have
\bea  
&&\Psi_1(\Y_{aa_1;i}^{a_1;(1)})( e^{\pi i L(0)}w_a; -z, \tau)  \nn
&&\hspace{1cm}= 
E\big( 
\tr_{W_{a_1}} \Y_{aa_1;i}^{a_1;(1)}(\mathcal{U}(e^{-2\pi iz})e^{\pi iL(0)}w_a, 
 e^{-2\pi iz}) q_{\tau}^{L(0)-\frac{c}{24}} \big) \nn
&&\hspace{1cm}= 
E\big( \tr_{(W_{a_1})'} \tilde{A}_0(\Y_{aa_1;i}^{a_1;(1)})(
e^{\pi iL(0)}\mathcal{U}(e^{2\pi iz}) w_a, 
e^{\pi i}  e^{2\pi iz}) q_{\tau}^{L(0)-\frac{c}{24}} \big) \nn
&&\hspace{1cm}= 
E\big( \tr_{(W_{a_1})'} \tilde{A}_0(\Y_{aa_1;i}^{a_1;(1)})(
\mathcal{U}(e^{2\pi iz}) w_a, e^{2\pi iz}) q_{\tau}^{L(0)-\frac{c}{24}} \big) \nn
&&\hspace{1cm}= \Psi_1(\tilde{A}_0(\Y_{aa_1;i}^{a_1;(1)}))(w_a; z, \tau)
  \label{SS-A-0-equ-2}
\eea
Therefore, combining (\ref{SS-equ})
and (\ref{SS-A-0-equ-2}), we obtain (\ref{SS-A-0}). 
\epf

The following lemma is proved in \cite{BK2}. 
\begin{lemma} {\rm 
Let $D^2=\sum_{a\in \I}\dim^2 a$. Then $D\neq 0$ and we have 
\beq  \label{BK-lemma}
\begin{picture}(14,2)
\put(3, 1){$\sum_{a\in \I} \frac{\dim a_2\dim a}{D^2}$}
\put(6.5,0){\resizebox{1.5cm}{2cm}
{\includegraphics{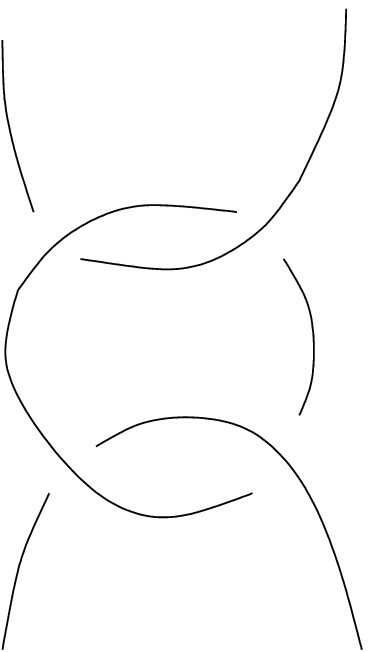}}}
\put(6.1,0){$a_1$}\put(8.1,0){$a'_1$}\put(8,1){$a$}
\put(8.1,1.8){$a_2'$}\put(6,1.8){$a_2$}

\put(9,1){$=$}
\put(9.5,1){$\delta_{a_1a_2}$}
\put(11,0){\resizebox{1cm}{2cm}
{\includegraphics{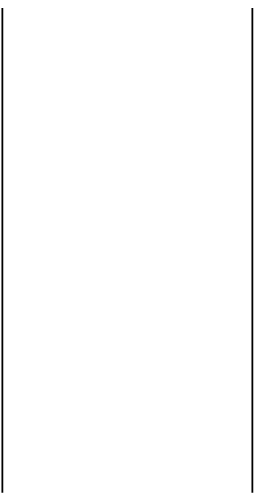}}}
\put(11.2,0){$a_1$}\put(12.1,0){$a'_1$}\put(13,0){.}
\end{picture}
\eeq
}
\end{lemma}

\begin{prop}  {\rm
\beq  \label{See-2-D-2}
(S_e^e)^2 = \frac{1}{D^2}. 
\eeq
}
\end{prop}
\pf
By (\ref{S-a-map}), we have 
\beq  \label{SS-a-map}
\begin{picture}(14,3)
\put(1, 1.4){$\ds S^2(a):$}
\put(3,0.7){\resizebox{1.5cm}{1.5cm}{\includegraphics{Y.eps}}}
\put(2.7, 0.7){$a$}\put(4.6, 0.7){$a_1$}\put(3.3, 2){$a_1$}
\put(3.7, 1){$i$}
\put(5,1.4){$\ds \mapsto$}\put(6,1.4){$\ds \sum_{a_2,a_3\in \I}$}
\put(6.9,1.4){$\ds (S_e^e)^2\dim a_2\dim a_3$}
\put(11,0){\resizebox{2.5cm}{3cm}{\includegraphics{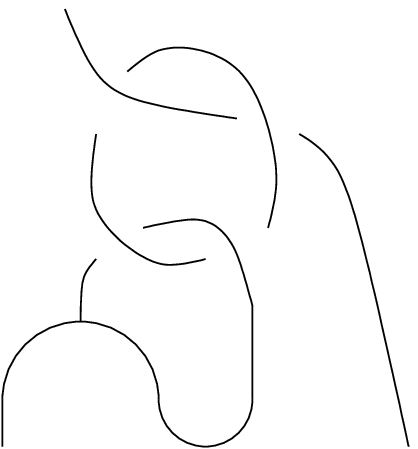}}}
\put(10.6, 0){$\ds a$}\put(11.5,0.5){$i$}
\put(11.1, 1.1){$a_1$}\put(11, 1.8){$a_2$}
\put(13.6, 0){$a_3$}\put(10.9, 2.8){$a_3$}
\put(14,0){.}
\end{picture}
\eeq
Apply (\ref{dual-twist}) to the graph in the right hand side of 
(\ref{SS-a-map}), we obtain
\beq  \label{SS-a-map-deform-1}
\begin{picture}(14,3)
\put(0.5,1){$\ds \sum_{a_2\in \I}\frac{\dim a_2\dim a_3}{D^2}$}
\put(4,0){\resizebox{2.5cm}{3cm}{\includegraphics{car-SS-map.eps}}}
\put(3.6, 0){$\ds a$}\put(4.5,0.5){$i$}
\put(4.1, 1.1){$a_1$}\put(4, 1.8){$a_2$}
\put(4, 2.8){$a_3$}\put(6.6, 0){$a_3$}

\put(7,1){$=$}
\put(7.5,1){$\ds \sum_{a_2\in \I}\frac{\dim a_2\dim a_3}{D^2}$}

\put(11,0){\resizebox{2.5cm}{3cm}{\includegraphics{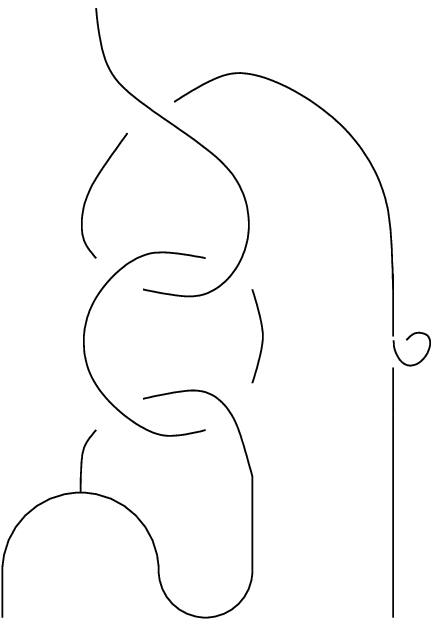}}}
\put(10.6,0){$a$}\put(11,0.8){$a_1$}\put(10.9, 1.4){$a_2$}
\put(11.1,2.8){$a_3$}\put(13.5,0){$a_3$}\put(11.4,0.1){$i$}
\put(14,0){.}
\end{picture}
\eeq
By (\ref{BK-lemma}), the right hands side of 
(\ref{SS-a-map-deform-1}) equals to 
\beq
\begin{picture}(14,2.5)
\put(1,1){$\delta_{a_1a_3}$}
\put(2.5,0){\resizebox{2.5cm}{2.5cm}
{\includegraphics{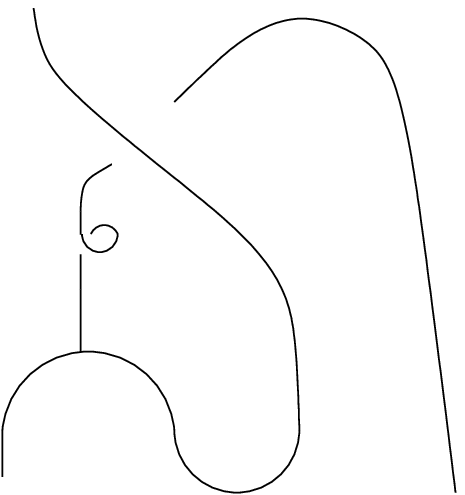}}}
\put(2.1,0){$a$}\put(2.2,2.3){$a'_1$}\put(5.1, 0){$a'_1$}
\put(2.9,0.3){$i$}

\put(6, 1){$=$}

\put(7,1){$\delta_{a_1a_3}$}
\put(8,0){\resizebox{2.5cm}{2.5cm}
{\includegraphics{Y-A--1.eps}}}

\put(7.7,0){$a$}\put(9.2,0){$a'_1$}\put(10.6, 2.3){$a'_1$}
\put(9.5,1.4){$i$}

\put(11,0){.}
\end{picture}
\eeq
By (\ref{SS-A-0}) and (\ref{tilde-A-0-graph}), 
we obtain (\ref{See-2-D-2}). 
\epf

So far, we have determined $S_e^e$ up to a sign. 
Now we consider the relation between $S$ and another generator
of modular group
$T: \tau \mapsto \tau +1$. 
We define a $T$-action on $\Psi_1(\Y_{aa_1;i}^{a_1})$ as follow: 
\beq
T(\Psi_1(\Y_{aa_1;i}^{a_1}))(w_a; z, \tau) :=  
\Psi_1(\Y_{aa_1;i}^{a_1})(w_a; z, \tau+1).
\eeq
It is clear that this action induces an action of $T$ on
$\V_{aa_1}^{a_1}$ for all $a, a_1\in \I$, given by
\beq
T|_{\V_{aa_1}^{a_1}} = e^{2\pi i (h_a-\frac{c}{24})}
\eeq
where $h_a$ is the lowest conformal weight of $W_a$. 

\begin{lemma}
{\rm $S$ and $T$ satisfy the following relation:
\beq  \label{S-T-rel}
(T^{-1}S)^3 = S^2 = T^{-1} S^2 T. 
\eeq
}
\end{lemma}
\pf
Let $w_a\in W_a$. We have
\bea
&& (T^{-1}S)^3((\Psi_1(\Y_{aa_1;i}^{a_1;(1)})))(w_a; z, \tau) \nn
&&=\hspace{1cm} 
\Psi_1(\Y_{aa_1;i}^{a_1;(1)}) \left( \tau^{L(0)} 
\left( \frac{-1}{\frac{-1}{\tau-1}-1}\right)^{L(0)} 
\left( \frac{-1}{\tau-1} \right)^{L(0)} w_a ,\, -z, \tau \right)
\eea
Keeping in mind our choice of branch cut. 
Then it is easy to show that
$$
\tau^{L(0)} 
\left( \frac{-1}{\frac{-1}{\tau-1}-1}\right)^{L(0)} 
\left( \frac{-1}{\tau-1} \right)^{L(0)} w_a = e^{\pi iL(0)} w_a.
$$
By (\ref{SS-equ}), we obtain the first equality of (\ref{S-T-rel}).
The proof of the second equality (\ref{S-T-rel}) is similar. 
\epf

\begin{prop}
{\rm Let $p_{\pm}= \sum_{a\in \I} e^{\pm 2\pi i h_a} \dim^2 a$. Then we have
\bea
S_e^e = \frac{1}{p_-} e^{2\pi ic/8} = \frac{1}{p_+} e^{-2\pi ic/8}.
\eea
}
\end{prop}
\pf
In the proof of the Theorem 3.1.16 in \cite{BK2}, 
Bakalov and Kirillov proved an identity, which, in our
own notation, can be written as follow: 
$$
\frac{1}{(S_e^e)^{2}D^{2}} e^{-2\pi i \frac{c}{24}}\, ST^{-1}S = 
 \frac{1}{S_e^e} \frac{p_+}{D^2} e^{2\pi i\frac{2c}{24}} \, TST.
$$
By (\ref{S-T-rel}) and the fact $p_-p_+=D^2$ which is proved
in \cite{BK2}, we simply obtain that 
$$
S_e^e =  \frac{1}{p_-} e^{2\pi ic/8}. 
$$
Using (\ref{See-2-D-2}) and $p_+p_-=D^2$, we also obtain the second
equality. 
\epf

We thus define 
\beq  \label{D-p-+}
D:= p_- e^{-2\pi ic/8}=p_+e^{2\pi ic/8}.
\eeq 
Notice that this notation is
compatible with the definition of $D^2$.  
Then the action of modular transformation $S(a)$ on $\oplus_{a_1\in \I} 
\hom_V(W_a\boxtimes W_{a_1}, W_{a_1})$
can be expressed graphically as follow: 
\beq   \label{S-a-map-1}
\begin{picture}(14,2)
\put(2, 1){$\ds S(a):$}
\put(4,0.2){\resizebox{1.5cm}{1.5cm}{\includegraphics{Y.eps}}}
\put(3.7, 0.2){$a$}\put(5.6, 0.2){$a_1$}\put(4.3, 1.5){$a_1$}
\put(4.7, 0.5){$i$}
\put(6,1){$\ds \mapsto$}\put(7,1){$\sum_{a_2\in \I}\frac{\dim a_2}{D}$}
\put(9.7,0){\resizebox{2.2cm}{2.2cm}
{\includegraphics{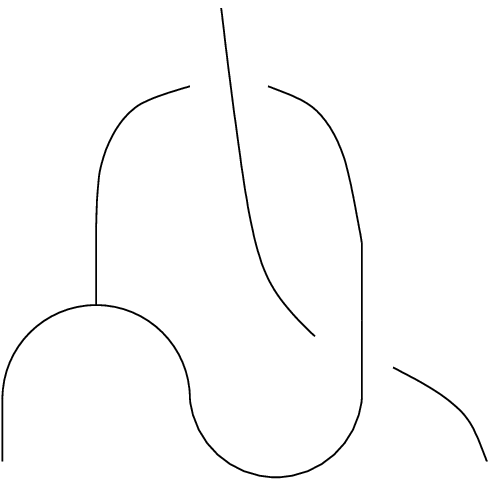}}}
\put(9.4, 0){$a$}\put(10,0.5){$i$}
\put(9.6, 1){$a_1$}\put(12.1, 0){$a_2$}
\put(12.7,0){.}
\end{picture}
\eeq

\begin{prop}
\beq   \label{S-a-map-2}
\begin{picture}(14,2)
\put(2, 1){$\ds S^{-1}(a):$}
\put(4,0.2){\resizebox{1.5cm}{1.5cm}{\includegraphics{Y.eps}}}
\put(3.7, 0.2){$a$}\put(5.6, 0.2){$a_1$}\put(4.3, 1.5){$a_1$}
\put(4.7, 0.5){$i$}
\put(6,1){$\ds \mapsto$}\put(7,1){$\sum_{a_2\in \I}\frac{\dim a_2}{D}$}
\put(9.7,0){\resizebox{2.2cm}{2.2cm}
{\includegraphics{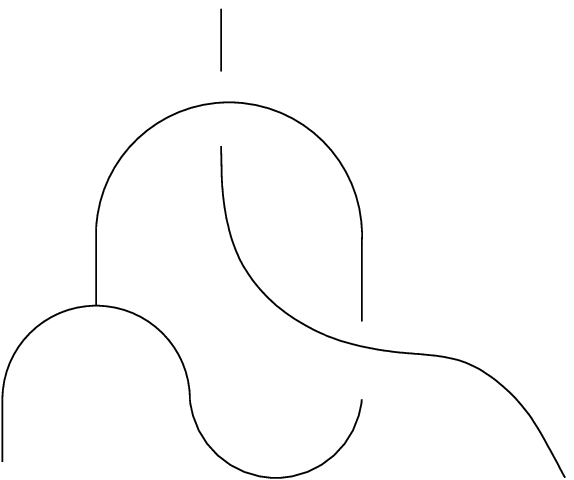}}}
\put(9.4, 0){$a$}\put(10,0.5){$i$}
\put(9.6, 1){$a_1$}\put(12.1, 0){$a_2$}
\put(12.7,0){.}
\end{picture}
\eeq
\end{prop}
\pf
Composing the map (\ref{S-a-map-1}) with (\ref{S-a-map-2}), 
we obtain a map given as follow:
\beq  \label{S-S-graph-1}
\begin{picture}(14,3)
\put(2,0.2){\resizebox{1.5cm}{1.5cm}{\includegraphics{Y.eps}}}
\put(1.7, 0.2){$a$}\put(3.6, 0.2){$a_1$}\put(2.3, 1.5){$a_1$}
\put(2.7, 0.5){$i$}
\put(4,1){$\ds \mapsto$}
\put(5,1){$\sum_{a_2, a_3\in \I}\frac{\dim a_2\dim a_3}{D^2}$}
\put(9,0){\resizebox{3cm}{3cm}
{\includegraphics{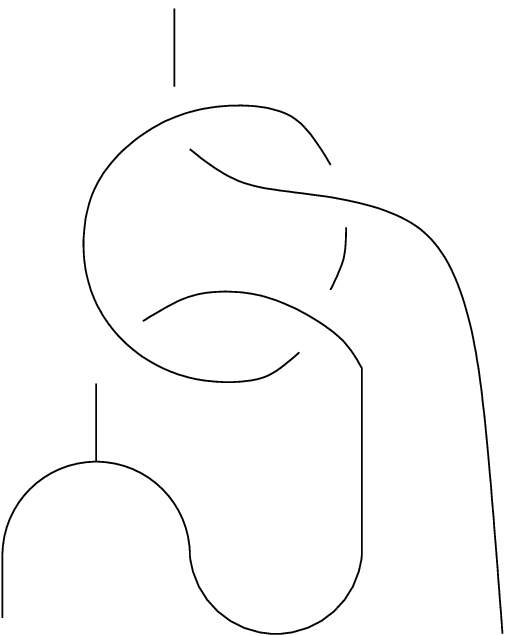}}}
\put(8.7, 0){$\ds a$}\put(9.5,0.5){$i$}
\put(9.1, 1){$a_1$}\put(9, 1.8){$a_2$}
\put(12.1, 0){$a_3$}\put(13,0){.}
\end{picture}
\eeq
Apply (\ref{BK-lemma}) to the graph in (\ref{S-S-graph-1}), 
it is easy to see that above
map is the identity map.
\epf

\begin{rema}
{\rm 
Bakalov and Kirillov obtained the same formula 
(\ref{S-a-map-1}) in \cite{BK2} 
by directly working with modular tensor category and 
solving equations obtained in the so-called Lego-Teichm\"{u}ller game
\cite{BK1}. In our approach, we see the direct link between
the modular transformations of $q$-traces of
the product (or iterate) of intertwining operators 
and their graphic representations in a modular tensor category. 
}
\end{rema}

\begin{prop}
\beq  \label{S-dual-graph-1}
\begin{picture}(14,2.2)
\put(1.5, 1){$\ds (S(a))^*:$}
\put(4,0.2){\resizebox{1.5cm}{1.5cm}{\includegraphics{dual-Y.eps}}}
\put(3.7, 1.3){$a$}\put(5.6, 1.3){$a_1$}\put(4.3, 0.2){$a_1$}
\put(4.7, 1.1){$j$}
\put(6,1){$\ds \mapsto$}\put(7,1){$\sum_{a_2\in \I}\frac{\dim a_2}{D}$}
\put(10.5,0){\resizebox{2.2cm}{2.2cm}
{\includegraphics{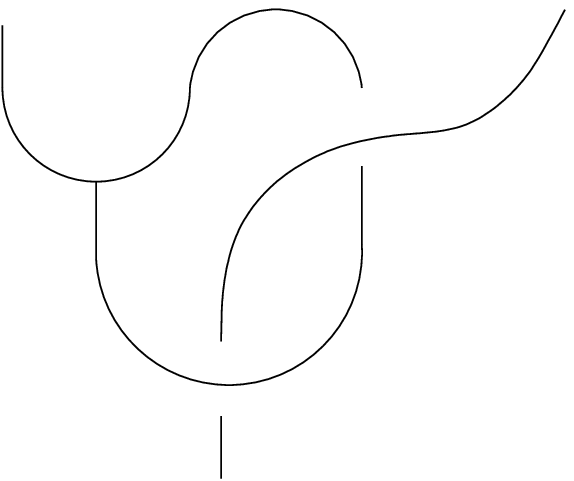}}}
\put(10.2, 2){$a$}\put(10.8,1.6){$j$}
\put(10.4, 1){$a_1$}\put(12.9, 2){$a_2$}
\put(13.5,0){,}
\end{picture}
\eeq
\beq  \label{S-dual-graph-2}
\begin{picture}(14,2.2)
\put(1.5, 1){$\ds (S^{-1}(a))^*:$}
\put(4,0.2){\resizebox{1.5cm}{1.5cm}{\includegraphics{dual-Y.eps}}}
\put(3.7, 1.3){$a$}\put(5.6, 1.3){$a_1$}\put(4.3, 0.2){$a_1$}
\put(4.7, 1.1){$j$}
\put(6,1){$\ds \mapsto$}\put(7,1){$\sum_{a_2\in \I}\frac{\dim a_2}{D}$}
\put(10.5,0){\resizebox{2.2cm}{2.2cm}
{\includegraphics{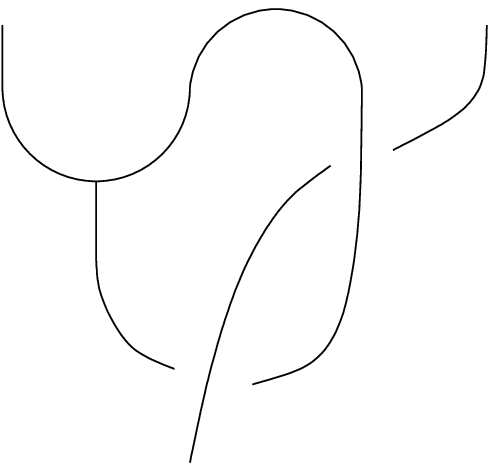}}}
\put(10.2, 2){$\ds a$}\put(10.8,1.6){$j$}
\put(10.4, 1){$a_1$}\put(12.9, 2){$a_2$}
\put(13.5,0){.}
\end{picture}
\eeq
\end{prop}
\pf
We only prove (\ref{S-dual-graph-1}). The proof 
of (\ref{S-dual-graph-2}) is analogous to that of 
(\ref{S-dual-graph-1}). 

It is enough to show that the pairing between the image of
(\ref{S-a-map}) and that of (\ref{S-dual-graph-1}) 
still gives $\delta_{ij}$. This can be proved as follow:
\beq  \label{S-dual-S-equ-1}
\begin{picture}(14,3)
\put(0.2, 1.5){$\sum_{a_2} \frac{\dim a_2}{D^2} $}
\put(2.5,0){\resizebox{3cm}{3cm}{\includegraphics{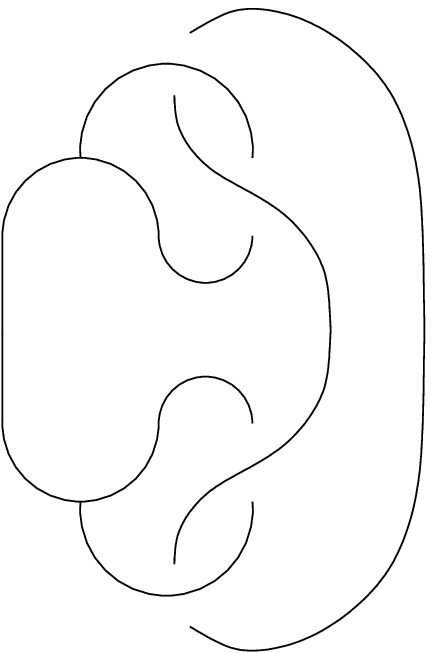}}}
\put(2.2,1.5){$a$}\put(2.7,2.5){$a_1$}\put(2.7,0.4){$a_1$}
\put(4.3,1.5){$a_2$}\put(3,1.9){$i$}\put(3,1){$j$}

\put(6.2,1.5){$=$}\put(6.7,1.5)
{$\sum_{a_2} \frac{\dim a_1\dim a_2}{\dim a_1D^2}$}  
\put(10,0)
{\resizebox{2.5cm}{3cm}{\includegraphics{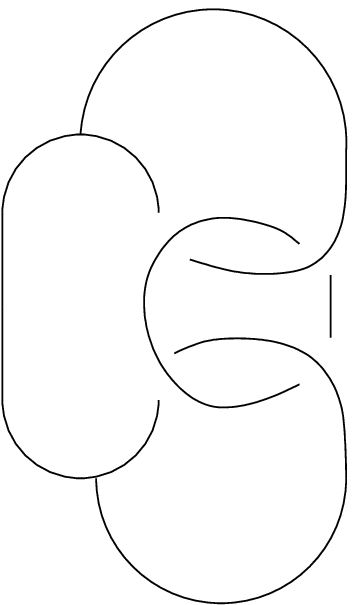}}}
\put(9.7,1.5){$a$}\put(10.1,2.6){$a_1$}\put(10.2,0.4){$a_1$}
\put(10.5,1.5){$a_2$}\put(10.5,1.9){$i$}\put(10.5,0.9){$j$}

\end{picture}
\eeq
By (\ref{BK-lemma}), the right hand side of 
(\ref{S-dual-S-equ-1}) equals to
\beq
\begin{picture}(14,2)
\put(2.5, 0.8){$\frac{1}{\dim a_1}$}
\put(4,0){\resizebox{2cm}{2cm}{\includegraphics{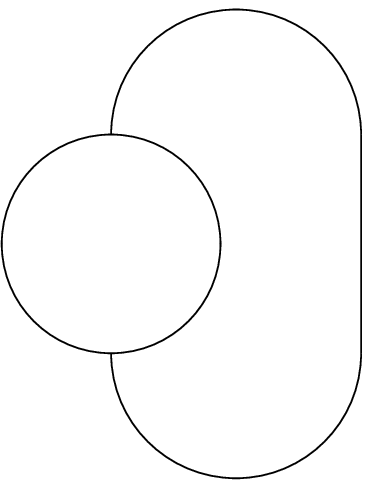}}}
\put(3.7,1){$a$}\put(5.3,1){$a_1$}
\put(4.2,1.7){$a_2$}\put(4.5,1.1){$i$}\put(4.5,0.7){$j$}

\put(6.7, 0.8){$=\, \, \delta_{ij}.$}
\end{picture}   \nonumber
\eeq
\epf

\renewcommand{\theequation}{\thesection.\arabic{equation}}
\renewcommand{\thethm}{\thesection.\arabic{thm}}
\setcounter{equation}{0}
\setcounter{thm}{0}

\section{Categorical formulations and constructions}

In this section, we give a categorical formulation of 
modular invariant conformal full field algebra over $V^L\otimes V^L$,
open-string vertex operator algebra over $V$ equipped with
nondegenerate invariant bilinear forms and Cardy condition. 
Then we introduce a notion called
Cardy $\mathcal{C}_V|\mathcal{C}_{V\otimes V}$-algebra. 
In the end, we give a categorical construction of such algebra
in the Cardy case \cite{FFFS2}.

\subsection{Modular invariant $\mathcal{C}_{V^L\otimes V^R}$-algebras}

We first recall the notion of coalgebra and 
Frobenius algebra (\cite{FS}) in a tensor category.  
\begin{defn}
{\rm 
A coalgebra $A$ in a tensor category 
$\mathcal{C}$ is an object with a
coproduct $\Delta\in \text{Mor}(A, A\otimes A)$ and a counit 
$\epsilon \in \text{Mor}(A, \one_{\mathcal{C}})$ such that 
\beq  \label{co-alg-def}
(\Delta \otimes \id_A) \circ \Delta = 
(\id_A \otimes \Delta) \circ \Delta, 
\hspace{1cm} (\epsilon \otimes \id_A)\circ \Delta 
= \id_A = (\id_A \otimes \epsilon)
\circ \Delta,
\eeq
which can also be expressed in term of the following
graphic equations:
\beq 
\begin{picture}(14,2)
\put(0.5,0){\resizebox{13cm}{2cm}{\includegraphics{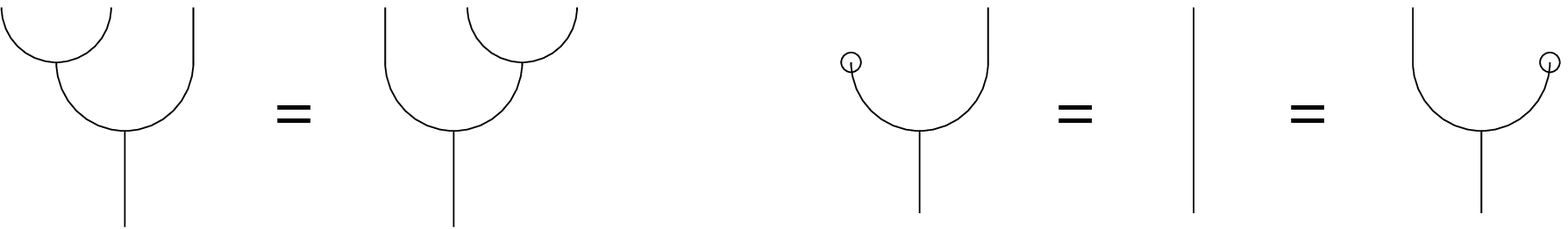}}}
\end{picture}  \nonumber
\eeq
}
\end{defn}

\begin{defn}
{\rm 
Frobenius algebra in $\mathcal{C}$ is 
an object that is both an algebra and a coalgebra
and for which the product and coproduct are related by 
\begin{equation}  \label{Frob-def-equ}
(\id_A \otimes m) \circ (\Delta \otimes \id_A) = \Delta \circ m = 
(m\otimes \id_A) \circ (\id_A \otimes \Delta),
\end{equation}
or as the following graphic equations, 
\beq \label{Frob-alg-def-fig}
\epsfxsize  0.7\textwidth
\epsfysize  0.15\textwidth
\epsfbox{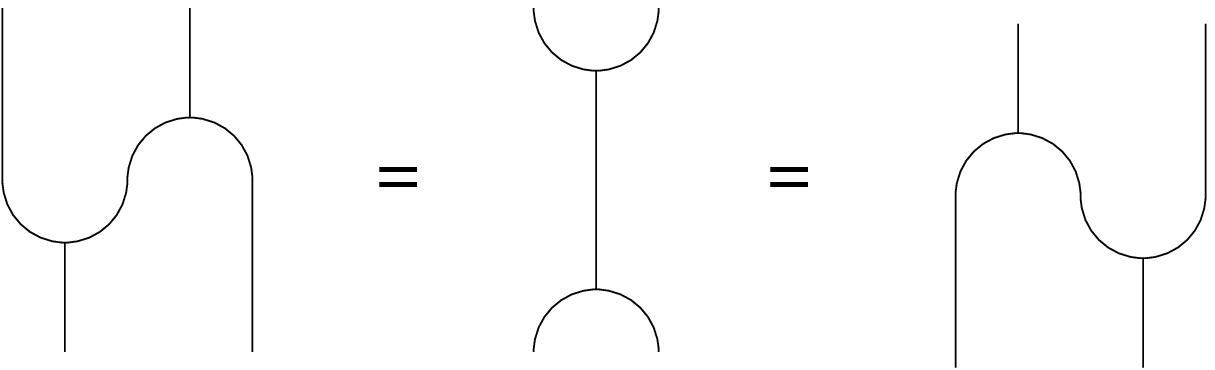}.
\eeq
A Frobenius algebra is called {\it symmetric} if the following
condition is satisfied. 
\beq  \label{symm-Frob-fig}
\begin{picture}(14,2)
\put(4.5,0){\resizebox{2cm}{2cm}{\includegraphics{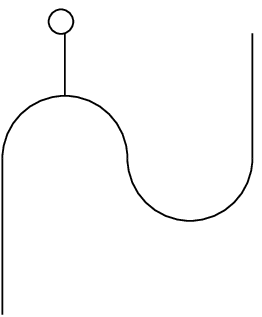}}}
\put(7.4,0.9){$=$}
\put(8.5,0){\resizebox{2cm}{2cm}{\includegraphics{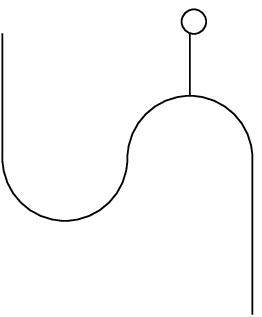}}}
\put(11,0){.}
\end{picture} 
\eeq
}
\end{defn}

Let $V^L$ and $V^R$ be vertex operator algebras satisfying
the conditions in Theorem \ref{ioa}. Then the vertex operator
algebra $V^L\otimes V^R$ also satisfies the conditions
in Theorem \ref{ioa} \cite{HKo1}. 
Thus $\mathcal{C}_{V^L\otimes V^R}$ also 
has a structure of modular tensor category. 
In particular, we choose the braiding structure on 
$\mathcal{C}_{V^L\otimes V^R}$ to be $\mathcal{R}_{+-}$ 
which is defined in \cite{Ko2}.
The twist $\theta_{+-}$ for each 
$V^L\otimes V^R$-module is defined by 
\beq
\theta_{+-} = e^{-2\pi iL^L(0)}\otimes e^{2\pi iL^R(0)}. 
\eeq
Duality maps are naturally induced from those 
of $\mathcal{C}_{V^L}$ and $\mathcal{C}_{V^R}$.

The following theorem is proved in \cite{Ko1}. 
\begin{thm}  \label{ffa-cat}
The category of conformal full field algebras 
over $V^L\otimes V^R$ 
equipped with nondegenerate invariant bilinear forms 
is isomorphic to the category of commutative Frobenius algebra
in $\mathcal{C}_{V^L\otimes V^R}$ with a trivial twist. 
\end{thm}

\begin{rema}
{\rm
In a ribbon category, it was proved in \cite{FFRS} 
that a commutative Frobenius algebra
with a trivial twist is equivalent to 
a commutative symmetric Frobenius algebra. 
}
\end{rema}

Let $\I^L$ and $\I^R$ 
denote the set of equivalent class of irreducible
$V^L$-modules and $V^R$-modules respectively. 
We use $a$ and $a_i$ for $i\in \N$ 
to denote elements in $\I^L$ and we use $e$ to
denote the equivalent class of $V^L$. 
We use $\bar{a}$ and $\bar{a}_i$ for $i\in \N$ 
to denote elements in $\I^R$, 
and $\bar{e}$ to denote the equivalent class of $V^R$. 
For each $a\in \I^L$ ($\bar{a}\in \I^R$), 
we choose a representative $W_{a}$ ($W_{\bar{a}}$). 
We denote the vector space of intertwining operators of type 
$\binom{W_{a_3}}{W_{a_1}W_{a_2}}$ and 
$\binom{W_{\bar{a}_3}}{W_{\bar{a}_1}W_{\bar{a}_2}}$
as $\V_{a_1a_2}^{a_3}$ and 
$\bar{\V}_{\bar{a}_1\bar{a}_2}^{\bar{a}_3}$ respectively,
the fusion rule as $N_{a_1a_2}^{a_3}$
and  $N_{\bar{a}_1\bar{a}_2}^{\bar{a}_3}$ respectively.

A conformal full field algebra over $V^L\otimes V^R$, 
denoted as $A_{cl}$, is
a direct sum of irreducible modules of $V^L\otimes V^R$, i.e.
\bea
A_{cl} = \oplus_{\alpha=1}^{N} W_{r(\alpha)} \otimes W_{\bar{r}(\alpha)},
\eea
where $r: \{1, \dots, N\} \rightarrow \I^L$ and
$\bar{r}: \{1, \dots, N\} \rightarrow \I^R$, for some $N\in \Z_+$. 
Let $\{ e_{ab;i}^{c}\}$ and $\{ \bar{e}_{\bar{a}\bar{b};j}^{\bar{c}}\}$ 
be basis for $\V_{ab}^c$ and $\bar{\V}_{\bar{a}\bar{b}}^{\bar{c}}$,
and $\{ f^{ab}_{c;i}\}$ and $\{ \bar{f}^{\bar{a}\bar{b}}_{\bar{c};j}\}$
be the dual basis respectively. Then the vertex operator $\mathbb{Y}$
can also be expanded as follow:
\beq
\mathbb{Y} = \sum_{\alpha, \beta, \gamma}\sum_{i,j} 
d_{\alpha\beta}^{\gamma}(f_{r(\gamma); i}^{r(\alpha)r(\beta)}, 
\bar{f}_{\bar{r}(\gamma); j}^{\bar{r}(\alpha)\bar{r}(\beta)} )
\, e^{r(\gamma); i}_{r(\alpha)r(\beta)} \otimes
\bar{e}^{\bar{r}(\gamma); j}_{\bar{r}(\alpha)\bar{r}(\beta)},
\eeq 
where $d_{\alpha\beta}^{\gamma}$ defines a bilinear map
$$
(\V^{r(\gamma)}_{r(\alpha)r(\beta)})^*\otimes 
(\bar{\V}^{\bar{r}(\gamma)}_{\bar{r}(\alpha)\bar{r}(\beta)})^* \rightarrow \C
$$
for all $\alpha, \beta, \gamma =1, \dots, N$ for some $N\in \N$. 

Since the trace function pick out $\gamma=\beta$ terms,
we define $\mathbb{Y}_{\alpha}^{diag}$ by
\beq \label{Y-diag-alpha}
\mathbb{Y}_{\alpha}^{diag} := \sum_{\beta}\sum_{i,j} 
d_{\alpha\beta}^{\beta}(f_{r(\beta); i}^{r(\alpha)r(\beta)}, 
\bar{f}_{\bar{r}(\beta); j}^{\bar{r}(\alpha)\bar{r}(\beta)} )
\, e^{r(\beta); i}_{r(\alpha)r(\beta)} \otimes
\bar{e}^{\bar{r}(\beta); j}_{\bar{r}(\alpha)\bar{r}(\beta)}.
\eeq
Let $\mathbb{Y}^{diag} := \sum_{\alpha} \mathbb{Y}_{\alpha}^{diag}$.
Of course, it is obvious to see that such defined $\mathbb{Y}^{diag}$
is independent of the choice of basis. We denote the representation of 
the modular transformation $S: \tau \mapsto \frac{-1}{\tau}$
on $\oplus_{b\in \I^L} \V_{ab}^b$ and 
$\oplus_{\bar{b}\in \I^R} \V_{\bar{a}\bar{b}}^{\bar{b}}$, 
by $S^L(a)$ and $S^R(\bar{a})$ respectively.

In \cite{HKo3}, we defined the notion of modular invariant
conformal full field algebra over $V^L\otimes V^R$ (see
\cite{HKo3} for the precise definition). It basically
means that $n$-point genus-one correlation functions built out of
$q$-$\bar{q}$-traces are invariant under the action of 
modular group $SL(2, \Z)$ for all $n\in \N$.  
Moreover, we proved the following results in \cite{HKo3}. 
\begin{prop}
$A_{cl}$, a conformal full field algebra over 
$V^L\otimes V^R$, is modular invariant if it 
satisfies $c^L-c^R=0 \,\,\, \text{\rm mod}\,\, 24$ and
\beq  \label{S-S-1-Y-diag}
S^L(r(\alpha)) \otimes (S^R(\bar{r}(\alpha)))^{-1}  
\mathbb{Y}_{\alpha}^{diag}  = \mathbb{Y}_{\alpha}^{diag}
\eeq
for all $\alpha=1, \dots, N$. 
\end{prop}

We denote the morphism in
$$
\text{Hom}_{C_{V^L\otimes V^R}} \big( (W_{r(\alpha)} \otimes W_{\bar{r}(\alpha)}) 
\boxtimes
(W_{r(\beta)} \otimes W_{\bar{r}(\beta)}), W_{r(\beta)} \otimes W_{\bar{r}(\beta)} 
\big)
$$
which corresponds to $\mathbb{Y}_{\alpha}^{diag}$ by 
$m_{\mathbb{Y}_{\alpha}^{diag}}$. 
Then (\ref{S-S-1-Y-diag}) is equivalent to the following categorical condition:
\beq  \label{S-S-2-Y-diag}
S^L(r(\alpha)) \otimes (S^R(\bar{r}(\alpha)))^{-1}  
m_{\mathbb{Y}_{\alpha}^{diag}}= m_{\mathbb{Y}_{\alpha}^{diag}}
\eeq
for all $\alpha=1, \dots, N$. 

Now we choose a basis $\{ b_{a\otimes \bar{a}; i}^{A_{cl}} \}$ 
of $\hom_{C_{V^L\otimes V^R}}(W_a\otimes W_{\bar{a}}, A_{cl})$
and a basis  $\{ b^{a\otimes \bar{a}; i}_{A_{cl} } \}$  
of $\hom_{C_{V^L\otimes V^R}}(A_{cl}, W_a\otimes W_{\bar{a}})$ as follow: 
\beq
\begin{picture}(14,2)
\put(3,0.9){$b_{a\otimes b; i}^{A_{cl}} :=$}
\put(4.8,0){\resizebox{0.3cm}{2cm}
{\includegraphics{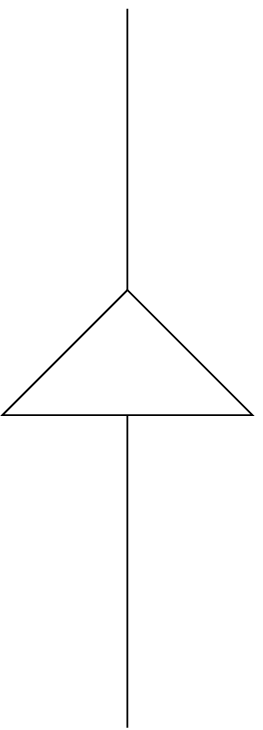}}}
\put(5.1,0){$a\otimes \bar{a}$}\put(5.1, 1.8){$A_{cl}$}\put(5.2, 0.9){$i$}
\put(6.3,0){,}

\put(7.5,0.9){$b^{a\otimes b; i}_{A_{cl}}:=$}
\put(9.3,0){\resizebox{0.3cm}{2cm}
{\includegraphics{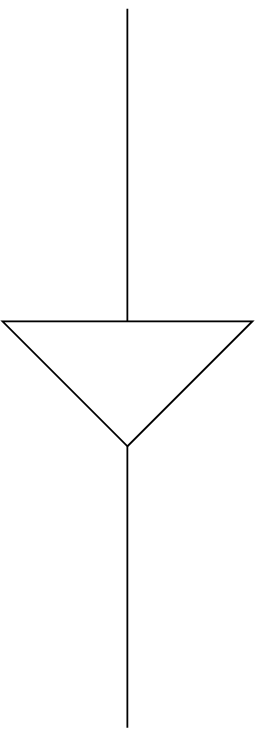}}}
\put(9.6,0){$A_{cl}$}\put(9.6, 1.8){$a\otimes \bar{a}$}\put(9.7, 0.9){$i$}
\put(10.3,0){,}
\end{picture}
\eeq
satisfying the following conditions: for all $a\in \I^L, \bar{a}\in \I^R$,  
\beq
\begin{picture}(14,3)
\put(3.5,0){\resizebox{0.4cm}{3cm}
{\includegraphics{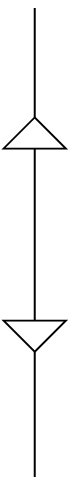}}}
\put(4.2, 1.4){$=\, \, \delta_{ij}$}
\put(2.6,2.8){$a\otimes \bar{a}$}
\put(2.6,0){$a\otimes \bar{a}$} \put(3, 1.4){$A_{cl}$}
\put(4, 2.1){$j$}\put(4,0.9){$i$}

\put(6.4,1.4){and}

\put(8,1.4){$\sum_{a, \bar{a},i}$}
\put(9.5,0){\resizebox{0.4cm}{3cm}
{\includegraphics{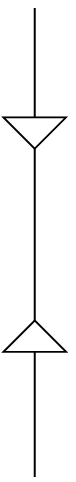}}}
\put(9.8,2.8){$A_{cl}$}\put(9.8,0){$A_{cl}$}
\put(9.8,1.4){$a\otimes \bar{a}$}
\put(9.2,0.8){$i$}\put(9.2,2.2){$i$}

\put(11.2,1.4){$=$} 
\put(12.2,0){\resizebox{0.1cm}{3cm}
{\includegraphics{id.eps}}}
\put(12.4,0){$A_{cl}$}
\put(12.4,2.8){$A_{cl}$}
\put(13,0){.}

\end{picture}
\eeq
Then the condition (\ref{S-S-2-Y-diag}) can be further expressed graphically 
as follow: 
\beq \label{mod-inv-A-cl}
\begin{picture}(14, 3)

\put(1, 1.4){$\frac{\dim a\otimes \bar{a}}{D^LD^R}$}

\put(3,0){\resizebox{3cm}{3cm}
{\includegraphics{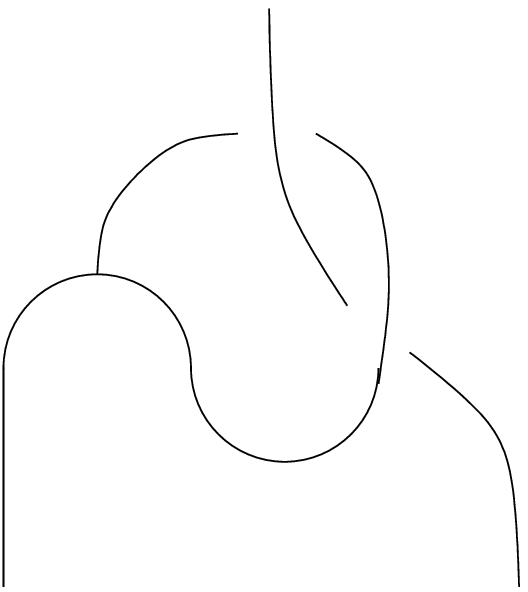}}}
\put(2.3, 0){$A_{cl}$}\put(3, 2){$A_{cl}$} \put(4.2, 1.3){$A_{cl}$}
\put(6.1,0){$a\otimes \bar{a}$}\put(4.7, 2.8){$a\otimes \bar{a}$}

\put(7,1.4){$=$}

\put(7.5, 1.4){$\sum_{r(k)=a\otimes \bar{a}}$}

\put(9.5,0){\resizebox{2cm}{3cm}
{\includegraphics{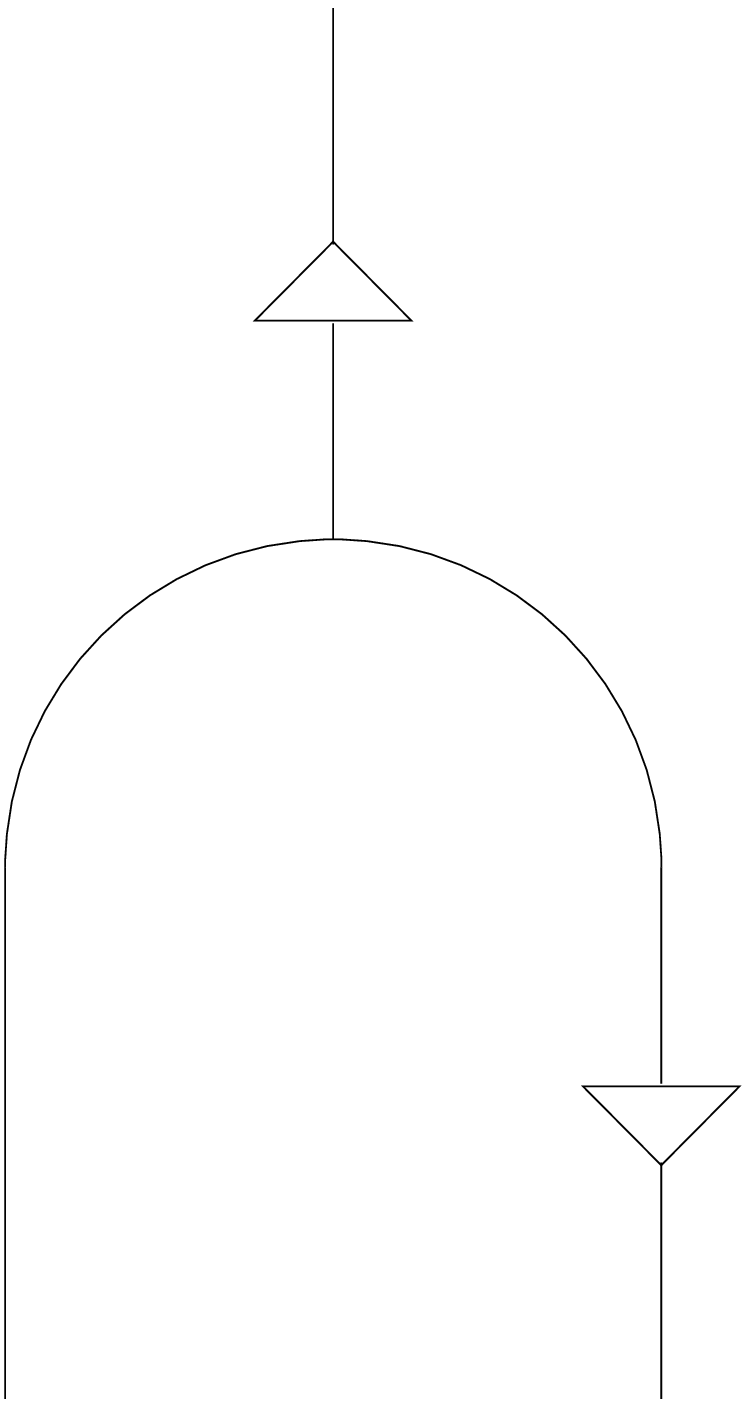}}}
\put(8.8,0){$A_{cl}$}
\put(11.4,0){$a\otimes \bar{a}$}\put(10.5,2.9){$a\otimes \bar{a}$}
\put(11.6,0.5){$k$}\put(10.7, 2.4){$k$}\put(9.7, 2){$A_{cl}$}
\put(11.4, 1.2){$A_{cl}$}

\end{picture}
\eeq
for all $a\in \I^L, \bar{a}\in \I^R$,  
where $D^L$ and $D^R$ are the $D$ defined by (\ref{D-p-+}) in
$\mathcal{C}_{V^L}$ and $\mathcal{C}_{V^R}$ respectively.

\begin{defn}
{\rm 
Let $V^L$ and $V^R$ be so that 
$c^L-c^R=0 \,\,\, \text{\rm mod}\,\, 24$. 
A modular invariant $\mathcal{C}_{V^L\otimes V^R}$-algebra
is an associative algebra $(A, \mu_A, \iota_A)$ satisfying
the condition (\ref{mod-inv-A-cl}).
}
\end{defn}

Some properties of modular invariant algebra follow immediately from
above definition, such as the following famous condition \cite{MSei3}: 
\beq
\sum_{a\in \I^L, \bar{a}\in \I^R} (S^L(e))_{ab} \,\,  N_{b\bar{b}}\,\,  
(S^R(\bar{e}))_{\bar{b}\bar{a}}^{-1} = 
N_{a\bar{a}}
\eeq
where $N_{a\bar{a}}$ is the multiplicity of $W_a\otimes W_{\bar{a}}$ in $A_{cl}$
for $a\in \I^L, \bar{a}\in \I^R$. We leave a 
more systematic study of modular invariant
$\mathcal{C}_{V^L\otimes V^R}$-algebras to elsewhere.

\begin{thm}
Let $V^L$ and $V^R$ be so that 
$c^L-c^R=0 \,\,\, \text{\rm mod}\,\, 24$. 
The following two notions are equivalent: 
\bnu
\item Modular invariant conformal
full field algebra over $V^L\otimes V^R$
equipped with a nondegenerate invariant bilinear form.   
\item Modular invariant 
commutative Frobenius algebra with a trivial twist. 
\enu
\end{thm}
\pf
The Theorem follows from Theorem \ref{ffa-cat} and 
the equivalence between 
(\ref{S-S-1-Y-diag}) and (\ref{mod-inv-A-cl}) immediately. 
\epf

\begin{rema} {\rm
In the case $V^L\cong \C \cong V^R$, a modular invariant 
commutative Frobenius algebra with a trivial twist in $\mathcal{C}_{V^L\otimes V^R}$
is simply a commutative Frobenius algebra over $\C$,
which is equivalent to a 2-dimensional topological field theory 
(see for example \cite{BK2}). 
In this case, the modular invariance condition holds automatically. 
}
\end{rema}

\subsection{Cardy $\mathcal{C}_V|\mathcal{C}_{V\otimes V}$-algebras}

For an open-string vertex operator algebra $V_{op}$ over $V$
equipped with a nondegenerate invariant bilinear form 
$(\cdot, \cdot)_{op}$, there is an isomorphism 
$\varphi_{op}: V_{op} \rightarrow V_{op}$ induce from 
$(\cdot, \cdot)_{op}$ (recall (\ref{varphi-form-op})).

In this case, $V_{op}$ is a $V$-module and $Y_{op}^f$ is 
an intertwining operator. 
By comparing (\ref{inv-form-osvoa-1}) with (\ref{omega-r-A-r}), 
and (\ref{inv-form-osvoa-2}) with (\ref{A-r-omega-r}), 
we see that the conditions (\ref{inv-form-osvoa-1}) and
(\ref{inv-form-osvoa-2}) can be rewritten as 
\bea
Y_{op}^f &=& \varphi_{op}^{-1} \circ
\sigma_{123}(Y_{op}^f) \circ (\varphi_{op} \otimes \id_{V_{op}})  
\label{sigma-123-varphi-op}  \\
&=& \varphi_{op}^{-1} \circ \sigma_{132}(Y_{op}^f) \circ 
(\id_{op} \otimes \varphi_{op}). \label{sigma-132-varphi-op}
\eea

\begin{rema}  {\rm
The representation theory of open-string vertex operator algebra
can be developed.  
In that context, $\sigma_{123}(Y_{op}^f)$ gives $V'_{op}$ a
right $V_{op}$-module structure and the equation 
(\ref{sigma-123-varphi-op}) is equivalent to the statement 
that $\varphi_{op}$ is an isomorphism between two 
right $V_{op}$-modules. Similarly, $\sigma_{132}(Y_{op}^f)$ 
gives $V'_{op}$ a left $V_{op}$-module structure and the equation 
(\ref{sigma-132-varphi-op}) is equivalent to the statement 
that $\varphi_{op}$ is an isomorphism between two 
left $V_{op}$-modules. But we do not need it in this work.
}
\end{rema}

\begin{thm}
The category of open-string vertex operator algebras over $V$ 
equipped with a nondegenerate invariant bilinear form is 
isomorphic to the category of symmetric 
Frobenius algebras in $\mathcal{C}_V$. 
\end{thm}
\pf
We have already shown in \cite{HKo1} that an open-string vertex 
operator algebra over $V$ is equivalent to an associative algebra
in $\mathcal{C}_V$. 

Let $V_{op}$ be an open-string vertex operator algebra over $V$.  
Giving a nondegenerate invariant bilinear form 
(recall (\ref{inv-form-osvoa-1}) and (\ref{inv-form-osvoa-2})) 
on $V_{op}$ is equivalent to give an isomorphism 
$\varphi_{op}: V_{op}\rightarrow V'_{op}$ satisfying the conditions 
(\ref{sigma-123-varphi-op}) and 
(\ref{sigma-132-varphi-op}). If we define 
\beq  \label{F-r-act-F-dual}
\begin{picture}(14,2)
\put(3,0){\resizebox{1.5cm}{2cm}
{\includegraphics{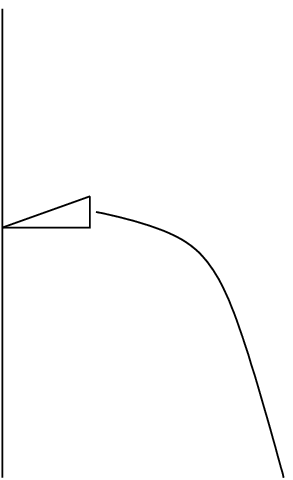}}}
\put(2.4,0){$V'_{op}$}\put(4.6,0){$V_{op}$}
\put(6,1){$:=$}
\put(7.5,0){\resizebox{2.5cm}{2cm}
{\includegraphics{sigma-123.eps}}}
\put(6.9,0){$V'_{op}$}\put(8.7,0){$V_{op}$}\put(10.1,1.8){$V'_{op}$}
\put(11,0){,}
\end{picture}
\eeq
\beq \label{F-l-act-F-dual}
\begin{picture}(14,2)
\put(3,0){\resizebox{1.5cm}{2cm}
{\includegraphics{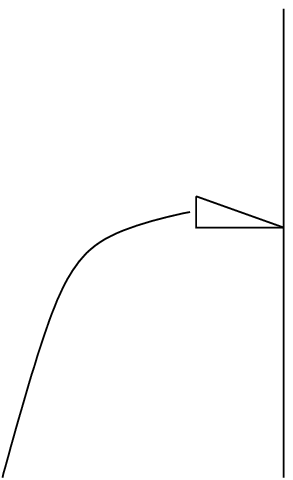}}}
\put(2.4,0){$V_{op}$}\put(4.6,0){$V'_{op}$}
\put(6,1){$:=$}
\put(7.5,0){\resizebox{2.5cm}{2cm}
{\includegraphics{sigma-132.eps}}}
\put(6.9,1.8){$V'_{op}$}\put(8.2,0){$V_{op}$}\put(10.1,0){$V'_{op}$}
\put(11,0){,}
\end{picture}
\eeq
then (\ref{sigma-123-varphi-op}) and 
(\ref{sigma-132-varphi-op}) can be rewritten as
\beq  \label{F-r-act-F-varphi}
\begin{picture}(14,2)
\put(3,0){\resizebox{2cm}{2cm}{\includegraphics{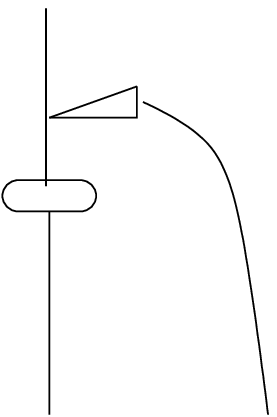}}}
\put(2.7,0){$V_{op}$}\put(5.1,0){$V_{op}$}\put(2.7,1.8){$V'_{op}$}
\put(2.3,1.1){$\varphi_{op}$}
\put(6,1){$=$}
\put(7.5,0){\resizebox{1.5cm}{2cm}
{\includegraphics{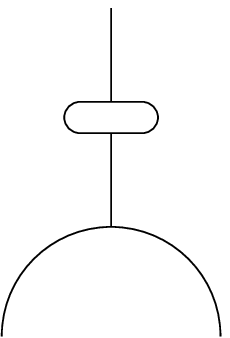}}}
\put(6.9,0){$V_{op}$}\put(9.1,0){$V_{op}$}\put(7.5,1.8){$V'_{op}$}
\put(8.7, 1.3){$\varphi_{op}$}
\put(10,0){.}
\end{picture}
\eeq
\beq  \label{F-l-act-F-varphi}
\begin{picture}(14,2)
\put(3,0){\resizebox{2cm}{2cm}{\includegraphics{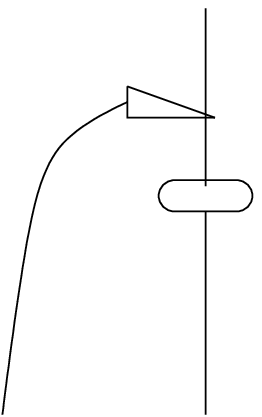}}}
\put(2.4,0){$V_{op}$}\put(4.8,0){$V_{op}$}\put(4.8,1.8){$V'_{op}$}
\put(5.1, 1){$\varphi_{op}$}
\put(6,1){$=$}
\put(7.5,0){\resizebox{1.5cm}{2cm}{\includegraphics{op-varphi-2.eps}}}
\put(6.9,0){$V_{op}$}\put(9,0){$V_{op}$}\put(7.5,1.8){$V'_{op}$}
\put(8.7, 1.3){$\varphi_{op}$}
\put(10,0){,}
\end{picture}
\eeq
Using the map $\varphi_{op}$ and its inverse, we can obtain
a natural coalgebra structure on $F$ defined as follow: 
\beq  \label{co-F}
\begin{picture}(14,2)
\put(2,0.7){$\Delta_{V_{op}}$}
\put(3,0.7){$=$}
\put(4,0){\resizebox{3cm}{2cm}{\includegraphics{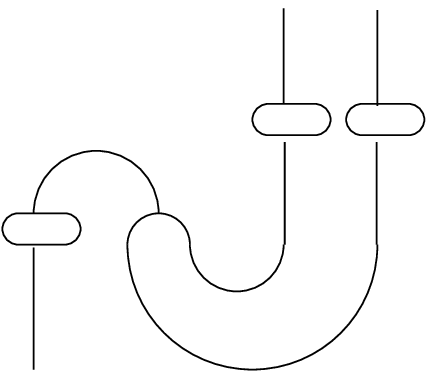}}}
\put(3.6,0){$V_{op}$}\put(5.3,1.8){$V_{op}$}\put(6.8,1.8){$V_{op}$}
\put(3.4,0.7){$\varphi_{op}$}
\put(5, 1.2){$\varphi_{op}^{-1}$}\put(7, 1.2){$\varphi_{op}^{-1}$}
\put(9,0.7){$\epsilon_{V_{op}}$}
\put(10,0.7){$=$}
\put(11,0){\resizebox{1cm}{1.5cm}{\includegraphics{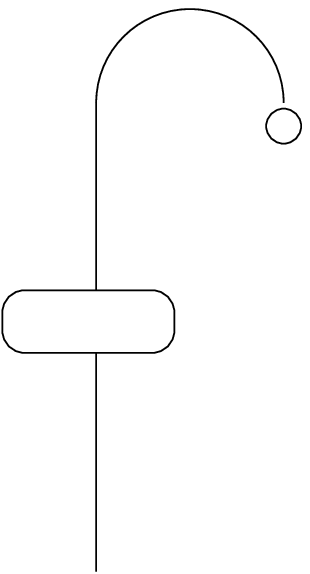}}}
\put(11.7, 0.6){$\varphi_{op}$}
\put(10.7,0){$V_{op}$}\put(13,0){.}
\end{picture}
\eeq
Similar to the proof of \cite[Thm.\, 4.15]{Ko1}, 
(\ref{F-r-act-F-varphi}) and (\ref{F-l-act-F-varphi}) imply
that such defined $\Delta_{V_{op}}$ and $\epsilon_{V_{op}}$
give $V_{op}$ a Frobenius algebra structure. 
Moreover, we also showed in \cite{Ko1} that (\ref{F-r-act-F-varphi}) 
implies the equality between $\varphi_{op}$ and
the left hand side of (\ref{symm-Frob-fig}).
Similarly, using (\ref{F-l-act-F-varphi}), we can show that
the right hand side of (\ref{symm-Frob-fig}) also 
equals to $\varphi_{op}$. Thus $V_{op}$ has a structure of 
symmetric Frobenius algebra.

We thus obtain a functor from the first category to the 
second category.

Conversely, given a symmetric Frobenius algebra in $\mathcal{C}_V$. 
In \cite{HKo1}, we showed that it gives an open-string vertex operator
algebra over $V$. It is shown in \cite{FRS2} that 
either side of (\ref{symm-Frob-fig}) is an isomorphism.
Take $\varphi_{op}$ to be either side of (\ref{symm-Frob-fig}).
Then (\ref{F-r-act-F-varphi}) and (\ref{F-l-act-F-varphi}) 
follow automatically from the
definition (\ref{F-r-act-F-dual}) and (\ref{F-l-act-F-dual}). 
They are nothing but the invariance properties
(recall (\ref{sigma-123-varphi-op})(\ref{sigma-132-varphi-op})) 
of the bilinear form associated with $\varphi_{op}$. 
Thus we obtain a functor from the second category to 
the first category. 

It is routine to check that these two functors are inverse to
each other. 
\epf

Now we consider an open-closed field algebra over $V$ given in 
(\ref{opcl-fa}) equipped with nondegenerate invariant bilinear
forms $(\cdot, \cdot)_{op}$ and $(\cdot, \cdot)_{cl}$. 
We assume that $V_{cl}$ and $V_{op}$ have the following decompositions:
$$
V_{cl} = \oplus_{i=1}^{N_{cl}}\, W_{r_L(i)} \otimes W_{r_R(i)}, \quad\quad\quad
V_{op} =\oplus_{i=1}^{N_{op}}\, W_{r(i)}
$$
where $r_L, r_R: \{ 1, \dots, N_{cl}\}\rightarrow \I$ 
and $r: \{1, \dots, N_{op}\}\rightarrow \I$. 
We denote the embedding 
$b_{(i)}^{V_{op}}:W_{r(i)} \hookrightarrow V_{op}$, 
the projection $b_{V_{op}}^{(i)}:V_{op} \rightarrow W_{r(i)}$ 
and $b_{(i)}^{V_{cl}}: W_{r_L(i)}\boxtimes W_{r_R(i)} \hookrightarrow T(V_{cl})$
by the following graphs: 
\beq
\begin{picture}(14,2)
\put(2,0.9){$b_{(i)}^{V_{op}}=$}
\put(3.5,0){\resizebox{0.3cm}{2cm}
{\includegraphics{embed.eps}}}
\put(3.8,0){$r(i)$}\put(3.8, 1.8){$V_{op}$}\put(3.9, 0.9){$i$}
\put(5,0){,}

\put(6.5,0.9){$b^{(i)}_{V_{op}}=$}
\put(8,0){\resizebox{0.3cm}{2cm}
{\includegraphics{proj.eps}}}
\put(8.3,0){$V_{op}$}\put(8.3, 1.8){$r(i)$}\put(8.4, 0.9){$i$}
\put(9,0){,}

\put(10.5,0.9){$b_{(i)}^{V_{cl}}=$}
\put(12.3,0){\resizebox{0.3cm}{1.6cm}
{\includegraphics{embed.eps}}}
\put(12.7,0){\resizebox{0.3cm}{1.6cm}
{\includegraphics{embed.eps}}}
\put(12.3,1.7){$T(V_{cl})$}\put(12.1,0.6){$i$}\put(13.1,0.6){$i$}
\put(11.5,0){$r_L(i)$}\put(12.9,0){$r_R(i)$}\put(13.9,0){.}
\end{picture}
\eeq
We denote the map 
$\iota_{cl-op}: T(V_{cl}) \rightarrow V_{op}$ \cite{Ko2} 
by the following graph 
\beq  \label{iota-*-def}
\begin{picture}(14,2.2)
\put(5.5, 0.9){$\iota_{cl-op} :=$}
\put(7.5, 0.4){\resizebox{1cm}{1.5cm}
{\includegraphics{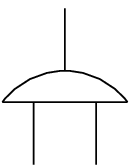}}}
\put(7.6,0){$T(V_{cl})$}\put(7.8,2){$V_{op}$}
\put(9.2,0){.}
\end{picture}
\eeq

Now we can express the Cardy condition (\ref{cardy-ioa}) in graphs. 
The left hand side of (\ref{cardy-ioa}) can be expressed by:
\beq  \label{cardy-L-graph-1}
\begin{picture}(14,3)
\put(1.5, 1.3){$\sum_{i=1}^{N_{cl}}$}

\put(3.8,0){\resizebox{3cm}{3cm}
{\includegraphics{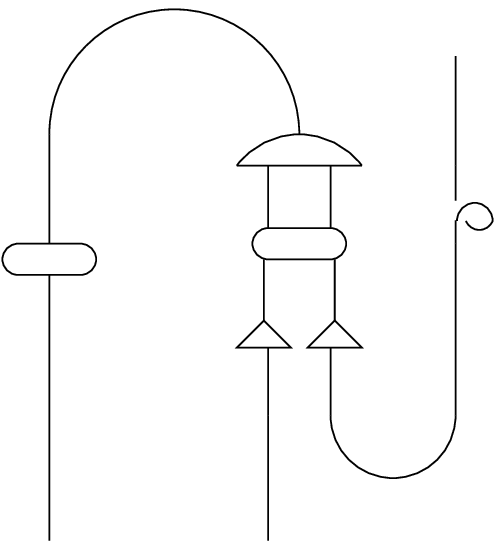}}}
\put(3.4,0){$V_{op}$} 
\put(3.1,1.6){$\varphi_{op}$}
\put(4.6, 1.6){$\varphi_{cl}^{-1}$}\put(5,1){$i$}
\put(6.1, 1){$i$}\put(6.7,2.5){$r_R(i)$}
\put(5.5,0){$(r_L(i))'$}

\put(7.3,1.5){$\otimes$}

\put(8.3,0){\resizebox{3cm}{3cm}
{\includegraphics{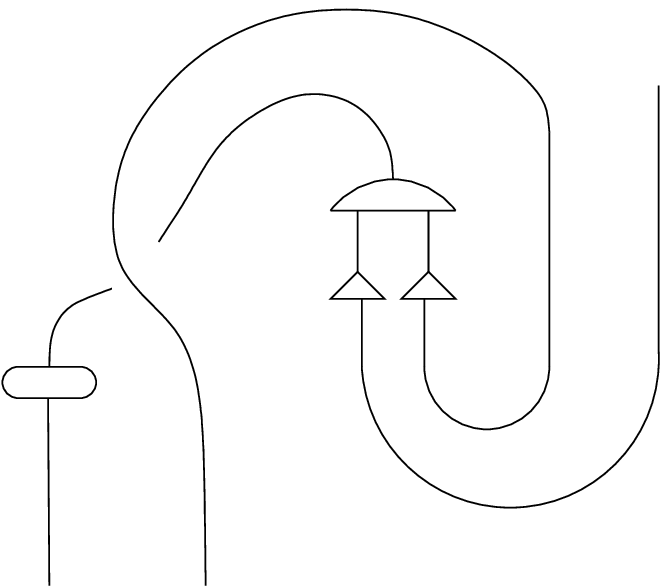}}}
\put(7.9,0){$V_{op}$} 
\put(7.6, 1){$\varphi_{op}$}
\put(9.6, 1.4){$i$}\put(10.4,1.4){$i$}
\put(9.3, 0){$r'_R(i)$}\put(11.4, 2.3){$(r_L(i))'$}

\put(13,0){.}

\end{picture}
\eeq

By the universal property of tensor product, for 
$a_i\in \I, i=1,\dots, 6$, we have a canonical isomorphism: 
\bea
\oplus_{a\in \I} \mathcal{V}_{a_1a}^{a_4}\otimes \mathcal{V}_{a_2a_3}^{a}
&\xrightarrow{\cong} &
\hom_V(W_{a_1}\boxtimes (W_{a_2} \boxtimes W_{a_3}), W_{a_4})  \nn
\Y_1 \otimes \Y_2 &\mapsto & 
m_{\Y_1} \circ (\id_{W_{a_1}} \boxtimes m_{\Y_2}).
\eea
Under this canonical isomorphism, the Cardy condition 
(\ref{cardy-ioa}) can be viewed as a condition on
two morphisms in 
$\hom_V(V_{op}\boxtimes (V_{op} \boxtimes W_{r_R(i)}), W_{r_R(i)})$.
In particular, the left hand side of (\ref{cardy-ioa}) viewed
as a morphism in 
$\hom_V(V_{op}\boxtimes (V_{op} \boxtimes W_{r_R(i)}), W_{r_R(i)})$
can be expressed as follow: 
\beq  \label{cardy-L-graph-2}
\begin{picture}(14,4)
\put(3, 1.8){$\sum_{i=1}^{N_{cl}}$}
\put(5,0){\resizebox{7cm}{4cm}
{\includegraphics{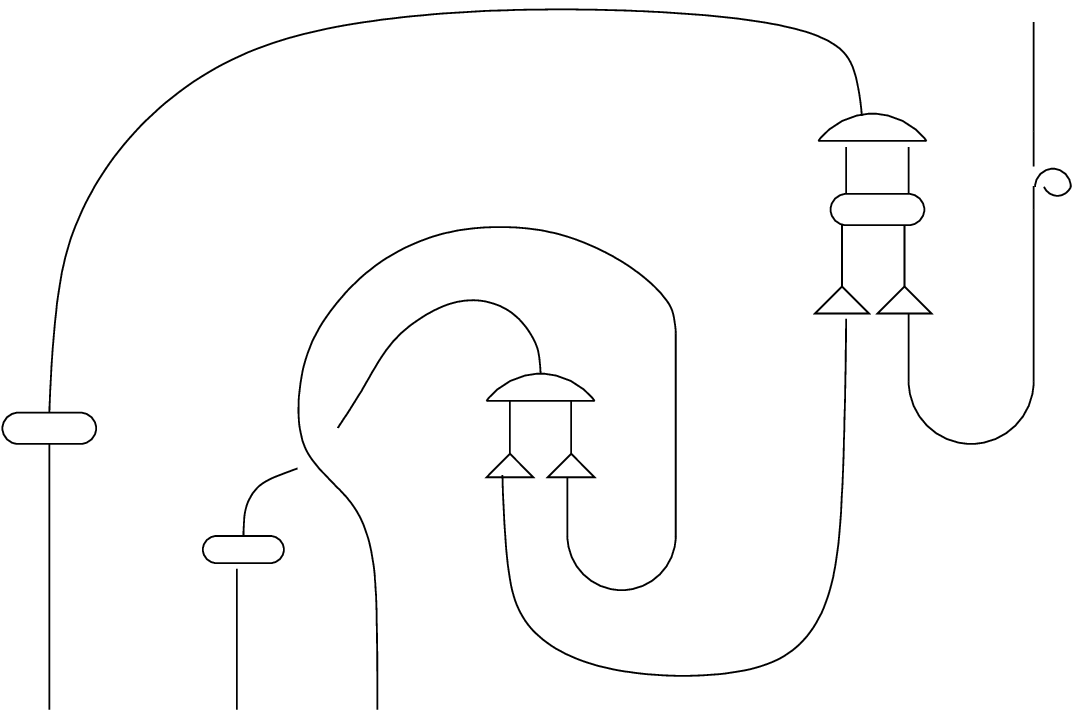}}}
\put(4.7,0){$V_{op}$} 
\put(4.3, 1.6){$\varphi_{op}$}\put(5.6, 0.9){$\varphi_{op}$}
\put(5.9,0){$V_{op}$} 
\put(7.6,0){$r_R(i)$}
\put(7.9,1.2){$i$}\put(9, 1.2){$i$}
\put(8.9, 2.8){$T(\varphi_{cl}^{-1})$}
\put(10.1, 2.2){$i$}\put(11.2,2.2){$i$}
\put(11.8,3.6){$r_R(i)$}\put(12.5,0){.}
\end{picture}
\eeq

We define a morphism $\iota_{cl-op}^*: V_{op} \rightarrow T(V_{cl})$ by
\beq  \label{iota-up-def}
\begin{picture}(14,2.2)
\put(3.1, 1){$\iota_{cl-op}^*$}\put(4.2, 1){$=$}
\put(5, 0.4){\resizebox{1cm}{1.5cm}
{\includegraphics{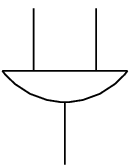}}}
\put(5.4,0){$V_{op}$}\put(5, 2){$T(V_{cl})$}

\put(6.5, 1){$:=$}

\put(8, 0){\resizebox{3.5cm}{2.2cm}
{\includegraphics{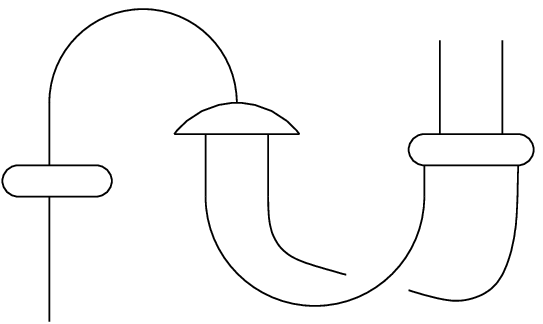}}}
\put(7.4,0.9){$\varphi_{op}$}
\put(11.6, 1.1){$T(\varphi_{cl}^{-1})$}
\put(12,0){.}

\end{picture}
\eeq
Then using the morphism $\iota_{cl-op}^*$, we can rewrite
the graph in (\ref{cardy-L-graph-2}) as follow: 
\beq  \label{cardy-L-graph-3-4}
\begin{picture}(14,4)
\put(1, 0){\resizebox{6cm}{4cm}
{\includegraphics{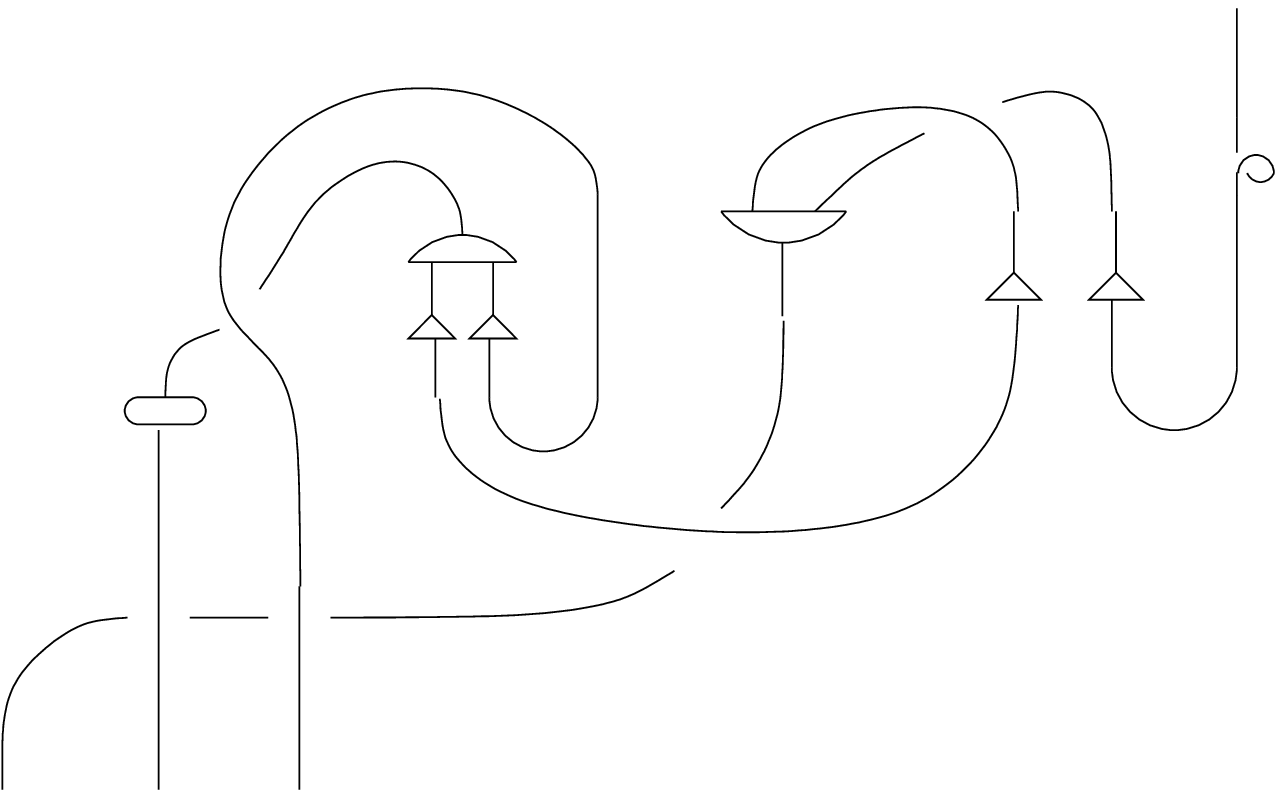}}}
\put(0.9, 1.9){$\varphi_{op}$}\put(0.3,0){$V_{op}$}
\put(1.1,0){$V_{op}$}
\put(2.7,2.2){$i$}\put(3.5,2.2){$i$}
\put(5.4, 2.5){$i$}\put(6.5,2.5){$i$}
\put(2.5,0){$r_R(i)$}\put(7,3.7){$r_R(i)$}

\put(7.9, 1.9){$=$}

\put(9.5, 0){\resizebox{4cm}{4cm}
{\includegraphics{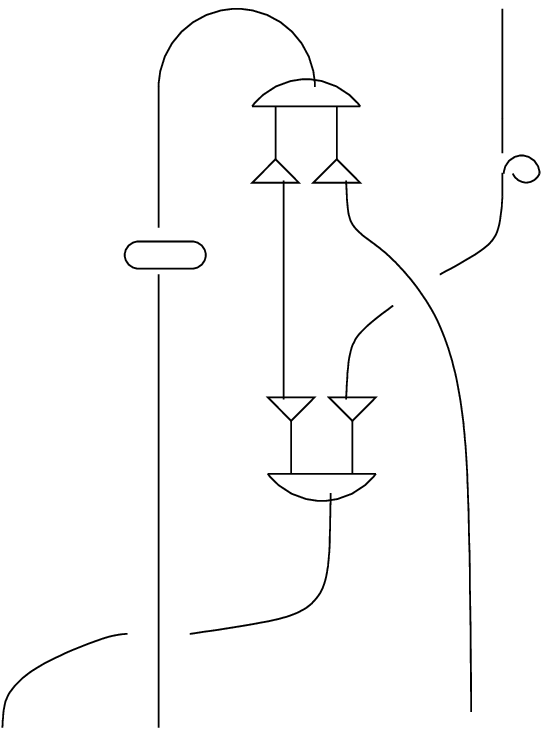}}}
\put(10.8,0){$V_{op}$}\put(8.9,0){$V_{op}$}
\put(13.1,0){$r_R(i)$}\put(9.7, 2.6){$\varphi_{op}$}
\put(13.4,3.7){$r_R(i)$}
\put(11.1,2.9){$i$}\put(12.3,2.9){$i$}
\put(11.2,1.7){$i$}\put(12.4,1.7){$i$}
\put(14,0){.}

\end{picture}
\eeq

Using (\ref{S-a-map-1}), (\ref{cardy-L-graph-2}) and 
(\ref{cardy-L-graph-3-4}), we obtain a graphic version of 
the Cardy condition (\ref{cardy-ioa}) as follow: 
\beq \label{cardy-cat}
\begin{picture}(14,4)

\put(2.2, 1.8){$\sum_{r_R(i)=a}$}
\put(4.2, 0){\resizebox{3.5cm}{4cm}
{\includegraphics{cardy-L-graph-4.eps}}}
\put(3.6,0){$V_{op}$}\put(5.3,0){$V_{op}$}
\put(4.3, 2.6){$\varphi_{op}$}
\put(5.6, 2.9){$i$}\put(6.6, 2.9){$i$}
\put(5.7, 1.7){$i$}\put(6.7, 1.7){$i$}
\put(7.6,3.8){$r_R(i)$}\put(7.3, 0){$r_R(i)$}

\put(8.4,1.9){$=$}

\put(9.1, 1.9){$\frac{\dim a}{D}$}

\put(10.6, 0){\resizebox{3.3cm}{4cm}
{\includegraphics{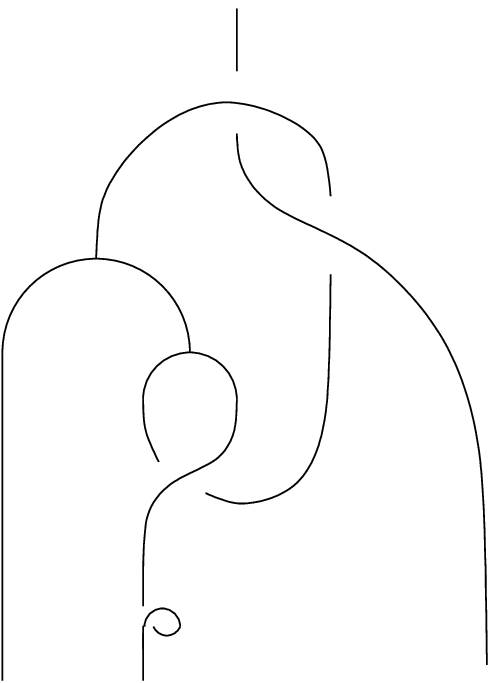}}}
\put(10,0){$V_{op}$}\put(10.9,0){$V_{op}$}\put(11, 1.3){$V_{op}$}
\put(10.7, 3){$V_{op}$}
\put(11.8, 2.3){$V_{op}$}\put(12.4, 3.8){$a$}
\put(14,0){$a$}

\end{picture}
\eeq
Using the Frobenius properties of $V_{op}$, one can show that 
(\ref{cardy-cat}) is equivalent to the following condition:
\beq   \label{cardy-cat-symm}
\begin{picture}(14,3)

\put(2.2, 1.4){$\sum_{r_R(i)=a'}$}
\put(4.5, 0){\resizebox{1.5cm}{3cm}
{\includegraphics{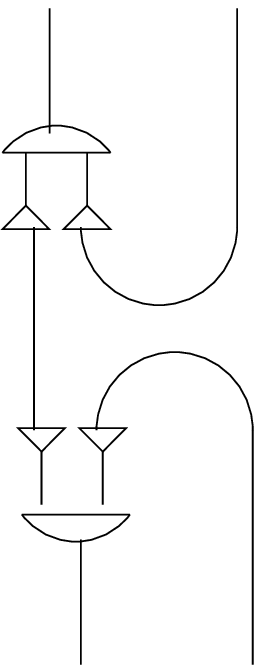}}}
\put(4.3,0){$V_{op}$}\put(4.1,2.8){$V_{op}$}
\put(4.3, 1.9){$i$}\put(5.2, 1.9){$i$}
\put(4.4, 1){$i$}\put(5.3,1){$i$}
\put(6,2.8){$a$}\put(6.1, 0){$a$}

\put(7.5,1.4){$=$}

\put(8.2, 1.4){$\frac{\dim a}{D}$}

\put(10, 0){\resizebox{3cm}{3cm}
{\includegraphics{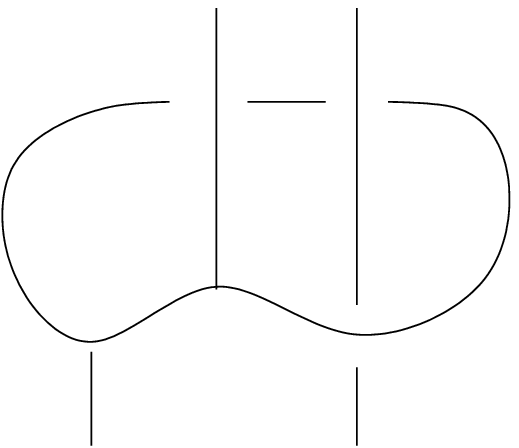}}}
\put(9.9,0){$V_{op}$}
\put(10.6, 2.8){$V_{op}$}
\put(12.2, 2.8){$a$}
\put(12.2,0){$a$}\put(13.5,0){.}

\end{picture}
\eeq
The asymmetry between chiral and antichiral parts in 
(\ref{cardy-cat})(\ref{cardy-cat-symm}) is superficial. 
Using (\ref{BK-lemma}), one can show that 
(\ref{cardy-cat-symm}) is further equivalent to the following condition: 
\beq  \label{cardy-cat-symm-1}
\begin{picture}(14,3)

\put(2.7, 1.4){$\frac{\dim a}{D}$}
\put(4.5, 0){\resizebox{2.3cm}{3cm}
{\includegraphics{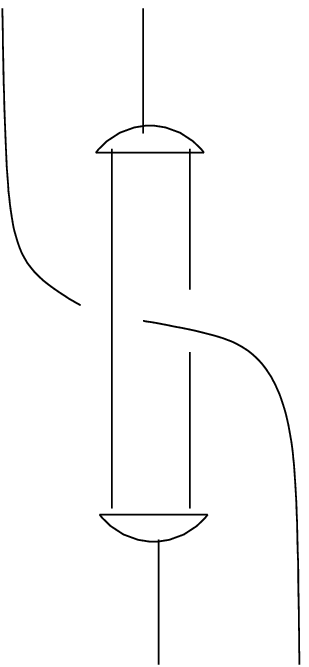}}}
\put(5.1,0){$V_{op}$}\put(5, 2.8){$V_{op}$}
\put(4.2,2.8){$a$}\put(6.9, 0){$a$}

\put(7.5,1.4){$=$}

\put(8.2, 1.4){$\sum_{r(k)=a}$}

\put(10.3, 0){\resizebox{2.5cm}{3cm}
{\includegraphics{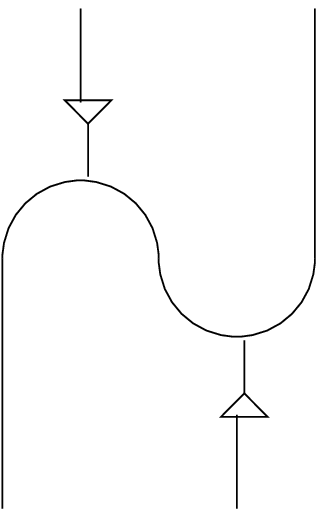}}}
\put(9.7,0){$V_{op}$}
\put(12.9, 2.8){$V_{op}$}
\put(10.6, 2.8){$a$}\put(10.5, 2.3){$k$}
\put(12.3,0){$a$}\put(12.5, 0.4){$k$}
\put(13.8,0){.}

\end{picture}
\eeq

We recall a definition in \cite{Ko2}. 
\begin{defn}  \label{def-opcl-C-alg}
{\rm 
An {\it open-closed $\mathcal{C}_V|\mathcal{C}_{V\otimes V}$-algebra},
denoted as $(A_{op}|A_{cl}, \iota_{cl-op})$, consists of
a commutative symmetric associative algebra $A_{cl}$
in $\mathcal{C}_{V\otimes V}$, an associative algebra $A_{op}$
in $\mathcal{C}_V$ and an associative algebra morphism
$\iota_{cl-op}: T(V_{cl}) \rightarrow V_{op}$,
satisfying the following condition: 
\beq \label{iota-m-comm-fig}
\begin{picture}(14,2.5)
\put(4,0.5){\resizebox{2cm}{2cm}
{\includegraphics{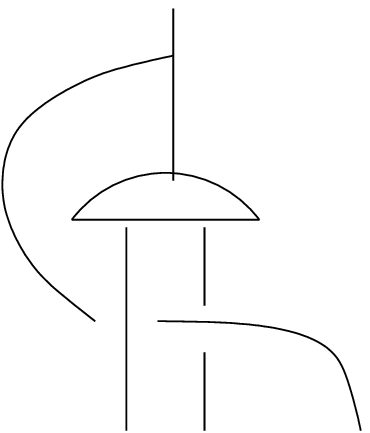}}}
\put(6, 0){$V_{op}$}\put(5.1, 2.4){$V_{op}$}
\put(4.5,0){$T(V_{cl})$}

\put(7,0.9){$=$}

\put(8.3,0.5){\resizebox{2cm}{2cm}
{\includegraphics{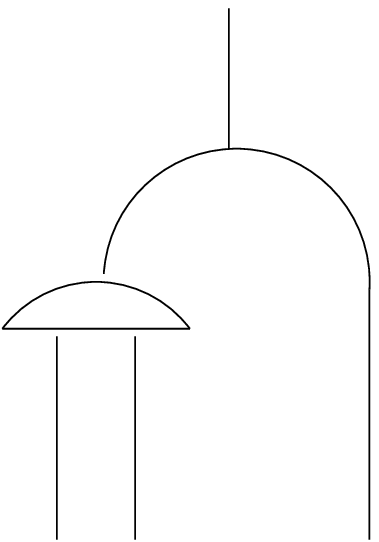}}}
\put(8.2,0){$T(V_{cl})$}
\put(10.1,0){$V_{op}$}\put(9.6,2.4){$V_{op}$}\put(11,0){.}

\end{picture}
\eeq
}
\end{defn}

The following Theorem is proved in \cite{Ko2}. 
\begin{thm}  \label{thm-opcl-C-alg}
The category of open-closed field algebras over $V$ is
isomorphic to the category of open-closed 
$\mathcal{C}_V|\mathcal{C}_{V\otimes V}$-algebras. 
\end{thm}

\begin{defn}  {\rm 
A {\it Cardy $\mathcal{C}_V|\mathcal{C}_{V\otimes V}$-algebra} is
an open-closed $\mathcal{C}_V|\mathcal{C}_{V\otimes V}$-algebra
$(A_{op}|A_{cl}, \iota_{cl-op})$ such that $A_{cl}$ is modular invariant 
commutative symmetric Frobenius algebra in $\mathcal{C}_{V\otimes V}$ 
and $A_{op}$ a symmetric Frobenius algebra
in $\mathcal{C}_V$ and the Cardy condition (\ref{cardy-cat-symm}) 
or (\ref{cardy-cat-symm-1}) hold. }
\end{defn}

\begin{rema}  {\rm
Notice that, in the case $V=\C$, 
Cardy $\mathcal{C}_V|\mathcal{C}_{V\otimes V}$-algebra 
(using (\ref{cardy-cat-symm-1}))
exactly coincides with the usual algebraic formulation of 
2-dimensional open-closed topological field theory 
\cite{La}\cite{Mo1}\cite{Mo2}\cite{Se2}\cite{MSeg}\cite{AN}\cite{LP}. 
As we discussed in the introduction, we believe that 
open-closed partial conformal field theories of all genus satisfying
the $V$-invariant boundary condition \cite{Ko2}
are classified by Cardy $\mathcal{C}_V|\mathcal{C}_{V\otimes V}$-algebras.
}
\end{rema}

The following result is clear. 
\begin{thm}
The category of open-closed field algebras over $V$ equipped with nondegenerate invariant bilinear forms and satisfying the modular invariance condition and the Cardy condition is isomorphic to the category of Cardy $\mathcal{C}_V|\mathcal{C}_{V\otimes V}$-algebras.
\end{thm}

\begin{defn} \label{def:dbrane}  {\rm 
A {\it $V$-invariant D-brane}
associated to a closed algebra $A_{cl}$ in 
$\mathcal{C}_{V\otimes V}$ is a pair $(A_{op}, \iota_{cl-op})$ 
such that the triple $(A_{op}|A_{cl}, \iota_{cl-op})$ gives a 
Cardy $\mathcal{C}_V|\mathcal{C}_{V\otimes V}$-algebra. }
\end{defn}

D-branes usually form a category as we will see in an example 
in the next subsection. 

\subsection{Constructions}

In this section, we give a categorical construction of 
Cardy $\mathcal{C}_V|\mathcal{C}_{V\otimes V}$-algebra. 
This construction is 
called Cardy case in physics literature \cite{FFFS2}.

Let us first recall the diagonal construction of the close algebra $V_{cl}$
\cite{FFRS}\cite{HKo2}\cite{Ko1} \cite{HKo3}. 
We will follow the categorical construction given in \cite{Ko1}.

Let $V_{cl}$ is be the object in $\mathcal{C}_{V\otimes V}$ 
given as follow
\beq  \label{cvoa-const-1}
V_{cl} = \oplus_{a\in \mathcal{I}} \,\, W_a\otimes W'_a.
\eeq
The decomposition of $V_{cl}$ as a direct sum gives 
a natural embedding $V\otimes V \hookrightarrow V_{cl}$. 
We denote this embedding as $\iota_{cl}$. 
We define a morphism  
$\mu_{cl} \in \hom_{V\otimes V}(V_{cl} \boxtimes V_{cl}, V_{cl})$ by
\beq  \label{cvoa-const-2}
\mu_{cl} = \sum_{a_1,a_2,a_3\in \mathcal{A}} \sum_{i,j=1}^{N_{a_1a_2}^{a_3}} 
\langle f^{a_1a_2}_{a_3;i}, f^{a'_1a'_2}_{a'_3;j}\rangle \, 
e_{a_1a_2;i}^{a_3} \otimes e_{a'_1a'_2; j}^{a'_3},
\eeq
where $e_{a_1a_2;i}^{a_3}$ and 
$f^{a_1a_2;j}_{a_3}$ are basis vectors given in (\ref{basis-fig}) 
and (\ref{dual-basis-pic}) and 
$\langle \cdot, \cdot\rangle$ is a 
bilinear pairing given by
\beq  \label{pairing-equ}
\begin{picture}(14, 2.5)
\put(6,0){\resizebox{3.5cm}{2.5cm}{
\includegraphics{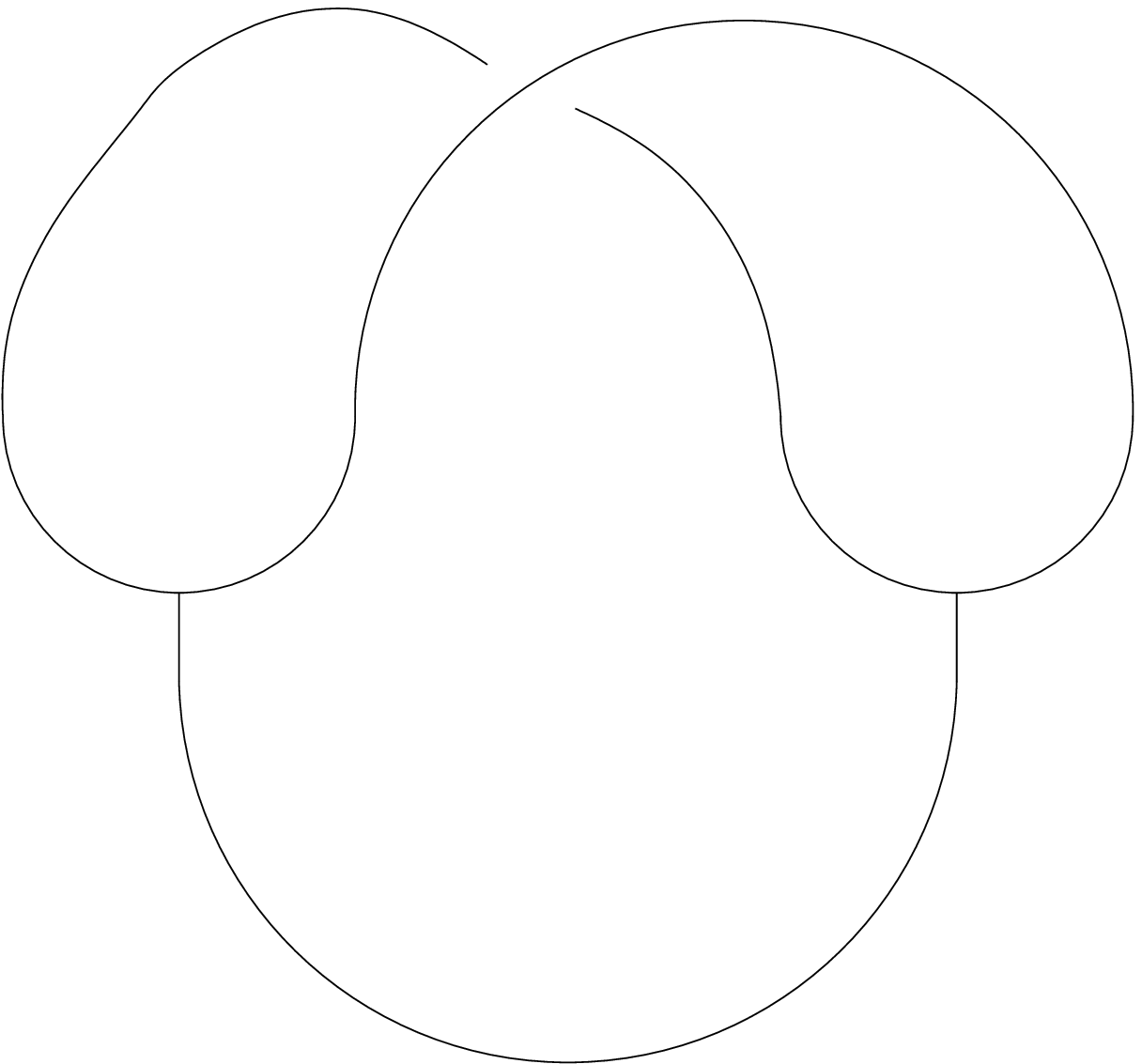}}}
\put(4, 1.2){$\frac{1}{\dim a_3}$}
\put(5.6, 1.6){$a_1$}
\put(6.5, 1.3){$i$}\put(6.1, 0.7){$a_3$}
\put(7.3, 1.6){$a_2$}
\put(7.9, 1.6){$a'_1$}\put(8.9, 1.3){$j$}
\put(9.6, 1.6){$a'_2$}\put(9, 0.7){$a'_3$}
\put(10.5,0){.}
\end{picture}
\eeq

Notice that $V_{cl}'$ has the same decomposition 
as $V_{cl}$ in (\ref{cvoa-const-1}). They
are isomorphic as $V\otimes V$-modules.  
There is, however, no canonical isomorphism. 
Now we choose a particular isomorphism  
$\varphi_{cl}: V_{cl}\rightarrow V_{cl}'$ given by 
\beq  \label{varphi-F}
\varphi_{cl}=\oplus_{a\in \I} \,\,
\frac{D}{\dim a} \, e^{-2\pi ih_{a}} \, \id_{W_{a}\otimes W'_a}.
\eeq
The isomorphism $\varphi_{cl}$ induces a 
nondegenerate invariant bilinear form on $V_{cl}$ viewed as
$V\otimes V$-module.

The following Theorem is a categorical version of Theorem 5.1
in \cite{HKo3}. We give a categorical proof here. 
\begin{thm}
$(V_{cl}, \mu_{cl}, \iota_{cl})$ together with
isomorphism $\varphi_{cl}$ gives a modular invariant 
commutative Frobenius algebra in $\mathcal{C}_{V\otimes V}$ 
with a trivial twist.
\end{thm}
\pf
It was proved in \cite{Ko1} that 
$(V_{cl}, \mu_{cl}, \iota_{cl})$ together with
isomorphism $\varphi_{cl}$ gives a commutative Frobenius algebra
with a trivial twist. It remains to show the modular invariance.

First, the bilinear pairing $\langle \cdot, \cdot\rangle$ given 
in  (\ref{pairing-equ}) can be naturally extended to  
a bilinear form, still denoted as
$\langle \cdot, \cdot\rangle$, 
on $\oplus_{a,a_1\in \I}(W_{a_1}, W_a\boxtimes W_{a_1})$ as follow: 
\beq
\langle f^{a a_1}_{a_1;i}, f^{bb_1}_{b_1;j}\rangle :=
\delta_{ab'}\,\, \delta_{a_1b'_1} \,
\langle f^{a a_1}_{a_1;i}, f^{a'a'_1}_{a'_1;j}\rangle.
\eeq

Then it is easy to see that to prove the modular invariance of $V_{cl}$  
is equivalent to prove that the bilinear form $\langle \cdot, \cdot\rangle$ 
on $\oplus_{a,a_1\in \I}(W_{a_1}, W_a\boxtimes W_{a_1})$ defined above
is invariant under the action of $(S^{-1})^* \otimes S^*$. 
Clearly, when $b\neq a'$, 
$\langle (S^{-1}(a))^* f_{a_1;i}^{aa_1},
(S(b))^* f_{b_1;j}^{bb_1}\rangle=0$.
When $b=a'$, we have (using (\ref{S-dual-graph-1}), 
(\ref{S-dual-graph-2}) and (\ref{BK-lemma}))
$$
\begin{picture}(14,4)
\put(0,2){$\ds \langle (S^{-1}(a))^* f_{a_1;i}^{aa_1}, 
             (S(a'))^* f_{b_1;j}^{a'b_1}\rangle$}
\put(5.5, 2){$=$}

\put(6, 2){$\ds \sum_{a_3\in \I} \frac{\dim a_3}{D^2}$} 

\put(8.5, 0){\resizebox{5cm}{4cm}{\includegraphics{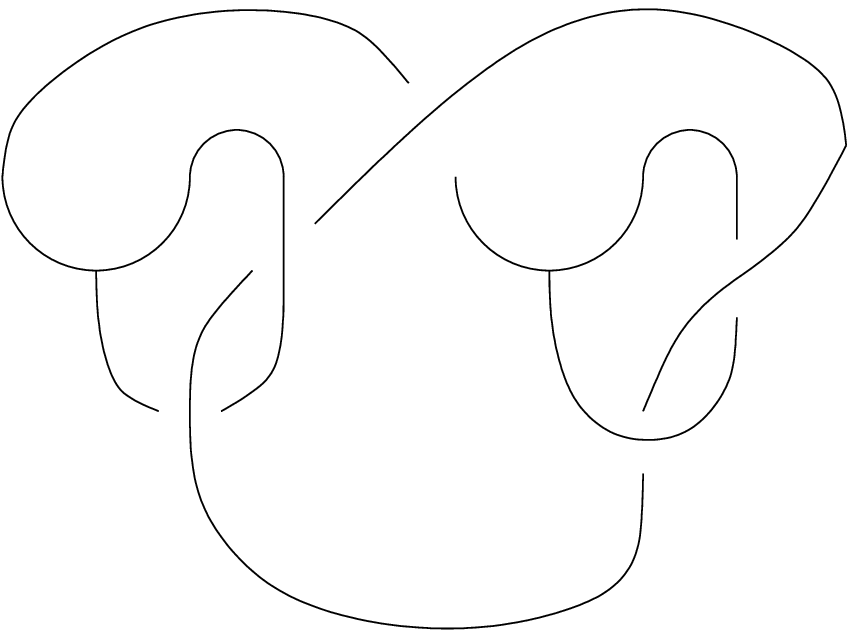}}}
\put(8.2, 3.2){$a$}\put(8.6, 1.8){$a_1$}\put(11.3, 1.8){$b_1$}
\put(9, 2.5){$i$}\put(11.7,2.6){$j$}\put(9.2, 0.5){$a_3$}
\end{picture}
$$
$$
\begin{picture}(14,4)
\put(0,2){$=$}
\put(0.5, 2){$\ds \sum_{a_3\in \I} \frac{\dim a_3}{D^2}$}

\put(2.8, 0){\resizebox{4cm}{4cm}{\includegraphics{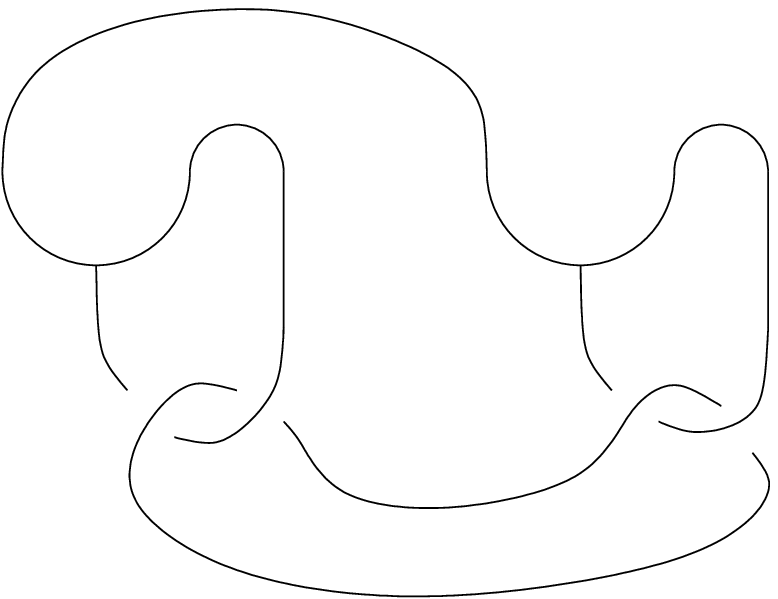}}}
\put(2.5, 3.2){$a$}\put(2.8, 1.7){$a_1$}\put(5.4,1.7){$b_1$}
\put(3.2, 2.4){$i$}\put(5.7,2.5){$j$}\put(2.9, 0.5){$a_3$}

\put(7.1,2){$=$}
\put(7.6, 2){$\ds \sum_{a_3\in \I} \frac{\dim a_3}{D^2}$}

\put(9.8, 0){\resizebox{4cm}{4cm}{\includegraphics{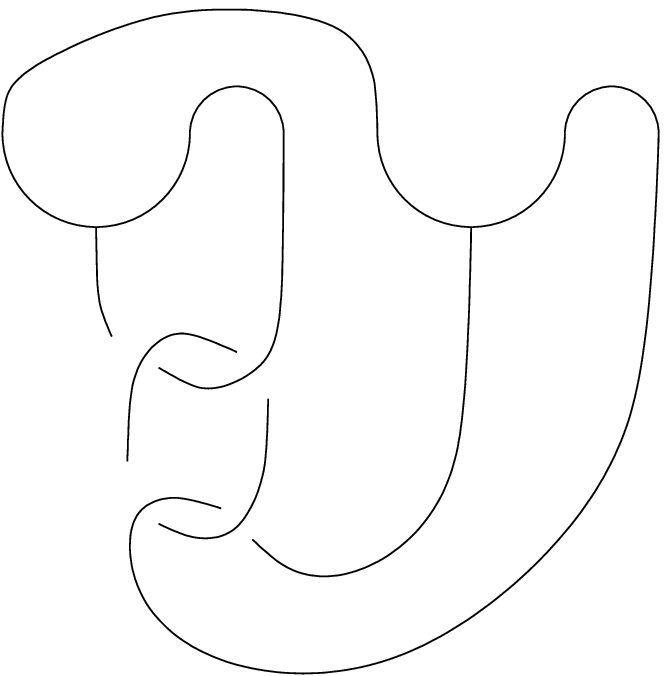}}}
\put(9.5, 3.2){$a$}\put(9.9, 2.3){$a_1$}\put(12.8, 2){$b_1$}
\put(10.3, 2.9){$i$}\put(12.6,2.9){$j$}\put(10, 1.4){$a_3$}
\end{picture}
$$
$$
\begin{picture}(14,4)
\put(0,2){$=$}
\put(0.5, 2){$\ds \sum_{a_3\in \I} \frac{\dim a_3}{D^2}$}

\put(2.8, 0){\resizebox{4cm}{4cm}{\includegraphics{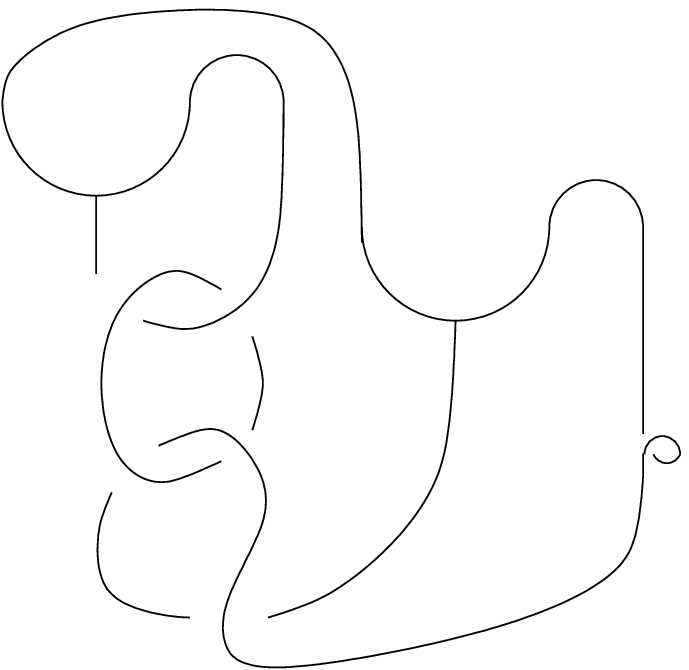}}}
\put(2.5, 3.2){$a$}\put(2.9, 2.5){$a_1$}\put(5.6,1.6){$b_1$}
\put(2.9, 0.6){$b'_1$}
\put(3.3, 3){$i$}\put(5.4,2.4){$j$}\put(2.9, 1.8){$a_3$}
\put(7.1,2){$=$}

\put(7.8, 2){$\ds \frac{\delta_{a_1b'_1}}{\dim a_1}$}
\put(9.5, 0){\resizebox{4cm}{4cm}{\includegraphics{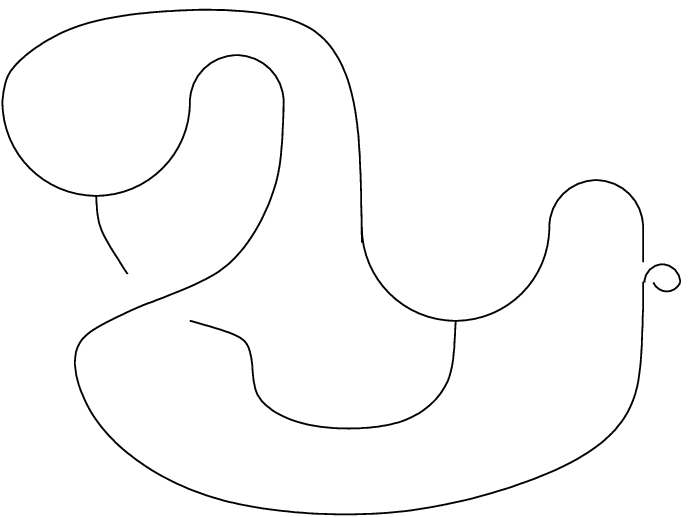}}}
\put(9.2, 3.2){$a$}\put(9.7, 2){$a_1$}\put(10, 2.7){$i$}
\put(12.1,1.8){$j$}\put(9.5, 1){$a'_1$}
\end{picture}
$$
$$
\begin{picture}(14,4)
\put(0,1.7){$=$}
\put(0.5, 1.7){$\ds \frac{\delta_{a_1b'_1}}{\dim a_1}$}

\put(2.8, 0){\resizebox{3.5cm}{3.5cm}{\includegraphics{pair-0--1.eps}}}
\put(2.5, 2.2){$a$}\put(4.1, 2.2){$a_1$}\put(3.3, 1.8){$i$}
\put(5.7,1.8){$j$}\put(2.8, 0.7){$a_1$}
\put(6.5,0){.}
\end{picture}
$$
\epf

Now we define $V_{op}$. Let $X$ be a $V$-module. 
Let $e_X: X' \boxtimes X \rightarrow V$ and 
$i_X: V \rightarrow X\boxtimes X'$ be the duality maps defined in \cite{Ko1}. 
$V_{op} := X \boxtimes X'$ has a natural structure of symmetric Frobenius algebra \cite{FS} with $\iota_{op}:= i_X$,  $\mu_{op}:= \id_X \boxtimes e_X \boxtimes \id_{X'}$, $\epsilon_{op}:= e_{X'}$ and $\Delta_{op}:= \id_X \boxtimes i_X \boxtimes \id_{X'}$. 

Now we define a map $\iota_{cl-op}: T(V_{cl}) \rightarrow V_{op}$ by 
\beq  \label{iota-cardy}
\begin{picture}(14,2)

\put(4, 0){\resizebox{1.2cm}{2cm}{
\includegraphics{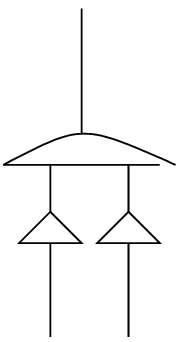}}}
\put(3.1,1.7){$X\boxtimes X'$}
\put(4,0){$a$}\put(5,0){$a'$}

\put(5.7, 0.9){$:=$}

\put(7,0){\resizebox{2cm}{2cm}{
\includegraphics{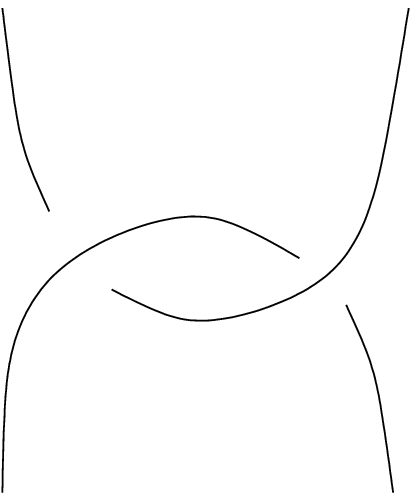}}}
\put(6.7,0){$a$}\put(9.1,0){$a'$}
\put(6.6,1.7){$X$}\put(9.1,1.7){$X'$}

\put(9.5,0){.}
\end{picture}
\eeq

\begin{lemma}
\beq  \label{iota-up-cardy-case}
\begin{picture}(14,2)

\put(4, 0){\resizebox{1.2cm}{2cm}{
\includegraphics{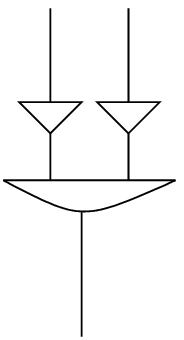}}}
\put(3.1,0){$X\boxtimes X'$}
\put(4,1.8){$a$}\put(5,1.8){$a'$}

\put(5.7, 0.9){$=$}

\put(6.5,0.9){$\, \frac{\dim a}{D} $}

\put(8.5,0){\resizebox{2cm}{2cm}{
\includegraphics{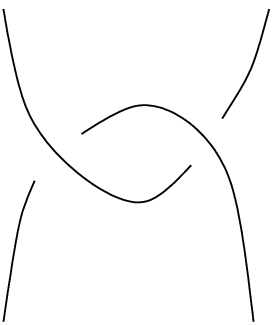}}}
\put(8.3,1.8){$a$}\put(10.6,1.8){$a'$}
\put(8.2,0){$X$}\put(10.6, 0){$X'$}

\put(11.5,0){.}
\end{picture}
\eeq
\end{lemma}
\pf
By (\ref{iota-up-def}) and (\ref{iota-cardy}), we have
\bea  \label{iota-up-cardy-proof}
\begin{picture}(14,2)

\put(0.5, 0){\resizebox{1cm}{2cm}{
\includegraphics{iota-up-def-2.eps}}}
\put(1.1,0){$X\boxtimes X'$}
\put(0.5,1.8){$a$}\put(1.3,1.8){$a'$}

\put(2.3,0.9){$=$}

\put(2.8,0.9){$ \frac{\dim a}{D}$}

\put(4.5, 0){\resizebox{3.5cm}{2cm}{
\includegraphics{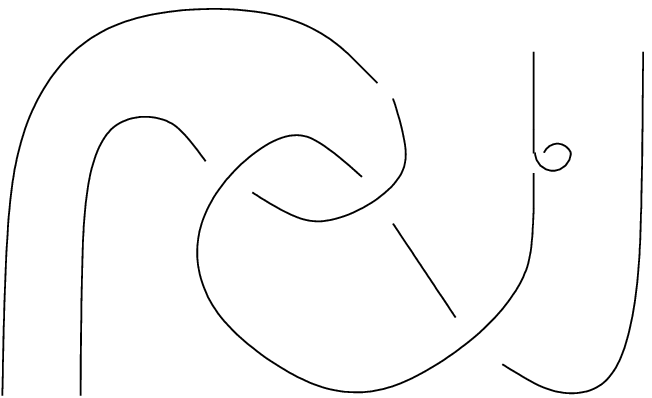}}}
\put(4.1,0){$X$}\put(5.1,0){$X'$}
\put(7.1,1.7){$a$}\put(8.1,1.7){$a'$}

\put(8.5,0.9){$= \,\, \frac{\dim a}{D}$}

\put(11, 0){\resizebox{3cm}{2cm}{
\includegraphics{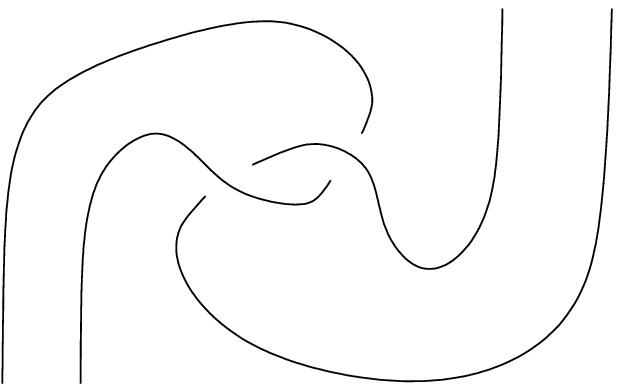}}}
\put(10.6,0){$X$}\put(11.5,0){$X'$}
\put(13.1,1.8){$a$}\put(14.1, 1.8){$a'$}

\end{picture}
\eea
It is easy to see that the last figure in 
(\ref{iota-up-cardy-proof}) can be
deformed to that on the right hand side of 
(\ref{iota-up-cardy-case}). 
\epf

\begin{thm} \label{thm-const}
$(V_{cl}, V_{op}, \iota_{cl-op})$ is a Cardy 
$\mathcal{C}_V|\mathcal{C}_{V\otimes V}$-algebra.
\end{thm}
\pf
Recall that $T(V_{cl})$ together with multiplication morphism 
$\mu_{T(V_{cl})} = T(\mu_{cl}) \circ \varphi_2$ and 
morphism $\iota_{T(V_{cl})} = T(\iota_{cl}) \circ \varphi_0$ is
an associative algebra. 
We first prove that $\iota_{cl-op}$ is an algebra morphism. 
It is clear that $\iota_{cl-op}\circ \iota_{T(V_{cl})}=\iota_{op}$. 
It remains to show the following identity
\beq  \label{iota-is-morphism}
\iota_{cl-op} \circ \mu_{T(V_{cl})} = 
\mu_{op} \circ (\iota_{cl-op} \boxtimes \iota_{cl-op}).
\eeq
By the definition of $\varphi_2$, that of $\mu_{cl}$ 
and (\ref{iota-cardy}), we obtain 
\beq  \label{iota-alg-morph-equ-1}
\begin{picture}(14,3)

\put(0.2, 1.4){$\iota_{cl-op} \circ \mu_{T(V_{cl})}$}

\put(2.7, 1.4){$=$}

\put(3.2, 1.4){$\sum_{a,b,c\in \I}\,\, \sum_{i, j} \,\,
\frac{1}{\dim c}$}

\put(7.7, 0.5){\resizebox{2cm}{2cm}{
\includegraphics{pair-0--1.eps}}}
\put(7.7, 0.9){$c$}\put(9.5,0.9){$c'$}
\put(7.4, 1.6){$a$}\put(8.5,1.6){$b$}
\put(8, 1.5){$i$}\put(9.3, 1.6){$j$}

\put(10.5, 0){\resizebox{3cm}{3cm}{
\includegraphics{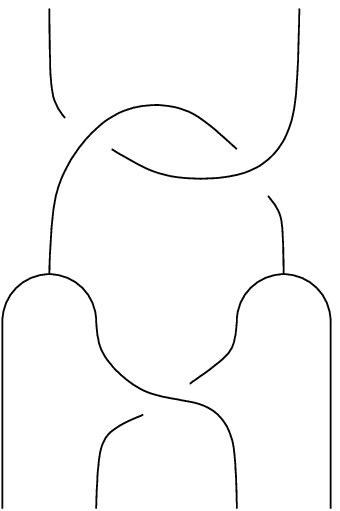}}}
\put(10.3, 0){$a$}\put(11.6,0){$a'$}
\put(12.4,0){$b$}\put(13.6,0){$b'$}
\put(10.9,1){$i$}\put(13,1){$j$}
\put(10.8, 1.7){$c$}\put(13.2, 1.7){$c'$}
\put(13.3, 2.8){$X'$}
\put(10.6,2.8){$X$}

\end{picture}
\eeq
It is easy to see that the right hand side of 
(\ref{iota-alg-morph-equ-1}) equals to 
\beq  \label{cardy-alg-equ-2}
\begin{picture}(14,4)

\put(2, 1.9){$\sum_{a,b,c\in \I}\,\, \sum_{i, j}\,\, \sum_{c_1\in \I} $}

\put(7.5, 0){\resizebox{4cm}{4cm}{
\includegraphics{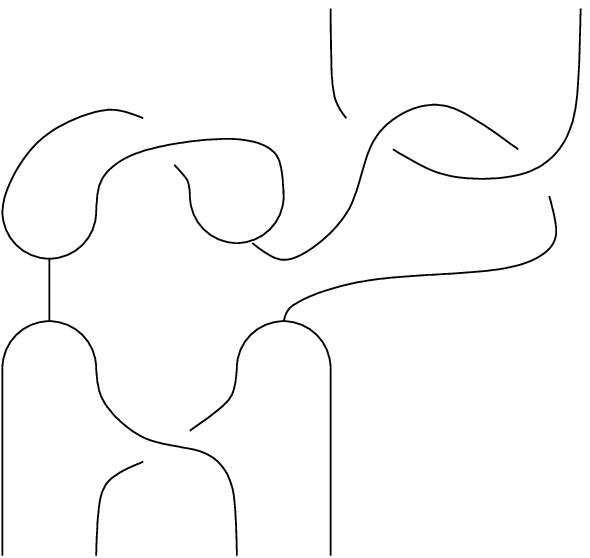}}}
\put(7.3,0){$a$}\put(8.3,0){$a'$}
\put(8.9,0){$b$}\put(9.9,0){$b'$}
\put(7.8,1.3){$i$}\put(9.4,1.3){$j$}
\put(7.5, 1.9){$c_1$}\put(11.4, 2.3){$c$}
\put(7.8, 2.3){$i$}\put(7.3,2.5){$a$}
\put(8.4, 2.5){$b$}\put(9.1,2.5){$j$}
\put(9.4,3.8){$X$}\put(11.6,3.8){$X'$}

\end{picture}
\eeq
Using (\ref{exp-id}) to sum up the indices $c_1$ and $i$,
we obtain that (\ref{cardy-alg-equ-2}) further equals to
\beq
\begin{picture}(14,3.5)

\put(3.5, 1.7){$\sum_{a,b,c\in \I}\,\, \sum_{j}$}

\put(7.5, 0){\resizebox{4cm}{3.5cm}{
\includegraphics{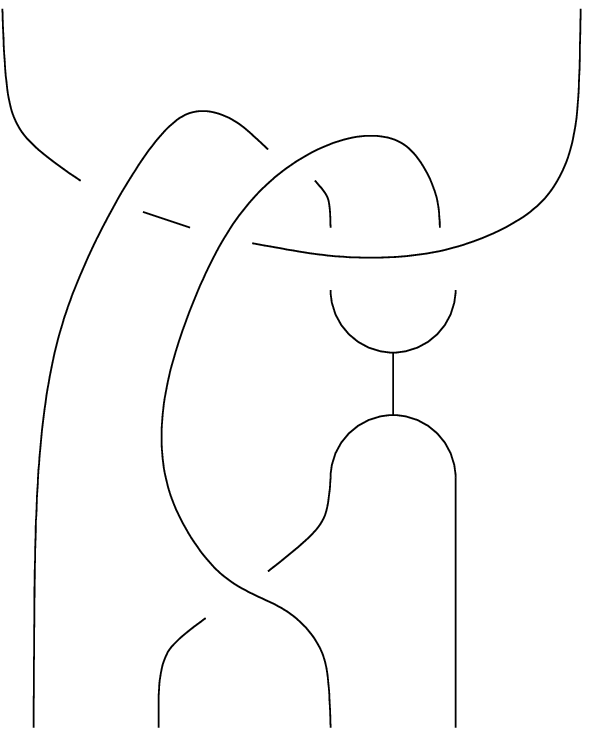}}}
\put(7.4,0){$a$}\put(8.3,0){$a'$}
\put(9.9,0){$b$}\put(10.7,0){$b'$}
\put(10.1, 1.1){$j$}\put(10.3, 1.6){$c$}
\put(10.1, 2){$j$}
\put(7.2,3.3){$X$}\put(11.6,3.3){$X'$}

\end{picture}
\eeq
Using (\ref{exp-id}) again, we obtain
\beq
\begin{picture}(14,3)

\put(1.5, 1.7){$\sum_{a,b\in \I} $}

\put(3.5, 0){\resizebox{3cm}{3cm}{
\includegraphics{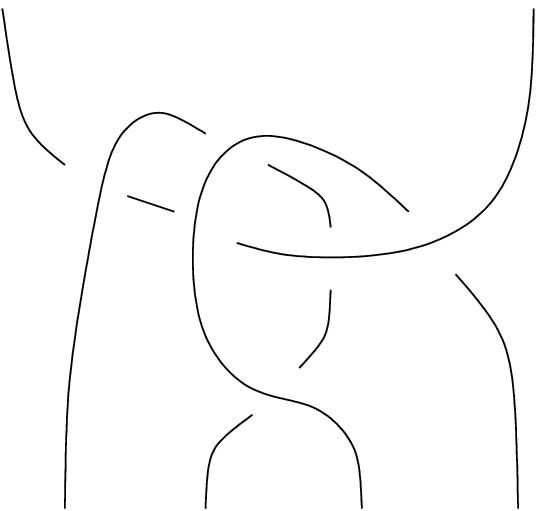}}}
\put(3.6,0){$a$}\put(4.3,0){$a'$}\put(5.7,0){$b$}
\put(6.5,0){$b'$}
\put(3.1,2.8){$X$}\put(6.6, 2.8){$X'$}

\put(7.3, 1.4){$=$} 

\put(8.2, 1.4){$\sum_{a,b\in \I}$}

\put(10.7, 0.25){\resizebox{3cm}{2.5cm}{
\includegraphics{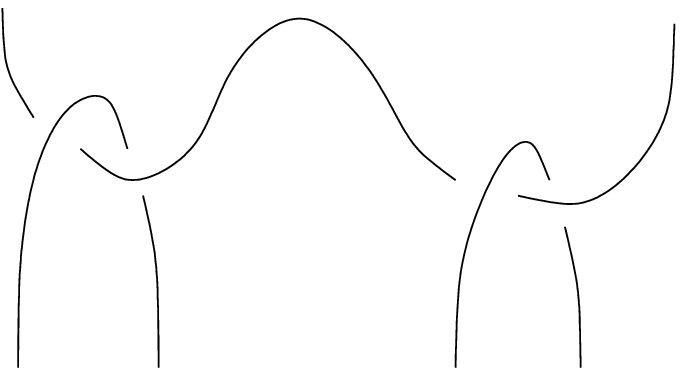}}}
\put(10.6,0.25){$a$}\put(11.5,0.25){$a'$}
\put(12.5,0.25){$b$}\put(13.4,0.25){$b'$}
\put(10.4, 2.55){$X$}\put(13.8,2.55){$X'$}

\put(14,0){,}

\end{picture}
\eeq
the right hand side of which is nothing but 
$\mu_{op} \circ (\iota_{cl-op} \boxtimes \iota_{cl-op})$.

The commutativity (\ref{iota-m-comm-fig}) follows from the following
identity:
\beq
\begin{picture}(14,3)

\put(3.5, 0.25){\resizebox{3cm}{2.5cm}{
\includegraphics{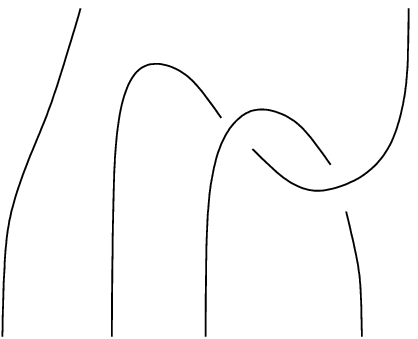}}}
\put(3.1,0.25){$X$}\put(3.8,0.25){$X'$} \put(5.1,0.25){$a$}
\put(6.3,0.25){$a'$}\put(6.6,2.55){$b'$}

\put(7.5, 1.4){$=$}

\put(8.5, 0){\resizebox{3.5cm}{3cm}{
\includegraphics{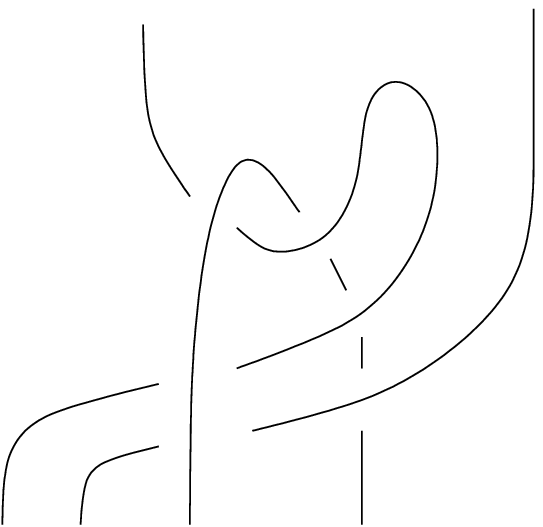}}}
\put(8.1,0){$X$}\put(8.6,0){$X'$}\put(9.9,0){$a$}
\put(11,0){$a'$}\put(12.1,2.8){$X'$}

\put(13,0){.}

\end{picture}
\eeq

In summary, we have proved that the triple 
$(V_{cl}, V_{op}, \iota_{cl-op})$ is an open-closed 
$\mathcal{C}_V|\mathcal{C}_{V\otimes V}$-algebra. 
It remains to show that Cardy condition (\ref{cardy-cat-symm}) holds. 
We use (\ref{iota-cardy}), (\ref{iota-up-cardy-case}) and the definition of $\mu_{op}$ and $\Delta_{op}$ to express both sides of (\ref{cardy-cat-symm}) graphically. Then it is easy to see that both sides are the deformations of each other. 
\epf

Theorem \ref{thm-const} reflects the general fact that consistent 
open theories (or {\rm D}-branes) for a given closed theory are not unique.
Instead they form a category. 
There are many good questions one can ask about 
Cardy $\mathcal{C}_V|\mathcal{C}_{V\otimes V}$-algebra, for example 
its relation to the works of Fuchs, Runkel, Schweigert and Fjelstad
\cite{FS}\cite{FRS1}-\cite{FRS4}\cite{FjFRS1}\cite{FjFRS2}. 
We leave such topics to \cite{cardysew} and future publications.

\appendix

\section{The Proof of Lemma \ref{b-F-map}}

\begin{lemma}
{\rm 
For $w_1\in W_{a_1}, w_a\in W_{a_2}$ and $w'_{a_3}\in (W_{a_3})'$
and $\Y_{a_1a_2}^{a_3}\in \V_{a_1a_2}^{a_3}$,  
we have  
\bea  \label{A-r-U}
&&\langle w'_{a_3}, \Y(\mathcal{U}(x)w_{a_1}, x)w_{a_2}\rangle  \nn
&&\hspace{1cm}=
\langle \tilde{A}_r(\Y)(\mathcal{U}(e^{(2r+1)\pi i}x^{-1})
e^{-(2r+1)\pi iL(0)}w_{a_1}, e^{(2r+1)\pi i}x^{-1})w'_{a_3}, w_{a_2}\rangle. 
\eea
}
\end{lemma}
\pf
Using the definition of $\hat{A}_r$, we see that the left
hand side of (\ref{A-r-U}) equals to 
\beq \label{A-r-U-equ-1}
\langle \tilde{A}_r(\Y)(e^{-xL(1)}x^{-2L(0)}
\mathcal{U}(x)w_{a_1},
e^{(2r+1)\pi i}x^{-1})w'_{a_3}, w_{a_2}\rangle.
\eeq
In \cite{H9}, the following formula
\beq  \label{A-r-U-equ-2}
e^{xL(1)}x^{-2L(0)} e^{(2r+1)\pi iL(0)} \mathcal{U}(x) e^{-(2r+1)\pi iL(0)}
=\mathcal{U}(x^{-1})
\eeq
is proved. Applying (\ref{A-r-U-equ-2}) to (\ref{A-r-U-equ-1}), 
we obtain (\ref{A-r-U}) immediately. 
\epf

Now we are ready to give a proof of Lemma \ref{b-F-lemma}. 

\pf
We have, for $w_{a_2}\in W_{a_2}, w_{a_3}\in W_{a_3}$, 
\bea  \label{b-equ-1}
&& ((\Psi_2(\Y_{aa_1; i}^{a_1;(1)}\otimes \Y_{a_2a_3;j}^{a;(2)})) 
(w_{a_2}\otimes w_{a_3})) (z_1,z_2+\tau,\tau)
 \nn
&&=E\bigg(\tr_{W_{a_{1}}}
\Y_{aa_{1}; i}^{a_{1};(1)}(\mathcal{U}(e^{2\pi i(z_{2}+\tau)})\cdot\nn
&&\quad\quad\quad \quad\quad\quad\quad\cdot
\Y_{a_{2}a_{3}; j}^{a;(2)}
(w_{a_{2}}, z_{1}-(z_{2}+\tau))w_{a_{3}}, e^{2\pi i (z_{2}+\tau)})
q_{\tau}^{L(0)-\frac{c}{24}}\bigg)\nn
&&=\sum_{b\in \I} \sum_{k,l} 
F^{-1}(\Y_{aa_1; i}^{a_1;(1)}\otimes \Y_{a_2a_3;j}^{a;(2)};
\Y_{a_{2}b; k}^{a_{1};(3)}\otimes \Y_{a_{3}a_{1}; l}^{b;(4)} ) \nn
&&\hspace{1.7cm} E\bigg(\tr_{W_{a_{1}}}
\Y_{a_{2}b; k}^{a_{1};(3)}(\mathcal{U}
(e^{2\pi iz_{1}})w_{a_{2}}, e^{2\pi i z_{1}})\cdot \nn
&&\hspace{3.3cm} 
\cdot \Y_{a_{3}a_{1}; l}^{b;(4)}(\mathcal{U}
(e^{2\pi i(z_{2}+\tau)})w_{a_{3}}, e^{2\pi i (z_{2}+\tau)})
q_{\tau}^{L(0)-\frac{c}{24}}\bigg).
\eea
Using the $L(0)$-conjugation formula, we can 
move $q_{\tau}$ from the right side of 
$\Y_{a_3a_1;l}^{b;(4)}$ to the left side of $\Y_{a_3a_1;l}^{b;(4)}$. Then
using the following property of trace: 
\beq  \label{tr-AB}
\tr_{W_{a_{1}}} (AB) = \tr_{W_{b}} (BA), 
\eeq
for all $A: W_{b} \rightarrow \overline{W_{a_1}}, 
B: W_{a_1}\rightarrow  \overline{W_b}$ whenever 
the multiple sums in either side of 
(\ref{tr-AB}) converge absolutely,
we obtain that the left hand side of (\ref{b-equ-1}) equals to 
\bea
&&\sum_{b\in \I}\sum_{k,l} 
F^{-1}(\Y_{aa_1; i}^{a_1;(1)}\otimes \Y_{a_2a_3;j}^{a;(2)};
\Y_{a_{2}b; k}^{a_{1};(3)}\otimes \Y_{a_{3}a_{1}; l}^{b;(4)} )  \nn
&&\hspace{0.2cm} E\bigg( \tr_{W_{b} }
\Y_{a_{3}a_{1}; l}^{b;(4)}(\mathcal{U}
(e^{2\pi iz_{2}})w_{a_{3}}, e^{2\pi i z_{2}})
\Y_{a_{2}b; k}^{a_{1};(3)}(\mathcal{U}
(e^{2\pi iz_{1}})w_{a_{2}}, e^{2\pi i z_{1}})q_{\tau}^{L(0)-\frac{c}{24}}\bigg).\nn
\label{beta-Y-equ-1}
\eea
Now apply (\ref{A-r-U}) to (\ref{beta-Y-equ-1}). We then
obtain that (\ref{beta-Y-equ-1}) equals to 
\bea
&&\sum_{b\in \I}\sum_{k,l} 
F^{-1}(\Y_{aa_1; i}^{a_1;(1)}\otimes \Y_{a_2a_3;j}^{a;(2)};
\Y_{a_{2}b; k}^{a_{1};(3)}\otimes \Y_{a_{3}a_{1}; l}^{b;(4)} )  \nn
&&\hspace{0.2cm} 
E\bigg(\tr_{(W_b)'} 
\tilde{A}_{r}(\Y_{a_{2}b; k}^{a_{1};(3)})(\mathcal{U}
(e^{(2r+1)\pi i}e^{-2\pi iz_{1}})e^{-(2r+1)\pi iL(0)} w_{a_{2}}, 
e^{(2r+1)\pi i}e^{-2\pi iz_{1}})  \nn
&&\hspace{0.5cm} \tilde{A}_{r}(\Y_{a_{3}a_{1}; l}^{b;(4)})(\mathcal{U}
(e^{(2r+1)\pi i}e^{-2\pi iz_{2}})e^{-(2r+1)\pi iL(0)} w_{a_{3}}, 
e^{(2r+1)\pi i}e^{-2\pi iz_{2}}) q_{\tau}^{L(0)-\frac{c}{24}}\bigg). \nonumber
\eea
Now apply the associativity again and be careful about the
branch cut as in \cite{H9}, then the 
left hand side of (\ref{b-equ-1}) further equals to 
\bea
&&\sum_{b\in \I}\sum_{k,l}\sum_{c\in \I} \sum_{p,q} 
F^{-1}(\Y_{aa_1; i}^{a_1;(1)}\otimes \Y_{a_2a_3;j}^{a;(2)};
\Y_{a_{2}b; k}^{a_{1};(3)}\otimes \Y_{a_{3}a_{1}; l}^{b;(4)} ) \nn
&&\hspace{2.5cm} 
F(\tilde{A}_{r}(\Y_{a_{2}b; k}^{a_{1};(3)})\otimes
  \tilde{A}_{r}(\Y_{a_{3}a_{1}; l}^{b;(4)}), 
  \Y_{cb';p}^{b';(5)} \otimes \Y_{a_2a_3;q}^{c;(6)} ) \nn
&&\hspace{0.2cm}
E\bigg(\tr_{(W_b)'} \, q_{\tau}^{L(0)-\frac{c}{24}}\,
\Y_{cb';p}^{b';(5)}(\mathcal{U}(e^{(2r+1)\pi i}e^{-2\pi iz_{2}}) \cdot \nn
&&\hspace{0.3cm}
\cdot \Y_{a_2a_3;q}^{c;(6)}( e^{-(2r+1)\pi iL(0)} w_{a_{2}}, e^{\pi i}(z_1-z_2))
e^{-(2r+1)\pi iL(0)}w_{a_{3}}, 
e^{(2r+1)\pi i}e^{-2\pi iz_{2}}) \bigg).  \nn
&&=\sum_{b\in \I}\sum_{k,l}\sum_{c\in \I} \sum_{p,q} 
F^{-1}(\Y_{aa_1; i}^{a_1;(1)}\otimes \Y_{a_2a_3;j}^{a;(2)};
\Y_{a_{2}b; k}^{a_{1};(3)}\otimes \Y_{a_{3}a_{1}; l}^{b;(4)} ) \nn
&&\hspace{2.5cm} 
F(\tilde{A}_{r}(\Y_{a_{2}b; k}^{a_{1};(3)})\otimes
  \tilde{A}_{r}(\Y_{a_{3}a_{1}; l}^{b;(4)}), 
  \Y_{cb';p}^{b';(5)} \otimes \Y_{a_2a_3;q}^{c;(6)} ) \nn
&&\hspace{0.2cm}
E\bigg(\tr_{W_b}  \, q_{\tau}^{L(0)-\frac{c}{24}}\,  
\hat{A}_r(\Y_{cb';p}^{b';(5)})(e^{(2r+1)\pi iL(0)} \nn
&&\hspace{1cm}
\Y_{a_2a_3;q}^{c;(6)}( e^{-(2r+1)\pi iL(0)} w_{a_{2}}, e^{\pi i}(z_1-z_2))
e^{-(2r+1)\pi iL(0)}w_{a_{3}}, e^{2\pi iz_{2}}) \bigg)
\eea
Choosing $r=0$ and using 
$\Y(\cdot, e^{2\pi i}x)\cdot = \Omega_0^2(\Y)(\cdot, x)\cdot$, 
we obtain
\bea  \label{b-F-equ}
&&((\Psi_2(\Y_{aa_1; i}^{a_1;(1)}\otimes \Y_{a_2a_3;j}^{a;(2)})) 
(w_{a_2}\otimes w_{a_3})) (z_1,z_2+\tau,\tau)  \nn
&&\hspace{0.5cm}=\sum_{b\in \I}\sum_{k,l}\sum_{c\in \I} \sum_{p,q} 
F^{-1}(\Y_{aa_1; i}^{a_1;(1)}\otimes \Y_{a_2a_3;j}^{a;(2)};
\Y_{a_{2}b; k}^{a_{1};(3)}\otimes \Y_{a_{3}a_{1}; l}^{b;(4)} ) \nn
&&\hspace{3cm} 
F(\tilde{A}_{0}(\Y_{a_{2}b; k}^{a_{1};(3)})\otimes
  \tilde{A}_{0}(\Y_{a_{3}a_{1}; l}^{b;(4)}), 
  \Y_{cb';p}^{b';(5)} \otimes \Y_{a_2a_3;q}^{c;(6)} ) \nn
&&\hspace{0.8cm}
E\bigg(\tr_{W_b} \hat{A}_0(\Y_{cb';p}^{b';(5)})(
\Omega_{0}^2(\Y_{a_2a_3;q}^{c;(6)})
(w_{a_{2}}, z_1-z_2)w_{a_{3}}, e^{2\pi iz_{2}}) q_{\tau}^{L(0)-\frac{c}{24}} 
\bigg).
\eea
By the linear independency proved in \cite{H11} 
of the last factor in each term of above sum,
it is clear that $\beta$ induces a map given by (\ref{b-F-map}).
\epf

\noindent {\small \sc Max-Planck-Institute for Mathematics
in the Sciences, Inselstrasse 22, D-04103, Leipzig, Germany} 

\noindent and 

\noindent {\small \sc 
Institut Des Hautes \'{E}tudes Scientifiques, 
Le Bois-Marie, 35, Route De Chartres,
F-91440 Bures-sur-Yvette, France} 

\noindent {\small \sc 
Max-Planck-Institute for Mathematics, Vivatsgasse 7, 
D-23111 Bonn, Germany}

\noindent {\em E-mail address}: kong@mpim-bonn.mpg.de

\end{document}